\documentclass{book}
\usepackage{latexsym}
\usepackage{amsmath,amssymb}
\usepackage{txfonts}
\usepackage[dvips]{graphicx, color}
\begin{document}
\title{\bf \Huge{Rank Two Non-Abelian Zeta and Its Zeros}}
\author{WENG Lin}
\date{December 4, 2004}
\maketitle
~\vskip 7.0cm
{
\centerline{\Huge Dedicated to Prof. Dr. Lion}
\vskip 0.50cm
\centerline{\huge of the $\frak{L}$ 
{\usefont{T1}{cmss}{bx}{n} Academy}}
\vskip 0.50cm
\centerline{\Huge for His Forty Years Birthday}}
\pagebreak

\noindent
{\Large Abstract.} In this paper, we first reveal an intrinsic
relation between non-abelian zeta functions and Epstein zeta functions
for algebraic number fields.
Then, we expose a fundamental relation between stability of
lattices and distance to cusps. Next,
using these two relations, we explicitly express  
rank two zeta functions in terms of the well-known 
Dedekind zeta functions. Finally, based on such an expression,
we show that all zeros of
rank two non-abelian zeta functions are entirely
sitting on the critical line whose real part equals to $\frac {1}{2}$.
\vskip 0.45cm
As such, this work is built up on the classics of number theory.
Many fine pieces of algebraic and analytic number theory 
are beautifully unified 
under the name of non-abelian zeta functions. To give the reader 
an idea on what we mean, let us only indicate the following aspects collected 
in the present work;

\noindent
1) New Geo-Arithmetic cohomology for lattices over number fields, by 
further developing Tate's fundamental work, known as Tate's Thesis,
 along with the line of Tate, Iwasawa, and van der Geer \& Schoof;

\noindent
2) A definition of non-abelian zeta functions for number fields, as a natural 
generalization (and hence offering a natural framework) for the classical 
Dedekind zeta functions;

\noindent
3) A relation between non-abelian zeta and Epstein type zeta 
functions, via the well-known Mellin transformation;

\noindent
4) A classification of lattices first according to their volumes and unit 
twists, in connection with an intrinsic relation between $GL_n$ and $SL_n$ 
over a number field $K$ using Dirichlet's Unit Theorem;
and hence a relation between isometry classes of rank two lattices over 
ring of integers and the upper half space model, 
as an extension of a discussion of Hayashi;
 
\noindent
5) A construction of a fundamental domain for the action of special 
automorphism group of rank two lattices on the associated upper half 
space using normalized Siegel type distances to cusps, 
by generalizing Siegel's original construction for totally real fields;

\noindent
6) An intrinsic relation between stability of lattices 
and  distances to cusps:
a lattice is semi-stable if and only if its distances to all cusps are
at least one, using a result of Hayashi;

\noindent 
7) A Fourier expansion for the associated Epstein zeta function, 
along with the classical line, in particular, that presented in the books of 
Kubota, Elstrodt-Grunewald-Mennicke, and Terras;

\noindent
8) An explicit expression of rank two non-abelian zeta in terms of 
 the associated Dedekind zeta function,  as an application of 
Rankin-Selberg \& Zagier method;

\noindent 
9) An analogue of the Riemann Hypothesis for  rank two non-abelian
zetas, using 
an argument of Suzuki and Lagarias, who first show that the rank two zeta 
for the field of rationals satisfies the generalized RH motivated by 
Titchmarsh and de Branges, respectively.
\vskip 2.0cm
\noindent
{\Large Acknowledgement \huge$\heartsuit$:}

\noindent
Morally, I want to express my thanks to  
XiaoYing, LerLer and AnAn, my family members,
for their constant support and understanding; to Deninger and Ueno 
for their encouragement and constant support; and to Elstrodt and 
Freitag for their interests.

\noindent
Mathematically, I would like to thank Deninger and Zagier for their
discussions; Suzuki and  Lagarias for sharing with me their discovery that 
rank two zeta for the field of rationals satisfies the generalized Riemann 
Hypothesis. 

\noindent 
Our geniune attempt to expose this beautiful part of mathematics 
began with our two weeks stay at Hokudai during
the 2004 Summer Vacation for Research of Kyudai after noticing 
that this part of the Project has been seriously delayed;
Key ingredients were found when we visited M\"unster almost immediately after; 
The whole picture emerged upon our short yet fine call at Heidelberg. 
Accordingly, 
I would like to thank Nakamura, Deninger and Freitag for their kind helps.

\noindent
The project is partially supported by JSPS.
\newpage
~
\vskip 1.0cm
\centerline{\Large Contents}
\vskip 0.45cm
\noindent
Dedication

\noindent
Abstract

\noindent
Acknowledgement
\vskip 0.5cm
\noindent
{\bf Ch.1. Non-Abelian Zeta Functions and Eisenstein Series}

\noindent
1.1. Projective $\mathcal O_K$-modules

\noindent
1.2 $\mathcal O_k$-Lattices: First Level

\noindent
1.3. Space of $\mathcal O_K$-Lattices: First Level

\noindent
1.4. Semi-Stable Lattices: First Level

\noindent
1.5. Minkowski Metric versus Canonical Metric: A Bridge

\noindent
1.6. Space of $\mathcal O_K$-Lattices via Special Linear Groups

\noindent
1.7. Structure of Moduli Space:  Action of $\mathcal O_K$-Units

\noindent
1.8. Non-Abelian Zeta Functions for Number Fields

\noindent
1.9. Non-Abelian Zeta Functions and Epstein Zeta Functions

\noindent
Appendix: Higher Dimensional Gamma Function
\vskip 0.5cm
\noindent
{\bf Ch.2. Rank Two $\mathcal O_K$-Lattices: Stability and Distance to Cusps}

\noindent
2.1. Upper Half Space Model

2.1.1. Upper Half Plane

2.1.2. Upper Half Space

2.1.3. Rank Two $\mathcal O_K$-Lattices

\noindent
2.2. Cusps

2.2.1. Upper Half Plane

2.2.2. Upper Half Space

2.2.3. Rank Two $\mathcal O_K$-Lattices

2.2.3.A. Totally Real Fields

2.2.3.B. Genaral Number Fields

\noindent
2.3. Stablizer Groups of Cusps

2.3.1. Upper Half Plane

2.3.2. Upper Half Space

2.3.3. Rank Two $\mathcal O_K$-Lattices

2.3.3.A. Totally Real Fields

2.3.3.A.1. Stablizer Groups for Cusps

2.3.3.A.2. Action of $\Gamma$

2.3.3.A.3. Fundamental Domain of $\Gamma_\lambda$

2.3.3.B. General Number Fields

\noindent
2.4. Fundamental Domain

2.4.1. Upper Half Plane

2.4.2. Upper Half Space

2.4.3. Rank 2 $\mathcal O_K$-Lattices

2.4.3.A. Totally Real Fields

2.4.3.A.1. Distance to Cusps

2.4.3.A.2. Fundamental Domain for $\Gamma$

2.4.3.B. General Number Fields

\noindent
2.5. Stability and Distance to Cusps

2.5.1. Upper Half Plane

2.5.2. Upper Half Space

2.5.3. Rank Two $\mathcal O_K$-Lattices: Level One

2.5.4. Rank Two $\mathcal O_K$-Lattices: Normalization or Convention

2.5.4.A. $SL(2,\mathbb Z)$ Acts on the Upper Half Plane

2.5.4.B. Identification Between 
$SL(2,\mathbb R)\Big/ SO(2)$ and $\mathcal H$

2.5.4.C. Metrized Structures

2.5.4.D. Automorphism Group $\mathrm{Aut}_{\mathcal O_K}(\mathcal O_K
\oplus\frak a)$

2.5.4.E. Cusp-Ideal Class Correspondence

2.5.4.F. Actions of $SL(\mathcal O_K\oplus\frak a)$ on Lattices and 
on $\mathcal H^{r_1}\times{\mathbb H}^{r_2}$

2.5.4.G. Stablizer Groups of Cusps

2.5.5.  Stability and Distance to Cusps

2.5.6. Moduli Space of Rank Two Semi-Stable $\mathcal O_K$-Lattices
\vskip 0.5cm
\noindent
{\bf Ch. 3. Epstein Zeta Functions and Their Fourier Expansions}

\noindent
3.1. Upper Half Plane

\noindent
3.2. Upper Half Space

\noindent
3.3. Rank Two $\mathcal O_K$-Lattices

3.3.1. Epstein Zeta Function and Eisenstein Series

3.3.2. Fourier Expansion: Constant Term
\vskip 0.5cm
\noindent
{\bf Ch. 4. Explicit Formula for Rank Two Zeta Functions: Rankin-Selberg 
\& Zagier Method}

\noindent
4.1. Upper Half Plane

4.1.1. Geometric Approach

Appendix: Rankin-Selberg \& Zagier Method (II)

4.1.2. Rank Two Non-Abelian Zeta Function For $\mathbb Q$

\noindent
4.2. Upper Half Space Model: Rankin-Selberg Method

\noindent
4.3. $\mathcal O_K$-Lattices of Rank Two

4.3.1. Rankin-Selberg Method

4.3.2. Explicit Formula for Rank Two Zeta
\vskip 0.5cm
\noindent
{\bf Ch. 5. Zeros of Rank Two Non-Abelian Zeta Functions for Number Fields}

\noindent
5.1. Zeros of Rank Two Non-Abelian Zeta Function of $\mathbb Q$

5.1.1. Product Formula for Entire Function of Order 1

5.1.2. Proof

5.1.3. A Simple Generalization

\noindent
5.2. Zeros of Rank Two Zetas for Number Fields: 
Generalized Riemann Hypothesis
\vskip 0.5cm
\noindent
{\bf References}
\chapter{Non-Abelian Zeta Functions and Eisenstein Series}
\section{Projective $\mathcal O_K$-modules}

Let $K$ be an algebraic number field, i.e., a finite field 
extension of  the field of rationals $\mathbb Q$.
 Denote by $\mathcal O_K$ the ring of integers of $K$. Then by definition,
an $\mathcal O_K$-module $M$ is called {\it projective} 
 if there exists an $\mathcal O_K$-module $N$ such that 
$M\oplus N$ is a free $\mathcal O_K$-module. For 
projective $\mathcal O_K$-modules, it is  well-known that
\vskip 0.20cm
\noindent
(i) If $M$ is a projective 
$\mathcal O_K$-module and is a part of a
 short exact sequence (of $\mathcal O_K$-modules) 
$0\to M_0\to M_1\to M\to 0$, then  $M_1\simeq M_0\oplus M$.

\noindent
(ii) All fractional $\mathcal O_K$-ideals are projective; and

\noindent
(iii) Rank 1 projective $\mathcal O_K$-submodules in 
$K$ are simply fractional $\mathcal O_K$-ideals.
\vskip 0.20cm
\noindent
Thus, by  finiteness of the ideal class group of $K$,
up to equivalence relation (defined by isomorphisms as 
$\mathcal O_K$-modules), there are only finitely many 
rank 1 projective $\mathcal O_K$-modules in $K$.
 Therefore, we may choose integral $\mathcal O_K$-ideals 
$\frak a_i, i=1,\cdots,h$ with $h=h(K)$, the
  class number of $K$, such that 
\vskip 0.20cm
\noindent
(a) Any rank 1 projective $\mathcal O_K$-module is 
isomorphic to one of the $\frak a_i$; while

\noindent
(b) None of the $\frak a_i$ and $\frak a_j$ are isomorphic
 to each other  if $i\not=j$, $i,j=1,\cdots,h$.
\vskip 0.20cm
\noindent
We will fix a choice of $\frak a_i, i=1,\cdots, h$ satisfying 
(a) and (b) above for the rest of this paper, and sometimes use 
$\frak a$ as a running symbol for them.
\vskip 0.30cm
 With this discussion of the rank 1  projective 
$\mathcal O_K$-modules, now let us consider higher rank 
 $\mathcal O_K$-modules. Clearly  for a fractional 
ideal $\frak a$,
 $$P_{\frak a}:=P_{r;\frak a}:=\mathcal O_K^{r-1}\oplus 
\frak a$$ is a rank $r$ projective 
 $\mathcal O_K$-module. The nice thing is that such 
types of projective $\mathcal O_K$-modules, up to
  isomorphism, give all rank $r$ projective 
$\mathcal O_K$-modules. Indeed, we have the following structural
 \vskip 0.30cm
\noindent
{\bf{\Large Proposition}.} 
 (1)  {\it For fractional ideals $\frak a$ and $\frak b$, 
$P_{r;\frak a}\simeq P_{r;\frak b}$ if and only if 
 $\frak a\simeq \frak b$;}
 
 \noindent
 (2) {\it For a rank $r$  projective $\mathcal O_K$-module $P$, 
there exists a  fractional ideal $\frak a$ such that
 $P\simeq P_{\frak a}.$}
 
 The proof of this proposition is based on an induction on the rank 
and the following basic relation about  
 fractional ideals: For two   fractional ideals $\frak a$ and $\frak b$,
as $\mathcal O_K$-modules, $\frak a\oplus\frak b
 \simeq \mathcal O_K\oplus\frak a\frak b.$ 
 (The reader can find a complete proof in [FT].) So we omit the details.
\vskip 0.30cm  
Now use the natural inclusion of fractional ideals in $K$ to embed
$P_{r;\frak a}$ into $K^r$. View an element in $K^r$ as a column
vector. In such a way, any $\mathcal O_K$-module
 morphism  $A:P_{r;\frak a}\to P_{r;\frak b}$ may 
be written down as an element in $A\in M_{r\times r}(K)$ so that
the image $A(x)$ of $x$ under $A$ becomes simply the matrix multiplication
 $A\cdot x$. In particular, one checks easily the following 
\vskip 0.30cm
\noindent
{\bf{\large Lemma}.} {\it If $A\in GL(r,F)$ 
defines an $\mathcal O_K$-isomorphism 
$A:P_{r;\frak a}\to P_{r;\frak b}$, then 
$$\frak b\simeq (\det A)\cdot\frak a.$$ 

\noindent
In particular, in the case when both $\frak a$ and $\frak b$ are integral 
$\mathcal O_K$-ideals,}

\noindent
(i) $\det A\in U_K$, {\it the group of units of $K$;}

\noindent
(ii)  $\frak a=\frak b$; {\it and} 

\noindent
(iii) $A\in \mathrm{Aut}_{\mathcal O_K}(P_{r;\frak a}).$

\section{$\mathcal O_K$-Lattices: First Level}

We start our first level discussion on lattices in the next three sections,
mainly follow [Gr] for the expositions.

Let $\sigma$ be an Archimedean place of $K$, and  
$K_\sigma$ be the $\sigma$-completion of $K$. It is 
well-known that the $\mathbb R$-algebra $K_\sigma$ is 
either equal to $\mathbb R$, or equal to $\mathbb C$.
 Accordingly, we call $\sigma$ (to be) real or complex, 
write sometimes in terms of $\sigma:\mathbb R$ or 
 $\sigma:\mathbb C$ accordingly.

If $V_\sigma$ is a finite dimensional $K_\sigma$-vector 
space, then an inner product on $V_\sigma$ is a 
positive definite bilinear form $V_\sigma\times V_\sigma
\to K_\sigma$ which is symmetric if $\sigma$ is real 
and is hermitian if $\sigma$ is complex. When equipped 
with an inner product, $V_\sigma$ is called a {\it metrized 
space}. 
\vskip 0.20cm
By definition, an $\mathcal O_K$-{\it lattice} $\Lambda$ 
consists of 

\noindent
(1) a projective $\mathcal O_K$-module 
$P=P(\Lambda)$ of finite rank; and 

\noindent
(2) an inner product on the vector space 
$V_\sigma:=P\otimes_{\mathcal O_K}K_\sigma$ for each of the
 Archmidean place $\sigma$ of $K$.

\noindent
Set $V=P\otimes_{\mathbb Z}\mathbb R$ 
so that $V=\prod_{\sigma\in S_\infty}V_\sigma$,
  where $S_\infty$ denotes the collection of all (inequivalent) 
Archimedean places of $K$. Indeed, this is a
  direct consequence of the fact that as a $\mathbb Z$-module,
 an $\mathcal O_K$-ideal is of rank $n=r_1+2r_2$ 
  where $n=[K:\mathbb Q]$, $r_1$ denotes the number of real places 
and $r_2$ denotes the number of complex places
   (in $S_\infty$).

\section{Space of $\mathcal O_K$-Lattices: First Level}

Let $P$ be a rank $r$ projective $\mathcal O_K$-module. Denote by 
$GL(P):=\mathrm{Aut}_{\mathcal O_K}(P)$. Let
 $\widetilde{\mathbf{\Lambda}}:=\widetilde{\mathbf{\Lambda}}(P)$ be 
the space of ($\mathcal O_K$-)lattices $\Lambda$ whose 
 underlying $\mathcal O_K$-module is $P$. For $\sigma\in S_\infty$, 
let $\widetilde{\mathbf{\Lambda}}_\sigma$ be 
 the space of inner products on $V_\sigma$; if a basis is chosen for
 $V_\sigma$ as a real or a complex vector space according to whether
$\sigma$ is real or complex,   
 $\widetilde{\mathbf{\Lambda}}_\sigma$ may be realized as
 an open set of a real or complex vector space. (See \S1.6 below
for details.) We have
  $\widetilde{\mathbf{\Lambda}}=\prod_{\sigma\in S_\infty}
\widetilde{\mathbf{\Lambda}}_\sigma$ and this provides 
  us with a natural topology on $\widetilde{\mathbf{\Lambda}}.$

As in \S1.1, consider $GL(P)$ to act on $P$ from the left. 
With the chosen basis, elements of 
$\widetilde{\mathbf{\Lambda}}$ will be thought of 
as column vectors, and matrices of linear maps will
 be written on the left.
\vskip 0.30cm
Given $\Lambda\in \widetilde{\mathbf{\Lambda}}$ and
$u,w\in V_\sigma$, let $\langle u,w\rangle_{\Lambda,\sigma}$ or
$\langle u,w\rangle_{\rho_\Lambda(\sigma)}$
denote the value of the inner product on the 
vectors $u$ and $w$ associated to the lattice $\Lambda$.
 
As such, if $A\in GL(P)$, we may define a new 
lattice $A\cdot\Lambda $ in 
$\widetilde{\mathbf{\Lambda}}$ by the
following formula $$\langle u,w
\rangle_{A\cdot\Lambda ,\sigma}:=
\langle A^{-1}\cdot u,A^{-1}\cdot w\rangle_{\Lambda,\sigma}.$$ 
This defines an action of $GL(P)$ on 
$\widetilde{\mathbf{\Lambda}}$ from the right. Clearly, 
then the map $v\mapsto Av$ gives an isometry
 $\Lambda\cong A\cdot\Lambda $ of the lattices. 
(By an isometry here, we mean an isomorphism of 
 $\mathcal O_K$-modules for the underlying 
$\mathcal O_K$-modules subjecting the condition that the
 isomorphism also keeps the inner product unchanged.) 
Conversely, suppose that $A:\Lambda_1\cong\Lambda_2$
  is an isometry of $\mathcal O_K$-lattices, each of 
which is in $\widetilde{\mathbf{\Lambda}}$. Then, $A$
   defines an element, also denoted by $A$, of $GL(P)$. 
Clearly $\Lambda_2\cong A\cdot \Lambda_1.$
(Here we follow a nice notation initiated by Miyaoka; 
an isomorphism of corresponding algebraic structures
 is only denoted by $\simeq$, while an
 isometry is denoted by $\cong$.)
\vskip 0.20cm 
Therefore, the orbit set $GL(P)\backslash 
\widetilde{\mathbf{\Lambda}}(P)$ can be regarded as the set of isometry
 classes of $\mathcal O_K$-lattices whose underlying 
$\mathcal O_K$-modules are isomorphic to $P$.
\vskip 0.30cm
We end this section by introducing an operation among
 the lattices in $\widetilde{\mathbf{\Lambda}}$. If $T$ 
is a positive real number, then from $\Lambda$, 
we can produce a new $\mathcal O_K$-lattice called
$\Lambda[T]$ by multiplying each of the inner products 
on $\Lambda$, or better, on $\Lambda_\sigma$ for 
$\sigma\in S_\infty$, by $T^2$. Let then 
${\mathbf{\Lambda}}={\mathbf{\Lambda}}(P)$ be the quotient of
 $\widetilde{\mathbf{\Lambda}}$ by the equivalence relation 
$\Lambda\sim\Lambda[T]$. As such, 
 ${\mathbf{\Lambda}}$ admits a natural topological
 structure as well. Furthermore, as it becomes clear later, 
the construction of ${\mathbf{\Lambda}}$
 from $\widetilde{\mathbf{\Lambda}}$  plays a key role when
we want to get the compactness statement for our moduli spaces.
(Indeed,  the $[T]$-construction naturally fixes 
a specific volume for a certain family of lattices,
while  does not really change the \lq essential' structures of 
lattices involved. Consequently, by reduction theory,
semi-stable lattices of a fixed volume form a compact family.)

\section{Semi-Stable Lattices: First Level}

Let $\Lambda$ be an $\mathcal O_K$-lattice with 
underlying $\mathcal O_K$-module $P$. Then any submodule
$P_1\subset P$ can be made into an 
$\mathcal O_K$-lattice by restricting the inner 
product on each $V_\sigma$ 
to the subspace $V_{1,\sigma}:=P_1\otimes_KK_\sigma$. 
Call the resulting $\mathcal O_K$-lattice $\Lambda_1:=\Lambda\cap P_1$
 and write $\Lambda_1\subset \Lambda$. If moreover,
 $P/P_1$ is  projective, we say that $\Lambda_1$
is a {\it sublattice} of $\Lambda$.

The orthogonal projections $\pi_\sigma:V_\sigma\to 
V_{1,\sigma}^\perp$ to the orthogonal complement $V_{1,\sigma}^\perp$
of $V_{1,\sigma}$ in $V_{\sigma}$
provide isomorphisms 
$(P/P_1)\otimes_{\mathcal O_K}K_\sigma\simeq V_{1,\sigma}^\perp$, 
which can be used to make $P/P_1$ into an
 $\mathcal O_K$-lattice. We call this resulting lattice the 
 {\it quotient lattice} of $\Lambda$ by $\Lambda_1$, and denote it by 
$\Lambda/\Lambda_1$. 
\vskip 0.20cm
There is a procedure called {\it restriction of scalars} 
which makes an $\mathcal O_K$-lattice into a standard 
$\mathbb Z$-lattice. Recall that $V=\Lambda
\otimes_{\mathbb Z}\mathbb R=\prod_{\sigma\in S_\infty}V_\sigma$. 
Define an inner product on the real vector space $V$ by
 $$\langle u,w\rangle_\infty:=
\sum_{\sigma:\mathbb R}\langle u_\sigma,w_\sigma\rangle_\sigma
+\sum_{\sigma:\mathbb C}\mathrm{Re}\,\langle u_\sigma,
w_\sigma\rangle_\sigma.$$ Let $\mathrm{Res}_{K/\mathbb Q}\Lambda$ 
denote the $\mathbb Z$-lattice obtained by 
equipped $P$, regarding as a $\mathbb Z$-module, with this 
inner product (at the unique infinite place 
$\infty$ of $\mathbb Q$).

We let $\mathrm{rk}(\Lambda)$ denote the $\mathcal O_K$-module
 rank of $P$ (or of $\Lambda$) and let 
$\mathrm{dim}(\Lambda)$  denote the rank of $P$ as
 $\mathbb Z$-module. Clearly, 
$$\mathrm{dim}(\Lambda)=\mathrm{rk}(\Lambda)\cdot 
\mathrm{dim}(\mathcal O_K)=\mathrm{rk}(\Lambda)
\cdot[F:\mathbb Q].$$
We define the {\it Lebesgue volume} of $\Lambda$, 
denoted by $\mathrm{Vol}_{\mathrm{Leb}}(\Lambda)$, to be
the (co)volume of the lattice $\mathrm{Res}_{K/\mathbb Q}
\Lambda$ inside its inner product space $V$. As such, this 
volume may be computed as $\Big|\mathrm{det}\,\langle l_i,e_j\rangle\Big|$,
 where $\{l_i\}$ is a $\mathbb Z$-basis of 
$\mathrm{Res}_{K/\mathbb Q}\Lambda$ and $\{e_j\}$ is 
an orthonormal basis of $\Lambda$ with respect to 
$\langle\cdot,\cdot\rangle_\infty$. For examples, 

\noindent
(a) If $\mathrm{dim}\Lambda=0$, then 
$\mathrm{Vol}_{\mathrm{Leb}}(\Lambda)=1;$

\noindent
(b) If $\mathrm{dim}\Lambda=1$, then 
$\mathrm{Vol}_{\mathrm{Leb}}(\Lambda)$ is the length of a generator of 
$\Lambda$;

\noindent
(c) If $\mathrm{dim}\Lambda=2$, then 
$\mathrm{Vol}_{\mathrm{Leb}}(\Lambda)$ is the area of a fundamental 
parallelopiped, and so on.

Clearly, if $P'$ is a submodule of 
finite index in $P$, then 
$$\mathrm{Vol}_{\mathrm{Leb}}(\Lambda')=
[P:P']\mathrm{Vol}_{\mathrm{Leb}}(\Lambda),$$ 
where $\Lambda'=\Lambda\cap P$ is the lattice induced from $P'$.
\vskip 0.20cm
\noindent 
{\bf{\large Examples}}: 1) Take $P=\mathcal O_K$ and for each place 
$\sigma$, let $\{1\}$ be an orthonormal basis of 
$V_\sigma=K_\sigma$, i.e., equipped 
$V_\sigma=\mathbb R$ or $\mathbb C$ with the standard Lebesgue
 measure. This makes $\mathcal O_K$ into an 
$\mathcal O_K$-lattice $\overline{\mathcal O_K}=(\mathcal O_K,\mathbf 1)$ 
in a natural  way. It is a well-known fact, see e.g.,
[L1], that $$\mathrm{Vol}_{\mathrm{Leb}}
\Big(\overline{\mathcal O_K}\Big)=2^{-r_2}\cdot \sqrt{\Delta_F},$$ where
  $\Delta_F$ denotes the absolute value of the discriminant of $K$.

 More generally, take $P=\frak a$ an fractional 
  idea of $K$ and equip the same inner product as above 
on $V_\sigma$. Then $\frak a$ becomes an
   $\mathcal O_K$-lattice $\overline {\frak a}=(\frak a,\mathbf 1)$ 
in a natural 
way with $\mathrm{rk}(\frak a)=1$.  It is a well-known fact, see e.g.,
[Neu], that
$$\mathrm{Vol}_{\mathrm{Leb}}\Big(\overline {\frak a}\Big)=
2^{-r_2}\cdot \Big(N(\frak a)\cdot\sqrt{\Delta_K}\Big),$$ 
where $N(\frak a)$ denote the norm of $\frak a$.
\vskip 0.30cm
Due to the appearence of the factor $2^{-r_2}$, we  also define 
the {\it canonical volume} of $\Lambda$, denoted
 by $\mathrm{Vol}_{\mathrm{can}}(\Lambda)$ or simply by
$\mathrm{Vol}(\Lambda)$, to be 
$2^{r_2\mathrm{rk}(\Lambda)}\mathrm{Vol}_{\mathrm{Leb}}(\Lambda).$ 
(This canonical volume is a theoretically correct one.  
In fact, in Arakelov theory, where the
 minus of the log of the canonical vomule
 is defined to be  the Arakelov-Euler characteristic $\chi(\Lambda)$
of $\Lambda$. That is, 
$$\chi(\Lambda):=-\log\Big(\mathrm{Vol}(\Lambda)\Big).$$
See also the justification  given in the next section.) So in 
particular, $$\mathrm{Vol}\Big(\overline 
{\frak a}\Big)=
N(\frak a)\cdot\sqrt{\Delta_K},$$ with $$\mathrm{Vol}
\Big(\overline {\mathcal O_K}\Big)
=\sqrt{\Delta_K}$$ as its special case.

\noindent 
2) The $[T]$-construction changes volumes of lattices in the following way
 $$\mathrm{Vol}(\Lambda[T])
=T^{\mathrm{dim}(\Lambda)}\cdot \mathrm{Vol}(\Lambda).$$

Now we are ready to introduce our first key definition.

\vskip 0.20cm
\noindent
{\bf \Large Definition.} An $\mathcal O_K$ lattice $\Lambda$ is called 
{\it semi-stable} (resp. {\it stable}) if for any proper 
sublattice $\Lambda_1$ of $\Lambda$,
 $$\mathrm{Vol}(\Lambda_1)^{\mathrm{rk}(\Lambda)}\geq\,(\mathrm{resp.}>)
\,\mathrm{Vol}(\Lambda)^{\mathrm{Vol}(\Lambda_1)}.$$ 

Clearly the last inequality is equivalent to 
$$\mathrm{Vol}_{\mathrm{Leb}}(\Lambda_1)^{\mathrm{rk}(\Lambda)}\geq
\mathrm{Vol}_{\mathrm{Leb}}
(\Lambda)^{\mathrm{Vol}(\Lambda_1)}.$$ So it does not matter which volume, 
the canonical one or the Lebesgue one, we use.
\vskip 0.20cm
\noindent
{\bf Remark.} Despite the fact that we introduce the stability for lattices
independently, many others, notably Stuhler,
introduced the stability earlier. (See e.g. [Gr], [St].) 

\section{Minkowski Metric versus Canonical Metric: A Bridge}

This section is specially introduced for the reader who wants to 
see clearly the relation between
Lebesgue  and  canonical volumes. We mainly 
follow [Neu] for the presentation. For the first reading, one can skip it.
\vskip 0.30cm
Minkowski's basic idea using Geometry of Numbers to study 
algebraic number field $K/\mathbb Q$ of degree $n$ is to 
 interpret its numbers as points in an $n$-dimensional space.
To view such points, let us consider the canonical map 
$$j:K\to K_{\mathbb C}:=\prod_\tau\mathbb C, a\mapsto j(a):=
(\tau a)$$ induced from the $n$ complex 
embeddings $\tau:K\to\mathbb C$. (Here, if $\sigma:K\to
 K_\sigma=\mathbb R$ is a real Archimedean place, 
we use the natural embedding 
 $\mathbb R\hookrightarrow\mathbb C=\mathbb R+i\mathbb R$ 
 so as to get a natural map $\tau:K\to\mathbb C$ ending 
with $\mathbb C$.) The $\mathbb C$-vector 
 space $K_\mathbb C$ is equipped with the hermitian 
scalar product $$\langle x,y\rangle=\sum_\tau x_\tau
  \bar y_\tau.$$
In the sequel, we always view $K_{\mathbb C}$ as the
 hermitian space with respect to this standard metric.
\vskip 0.20cm
Let $G(\mathbb C|\mathbb R)$ be the Galois group generated by 
complex conjugation $F:z\mapsto\bar z$. Then $F$ acts
 on both   the factors $\mathbb C$ of the product $\prod_\tau\mathbb C$
 and  on the index set of $\tau$'s at the same time; 
 that is to say, $a\mapsto\bar a$ for $a\in\mathbb C$, while $\tau\mapsto
\bar\tau$ for each embedding $\tau:K\to \mathbb C$ with associated complex
  conjugate $\bar\tau:K\to\mathbb C$. Altogether,
 this defines an involution $F:K_{\mathbb C}\to
   K_{\mathbb C}$, which, in terms of points 
$z=(z_\tau)\in K_{\mathbb C}$, is given by
   $(Fz)_\tau=\bar z_{\bar\tau}.$ Clearly,
the scalar product $\langle\cdot,\cdot\rangle$ is equivariant under $F$. 
 That is to say, $\langle Fx,Fy\rangle=\langle x,y\rangle.$

We now concentrate on the $\mathbb R$-vector space 
$$K_{\mathbb R}:=K_{\mathbb C}^+:=\Big[\prod_\tau\mathbb C\Big]^+$$
 consisting of the $G(\mathbb C|\mathbb R)$-invariant, 
i.e., $F$-invariant, points of $K_\mathbb C$. Easily,
the $F$-invariant points of $K_\mathbb C$  are exactly these 
points  $(z_\tau)$ such that 
$z_{\bar\tau}=\bar z_\tau$. In particular, then we have
$$K_{\mathbb R}\simeq \mathbb R^{r_1}\times \mathbb C^{r_2}.$$ via the map 
$$(x_{\sigma_1},\dots, x_{\sigma_{r_1}};
z_{\sigma_{r_1+1}},\bar z_{\sigma_{r_1+1}},
\cdots,z_{\sigma_{r_1+r_2}},\bar z_{\sigma_{r_1+r_2}})\mapsto
(x_{\sigma_1},\dots, x_{\sigma_{r_1}};
z_{\sigma_{r_1+1}},,\cdots,z_{\sigma_{r_1+r_2}}).$$ 
Since $\overline{\tau}(a)=\overline{\tau a}$ for $a\in K$, we have
$F(j(a))=F(a)$ for all $a\in K$. This then further induces a map 
$j:K\to K_{\mathbb R}$. 

The restriction of the hermitian scalar product 
$\langle\cdot,\cdot\rangle$ from $K_{\mathbb C}$ to 
$K_{\mathbb R}$ gives a scalar product 
$$\langle\cdot,\cdot\rangle:K_{\mathbb R}\times K_{\mathbb R}\to
 \mathbb R$$ on the
$\mathbb R$-vector space $K_\mathbb R$. We call the
 Euclidean vector space 
$K_\mathbb R=[\prod_\tau \mathbb C]^+$ the {\it Minkowski 
space}, its scalar product $\langle\cdot,\cdot\rangle$
 the {\it canonical metric}, and the association
Haar measure the canonical measure. 

Moreover, the map $j:K\to K_\mathbb R$ identifies the 
vector space $K_\mathbb R$ with the tensor product
 $K\otimes_\mathbb Q\mathbb R.$ That is, we have  a natural identification
 $$K\otimes_\mathbb Q
\mathbb R\simeq K_\mathbb R,\qquad 
 a\otimes x\mapsto (j(a))\cdot x.$$ Likewise,
$K\otimes_\mathbb Q\mathbb C\simeq K_\mathbb C.$ With this said, 
then the inclusion
 $K_\mathbb R\subset K_\mathbb C$ corresponds exactly to the canonical map 
 $K\otimes_\mathbb Q\mathbb R\to K
\otimes_\mathbb Q\mathbb C$  induced from the natural inclusion 
 $\mathbb R\hookrightarrow \mathbb C$. And simply, $F$ 
corresponds to $a\otimes z\mapsto a\otimes \bar z$.
 \vskip 0.30cm
On the other hand, the Lebesgue volume is associated with the following
alternative  explicit descripion of the Minkowski space $K_\mathbb R$. 
When the  embeddings $\tau:K\to\mathbb C$ are real,  their images 
land already in $\mathbb R$. As for
complexes, i.e., these which are not real, by writing them in pairs, 
we may list (all real embeddings as 
$\sigma_1,\cdots,\sigma_{r_1}:K\to\mathbb R$, and) all complex (conjugate)
embeddings as 
$\tau_1,\bar\tau_1,\cdots,\tau_{r_2},\bar\tau_{r_2}:K\to\mathbb C$. 
(Here $[K:\mathbb Q]=n=r_1+2r_2.$) Choose from each pair $\tau_j,\,\bar\tau_j$
a certain 
fixed complex embedding $\tau_j$. With this done, then let $\sigma$ vary 
over the family of real embeddings and $\tau$ over the 
family of chosen complex embeddings.
Write $\frak p$ as a (common) running symbol for $\sigma_i$ and 
$\{\tau_j,\,\bar\tau_j\}$. Since $F$ leaves 
the $\sigma$ invariant, but exchanges $\tau$ and $\bar\tau$, we have
 $$K_{\mathbb R}=\Big\{(z_\frak p)\in\prod_\frak p\mathbb C:
z_\sigma\in\mathbb R,z_{\bar\tau}=\overline
  {z_\tau}\in \mathbb C\Big\}.$$ As a direct consequence, 
we then obtain a natural  
isomorphism $$f:K_\mathbb R\to \prod_\frak p
  \mathbb R=\mathbb R^{r_1+2r_2=n=[K:\mathbb Q]},\qquad(z_\frak p)
\mapsto(x_\frak p),$$ where for the reals 
$x_\sigma=z_\sigma$ while for the complexes $x_\tau=\Re(z_\tau), 
x_{\bar\tau}=\Im(z_\tau).$ (So in particular
for complex $\tau$, 
$$(z_\tau,z_{\bar\tau})=(x_\tau+iy_\tau,x_\tau-iy_\tau)\mapsto (x_\tau,y_\tau),\qquad\mathrm{and}\qquad
(x_\tau^2+y_\tau^2)+(x_\tau^2+y_\tau^2)\mapsto 2(x_\tau^2+y_\tau^2).)$$

This isomorphism transforms the canonical metric 
$\langle\cdot,\cdot\rangle$ on $K_{\mathbb R}$ into the scalar product
 $(x,y)=\sum_\frak pN_\frak p x_\frak p y_\frak p$ 
where $N_\frak p=1$ resp. $N_\frak p=2$ if $\frak p$ 
 is real, resp. complex. The scalar product $(x,y)
=\sum_\frak pN_\frak p x_\frak p y_\frak p$ transforms 
 the canonical measure from $K_\mathbb R$ to a measure on  
$\mathbb R^{r_1+2r_2}$. It obviously differes from the standard
  Lebesgue measure by $$\mathrm{Vol}_{\mathrm{can}}(X)
=2^{r_2}\mathrm{Vol}_{\mathrm{Leb}}(f(X)).$$
   Minkowski himself worked with the Lebesgue measure on 
$\mathbb R^{r_1+2r_2}$, and afterwards, most works, e.g., 
    research papers and textbooks, follow suit. As we said before, 
we will write $\mathrm{Vol}_{\mathrm{can}}$ 
   simply as $\mathrm{Vol}$. Also in the calculations below, 
we use both canonical and Lebesgue measures without any clear 
indication on which one we really use for our convenience. 
\vskip 0.20cm
As said too,  the canonical measures has an 
advantage theoretically. For example, we have
    the following

\noindent
{\bf {\large Arakelov-Riemann-Roch Formula}:} 
For an  $\mathcal O_K$-lattice $\Lambda$ of rank $r$, 
$$-\log \Big(\mathrm{Vol}(\Lambda)\Big)=
\mathrm{deg}(\Lambda)-
\frac {r}{2}\log\Delta_K.$$ 
(For the reader who does not know the definition of the Arakelov 
degree, she or he may simply take this 
relation as the definition.)

\section{Space of $\mathcal O_K$-Lattices via Special Linear Groups}

Recall that, by definition, an $\mathcal O_K$-lattice 
$\Lambda$ consists of two aspects, i.e., a underlying projective 
$\mathcal O_K$-module $P$ and a metric structure on the space 
$V=\Lambda\otimes_{\mathbb Z}\mathbb R=
\prod_{\sigma\in S_\infty}V_\sigma.$ Moreover, for 
the projective $\mathcal O_K$-module $P$, in assuming 
that the $\mathcal O_K$-rank of $P$ is $r$, we can identify 
$P$ with one of the $P_i:=P_{r;\frak a_i}:=
\mathcal O_K^{(r-1)}\oplus\frak a_i$, where
 $\frak a_i, i=1,\cdots,h,$ are chosen integral $\mathcal O_K$-ideals 
of \S1.1
so that $\Big\{[\frak a_1],\,[\frak a_2],\,\ldots,\,[\frak a_h]\Big\}=CL(K)$
 the class group of $K$. 
In the sequel, we often use $P$ as a running symbol for the $P_i$'s.

With this said,  via the Minkowski embedding 
$K\hookrightarrow \mathbb R^{r_1}\times\mathbb C^{r_2}$,  
we obtain a natural embedding for $P$:
$$P:=\mathcal O_K^{(r-1)}\oplus\frak a\hookrightarrow 
K^{(r)}\hookrightarrow \Big(\mathbb R^{r_1}\times
\mathbb C^{r_2}\Big)^{r}\cong \Big(\mathbb R^r\Big)^{r_1}\times
\Big(\mathbb C^r\Big)^{r_2},$$ which is simply the space 
$V=\Lambda\otimes_{\mathbb Z}\mathbb R$ above. 
As a direct consequence,  our lattice $\Lambda$ then is determined
by a metric structure on $V=\prod_{\sigma\in S_\infty}V_\sigma$, or better, on 
$\big(\mathbb R^r\big)^{r_1}\times\big(\mathbb C^r\big)^{r_2}$.
Hence, we need to determine all metrized structures on 
$\big(\mathbb R^r\big)^{r_1}\times\big(\mathbb C^r\big)^{r_2}$.

For doing so, let us start with each component  $\mathbb R^{r}$ 
(resp. $\mathbb C^{r}$). From linear algebra, metrized structures
are characterized by the following two results:

\noindent
(i) For any $g\in GL(r,\mathbb R)$ (resp. 
$g\in GL(r,\mathbb C)$), there is an associated metric structure
$\rho(g)$ or simply $g$ on 
$\mathbb R^{r}$ (resp. on $\mathbb C^{r}$)  defined by the matrix 
$g\cdot g^t$ (resp. $g\cdot\overline{g^t}$). More precisely, 
for $x,\,y\in \mathbb R^{r}$ 
(resp. $\mathbb C^{r}$)
$$\langle x,y\rangle_{g}:=\langle x,y\rangle_{\rho(g)}:=x\cdot (g \overline{g^t})\cdot y^t=(xg)\cdot \overline{(yg)^t};$$

\noindent
(ii) Two matrices $g$ and $g'$ in $GL(r,\mathbb R)$ (resp. 
in $GL(r,\mathbb C)$) correspond to the same metrized 
structure on $\mathbb R^{r}$ (resp. $\mathbb C^{r}$) if and only 
if there is a matrix $A\in GL(r,\mathbb R)$ (resp. $GL(r,\mathbb C)$) such that
$g'=g\cdot A$ and $A\cdot A^t=E_r$ 
(resp. $A\cdot\overline{A^t}=E_r$). That is to say, $g$ and $g'$ differ
from each other by a matrix $A$ from the orthogonal 
group $O(r)$ (resp. from the unitary group 
$U(r)$). 
 
Therefore, all metrized structures on  $\mathbb R^{r}$ 
(resp. on $\mathbb C^{r}$) are parametrized by the quotient space
$GL(r,\mathbb R)/O(r)$ (resp. $GL(r,\mathbb C)/U(r)$).
Consequently, metrized structures on 
$\big(\mathbb R^r\big)^{r_1}\times\big(\mathbb C^r\big)^{r_2}$ 
are parametrized by the space 
$$\Big(GL(r,\mathbb R)\Big/O(r)\Big)^{r_1}\times 
\Big(GL(r,\mathbb C)\Big/U(r)\Big)^{r_2}.$$
 
Now what comes into our discussion is the construction of a new lattice 
$\Lambda[T]$ from $\Lambda$ associated to a fixed 
 positive real number $T$. In essence, such a construction makes 
it possible for us to concentrate only on a fixed \lq level' of the 
volumes for the lattices involved. Indeed,   the set of volumes of all
lattices in $\bold{\Lambda}(P)$ can be easily seen to be
coincided with the set of positive
real numbers $\mathbb R_+^*$.  And the $[T]$-construction, scaling only the 
metric by a constant factor,
does not change any other structures of the lattices. 
So we can effectively focus our attention only to the lattices with 
a fixed volume, say, 1 in the case of the field 
of rationals, or better,
 $N(\frak a)\cdot \Delta_K^{\frac{r}{2}}$ in the case of 
general number fields $K$.) 
 
Motivated by such a $[T]$-construction for lattices in mind, naturally
at the group level, we need to shift our discussion from the
general linear group $GL$ to the special linear group $SL$.
For this, let us start with a local discussion on $\mathcal O_K$-lattice 
structures. 
 
First, look at complex places $\tau$, which are easier. Since 
we are working over $\mathbb C^r$, so
the metric structures are parametrized by
  $GL(r,\mathbb C)/U(r)$. Clearly, by fixing a branch of the 
$n$-th root, we get  natural identifications
  $$\begin{matrix}GL(r,\mathbb C)&\to& SL(r,\mathbb C)\times \mathbb C^*&\to& 
 SL(r,\mathbb C)\times S^1 \times \mathbb R_+^*\\
 g&\mapsto&(\frac{1}{\root{r}\of{\det g}}g,\det g)&\mapsto& 
(\frac{1}{\root{r}\of{\det g}}g,\frac{\det g}{|\det g|},|\det g|)\end{matrix}$$ and $$U(r)\to SU(r)\times S^1,\qquad  
U\mapsto(\frac{1}{\root{r}\of{\det U}}U,\det U),$$ where 
$SL$ (resp. $SU$) denotes the special linear group
 (resp. the special unitary group) and $S^1$ denotes the 
unit circle $\{z\in\mathbb C: |z|=1\}$ in 
 $\mathbb C^*$. As a direct consequence,  
we obtain the following natural identification $$GL(r,\mathbb C)\Big/U(r)\cong 
\Big(SL(r,\mathbb C)\Big/SU(r)\Big)\times \mathbb R_+^*.$$
In particular, we can then use this latest quotient space 
$\Big(SL(r,\mathbb C)\Big/SU(r)\Big)\times \mathbb R_+^*$ induced 
from the spacial
linear group $SL$ to parametrize all metric structures on $\mathbb C^r$.
 
Then, let us trun to real places $\sigma$, which are slightly complicated. 
Since we are working over $\mathbb R^r$ now,
 the metric structures are  parametrized by $GL(r,\mathbb R)\Big/O(r)$. 
Here, one might try to use the same approach for $\mathbb C$ above 
 for the reals as well. However, this does not work directly, simply because 
$\root{r}\of{\det g}$ is not always well-defined in the reals (say, when
  $\det g$ is not positive). Thus, alternatively,  
as an intermediate step, we use the subgroups 
$$GL^+(r,\mathbb R):=
  \{g\in GL(r,\mathbb R):\det g>0\}\quad 
\mathrm{and}\quad O^+(r):=\{A\in O(r,\mathbb R):\det g>0\}.$$ 
Clearly, we have the following relations;

\noindent
(i) $O^+(r)=SO(r)$, the special orthogonal group consisting of 
these $A$'s in $O(r)$ whose determinants are exactly 1; and 

\noindent
(ii) $GL(r,\mathbb R)\Big/O(r)\cong GL^+(r,\mathbb R)\Big/SO(r);$ moreover

\noindent
(iii) There is an identification  
$$GL^+(r,\mathbb R)\to  SL(r,\mathbb R)\times \mathbb R_+^*, \qquad
  g\mapsto
  (\frac{1}{\root{r}\of{\det g}}g,\det g).$$ 

\noindent
As a direct consequence,  we obtain a natural identification
 $$GL(r,\mathbb R)\Big/O(r)\cong \Big(SL(r,\mathbb R)\Big/SO(r)\Big)
\times \mathbb R_+^*.$$ In particular, we can then use this latest quotient space 
$\big(SL(r,\mathbb R)/SO(r)\big)\times \mathbb R_+^*$ induced from the spacial
linear group $SL$ to parametrize all metric structures on $\mathbb R^r$.
\vskip 0.30cm 
Now we are ready to resume our global 
discussion on $\mathcal O_K$-lattices of rank $r$.
From  above, the metrized structures on 
$V=\prod_{\sigma\in S_\infty} V_\sigma\simeq(\mathbb R^r)^{r_1}\times
(\mathbb C^r)^{r_2}$ are  
parametrized by the space $$\Big(\Big(SL(r,\mathbb R)\Big/SO(r)\Big)^{r_1}
\times \Big(SL(r,\mathbb C)\Big/SU(r)\Big)^{r_2}\Big)\times 
(\mathbb R_+^*)^{r_1+r_2}.$$ 
Furthermore, when we really work with $\mathcal O_K$-lattice strucures on $P$, 
i.e., with the space $\bold\Lambda=\bold\Lambda(P)$, from 
the above parametrized space of metric 
structures on $V=\prod_{\sigma\in S_\infty} V_\sigma$, we need 
to further factor out $GL(P)$, i.e., the automorphism  group 
$\mathrm{Aut}_{\mathcal O_K}(\mathcal O_K^{(r-1)}\oplus\frak a)$ of 
$\mathcal O_K^{(r-1)}\oplus\frak a$ as $\mathcal O_K$-modules. So our next aim
is to use $SL$ to understand the quotient space
$$GL(P)\Big\backslash\Bigg(\Big(\Big(SL(r,\mathbb R)\Big/SO(r)\Big)^{r_1}
\times \Big(SL(r,\mathbb C)\Big/SU(r)\Big)^{r_2}\Big)\times 
(\mathbb R_+^*)^{r_1+r_2}\Bigg).$$
  
As such, naturally, now we want  

\noindent
(a) To study the structure of the group 
$\mathrm{Aut}_{\mathcal O_K}(\mathcal O_K^{(r-1)}\oplus\frak a)$ in terms
of $SL$ and units; 
and 

\noindent
(b) To see how this group acts on the space of metrized structures 
$$\Big(\big(SL(r,\mathbb R)/SO(r)\big)^{r_1}
\times \big(SL(r,\mathbb C)/SU(r)\big)^{r_2}\Big)\times 
(\mathbb R_+^*)^{r_1+r_2}.$$
\vskip 0.20cm
View $\mathrm{Aut}_{\mathcal O_K}(\mathcal O_K^{(r-1)}\oplus\frak a)$ 
as a subgroup of $GL(r,K)$. Easily, for
 an element $A=(a_{ij})\in \mathrm{Aut}_{\mathcal O_K}
(\mathcal O_K^{(r-1)}\oplus\frak a)$,  $\det A\in U_K$. 
That is to say, the determinant of an automorphism $A$  has 
to be a unit of $K$. Hence, by the facts 
  that $$A\big(\mathcal O_K^{(r-1)}\oplus\frak a\big)\subset
 \mathcal O_K^{(r-1)}\oplus\frak a, \qquad
  A^{-1}\big(\mathcal O_K^{(r-1)}\oplus\frak a\big)\subset 
\mathcal O_K^{(r-1)}\oplus\frak a,$$ in particular,
   by looking at how the entires $a_{ij}$ of $A$ play 
between $\mathcal O_K$ and $\frak a$, one checks 
   without too much difficulty that 
$$\begin{aligned}~&\mathrm{Aut}_{\mathcal O_K}
(\mathcal O_K^{(r-1)}\oplus\frak a)
=GL(r,\mathcal O_K^{(r-1)}\oplus\frak a)\\
:=&\Bigg\{(a_{ij})\in GL(r,K):\begin{aligned} &a_{rr} \&
a_{ij}\in \mathcal O_K,\\
& a_{ir}\in\frak a,\ a_{rj}\in\frak a^{-1},\end{aligned}  i,j=1,\cdots,r-1;
\qquad \det(a_{ij})\in U_K\Bigg\}.\end{aligned}$$
In other words, $$\mathrm{Aut}_{\mathcal O_K}
(\mathcal O_K^{(r-1)}\oplus\frak a)=\Bigg\{A\in GL(r,K)\cap\begin{pmatrix}
&&&\frak a\\
&\mathcal O_K&&\vdots\\
&&&\frak a\\
\frak a^{-1}&\ldots&\frak a^{-1}&\mathcal O_K\end{pmatrix}:\det A
\in U_K\Bigg\}.$$

To go further, we still need to see how  $\mathrm{Aut}_{\mathcal O_K}
(\mathcal O_K^{(r-1)}\oplus\frak a)$ decomposes
with respect to the shift from $GL$ to $SL$ adopted in the 
discussion on metrized structures. For this purpose, we first introduce 
the subgroup
$\mathrm{Aut}_{\mathcal O_K}^+(\mathcal O_K^{(r-1)}\oplus\frak a)$ 
of $\mathrm{Aut}_{\mathcal O_K}
(\mathcal O_K^{(r-1)}\oplus\frak a)$ consisting of these elements 
whose local determinants at real places are all 
positive. 
Clearly, $\mathrm{diag}(-1,1,\cdots,1)$ 
is an element of $O(r)$, which is supposed to be factored out in our final
discussion. Note also that $$GL(r,\mathbb R)\Big/O(r)\simeq
GL^+(r,\mathbb R)\Big/O^+(r)\qquad\mathrm{and}\quad O^+(r)=SO(r).$$
Consequenly,  we obtain a natural 
identification of quotient spaces between 
$$\mathrm{Aut}_{\mathcal O_K}(\mathcal O_K^{(r-1)}\oplus\frak a)
\Big\backslash\Big(\big(GL(r,\mathbb R)/O(r)\big)^{r_1}\times 
\big(GL(r,\mathbb C)/U(r)\big)^{r_2}\Big)$$ 
and  
$$\mathrm{Aut}_{\mathcal O_K}^+(\mathcal O_K^{(r-1)}\oplus\frak a)
\Big\backslash 
\Big(\big(GL^+(r,\mathbb R)/O^+(r)\big)^{r_1}\times
\big(GL(r,\mathbb C)/U(r)\big)^{r_2}\Big).$$ 

Now we are ready  to shift further to the special linear group $SL$. 
It is here that Dirichlet's Unit Theorem, i.e., 
finiteness of the group of units plays a key role.
As to be expected, the discussion here is a bit involved, for the 
reason that when dealing with metric 
 structures, locally, the genuine realizations are precisely given by 
the following identifications:
$$\begin{aligned}GL(r,\mathbb R)/O(r)\to& GL^+(r,\mathbb R)/SO(r)\\
\to& \Big(SL(r,\mathbb R)
/SO(r)\Big)\times \Big(\mathbb R_+^*
\cdot\mathrm{diag}(1,\cdots,1)\Big)\simeq \Big(SL(r,\mathbb R)/SO(r)\Big)
\times \mathbb R_+^*\end{aligned}$$ via
$$\begin{aligned} [A]&\mapsto[A^+]\\
&\mapsto \Big(\frac{1}{\root{r}\of {\det A^+}}A^+,
\mathrm{diag}\big(\root{r}\of{\det A^+},
 \cdots,\root{r}\of{\det A^+}\big)\Big)\mapsto 
\Big(\frac{1}{\root{r}\of {\det A^+}}A^+,\root{r}\of{\det A^+}\Big),
\end{aligned}$$ 
for real places, and 
$$\begin{aligned}GL(r,\mathbb C)/U(r)\to &\Big(SL(r,\mathbb C)
\times \mathbb C\Big)/\Big(SU(r)\times S^1\Big)\\
\to& 
\Big(SL(r,\mathbb C)/SU(r)\Big)\times \Big(\mathbb R_+^*\cdot
\mathrm{diag}(1,\cdots,1)\Big)\simeq \Big(SL(r,\mathbb C)/SU(r)\Big)
\times \mathbb R_+^*\end{aligned}$$ via
$$\begin{aligned} [A]&\mapsto[A]\\
&\mapsto \Big(\frac{1}{\root{r}\of {\det A}}A,\mathrm{diag}\big(
\root{r}\of{\det A}, \cdots,\root{r}\of{\det A}\big)\Big)\to 
\Big(\frac{1}{\root{r}\of {\det A}}A,\root{r}\of{\det A}\Big),\end{aligned}$$ 
for complex places.
Ideally, we want to have corresponding identifications 
for elements in $\mathrm{Aut}_{\mathcal O_K} 
(\mathcal O_K^{(r-1)}\oplus\frak a)$.   However, this cannot be achieved
in general, due to the fact that,  the $r$-th roots of a unit in $K$ 
lie only in a finite extension of $K$, which usually does not coincide with 
$K$ itself. So suitable modifications have to be made. More precisely,
we go as follows:

Recall that for a unit $\varepsilon \in U_K$,  

\noindent
(a) $\mathrm{diag}(\varepsilon,\cdots,\varepsilon)\in 
\mathrm{Aut}_{\mathcal O_K}(\mathcal O_K^{(r-1)}\oplus\frak a);$ and

\noindent
(b) $\det\mathrm{diag}(\varepsilon,\cdots,\varepsilon)
=\varepsilon^r\in U_F^r:=\{\varepsilon^r:\varepsilon\in U_K\}.$

\noindent
So to begin with, note that to pass from $GL$ to $SL$ over $K$, 
we need to use the intermediate subgroup $GL^+$. Consequently, 
we introduce a subgroup $U_K^+$ of $U_K$ by setting
$$U_K^+:=\{\varepsilon\in U_K:\varepsilon_\sigma>0,
\forall\sigma\ \mathrm{real}\}$$ so as to get
a well-controlled subgroup $U_K^{r,+}:=U_K^+\cap U_K^r$. Indeed, by 
Dirichlet's Unit Theorem, the quotient group 
$U_K^+/(U_K^+\cap U_K^r)$ is  finite. (See the next section for details.)
 
With this said, next we use $U_K^+\cap U_K^r$ to decomposite the automorphism
group  $\mathrm{Aut}_{\mathcal O_K}(\mathcal O_K^{(r-1)}\oplus\frak a).$ 
Thus, choose elements $$u_1,\cdots,
 u_{\mu(r,F)}\in U_K^+$$ such that 
$\Big\{[u_1],\cdots, [u_{\mu(r,F)}]\Big\}$ gives a complete 
representatives of the finite quotient group
$U_K^+\big/\Big(U_K^+\cap U_K^r\Big)$, where $\mu(r,K)$ denotes
 the cardinality of the group $U_K^+/(U_K^+\cap 
U_K^r)$. Set also  $$SL(\mathcal O_K^{(r-1)}\oplus\frak a)
:=SL(r,K)\cap GL(\mathcal O_K^{(r-1)}\oplus\frak a).$$  
\vskip 0.20cm
\noindent
{\bf{\large Lemma}:} {\it There exist elements $A_1, \ldots, A_{\mu(r,K)}$ 
in $GL^+(\mathcal O_K^{(r-1)}\oplus\frak a)$ such that}

\noindent 
(i) $\det A_i=u_i, i=1,\ldots,\mu(r,K)$;

\noindent
(ii) $A_1,\cdots,A_{\mu(r,K)}$ {\it consist of a completed representatives 
of the natural quotient 
$\mathrm{Aut}_{\mathcal O_K}^+(\mathcal O_K^{(r-1)}\oplus\frak a)$ 
by} $SL(\mathcal O_K^{(r-1)}
\oplus\frak a)\times \Big(U_K^{r,+}\cdot\mathrm{diag}(1,\cdots,1)\Big)$.
\vskip 0.30cm
\noindent
That is to say, for  automorphism groups, 

\noindent
(a) $\mathrm{Aut}_{\mathcal O_K}^+(\mathcal O_K^{(r-1)}\oplus\frak a)$
is naturally identified with the disjoint union $$\cup_{i=1}^{\mu(r,K)}
A_i\cdot\Bigg(SL(\mathcal O_K^{(r-1)}\oplus\frak a)\times \Big(U_K^{r,+}
\cdot\mathrm{diag}(1,\cdots,1)\Big)\Bigg);$$ and, consequently,

\noindent 
(b) The  $\mathcal O_K$-lattice 
structures $\bold\Lambda(P)$ on the projective $\mathcal O_K$-module 
$P=\mathcal O_K^{(r-1)}\oplus \frak a$ are parametrized by the disjoint union 
 $$\begin{aligned}\cup_{i=1}^{\mu(r,K)}&A_i\Big\backslash 
\Bigg(
\bigg(SL(\mathcal O_K^{(r-1)}\oplus\frak a)\Big\backslash 
 \Big(\Big(SL(r,\mathbb R)/SO(r)\Big)^{r_1}\times 
\Big(SL(r,\mathbb C)/SU(r)\Big)^{r_2}\Big)
\bigg)\\
&\qquad\times
 \Big(\Big|U_{K}^r\cap U_K^+\Big|\Big\backslash 
\Big(\mathbb R_+^*\Big)^{r_1+r_2}\Big)\Bigg).\end{aligned}$$

\noindent
Proof. With all said above, the proof of this proposition becomes
quite easy now -- It is a direct consequence of the follows:

\noindent
(1) For all $\varepsilon\in U_{K}^+$,
$\mathrm{diag}(\varepsilon,\cdots, \varepsilon)\in 
\mathrm{Aut}_{\mathcal O_K}^+
(\mathcal O_K^{(r-1)}\oplus\frak a)$ and its determinant 
belongs to $U_K^+\cap U_K^r$;
 
\noindent
(2) For $A\in \mathrm{Aut}_{\mathcal O_K}^+
(\mathcal O_K^{(r-1)}\oplus\frak a)$, by definition, 
$\det A\in U_{K}^+$. completed.
\vskip 0.30cm
Therefore, to understand the space of 
$\mathcal O_K$-lattice structures, beyond the spaces
$SL(r,\mathbb R)/SO(r)$ and $SL(r,\mathbb C)/SU(r)$, 
we further need to study 
 
\noindent
(i) the quotient space 
$\big|U_{K}^r\cap U_K^+\big|\Big\backslash 
 \Big(\mathbb R_+^*\Big)^{r_1+r_2}$; and more importantly,

\noindent
(ii) the (modular) space $$SL(\mathcal O_K^{(r-1)}
\oplus\frak a)\Big\backslash \Bigg(\Big(SL(r,\mathbb R)/SO(r)\Big)^{r_1}
\times \Big(SL(r,\mathbb C)/SU(r)\Big)^{r_2}\Bigg).$$ 
\vskip 0.30cm
Now denote by $\widetilde{\mathcal M}_{K,r}(\frak a)$ the moduli space of rank 
$r$ semi-stable $\mathcal O_K$-lattices
with underlying projective module $\mathcal O_K^{(r-1)}\oplus\frak a$. 
For our own convenience, for a set $X$ of (isometry classes of)
lattices, we use the notation $X_{\mathrm{ss}}$ 
to denote the subset of $X$ consisting of lattices
which are semi-stable. 
As such, then we have the following variation of the previous lemma.
\vskip 0.30cm
\noindent
{\bf{\Large Proposition}}. {\it There is a natural identification
between the moduli space $\widetilde{\mathcal M}_{K,r}(\frak a)$ of rank 
$r$ semi-stable $\mathcal O_K$-lattices
on the projective module $\mathcal O_K^{(r-1)}\oplus\frak a$
and the disjoint union of (the ss part of) the quotient spaces}
$$\begin{aligned}\cup_{i=1}^{\mu(r,K)}&A_i\Big\backslash 
\Bigg(
\bigg(SL(\mathcal O_K^{(r-1)}\oplus\frak a)\Big\backslash 
 \Big(\Big(SL(r,\mathbb R)/SO(r)\Big)^{r_1}\times 
\Big(SL(r,\mathbb C)/SU(r)\Big)^{r_2}\Big)
\bigg)_{\mathrm{ss}}\\
&\qquad\times
 \Big(|U_{K}^r\cap U_K^+|\backslash 
(\mathbb R_+^*)^{r_1+r_2}\Big)\Bigg).\end{aligned}$$

\noindent
Proof. By definition and the previous lemma, 
$$\begin{aligned}\widetilde{\mathcal M}_{K,r}(\frak a)\cong&
\bigg[\cup_{i=1}^{\mu(r,K)}A_i\Big\backslash 
\Bigg(
\bigg(SL(\mathcal O_K^{(r-1)}\oplus\frak a)\Big\backslash 
 \Big(\Big(SL(r,\mathbb R)/SO(r)\Big)^{r_1}\times 
\Big(SL(r,\mathbb C)/SU(r)\Big)^{r_2}\Big)
\bigg)\\
&\qquad\times
 \Big(|U_{K}^r\cap U_K^+|\backslash 
(\mathbb R_+^*)^{r_1+r_2}\Big)\Bigg)\bigg]_{\mathrm{ss}}.\end{aligned}$$
Moreover, by definition, we can interchange the subindex ss with the 
disjoint union symbol. With this said, it is sufficient to show that
$$\begin{aligned}&\Bigg[A_i\Big\backslash 
\Bigg(
\bigg(SL(\mathcal O_K^{(r-1)}\oplus\frak a)\Big\backslash 
 \Big(\Big(SL(r,\mathbb R)/SO(r)\Big)^{r_1}\times 
\Big(SL(r,\mathbb C)/SU(r)\Big)^{r_2}\Big)
\bigg)\\
&\hskip 5.0cm\times
 \Big(|U_{K}^r\cap U_K^+|\backslash 
(\mathbb R_+^*)^{r_1+r_2}\Big)\Bigg)\Bigg]_{\mathrm{ss}}\\
=&A_i\Big\backslash 
\Bigg(
\bigg(SL(\mathcal O_K^{(r-1)}\oplus\frak a)\Big\backslash 
 \Big(\Big(SL(r,\mathbb R)/SO(r)\Big)^{r_1}\times 
\Big(SL(r,\mathbb C)/SU(r)\Big)^{r_2}\Big)
\bigg)_{\mathrm{ss}}\\
&\hskip 5.0cm\times
 \Big(|U_{K}^r\cap U_K^+|\backslash 
(\mathbb R_+^*)^{r_1+r_2}\Big)\Bigg).\end{aligned}$$ 
Clearly, an action of an automorphism of a lattice 
does not change the semi-stability. Hence
we need to check whether 
$$\begin{aligned}&\Big[ 
\bigg(SL(\mathcal O_K^{(r-1)}\oplus\frak a)\Big\backslash 
 \Big(\Big(SL(r,\mathbb R)/SO(r)\Big)^{r_1}\times 
\Big(SL(r,\mathbb C)/SU(r)\Big)^{r_2}\Big)
\bigg)\\
&\hskip 5.0cm\times
 \Big(|U_{K}^r\cap U_K^+|\backslash 
(\mathbb R_+^*)^{r_1+r_2}\Big)\Big]_{\mathrm{ss}}\\
=&\bigg[SL(\mathcal O_K^{(r-1)}\oplus\frak a)\Big\backslash 
 \bigg(\Big(SL(r,\mathbb R)/SO(r)\Big)^{r_1}\times 
\Big(SL(r,\mathbb C)/SU(r)\Big)^{r_2}\bigg)
\bigg]_{\mathrm{ss}}\\
&\hskip 5.0cm\times
 \Big(|U_{K}^r\cap U_K^+|\backslash 
(\mathbb R_+^*)^{r_1+r_2}\Big).\end{aligned}$$ 
This is simple since a lattice $\Lambda$ is semi-stable if and only if
its $[T]$-modifications $\Lambda[T]$ are semi-stable for all $T>0$.

\section{Structure of Moduli Space:  Action of $\mathcal O_K$-Units}

To further understand the structure of moduli space of semi-stable 
$\mathcal O_K$-lattices, let us consider the quotient space 
$\big|U_{K}^r\cap U_K^+
\big|\big\backslash \big(\mathbb R_+^*\big)^{r_1+r_2}$. 
 
We start with $U_{K}^r\cap U_K^+$. Clearly,  $U_K^2\subset U_K^+$. 
On the other hand, by  Dirichlet's Unit 
Theorem, up to a finite torsion subgroup consisting of the roots of unity 
in $K$, the image $|U_K|$ of $U_K$ (under the natural logarithm map) is a 
$\mathbb Z$-lattice of rank $r_1+r_2-1$ in $\mathbb R^{r_1+r_2}$. 
As such, the image $|U_K^r|$ of $U_K^r$  corresponds 
simply to the sublattice $r|U_K|$, i.e., the one consists of all elements 
in the lattice $|U_K|$ which are $r$-times of elements in $|U_K|$.
Consequenly, $U_K^+$ as well as $U_K^+\cap U_K^r$ are all finite 
index subgroups of $U_K$. 

Next, let us look at the quotient $\big|U_{K}^r\cap U_K^+
\big|\big\backslash \big(\mathbb R_+^*\big)^{r_1+r_2}$. For this,
we adopt Neukirch's [Neu] presentation.

Let $X$ be a finite $G(\mathbb C|\mathbb R)$-set, 
i.e., a finite set with an involution 
 $\tau\mapsto\bar\tau, \forall\tau\in X$, and let $n=\#X$. 
Consider the $n$-dimensional $\mathbb C$-algebra
  $\bold C:=\prod_{\tau\in X}\mathbb C$ of all tuples 
$z:=(z_\tau)_{\tau\in X}, z_\tau\in \mathbb C$, with 
  componentwise addition and multiplication. Set involutions 
$z\mapsto\bar z\in\bold C$ (resp. $z\mapsto z^*$,
   resp., $z\mapsto\,^*z$) as follows;
for $z=(z_\tau)\in \bold C$, the element $\bar z\in\bold C$ 
(resp. $z^*\in\bold C$, resp.  $\,^*z\in\bold C$)
 is defined to be the element of $\bold C$ having the following components: 
$(\bar z)_\tau=\bar z_{\bar\tau}$, (resp. 
 $z^*_\tau=z_{\bar\tau}$, resp. $\,^*z_\tau=\overline{z_\tau}$). 
Clearly, $\bar z=\,^*z^*$. As such, the invariant subset 
 $\bold R:=[\prod_{\tau\in X}\mathbb C]^+:=\{z\in\bold C:z=\bar z\}$ 
forms an $n$-dimensional commutative
  $\mathbb R$-algebra, and $\bold C=\bold R\otimes
_{\mathbb R}\mathbb C.$
For example, for a number field $K$ of degree $n$ with 
$X=\mathrm{Hom}(K,\mathbb C)$, $\bold R$ coincides with the Minkowski
space $K_\mathbb R:=K\otimes_{\mathbb Q}\mathbb R$.

For the additive, resp. multiplicative group $\bold C$, 
resp. $\bold C^*$, we have the homomorphism
$$\mathrm{Tr}:\bold C\to\mathbb C,\qquad z\mapsto\sum_\tau z_\tau,$$ resp.
 $$N:\bold C^*\to\mathbb C^*,\qquad 
z\mapsto\prod_\tau z_\tau.$$
In other words, $\mathrm{Tr}(z)$ and $N(z)$ is the trace and the determinant 
of the endomorphism $\bold C\to\bold C$ 
defined by $x\mapsto z\cdot x$ respectively. Furthermore,
 we have on $\bold C$ the hermitian scalar product 
$$\langle x,y\rangle:=\sum_\tau x_\tau \overline{y_\tau}=
\mathrm{Tr}(x\cdot\,^*y)$$
which is invariant under conjugation, i.e., 
$\overline{\langle x,y\rangle}=\langle \bar x,\bar y\rangle.$ 
Thus, by restricting it to $\bold R$, we get a scalar product 
$\langle\cdot,\cdot\rangle$, i.e., an Euclidean metric, on the 
$\mathbb R$-vector space $\bold R$.

In $\bold R$, consider the subspace $$\bold R_{\pm}:=
\Big\{x\in\bold R:x=x^*\Big\}=\Big[\prod_\tau\bold R\Big]^+.$$ 
Clearly, for $x=(x_\tau)\in \bold R_\pm$, its components satisfy
$x_{\bar\tau}=x_\tau\in\bold R$. For our convenience, for
  $\delta\in\mathbb R$, we simply write $x>\delta$ to 
signify that $x_\tau>\sigma$ for all $\tau$. With this, then
we introduce the  multiplicative group
$$\bold R_+^*:=\Big\{x\in\bold R_\pm:x>0\Big\}=\Big[\prod_\tau \bold R_+^*
\Big]^+.$$ Clearly, 
$\bold R_+^*$ consists of the tuples $x=(x_\tau)$ of positive real numbers 
$x_\tau$ such that $x_{\bar \tau}=x_\tau$, 
and admits two homomorphisms:
$$|\ |:\bold R^*\to \bold R_+^*,\qquad x=(x_\tau)\mapsto |x|=(|x_\tau|),$$
and
$$\log: \bold R_+^*\to \bold R_\pm,\qquad x=(x_\tau)\mapsto 
\log x=(\log x_\tau).$$ 
For example, when $X=\mathrm{Hom}(K,\mathbb C)$, 
$\bold R_+^*=\mathbb R_{>0}^{r_1+r_2}$ is exactly the 
(unit) factor appeared in our description of the moduli space 
of semi-stable lattices above. Moreover, 
the $G(\mathbb C|\mathbb R)$-set $X=\mathrm{Hom}(K,\mathbb C)$ then
corresponds to the Minkowski space
 $K_\mathbb R=\bold R=[\prod_\tau\mathbb bC]^+,$ in which the field 
$K$ may be naturally embedded. In particular, 
$N((a))=|N_{K/\mathbb Q}(a)|=|N(a)|,$ where 
$N$ denotes the norm on $\bold R^*$.

Now let $\frak p=\{\tau,\bar\tau\}$ be a conjugation
class in $X$. We call $\frak p$ real or complex
 according to $\#\frak p=1$ or 2. Accordingly,
$\bold R_+^*=\prod_\frak p \bold R_{+\frak p}^*$ with
 $\bold R_{+\frak p}^*=\mathbb R_{+}^*$ when $\frak p$ 
is real and $\bold R_{+\frak p}^*=(\mathbb R_+^*\times 
\mathbb R_+^*)^+=\{(y,y):y\in \mathbb R_+^*\}$. Further,
define isomorphisms $\bold R_{+\frak p}^*\simeq 
\mathbb R_{+}^*$ by $y\mapsto y$ resp. $(y,y)\mapsto y^2$
for $\frak p$ real resp. complex, so as to
obtain a natural isomorphism $$\alpha:\bold R_{+}^*
\simeq\prod_\frak p \mathbb R_+^*.$$
 
With this, by $\frac{dy}{y}$ the Haar measure on $\bold R_{+}^*$, 
we mean that one corresponding to the product 
 measure $\prod_\frak p\frac{dt}{t}$, where $\frac{dt}{t}$ 
is the usual Haar measure on $\mathbb R_{+}^*$. 
We call the Haar measure thus defined  the {\it canonical 
measure} on $\bold R_{+}^*.$ Under the logarithm map
  $\log:\mathbb R_{+}^*
 \to \bold R_{\pm}$, it is mapped to the Haar measure $dx$ 
on $\bold R_{\pm}$ which under the isomorphism
$\bold R_{\pm}=\prod_\frak p\mathbb R_{\pm p}\to
\prod_\frak p\mathbb R$ (componentwisely given by
$x_\frak p\mapsto x_\frak p$  resp. $(x_\frak p,x_\frak p)\mapsto 2x_\frak p$
for $\frak p$ real resp. complex) corresponds to the standard
Lebesgue measure on $\prod_\frak p \mathbb R$.

Obviously, for a unit $\varepsilon$ in $U_K$, its $K/\mathbb Q$-norm 
gives $\pm 1$ (in $\mathbb Q$). Hence, the image $|U_K|$
 of the unit group $U_K$ under the map $|\ |:\bold R^*\to \bold R_+^*$ 
is contained in the norm-one
  hypersurface $$\mathbf{S}:=\Big\{x\in \bold R_+^*:N(x)=1\Big\}.$$ 
Write every $y\in \bold R_+^*$ in the form
   $$y=xt^{\frac{1}{n}},\qquad\mathrm{where}\ x=\frac{y}{N(y)^{\frac{1}{n}}},
\ t=N(y).$$
 We then obtain a direct decomposition 
   $\bold R^*_+=\mathbf{S}\times \mathbb R_+^*.$ Let $d^*x$ be the 
unique Haar measure on the mulitplicative
    group $\mathbf{S}$ such that the canonical Haar measure 
$\frac {dy}{y}$ on $\bold R^*_+$ becomes the
product measure $\frac{dy}{y}=d^*x\times\frac{dt}{t}.$
\vskip 0.20cm
The logarithm map $\log$ takes 
 $\mathbf{S}$ to the trace-zero space $$\mathbf{H}:=\Big\{x\in\bold R_\pm: 
\mathrm{Tr}(x)=0\Big\}$$ and the group $|U_K|$ is taken to a
  full ($\mathbb Z$-)latice $G=G_K$ in $\mathbf{H}$ 
(Dirichlet's Unit Theorem). 
We claim that the group $|U_K^+|$ of $U_K^+$ is
   also a full lattice $G^+=G_K^+$ in $\mathbf{H}$. Indeed, it is 
clear that $|U_K^2|\subset|U_K^+|\subset|U_K|$. But
$|U_K^2|=2|U_K|$ is a finite index subgroup of $|U_K|$. 
Thus, $[G:G^+]$ is finite, and
 being a subgroup of $G$, a full rank lattice, of finite index, $G^+$ has to 
be a full rank lattice. Similarly, one sees that
 the image $G_{K,r}^+$ of the group $|U_K^r\cap U_K^+|$ 
is a full rank lattice in $\mathbf{H}$ as well.

Choose now $F_{K,r}^+$ to be the preimage of an arbitrary fundamental
parallelopiped $\frak D_{K,r}^+$ of the lattice $G_{K,r}^+$ in $\mathbf{H}$, 
then the fundamental domain $F_{K,r}^+$ cuts up the 
norm-one hypersurface $\mathbf{S}$ into the disjoint union 
 $$\mathbf{S}=\cup_{\eta\in U_F^+}\eta^r F_{r,K}^+.$$

\noindent
{\bf{\large Lemma}}. {\it The fundamental domain $F_{r,K}^+$ of
 $U_K^r\cap U_K^+$ in $\mathbf{S}$ has the following volume with 
 respect to $d^*x$; $$\mathrm{Vol}(F_{r,K}^+)=
r^{r_1+r_2-1}R^+_K$$ where $R_K^+$ is the narrow regulator of $K$.}
 
\noindent
Proof. Since $I:=\{t\in \mathbb R_+^*:1\leq t\leq e\}$ 
has measure 1 with respect to $\frac{dt}{t}$, the 
 quantity $\mathrm{Vol}(F_{r,K}^+)$ is also the volume of
 $F_{r,K}^+\times I$ with respect to $d^*x\times
  \frac{dt}{t}$, i.e., the volume of $\alpha(F_{r,K}^+\times I)$ 
with respect to $\frac{dy}{y}$. The 
  composition $\psi$ of the isomorphisms
$$\bold R_+^*\buildrel\log\over\to \bold R_\pm\buildrel\phi\over
  \to\prod_{\frak p|\infty}\mathbb R=\mathbb R^{r_1+r_2}$$ 
transforms $\frac{dy}{y}$ into the Lebesgue 
  measure of $\mathbb R^{r_1+r_2}$, $$\mathrm{Vol}(F_{r,K}^+)
=\mathrm{Vol}_{\mathbb R^{r_1+r_2}}
  \Big((\psi\circ\alpha)(F_{r,K}^+\times I)\Big).$$ Let us 
compute the image $(\psi\circ\alpha)(F_{r,K}^+
  \times I)$. Let $\bold 1:=(1,1,\cdots,1)\in \mathbf{S}$. Then 
we find $(\psi\circ\alpha)(F_{r,K}^+\times I)=
  \frak e\cdot\log t^{1/n}=\frac{1}{n}\frak e\log t$ with the 
vector $\frak e=(e_{\frak p_1},\cdots,
   e_{\frak p_{r_1+r_2}})\in \mathbb R^{r_1+r_2}, e_{\frak p_i}=1$ 
or 2 depending whether $\frak p_i$ is 
   real or complex. By definition of $F_{r,K}^+$, we also have 
$(\psi\circ\alpha)(F_{r,K}^+\times \{1\})=
r\Phi_K^+$ where $\Phi_K^+$ denotes the fundamental parallelopiped of 
the totally positive unit lattice $G_K^+$ 
in the trace-zero space $\mathbf{H}$. 
This gives $$(\psi\circ\alpha)(F_{r,K}^+\times I)=
r\Phi_K^+\times [0,\frac{1}{n}]\frak e,$$ the parallelopiped 
spanned by the vectors $\frak e_1,\cdots, 
\frak e_{r_1+r_2-1},\frac {1}{n}\frak e$ if $\frak e_1,\cdots, 
\frak e_{r_1+r_2-1}$ span the fundamental 
domain $\Phi_K^+$. Its volume is $\frac{1}{n}r^{r_1+r_2-1}$ 
times the absolute value of the determinant 
$$\det\left(\begin{matrix}
\frak e_{1,1}&\cdots&\frak e_{r_1+r_2-1,1}&\frak e_{\frak p_1}\\
\cdots&\cdots&\cdots&\cdots\\
\frak e_{1,r_1+r_2}&\cdots&\frak e_{r_1+r_2-1,r_1+r_2}&
\frak e_{\frak p_{r_1+r_2-1}}\end{matrix}\right).$$
Adding the first $r_1+r_2-1$ lines to the last one, all entries 
of the last line becomes 0 except the last one,
 which is $\sum_i \frak e_{\frak p_i}=n$. By defintion,
the absolute value of the determinant of the matrix above 
these zeros is equal to the narrow regulator 
$R_K^+$. Thus we get 
 $\mathrm{Vol}(F_{r,K}^+)=r^{r_1+r_2-1}R^+_K$. This completes the proof.
\vskip 0.20cm
In summary, for the (unit) factor $|U_{K}^r\cap U_K^+
|\big\backslash \big(\mathbb R_+^*\big)^{r_1+r_2}$, 
its structure may be understood via a natural decomposition
$\Big(\big(U_{K}^r\cap U_K^+
\big)\big\backslash \bold S\Big)\times \mathbb R_+^*,$
where $\bold S$ denotes the norm-one
  hypersurface $\mathbf{S}:=\{x\in \bold R_+^*:N(x)=1\};$
together with a  disjoint union 
 $\mathbf{S}=\cup_{\eta\in U_F^+}\eta^r F_{r,K}^+,$
where  $F_{r,K}^+$ denotes a \lq fundamental parallogram' of
 $U_K^r\cap U_K^+$ in $\mathbf{S}$ with
$r^{r_1+r_2-1}R^+_K$ as its volume.

\section{Non-Abelian Zeta Functions for Number Fields}

Let $K$ be an algebraic number field (of finite degree $n$) with 
$\Delta_K$ the absolute value of its discriminant. Denote by 
$\mathcal O_K$ the ring of integers as usual. For a 
fixed positive integer $r\in\mathbb N$, denote by
${\mathcal M}_{F,r}$ the moduli space of semi-stable 
$\mathcal O_K$-lattices of rank $r$. Denote by $d\mu$ 
the natural associated (Tamagawa type) measure
(induced from that on $GL$). For each 
$\Lambda\in{\mathcal M}_{K,r}$, following van der Geer 
and Schoof, (see e.g. [GS], [Bor] and [We1,2],)
define the associated {\it 0-th geo-arithmetical cohomology} 
$h^0(K,\Lambda)$ by $$h^0(K,\Lambda):=
\log\Big(\sum_{x\in \Lambda}\exp\Big(-\pi
\sum_{\sigma:\mathbb R}\|x_\sigma\|_{\rho_\sigma}^2-2\pi
\sum_{\sigma:\mathbb C}\|x_\sigma\|_{\rho_\sigma}^2\Big)\Big)$$ 
where  $x=(x_\sigma)_{\sigma\in S_\infty}$ and 
$(\rho_\sigma)_{\sigma\in S_\infty}$ denote the 
$\sigma$-component of the metric $\rho=\rho_\Lambda$ determinet by 
the lattice $\Lambda$ with $S_\infty$  a collection of inequivalent 
Archimedean places of $K$.

Following [We1,2], we introduce the following
\vskip 0.20cm
\noindent
{\bf \Large Definition.} Define the {\it rank $r$ (non-abelian) zeta function 
$\xi_{K,r}(s)$ of a number field $K$} to 
be the integration $$\xi_{K,r}(s):=\int_{\Lambda\in
{\mathcal M}_{K,r}}\Big(e^{h^0(K,\Lambda)-1}\Big)\cdot
 \Big(e^{-s}\Big)^{-\log\mathrm{Vol}(\Lambda)}\,d\mu(\Lambda),\qquad\Re(s)>1.$$

By the Arakelov-Riemann-Roch Formula, one can write the 
non-abelian zeta function in the following form 
which fits more for practical purpose
$$\xi_{K,r}(s):=\Big(\Delta_K^{\frac{r}{2}}\Big)^s
\cdot \int_{\Lambda\in{\mathcal M}_{K,r}}
\Big(e^{h^0(K,\Lambda)-1}\Big)\cdot \Big(e^{-s}\Big)^{\mathrm{deg}
(\Lambda)}\,d\mu(\Lambda),\qquad\Re(s)>1.$$
 
In [We1,2], we, for an $\mathcal O_K$-lattice $\Lambda$, construct 
two geo-arithmetical cohomology groups $$H^0(K,\Lambda):=\Lambda,
\qquad\mathrm{and} \qquad H^1(K,\Lambda):=V(\Lambda)/\Lambda$$ 
where $V(\Lambda)$ denotes $V:=\prod_{\sigma\in S_\infty}V_\sigma$ 
for $V_\sigma:=\Lambda\otimes_{\mathcal O_K}K_\sigma$ equipped with the 
canonical measures. In such a way, both 
$H^0(K,\Lambda)$ and $H^1(K,\Lambda)$ are topological 
groups. More precisely, $H^0$ is discrete, while $H^1$ 
is compact. As a direct consequence, then the corresponding geo-arithmetical 
counts for these locally compact groups can be done
 by using Fourier analysis on them so as to naturally 
get not only the above $h^0$ but also a new $h^1$
in a very natural way for lattices. Moreover,  fundamental results 
corresponding to the Serre duality and Riemann-Roch Theorem
hold for these newly defined $h^i, i=0,1$ as well. 
To state them more clearly, as usual, introduce 
the {\it dualizing lattice} $\mathcal K_K$ of $K$ as the
 dual of the so-called different lattice $\overline{\frak D_K}$ of 
$K$. (Here by the different lattice $\overline{\frak D_K}$, we
  mean the rank one $\mathcal O_K$-lattice whose underlying 
module is given by the different $\frak d_K$ of $K$ and
   whose metric is induced from the canonical one via the 
natural embedding $\frak d_K\hookrightarrow
    K_{\mathbb R}$, the Minkowski space.) Also as usual, denote the 
(Arakelov) dual lattice of $\Lambda$ by $\Lambda^\vee$. 
Then we have the following
\vskip 0.20cm
\noindent
(1) ({\bf \large{Serre Duality=Pontragin Duality}})

\noindent 
(a) ({\bf Topologically}) 
$$\widehat{H^1(K,\Lambda)}\cong 
H^0(K,\mathcal K_K\otimes \Lambda^{\vee}),$$ where $\,\widehat ~\,$ 
denotes the Pontragin dual of a topological group;

\noindent
(b) ({\bf Analytically}) $$h^1(K,\Lambda)=h^0(K,\mathcal K_K\otimes \Lambda^{\vee});$$

\noindent
(2) ({\bf \large{Riemann-Roch Theorem}})
$$h^0(K,\Lambda)-h^1(K,\Lambda)=:\chi(K,\Lambda)=\mathrm{deg}(\Lambda)
-\frac {r}{2}\log\Delta_K.$$

\noindent
{\bf Remarks.} (1) While $H^0$ and $H^1$ together with $h^0$ and 
$h^1$ are quite similar to those for function
 fields via an adelic approach (see e.g., [Iw2], [Se] or [W]), 
two major differences should be noticed.

\noindent
(a) For number fields, $H^0$ is discrete and $H^1$ is compact, 
while for function fields, both $H^0$ and 
$H^1$ are linearly compact, i.e., are finite dimensional vector 
spaces over the base field;

\noindent
(b) For number fields, $h^i$ are defined using Fourier analysis,
 say, a weight of Gauss distribution is
attached to each element of $H^0$ in defining $h^0$.
But for function fields,
 $h^i$ are defined using a much simpler count. 
Say, when the base fields are finite, the counts are 
carried out by a direct counting process, i.e., 
every element in $H^i$ is counted with the naive weight 1.

\noindent
(2) It is remarkable to see that the analogue of 
Serre Duality has a certain topological counterpart via Pontragin Duality
for topological groups 
and an analytic counterpart via the Plancherel Formula, a special
kind of Fourier Inversion Formula.

\noindent
(3) The Riemann-Roch Theorem is a direct consequence of 
the Serre Duality and the Poisson Summation Formula. So
the above constructions and results are almost in Tate's 
Thesis, but not quite yet there.

\noindent
(4) A two dimensional analogue of such a theory seems to be 
very much in demanding -- Such a two dimensional 
theory is closely related with the Riemann Hypothesis via an 
intersection approach proposed in [We1].

\noindent
(5) The reader may learn how to appreciate the treatment here for $H^i$'s 
and $h^i$'s by consulting Weil's {\it Basic Number
 Theory} and Neukirch's {\it Algebraic Number Theory}. For the first 
one, mainly due to the lake of the construction 
 above, Weil, unlike in the rest of his book, treated  
zeta functions for number fields separately from that for function
fields, while for the second, Neukirch introduced
a different type of $h^i$ for which no duality is satisfied. 
\vskip 0.30cm
With all this well-prepared cohomology theory, standard yet
fundamental properties for non-abelian zeta functions can be easily deduced.
It works exactly as that for Artin zeta functions for curves over finite 
fields, as done by H. L. Schmid. Indeed, it is now a standard 
procedure to deduce the meromorphic continuation from the
 Riemann-Roch, to establish the functional equation from the Serre Duality 
and to locate the singularities from both Riemann-Roch and Serre Duality.
(For details, please see Moreno [Mo] and/or Weil [W] and/or [We1,2].)
That is to say, we have the following
\vskip 0.20cm
\noindent
{\bf{\Large Facts}.} (I) ({\bf Meromorphic Continuation})
{\it The rank $r$ non-abelian zeta function 
$\xi_{K,r}(s)$  is well-defined when $\Re(s)>1$ and admits a
 meromorphic continuation, denoted also by 
$\xi_{K,r}(s)$, to the whole complex $s$-plane;}

\noindent
(II) ({\bf Functional Equation}) $\xi_{K,r}(1-s)=\xi_{K,r}(s)$;

\noindent
(III) ({\bf Singularities \& Residues}) {\it $\xi_{K,r}(s)$ has only two
 singularities, all are simple poles, at $s=0,1$,
 with the same residues $\mathrm{Vol}\Big({\mathcal M}_{K,r}
\big([\Delta_K^{\frac{r}{2}}]\big)\Big)$,
where ${\mathcal M}_{K,r}\big([\Delta_K^{\frac{r}{2}}]\big)$ denotes 
the moduli space of rank $r$ semi-stable 
$\mathcal O_K$-lattices whose  volumes are fixed to be
$\Delta_K^{\frac{r}{2}}$.}

\noindent
{\bf Remarks}. (1) Due to the fact that the volumes of lattices are fixed, 
the semi-stable condition implies that the first Minkowski
successive minimums of the lattices involved
admit a natural lower bound away 
from 0 (depending only on $r$). Hence by the standard 
reduction theory, see e.g., Borel [Bo1,2], 
${\mathcal M}_{K,r}\big([\Delta_K^{\frac{r}{2}}]\big)$ is compact. 
Consequently, the volume  
$\mathrm{Vol}({\mathcal M}_{K,r}([\Delta_K^{\frac{r}{2}}]))$ appeared 
 above does make sense.

\noindent
(2) The Tamagawa type of volume $\mathrm{Vol}\Big({\mathcal M}_{K,r}
\big([\Delta_K^{\frac{r}{2}}]\big)\Big)$ is a new intrinsic 
 non-abelian invariant for the number field $K$.

\section{Non-Abelian Zeta Functions and Epstein Zeta Functions}

Recall that we can choose integral $\mathcal O_K$-ideals 
$\frak a_1=\mathcal O_K,\frak a_2,\cdots,
\frak a_h$ such that the ideal class group $CL(K)$ is given by
$\Big\{[\frak a_1],\cdots,[\frak a_h]\Big\},$
and that any rank $r$ projective $\mathcal O_K$-module $P$ is 
isomorphic to $P_{\frak a_i}$ for a certain $i, 1\leq i\leq h$.
Here, $P_\frak a:=P_{r,\frak a}:=\mathcal O_K^{(r-1)}\oplus\frak a$ 
for a fractional $\mathcal O_K$-ideal $\frak a$.  (Quite often,  
we use $\frak a$ as a running symbol for  
$\frak a_1,\frak a_2,\cdots,\frak a_h$.) Consequently,
 $${\mathcal M}_{K,r}=\cup_{i=1}^h\widetilde{\mathcal M}_{K,r}(\frak a_i)$$ 
with $\widetilde{{\mathcal M}}_{K,r}(\frak a_i)=:
\Big(\widetilde{\bold\Lambda}(P_{\frak a_i})\Big)_{\mathrm{ss}},$ the part of 
$\widetilde{\bold\Lambda}(P_{\frak a_i})$ consisting of only semi-stable
$\mathcal O_K$-lattices.

As such, introduce an intermidiate partial non-abelian 
zeta function $\widetilde\xi_{K,r;\frak a}(s)$ by setting
$$\widetilde\xi_{K,r;\frak a}(s):=\int_{\widetilde{\mathcal M}_{K,r}(\frak a)}
\big(e^{h^0(K,\Lambda)}-1\big)\cdot
 (e^{-s})^{-\log\mathrm {Vol}(\Lambda)}d\mu(\Lambda),\qquad\Re(s)>1.$$

\noindent
{\bf Remark}. It is interesting to see  functional 
equations among $\widetilde\xi_{K,r;\frak a}(s)$'s since 
the dual lattices for the lattices involved have the 
underlying projective module $\frak d_K^{(r-1)}\oplus
(\frak d_K\cdot \frak a^{-1})$, while our $\frak a$ which 
is one of the $\frak a_i$ has already been fixed (say, to be 
intergal). So it appears that it is better to leave such
 a matter untouched. However, it is not really that
 bad, as one can check easily that for any two fractional 
$\mathcal O_K$-ideals $\frak a,\frak b$, if 
 $[\frak a]=[\frak b]$ as ideal classes, then
$\widetilde\xi_{K,r;\frak a}(s)=\widetilde\xi_{K,r;\frak b}(s).$ 
Therefore we have indeed 
$$\widetilde\xi_{K,r;\frak a}(1-s)=
\widetilde\xi_{K,r;\frak d^{r}\frak a^{-1}}(s).$$ We leave the details to the
 reader for the reasons that 
only after checking this, he or she will get things right to carry on.

As such, we, after using the Proposition in \S6 and 
the notation there, to get
$$\widetilde\xi_{K,r;\frak a}(s)=\sum_{j=1}^{\mu(r,K)}
\xi_{K,r;\frak a;A_j}(s)$$ where
$$\xi_{K,r;\frak a;A_j}(s):=\int_{\Lambda\in\mathcal 
M_{K,r;A_j}(\frak a)}\Big(e^{h^0(K,\Lambda)}-1\Big)\cdot
 \Big(e^{-s}\Big)^{-\log\mathrm {Vol}(\Lambda)}d\mu(\Lambda),\qquad \Re(s)>1$$
with ${\mathcal M}_{K,r;A_j}(\frak a)$ the component of the moduli 
space of semi-stable $\mathcal O_K$-lattices whose
 points corresponding to these in
$$\begin{aligned}&\Big[A_i\Big\backslash 
\Bigg(
\bigg(SL(\mathcal O_K^{(r-1)}\oplus\frak a)\Big\backslash 
 \Big(\Big(SL(r,\mathbb R)/SO(r)\Big)^{r_1}\times 
\Big(SL(r,\mathbb C)/SU(r)\Big)^{r_2}\Big)
\bigg)\\
&\hskip 5.0cm\times
 \Big(|U_{K}^r\cap U_K^+|\backslash 
(\mathbb R_+^*)^{r_1+r_2}\Big)\Bigg)\Big]_{\mathrm{ss}}\\
=&A_i\Big\backslash 
\Bigg(
\bigg(SL(\mathcal O_K^{(r-1)}\oplus\frak a)\Big\backslash 
 \Big(\Big(SL(r,\mathbb R)/SO(r)\Big)^{r_1}\times 
\Big(SL(r,\mathbb C)/SU(r)\Big)^{r_2}\Big)
\bigg)_{\mathrm{ss}}\\
&\hskip 5.0cm\times
 \Big(|U_{K}^r\cap U_K^+|\backslash 
(\mathbb R_+^*)^{r_1+r_2}\Big)\Bigg),\end{aligned}$$
 under the 
natural identification set in \S6.
Moreover, since $A_i$ is simply an automorphism,  its action 
does not change the total volumes  as well as
 the $h^0$ of the lattices. Therefore, if we introduce further 
the (genuine) partial non-abelian zeta function 
$\xi_{K,r;\frak a}(s)$ by setting
$$\xi_{K,r;\frak a}(s):=\int_{{\mathcal M}_{K,r}(\frak a)}
\Big(e^{h^0(K,\Lambda)}-1\Big)\cdot
 \Big(e^{-s}\Big)^{-\log\mathrm {Vol}(\Lambda)}d\mu(\Lambda),\qquad \Re(s)>1$$ 
where ${\mathcal M}_{K,r}(\frak a)$ denotes the part of the moduli 
space of semi-stable $\mathcal O_K$-lattices
 whose points corresponding to these in
$$\begin{aligned}&\Bigg(
\bigg(SL(\mathcal O_K^{(r-1)}\oplus\frak a)\Big\backslash 
 \Big(\Big(SL(r,\mathbb R)/SO(r)\Big)^{r_1}\times 
\Big(SL(r,\mathbb C)/SU(r)\Big)^{r_2}\Big)
\bigg)_{\mathrm{ss}}\\
&\hskip 5.0cm\times
 \Big(|U_{K}^r\cap U_K^+|\backslash 
(\mathbb R_+^*)^{r_1+r_2}\Big)\Bigg).\end{aligned}$$ 
This then completes the proof of  the following

\noindent
{\bf \Large{Proposition}}.  {\it With the same notation as above,
$$\xi_{K,r;\frak a;A_j}(s)=\xi_{K,r;\frak a}(s),\qquad\forall j=1,\cdots,\mu(r,K).$$
In particular, $$\xi_{K,r}(s)=\mu(r,K)\cdot \sum_{i=1}^h
\xi_{K,r;\frak a_i}(s).$$}

This been said, to further understand the structure of non-abelian zeta 
function $\xi_{K,r}(s)$, we next investigate how the integrand 
$$\Big(e^{h^0(K,\Lambda)}-1\Big)\cdot
 (e^{-s})^{-\log\mathrm {Vol}(\Lambda)}d\mu(\Lambda)$$
behaves over the space $$\begin{aligned}&\Bigg(
\bigg(SL(\mathcal O_K^{(r-1)}\oplus\frak a)\Big\backslash 
 \Big(\Big(SL(r,\mathbb R)/SO(r)\Big)^{r_1}\times 
\Big(SL(r,\mathbb C)/SU(r)\Big)^{r_2}\Big)
\bigg)_{\mathrm{ss}}\\
&\hskip 5.0cm\times
 \Big(|U_{K}^r\cap U_K^+|\backslash 
(\mathbb R_+^*)^{r_1+r_2}\Big)\Bigg).\end{aligned}$$
  
By definition, $$e^{h^0(K,\Lambda)}-1=
\sum_{x\in\Lambda\backslash\{0\}}\exp\Big(-\pi
  \sum_{\sigma:\mathbb R}\|x_\sigma\|_{\rho_\sigma}-
2\pi\sum_{\sigma:\mathbb C}
  \|x_\sigma\|_{\rho_\sigma}\Big).$$ Thus, in terms 
of the embedding $$z\in\Lambda=\mathcal O_K^{(r-1)}
  \oplus\frak a\hookrightarrow K^{(r)}\hookrightarrow 
\big(\mathbb R^{r_1}\times\mathbb C^{r_2}\big)^r
  \simeq (\mathbb R^r)^{r_1}\times(\mathbb C^r)^{r_2},$$ 
$z$ maps to the corresponding point $(z_\sigma)$ and 
$\|z_\sigma\|_{\rho_\sigma}=\|g_\sigma z_\sigma\|$, 
where  the metric $\rho_\sigma$ is defined by
$g_\sigma\cdot g_\sigma^t$ for certain 
$g_\sigma\in GL(r, \mathbb R)$ when $\sigma$ is real, and
by $g_\sigma\cdot \bar g_\sigma^t$ for certain 
$g_\sigma\in GL(r, \mathbb C)$ when $\sigma$ is complex.

Recall that $\|g_\sigma z_\sigma\|$ is $O(r)$ resp. 
$U(r)$ invariant when $\sigma$ is real resp. 
complex. Similarly, $\mathrm{Vol}(\Lambda)$ is invariant.
Consequently, $\Big(e^{h^0(K,\Lambda)}-1\Big)\cdot
 (e^{-s})^{-\log\mathrm {Vol}(\Lambda)}$ is well-defined over
 $$\Big(GL(r,\mathbb R)/O(r)\Big)^{r_1}\times
\Big(GL(r,\mathbb C)/U(r)\Big)^{r_2}.$$

To go further, we next study how  $\Big(e^{h^0(K,\Lambda)}-
1\Big)\cdot (e^{-s})^{-\log\mathrm {Vol}(\Lambda)}$ changes 
when we apply the operation $\Lambda\mapsto\Lambda[t]$ for $t>0$. 
Clearly, in terms of each local component, 
  $\rho_\sigma\mapsto t_\sigma\rho_\sigma$ with 
$t_\sigma\in \mathbb R_+^*$, we have $\|x_\sigma\|_{t_\sigma\rho_\sigma}^2=t_\sigma^2
  \cdot\|x_\sigma\|_{\rho_\sigma}^2$. Hence 
$\Big(e^{h^0(K,\Lambda[t])}-1\Big)$ changes to
  $$\sum_{x\in\Lambda\backslash\{0\}}\exp\Big
(-\pi\sum_{\sigma:\mathbb R}\|x_\sigma\|_{\rho_\sigma}\cdot 
  t_\sigma^{\frac{r}{2}}-2\pi\sum_{\sigma:\mathbb C}
\|x_\sigma\|_{\rho_\sigma}\cdot t_\sigma^{\frac{r}{2}}
  \Big),$$
  while $\mathrm{Vol}(\Lambda[t])$ decomposes to 
$\mathrm{Vol}(\Lambda)\cdot\prod_{\sigma\in S_\infty}
  t_\sigma^r$ for $t=(t_\sigma)$. On the other hand, 
by changing the volume in such a way, $d\mu(\Lambda)$ becomes 
$\prod_{\sigma\in S_\infty} \frac{dt_\sigma}{t_\sigma}\cdot d\mu_1(\Lambda_1)$,  where $d\mu_1(\Lambda_1)$ denotes 
the corresponding volume form on the space of   semi-stable lattices
  corresponding to the points in 
$${\mathcal M}_{F,r;\frak a}\Big[N(\frak a)\cdot\Delta_K^{\frac{r}{2}}\Big]:=
\bigg(SL(\mathcal O_K^{(r-1)}\oplus\frak a)\Big\backslash 
  \Big(\Big(SL(r,\mathbb R)/SO(r)\Big)^{r_1}\times
 \Big(SL(r,\mathbb C)/SU(r)\Big)^{r_2}\Big)\bigg)_{\mathrm{ss}},$$
 due to the fact that $$\mathrm{Vol}\Big(\overline{\mathcal O_K^{(r-1)}
\oplus\frak a}\Big)=\Delta_K^{\frac{r-1}{2}}\cdot
 \big(N(\frak a)\cdot \Delta_K^{\frac{1}{2}}\big)=N(\frak a)\cdot
\Delta_K^{\frac{r}{2}}.$$ 
(As we are going to identify the moduli space of lattices with its realization
in terms of $SL$,  from now on we make no distinction between them.) 
Moreover, note that the $\mathcal O_K$-units have their 
  (total rational) norm 1, hence $\mathcal O_K$-units 
do not really change the total
 volume of the lattice. All in all, then we get for $\Re(s)>1$,
 $$\begin{aligned}~&\xi_{F,r;\frak a}(s)\\
=&\big(N(\frak a)\cdot\Delta_K^{\frac {r}{2}}\big)^s
\cdot\int_{{\mathbb R}_{>0}^{r_1+r_2}}
 t_\sigma^{rs}\prod_{\sigma\in S_\infty}\frac{dt_\sigma}{t_\sigma}\\
&\times\int_{\Lambda\in {\mathcal M}_{F,r;\frak a}\Big[N(\frak a)\cdot\Delta_K^{\frac{r}{2}}\Big]}
\sum_{x\in (\Lambda\backslash\{0\})/U_{r,F}^+}\exp\Big(-\pi\sum_{\sigma:\mathbb R}\|x_\sigma\|_{\rho_\sigma}
\cdot t_\sigma^{\frac{r}{2}}-2\pi\sum_{\sigma:
 \mathbb C}\|x_\sigma\|_{\rho_\sigma}\cdot 
t_\sigma^{\frac{r}{2}}\Big)\,d\mu_1(\Lambda).\end{aligned}$$
 
Therefore, by applying the Mellin transform and using the formula
 $$\int_0^\infty e^{-At^B}t^s\frac{dt}{t}=
\frac{1}{B}\cdot A^{-\frac{s}{B}}\cdot\Gamma(\frac{s}{B})$$ (whenever 
 both sides make sense,) we obtain that
  $$\begin{aligned}~&\xi_{F,r;\frak a}(s)\\
=&\big(N(\frak a)\cdot 
\Delta_K^{\frac {r}{2}}\big)^s\cdot
\int_{\Lambda\in {\mathcal M}_{F,r;\frak a}\Big[N(\frak a)\cdot\Delta_K^{\frac{r}{2}}\Big]}
\sum_{x\in (\Lambda\backslash\{0\})/U_{r,F}^+}\\
&\qquad \Bigg(
\prod_{\sigma:\mathbb R}\bigg(\frac{r}{2}\cdot 
\Big(\pi\|x_\sigma\|_{\rho_\sigma}\Big)^{-\frac{rs}{2}}
 \Gamma(\frac{rs}{2})\bigg)\cdot 
\prod_{\sigma:\mathbb C}\bigg(
\frac{r}{2}\cdot \Big(2\pi\|x_\sigma\|_{\rho_\sigma}
 \Big)^{-\frac{rs}{2}}\Gamma(rs)\bigg)\Bigg)\,d\mu_1(\Lambda)\\
=&\Big(\frac{r}{2}\Big)^{r_1+r_2}\cdot \Big(\pi^{-\frac{rs}{2}}
\Gamma(\frac{rs}{2})\Big)^{r_1}\cdot
 \Big(\big(2\pi\big)^{-rs}\Gamma(rs)\Big)^{r_2}\\
&\qquad\times
 \Big(N(\frak a)\cdot \Delta_K^{\frac{r}{2}}\Big)^s\cdot
\int_{\Lambda\in {\mathcal M}_{F,r;\frak a}\Big[N(\frak a)
\cdot\Delta_K^{\frac{r}{2}}\Big]
}\Big(\sum_{x\in (\Lambda\backslash\{0\})/
 U_{r,F}^+}\frac{1}{\|x\|_\Lambda^{rs}}\Big)
\,d\mu(\Lambda),\qquad \Re(s)>1\end{aligned}$$ (Here, in the last step, 
we also change the notation from $d\mu_1(\Lambda)$ to $d\mu(\Lambda)$
 for our own convenience.)

Accordingly, for $\Re(s)>1$, define the {\it completed Epstein zeta function} 
$\hat E_{K,r;\frak a}(s)$ by
 $$\hat E_{K,r;\frak a}(s):=\Big(\pi^{-\frac{rs}{2}}
\Gamma(\frac{rs}{2})\Big)^{r_1}\cdot
 \Big(\big(2\pi\big)^{-rs}\Gamma(rs)\Big)^{r_2}\cdot\bigg[
 \Big(N(\frak a)\cdot \Delta_K^{\frac{r}{2}}\Big)^s\cdot
\sum_{x\in (\Lambda\backslash\{0\})/
 U_{r,F}^+}\frac{1}{\|x\|_\Lambda^{rs}}\bigg].$$
All in all, what we have just said exposes the following
\vskip 0.30cm
\noindent
{\bf{\Large Facts.}}  (IV) ({\bf Decomposition})
{\it The rank $r$ non-abelian zeta funtion of $K$
admits a natural decomposition} 
$$\xi_{K,r}(s)=\mu(r,K)\cdot \sum_{i=1}^h
\xi_{K,r;\frak a_i}(s);$$

\noindent
(V) ({\bf Non-Abelian Zeta = Integration of Epstein Zeta})
{\it The partial rank $r$ non-abelian zeta function
 $\xi_{F,r;\frak a}(s)$ of $K$ associated to $\frak a$
is given by an integration of a completed Epstein
type zeta function:} 
$$\xi_{F,r;\frak a}(s)=\big(\frac{r}{2}\big)^{r_1+r_2}
\cdot\int_{
{\mathcal M}_{F,r;\frak a}[N(\frak a)\cdot\Delta_K^{\frac{r}{2}}]}
\hat E_{K,r;\frak a}(s)\,d\mu,\qquad\Re(s)>1.$$
\vskip 0.30cm
{\it Remark.} The relation between non-abelian zeta and Epstein zeta
was first established for $\mathbb Q$. (See my paper on
\lq Analytic truncation and Rankin-Selberg versus algebraic 
truncation and non-abelian zeta',
{\it Algebraic Number Theory and Related Topics}, RIMS Kokyuroku, 
No.1324 (2003).)  Consequently, in [We2,3],
we develop a general theory of non-abelian $L$-functions for global fields, 
using Langlands' theory of Eisenstein series. 
\vskip 0.45cm
\noindent
{{\Large Appendix:} {\Large Higher Dimensional Gamma Function}}
\vskip 0.30cm
For the reader who wants to know why $\Gamma$-factor 
appears in such a way in the above discussion, we now follow  Neukirch
to explain it in more professional way following [Neu].

For two tuples $z=(z_\tau),\ p=(p_\tau)\in \mathbb C$, 
we define the power $z^p:=(z_\tau^{p_\tau})\in \mathbb C$  
 by using $z_\tau^{p_\tau}:=e^{p_\tau\log z_\tau}$. Here, to
make it well-defined, we choose the principal branch for the logarithm 
and assume that $z_\tau$'s move only in the place cut 
along the negative real axis. Then, for 
$\mathbf{s}=(s_\tau)\in\bold C$ such that $\Re(s_\tau)>0$, 
we introduce the Gamma function associated to the 
$G(\mathbb C|\mathbb R)$-set $X$ by $$\Gamma_X(\mathbf{s}):=
\int_{\bold R_+^*}N(e^{-y}y^\mathbf{s})\frac{dy}{y}.$$
 This integral is then well-defined as well, according to our convention 
above. Indeed, the convergence of the integral can be
  reduced to the one for ordinary Gamma function as follows.

\noindent
{\bf{\large Lemma}.} {\it According to the decomposition of the 
$G(\mathbb C|\mathbb R)$-set $X$ into its conjugation classes $\frak p$,  
$$\Gamma_X(\mathbf{s})=\prod_\frak p\Gamma_\frak p(\mathbf{s}_\frak p),$$ 
where for $\frak p=\{\tau\}$, $\mathbf{s}_\frak p=s_\tau$ and the local factor
is simply $\Gamma(\mathbf{s}_\frak p)$, while for  $\frak p=\{\tau,\bar\tau\}$,
$\tau\not=\bar\tau$, $\mathbf{s}_\frak p=(s_\tau,s_{\bar\tau})$, and the 
local factor becomes
$2^{1-\mathrm{Tr}(\mathbf{s}_\frak p)}\Gamma(\mathrm{Tr}(\mathbf{s}_\frak p)),$
with $\mathrm{Tr}(\mathbf{s}_\frak p):=s_\tau+s_{\bar\tau}$.}

\noindent
Proof. The first statement is clear in view of the product 
decomposition $$\Big(\bold R_+^*,\frac{dy}{y}\Big)=
\Big(\prod_\frak p \bold R_{+;\frak p}^*,\prod_\frak p 
\frac{dy_\frak p}{y_\frak p}\Big).$$ The second is relative to a 
 $G(\mathbb C|\mathbb R)$-set $X$ which consists of only one 
conjugation class. If $\#X=1$,  trivially
  $\Gamma_X(\mathbf{s})=\Gamma(s)$. So let $X=\{\tau,\bar\tau\},
\tau\not=\bar\tau$. Mapping 
  $\psi:\mathbb R_+^*\to 
\bold R_+^*,\ t\mapsto (\sqrt t,\sqrt t)$, from definition,
$$\int_{\bold R_+^*}N(e^{-y}y^\mathbf{s})\frac{dy}{y}=
\int_{\bold R_+^*}N(e^{-(\sqrt t,\sqrt t)}
(\sqrt t,\sqrt t)^{(s_\tau,s_{\bar\tau})})\frac{dt}{t}=\int_0^\infty 
e^{-2\sqrt t}\sqrt t^{\mathrm{Tr}(\mathbf{s})}\frac{dt}{t}.$$
Note that $d(\frac{t}{2})^2=2\frac{dt}{t}$, the substitution $t\mapsto 
(\frac{t}{2})^2$ yields $$\int_{\bold R_+^*}N(e^{-y}y^\mathbf{s})
\frac{dy}{y}=2^{1-\mathrm{Tr}(\mathbf{s}_\frak p)}\Gamma
(\mathrm{Tr}(\mathbf{s}_\frak p)),$$
 as desired.
\vskip 0.20cm
Neukirch called the function $L_X(\mathbf{s})=N(\pi^{-\mathbf{s}/2})
\Gamma_X(\mathbf{s}/2)$ the $L$-function of the  
$G(\mathbb C|\mathbb R)$-set $X$. Accordingly, 
$$L_X(\mathbf{s})=\prod_\frak p L_\frak p(\mathbf{s}_\frak p)\qquad
\mathrm{with}\qquad L_\frak p(\mathbf{s}_\frak p)
=\begin{cases} 
\pi^{-\mathbf{s}_\frak p/2}\Gamma(\mathbf{s}_\frak p/2)
& \frak p\ \mathrm{real}\\
2(2\pi)^{-\mathrm{Tr}(-\mathbf{s}_\frak p)/2}
\Gamma(\mathrm{Tr}(\mathbf{s}_\frak p/2))
& \frak p\ \mathrm{complex}.\end{cases}$$
For a  complex number $s\in \mathbb C$, we put
 $\Gamma_X(s)=\Gamma(s\bold 1),$ where 
$\bold 1=(1,\cdots,1)$ is the unit element of $\bold C$. 
Denote by $r_1$, resp. $r_2$ the number of real, 
resp., complex, conjugation classes of $X$. Then 
$$\Gamma_X(s)=2^{(1-2s)r_2}\Gamma(s)^{r_1}\Gamma(2s)^{r_2}.$$
Similarly, we put $L_X(s)=L_X(s\bold 1)=\pi^{-ns/2}
\Gamma_X(s/2)$ with $n:=\#X$. Then, 
$$\begin{cases}L_\mathbb R(s)=L_X(s)=\pi^{-s/2}\Gamma(s/2), 
\qquad&\mathrm{if}\ X=\{\tau\}\\
L_\mathbb C(s)=L_X(s)=2(2\pi)^{-s}\Gamma(s), 
\qquad&\mathrm{if}\ X=\{\tau,\bar\tau\},\ \tau\not=\bar\tau,\end{cases}$$
and $$L_X(s)=L_\mathbb R(s)^{r_1}\cdot L_\mathbb
 C(s)^{r_2}.$$ Moreover, from the standard facts that
$$\Gamma(s+1)=s\,\Gamma(s),\quad \Gamma(s)\,\Gamma(1-s)
=\frac{\pi}{\sin\,\pi s},
\quad\mathrm{and}\ \ \frac{2^s}{2}\Gamma(\frac{s}{2})\,\Gamma(\frac{s+1}{2})=
\sqrt{\pi}\,\Gamma(s),$$ easily, we obtain the following 
\vskip 0.30cm
\noindent
{\bf Basic Relations.}

\noindent
(1) $L_\mathbb R(1)=1,\ L_\mathbb C(1)=\frac{1}{\pi};$

\noindent
(2) $L_\mathbb R(s+2)=\frac{s}{2\pi}L_\mathbb R(s),\hskip 1.72cm
 L_\mathbb C(s+1)=\frac{s}{2\pi}L_\mathbb C(s);$

\noindent
(3) $L_\mathbb R(1-s)L_\mathbb R(1+s)=\frac{1}
{\cos(\pi\frac {s}{2})},\qquad L_\mathbb C(s)L_\mathbb C(1-s)
=\frac{2}{\sin(\pi s)};$

\noindent
(4) $L_\mathbb R(s)L_\mathbb R(1+s)=L_\mathbb C(s)\qquad$ 
({\large Legendre's Duplication Formula}).

\chapter{Rank Two $\mathcal O_K$-Lattices: Stability and Distance to Cusps}

Typically, each section of this chapter consists of three parts: 
1) upper half plane, 2) upper half space and 3) moduli 
spaces of semi-stable $\mathcal O_K$-lattices of rank 2.
Parts 1) and 2) are preperations for Part 3), the central one,
and are for reader's convenience. As such, no originalities
in any sense from us in subsections 1) and 2). 
In fact, even the presentations mainly  
follow the classics such as Kubota [Kub],  
Elstrodt et al [EGM] and Siegel [S]. 

\section{Upper Half Space Model}

\subsection{Upper Half Plane}

The {\it upper half plane} $\mathcal H$ in complex plane $\mathbb C$ 
is defined to be 
$${\mathcal H}:=\{z=x+iy\in \mathbb C, x\in\mathbb R,
 y\in\mathbb R_+^*\}.$$
On ${\mathcal H}$, the natural hyperbolic metric is given by 
the line element $$ds^2:=\frac{dx^2+dy^2}{y^2}.$$
It is well-known that the geodesics with respect to this  
hyperbolic metric, which are sometimes called {\it hyperbolic lines}, are half 
 circles or half lines in ${\mathcal H}$ which are orthogonal to the 
boundary line $\mathbb R$ in the Euclidean sense. Moreover, 
the volume form of hyperbolic metric is given by 
$$d\mu:=\frac {dx\wedge dy}{y^2}.$$
Consequently, the associated hyperbolic Laplace-Beltrami operator can be 
written as  $$\Delta:=y^2\big(\frac{\partial^2}{\partial x^2}+
\frac{\partial^2}{\partial y^2}\big).$$ 

The natural action on $\mathcal H$ of the group $SL(2,\mathbb R)$ of 
real $2\times 2$ metrices with determinant one is given by: 
$$M\,z:=\frac{az+b}{cz+d},\qquad\forall 
M=\left(\begin{matrix} a&b\\ c&d\end{matrix}\right)\in 
SL(2,\mathbb R),\ \ z\in {\mathcal H}.$$
Easily, if we write $Mz=x^*+iy^*$ with $x^*,\,y^*\in\mathbb R$,
then $$x^*=\frac{(ax+b)(cx+d)+acy^2}{(cx+d)^2+c^2y^2},\qquad
y^*=\frac {y}{(cx+d)^2+c^2y^2}>0.$$
In particular, $y^*$ depends only on $z$ and the second row of $M$.

As said, ${\mathcal H}$ admits the real line $\mathbb R$ as its boundary. 
Consequently, to compactify it, we add on it the real projective line 
 $\mathbb P^1(\mathbb R)$ with $\infty=\left[\begin{matrix} 1\\ 
0\end{matrix}\right]$. Naturally, the above action of $SL(2,\mathbb R)$ 
also extends to $\mathbb P^1(\mathbb R)$ via 
 $$\left(\begin{matrix} a&b\\ c&d\end{matrix}\right)\left[\begin{matrix} x\\ 
y\end{matrix}\right]=
 \left[\begin{matrix} ax+by\\ cx+dy\end{matrix}\right].$$ 

Back to $\mathcal H$ itself. The stablizer of 
$i=(0,1)\in {\mathcal H}$ with respect to
the action of $SL(2,\mathbb R)$ on ${\mathcal H}$ is equal to 
$SO(2):=\{A\in O(2):\det A=1\}$. Since the action 
of $SL(2,\mathbb R)$ on ${\mathcal H}$ is transitive,  we can identify
the quotient $SL(2,\mathbb R)/SO(2)$ with $\mathcal H$ given by the 
quotient map induced from $$SL(2,\mathbb R)\to {\mathcal H},\qquad
g\mapsto g\cdot i.$$ 

Associated with the coset $M\cdot SO(2)$, $M\in SL(2,\mathbb R)$ is
the positive definite matrix $Y=M\cdot M^t$ of size $2\times 2$. 
This gives an injection of 
$SL(2,\mathbb R)\Big/SO(2)$ into the set $\mathcal P^+$ of positive 
definite matrices of size 2, and  hence an injection to the space of inner 
products on $\mathbb R^2$. If we choose by Iwasawa decomposition 
in each coset $M\cdot SO(2)$ the uniquely determined representative 
$\left(\begin{matrix} 1&x\\ 0&1\end{matrix}\right)\cdot \left(\begin{matrix} \sqrt{y}&0\\ 0&\frac{1}{\sqrt{y}}\end{matrix}\right)=\left(\begin{matrix} 
\sqrt y&\frac{x}{\sqrt y}\\ 0&\frac{1}{\sqrt y}
\end{matrix}\right),\ x\in \mathbb R,\ y>0$, we can parametrize the image 
of $SL(2,\mathbb R)\Big/SO(2)$ in 
$\mathcal P^+$ by means of $Y=M\cdot M^t=\left(\begin{matrix} 
\frac{x^2+y^2}{y}&\frac{x}{y}\\ \frac{x}{y}&
\frac{1}{y}\end{matrix}\right).$ This leads us to introduce 
$$z=\left(\begin{matrix} 1&x\\ 0&1\end{matrix}\right)\cdot 
\left(\begin{matrix} \sqrt{y}&0\\ 0&\frac{1}{\sqrt{y}}\end{matrix}\right)\, i
=\left(\begin{matrix} 1&x\\ 0&1\end{matrix}\right)(iy)=x+iy\in {\mathcal H}$$ 
as a coordinate for $Y$ and hence 
for $M\cdot SO(2)$ as well. Indeed, a computation certainly shows that 
the natural action of  $SL(2,\mathbb R)$ on
 $SL(2,\mathbb R)/SO(2)$ (induced from  the multiplication of cosets 
from the left) is  expressed in terms of the coordinate 
 $z=x+iy$ exactly by the formula above.
 
\subsection{Upper Half Space}

The {\it upper half space}   ${\mathbb H}$ in Euclidean 3-space $\mathbb R^3$ 
gives a convenient model of 3-dimensional 
hyperbolic space with in its properties closely resembles the upper 
half plane as a model of plane hyperbolic geometry. We use the following 
coordinates
$$\begin{aligned}{\mathbb H}:=&\mathbb C\times]0,\infty[\,=\,\Big\{(z,r):z=x+iy\in \mathbb C, 
r\in\mathbb R_+^*\Big\}\\
=&\Big\{(x,y,r):x,y
\in\mathbb R, r\in\mathbb R_+^*\Big\}.\end{aligned}$$
To facilitate computation, we will think of ${\mathbb H}$ as a subset of 
Hamilton's quaternions. As usual, if we write $1,\,i,\,j,\,k$ for the standard 
$\mathbb R$-basis of the quaternions, we may write points $P$ in 
 ${\mathbb H}$ as $$P=(z,r)=(x,y,r)=z+rj\qquad\mathrm{where}\  \
z=x+iy,\  j=(0,0,1).$$

We equip ${\mathbb H}$ with the hyperbolic metric coming from the line
element $$ds^2:=\frac{dx^2+dy^2+dr^2}{r^2}.$$ The geodesics 
with respect to the 
hyperbolic metric, which are sometimes 
called {\it hyperbolic lines}, are half circles or half lines in 
${\mathbb H}$ which are orthogonal to the
 boundary plane $\mathbb C$ in the Euclidean sense. The {\it hyperbolic 
planes} (also called geodesic hyperplanes),
  that is, the isometrically embedded copies of 2-dimensional 
hyperbolic space, are Euclidean hemispheres or
   half-planes which are perpendicular to the boundary 
$\mathbb C$ of ${\mathbb H}$ (in the Euclidean sense).

Moreover, the associated  hyperbolic volume form is given by 
$$d\mu:=\frac {dx\wedge dy\wedge dr}{r^3}.$$
Consequently, the hyperbolic Laplace-Beltrami operator associated to 
the hyperbolic metric $ds^2$ is simply 
$$\Delta:=r^2\Big(\frac{\partial^2}{\partial x^2}+\frac{\partial^2}
{\partial y^2}+\frac{\partial^2}
{\partial r^2}\Big)-r\frac{\partial}{\partial  r}.$$ 

The natural action of $SL(2,\mathbb C)$ on ${\mathbb H}$ and on its boundary 
$\mathbb P^1(\mathbb C)$ may be described as follows:
We represent an element of $\mathbb P^1(\mathbb C)$ by 
$\left[\begin{matrix} x\\ y\end{matrix}\right]$ where
 $x,y\in\mathbb C$ with $(x,y)\not=(0,0)$. Then the action of the matrix 
  $M=\left(\begin{matrix} a&b\\ c&d\end{matrix}\right)\in SL(2,\mathbb C)$ 
on $\mathbb P^1(\mathbb C)$ is defined to be
$$\left[\begin{matrix} x\\ y\end{matrix}\right]\mapsto \left(\begin{matrix} a&b\\ c&d\end{matrix}\right)\left[\begin{matrix} x\\ 
y\end{matrix}\right]:=
\left[\begin{matrix} ax+by\\ cx+dy\end{matrix}\right].$$ Moreover, if we 
represent points $P\in {\mathbb H}$ as 
quaternions whose fourth component equals zero, then the action of 
$M$ on ${\mathbb H}$ is defined to be 
$$P\mapsto M\,P:=(aP+b)(cP+d)^{-1},$$
where the inverse on the right is taken in the skew field of quaternions. 
Indeed, if we set $M(z+rj)=z^*+r^*j$ with $z^*\in\mathbb C, r^*\in\mathbb R$, 
then an obvious computation shows that 
 $$z^*:=\frac{(az+b)(\bar c\bar z+\bar d)+
a\bar c r^2}{|cz+d|^2+|c|^2r^2}, \qquad
r^*:=\frac{r}{|cz+d|^2+|c|^2r^2}=\frac{r}{\|cP+d\|^2}.$$  In particular, 
$r^*$  depends only on $P$ and the second row of $M$. Moreover
$r^*>0$, so  $M(z+rj)\in\mathbb H$ as well. (Here we 
have set $P=z+rj$ and used $\|cP+d\|$ to denote 
the Euclidean norm of the vector $cP+d\in \mathbb R^4$, which is 
indeed also the square root of the norm of $cP+d$ in the quaternions.)

Furthermore, with this action, the stablizer of $j=(0,0,1)\in 
{\mathbb H}$ in $SL(2,\mathbb C)$  is equal to 
$SU(2):=\{A\in U(2):\det A=1\}$. Since the 
action of $SL(2,\mathbb C)$ on ${\mathbb H}$ is transitive, 
we obtain also a natural identification ${\mathbb H}\simeq
SL(2,\mathbb C)/SU(2)$ via the quotient map induced from 
$SL(2,\mathbb C)\to {\mathbb H},\ \ g\mapsto g\cdot j.$ 

Associated with the coset $M\cdot SU(2)$, $M\in SL(2,\mathbb C)$ is
the positive definite hermitian matrix 
$Y=M\cdot\overline{M^t}$ of size $2\times 2$. This gives an injection of 
$SL(2,\mathbb C)/SU(2)$ into the set $\mathcal P_{\mathbb C}^+$ of 
positive definite hermitian matrices of size 2, and hence an 
injection to the space of hermitian inner products 
on $\mathbb C^2$. If we choose by Iwasawa decomposition in each coset
 $M\cdot SU(2)$ the uniquely determined representative 
$\left(\begin{matrix} 1&z\\ 0&1\end{matrix}\right)
\left(\begin{matrix} \sqrt r&0\\ 
0&\frac{1}{\sqrt r}\end{matrix}\right)=
\left(\begin{matrix} \sqrt r&\frac{z}{\sqrt r}\\ 
0&\frac{1}{\sqrt r}\end{matrix}\right),\ 
 z\in \mathbb C,\ r>0$, we can parametrize the image of 
$SL(2,\mathbb C)\Big/SU(2)$ in $\mathcal P_\mathbb C^+$ by
  means of $Y=M\cdot\bar M^t=\left(\begin{matrix} 
\frac{|z|^2+r^2}{r}&\frac{z}{r}\\ \frac{\bar z}{r}&\frac{1}{r}
  \end{matrix}\right).$ This leads us to introduce the quoternion 
$P=\left(\begin{matrix} 1&z\\ 0&1\end{matrix}\right)
\left(\begin{matrix} \sqrt r&0\\ 
0&\frac{1}{\sqrt r}\end{matrix}\right)\,j=z+rj\in{\mathbb H}$ as a coordinate for
$Y$ and hence for $M\cdot SU(2)$. Indeed, an obvious computation 
shows that the natural action of  $SL(2,\mathbb C)$ on
$SL(2,\mathbb C)\Big/SU(2)$ (induced from the multiplication of cosets 
from the left) is expressed in terms of the coordinate 
$P$ exactly by the formula above. We leave the details to the reader.

\subsection{Rank Two $\mathcal O_K$-Lattices: Upper Half Space Model}

With above discussion, we see that after identifying
${\mathcal H}$ with $SL(2,\mathbb R)/SO(2)$ and ${\mathbb H}$ with 
$SL(2,\mathbb C)/SU(2)$, using the discussion in Chapter 1, in particular,
\S1.9, we conclude that
$$\mathcal M_{K,2;\frak a}[N(\frak a)\cdot \Delta_K]\simeq 
\bigg(SL(\mathcal O_K\oplus\frak a)\Big\backslash 
\Big({\mathcal H}^{r_1}\times{\mathbb H}^{r_2}\Big)\bigg)_{\mathrm {ss}},$$ 
where as before ss 
means the subset consisting of points 
corresponding to rank two semi-stable $\mathcal O_K$-lattices
in $SL(\mathcal O_K\oplus\frak a)\Big\backslash 
\bigg(\Big(SL(2,\mathbb R)/SO(2)\Big)^{r_1}\times \Big(SL(2,\mathbb C)/SU(2)
\Big)^{r_2}\bigg).$

Put this in a more concrete term, if the metric on $\mathcal O_K\oplus\frak a$ 
is given by matrices 
$g=(g_\sigma)_{\sigma\in S_\infty}$ with $g_\sigma\in SL(2,K_\sigma)$, 
then the corresponding points on the 
right hand side is $g(\mathrm{ImJ})$ with 
$\mathrm {ImJ}:=(i^{(r_1)},j^{(r_2)})$, i.e., the point given by 
$(g_\sigma\bold\tau_\sigma)_{\sigma\in S_\infty}$ where 
$\tau_\sigma=i_\sigma:=(0,1)$ if $\sigma$ is real and 
 $\tau_\sigma=j_\sigma:=(0,0,1)$ if $\sigma$ is complex. 
As before, $SL(\mathcal O_K\oplus\frak a)$ 
 denotes elements in $GL(\mathcal O_K\oplus\frak a)$ with 
determinant 1, so that, in particular, 
$$SL(\mathcal O_K\oplus\frak a)=SL(2,K)\cap 
\left(\begin{matrix} \mathcal O_K&\frak a\\ \frak a^{-1}&\mathcal 
O_K\end{matrix}\right).$$ In other words, if
$\left(\begin{matrix} a&b\\ c&d\end{matrix}\right)\in 
SL(\mathcal O_K\oplus\frak a)$, then $ad-bc=1$ and $a, d\in
\mathcal O_K, b\in\frak a,$ and $c\in \frak a^{-1}$.
 
\section{Cusps}

\subsection{Upper Half Plane}

By definition, a subgroup $\Gamma\subset SL(2,\mathbb R)$ 
is called a {\it discontinuous group} if  for 
every $z\in {\mathcal H}$ and for every sequence $(T_n)_{n\geq 1}$ 
of distinct elements of $\Gamma$, the 
sequence $(T_nz)_{n\geq 1}$ has no accumulation point in 
${\mathcal H}$; and a subgroup  
$\Gamma\subset SL(2,\mathbb R)$ is called a {\it discrete group} 
if its image in $SL(2,\mathbb R)\subset \mathbb R^4$ 
is discrete with respect to the topology induced from $\mathbb R^4$. 
It is well-known that a subgroup   $\Gamma\subset SL(2,\mathbb R)$ 
is discontinuous if and only if it is discrete. We call a discrete 
subgroup of $SL(2,\mathbb R)$ a {\it Fuchsian group}.

Recall also that an element $\gamma\in SL(2,\mathbb R)$ is called 
{\it parabolic} if it conjugates in 
$SL(2,\mathbb R)$ to $\left(\begin{matrix} 1&1\\ 0&1\end{matrix}\right)$, 
or better if  $\mathrm{Tr}(\gamma)=2$. 
One knows that
an element $\gamma\not=I$ in $SL(2,\mathbb R)$ is parabolic 
if and only if it has exactly one fixed 
point on the boundary $\mathbb P^1(\mathbb R)$ of $\mathcal H$.

For a point $P\in \widehat{\mathcal H}:=
{\mathcal H}\cup\mathbb P^1(\mathbb R)$, as usual, define its 
{\it stablizer group} 
$\Gamma_P$ in $\Gamma$ by $\Gamma_P:=\{\gamma\in\Gamma:\gamma P=P\}$. 
Then for $\zeta\in \mathbb P^1
(\mathbb R)=\mathbb R\cup\{\infty\}$, we want to describe the stablizer 
group $\Gamma_\zeta$ of $\zeta$. 
This may be done by transforming $\zeta$ to 
$\infty$ and assuming without loss of generality that 
$\zeta=\infty$. Indeed, we can use the following two groups to proceed: 
$$B(\mathbb R):=\Bigg\{\left(\begin{matrix} a&b\\ 0&a^{-1}\end{matrix}\right):0
\not=a\in \mathbb R, b\in\mathbb R\Bigg\},\qquad
 N(\mathbb R):=\Bigg\{\left(\begin{matrix} 1&b\\ 0&1\end{matrix}\right): 
b\in\mathbb R\Bigg\}.$$
(The group $B(\mathbb R)$ is the Borel subgroup of $SL(2,\mathbb R)$ 
and $N(\mathbb R)$ is its unipotent radical. As abstract groups, 
$N(\mathbb R)$ is isomorphic to the additive group $\mathbb R^+$, while
  $B(\mathbb R)$ is isomorphic to the semi-direct product 
$\mathbb R^*$ by $\mathbb R^+$.)
Clearly, $SL(2,\mathbb R)_\infty=B(\mathbb R)$; and  if  
$\zeta\in \mathbb P^1(\mathbb R)$, then there is an
$A:=A_\zeta\in SL(2,\mathbb R)$ such that $A\cdot\infty=\zeta$. Consequently,
we have $\Gamma_\zeta=\Gamma\cap \Big(A\cdot B(\mathbb R)\cdot A^{-1}\Big)$. 
 
Put now $Z(\Gamma):=\Gamma\cap\{\pm I\}$. When 
$z\in \mathbb P^1(\mathbb R)$ is a fixed point of a parabolic 
element of $\Gamma$, we call $z$  a {\it cusp} of $\Gamma$. 
It is well-known that for a cusp $z$, 
$$\Gamma_z/Z(\Gamma)\simeq\mathbb Z,\qquad\mathrm{and}\qquad 
A_x^{-1}\cdot \Gamma_x \cdot A_x=\Bigg\{\pm\left(\begin{matrix} 
1&h\\ 0&1\end{matrix}\right)^m:m\in\mathbb Z\Bigg\}$$ for a certain $h>0$. Let 
$\mathcal C_\Gamma$ denote the collection of all cusps of
 $\Gamma$ and put ${\mathcal H}^*:={\mathcal H}^*_\Gamma:=
{\mathcal H}\cup \mathcal C_\Gamma$, then by definition, 
 $\Gamma$ is called a {\it Fuchsian group of the first kind} if 
$\Gamma\backslash {\mathcal H}^*$ is compact. It is also well-known that 
this definition is equivalent to the condition that 
$\Gamma$ is Fuchsian and the hyperbolic 
volume of $\Gamma\backslash {\mathcal H}^*$ is finite.

\subsection{Upper Half Space}

A subgroup $\Gamma\subset SL(2,\mathbb C)$ is called a {\it discontinuous 
group} if  for every 
$P\in {\mathbb H}$ and for every sequence $(T_n)_{n\geq 1}$ of 
distinct elements of $\Gamma$, the
 sequence $(T_nP)_{n\geq 1}$ has no accumulation point in 
${\mathbb H}$; while a subgroup  
 $\Gamma\subset SL(2,\mathbb C)$ is called a {\it discrete subgroup} 
if its image in $SL(2,\mathbb C)\subset \mathbb C^4$
is discrete with respect to the topology induced from that of $\mathbb C^4$. 
It is well-known  that 
a subgroup   $\Gamma\subset SL(2,\mathbb C)$ is 
 discontinuous if and only if it is discrete.

An element $\gamma\in SL(2,\mathbb C), \gamma\not=\{\pm I\}$ 
is called {\it parabolic} if $|\mathrm{Tr}(\gamma)|=2$. One 
checks that an element $\gamma\in SL(2,\mathbb C)$, $\gamma\not=I$ 
is parabolic if and only if it is 
conjugate in $SL(2,\mathbb C)$ to $\left(\begin{matrix} 1&1\\ 
0&1\end{matrix}\right)$, if and only if
it has exactly one fixed point in the boundary $\mathbb P^1(\mathbb C)$
of $\mathbb H$.

For a point $P\in \widehat{\mathbb H}:=\mathbb H \cup
\mathbb P^1(\mathbb C)$, as usual, define its {\it stablizer group} 
$\Gamma_P$ in $\Gamma$ by $\Gamma_P:=\{\gamma\in\Gamma:\gamma P=P\}$. 
We shall now give a description of the stablizer group $\Gamma_\zeta$ of
 $\zeta\in \mathbb P^1(\mathbb C)=\mathbb C\cup\{\infty\}$. As above,
we  do this by transforming $\zeta$ to 
 $\infty$ and assuming without loss of generality that 
$\zeta=\infty$. Similarly, let us first introduce the following 
 two groups: $$B(\mathbb C):=\Bigg\{\left(\begin{matrix} a&b\\ 
0&a^{-1}\end{matrix}\right):0\not=a\in \mathbb C, b
 \in\mathbb C\Bigg\}, N(\mathbb C):=\Bigg\{\left(\begin{matrix} 1&b\\ 
0&1\end{matrix}\right): b\in\mathbb C\Bigg\}\subset 
 B(\mathbb C).$$ (The group $B(\mathbb C)$ is the Borel 
subgroup of $SL(2,\mathbb C)$ and $N(\mathbb C)$ is 
 its unipotent radical. As abstract groups, 
$N(\mathbb C)$ is isomorphic to the additive group $\mathbb C^+$, and
$B(\mathbb C)$ is isomorphic to the semi-direct 
product of $\mathbb C^*$ by $\mathbb C^+$, defined by $\sigma_a(b)=a^2b$.)
Clearly, $SL(2,\mathbb C)_\infty=B(\mathbb C)$; and  if  
$\zeta\in \mathbb P^1(\mathbb C)$, then there is an 
$A=A_\zeta\in SL(2,\mathbb C)$ such that 
$A\cdot \infty=\zeta$. Consequently, we have $\Gamma_\zeta=\Gamma\cap 
\Big(A\cdot B(\mathbb C)\cdot A^{-1}\Big)$. 

For the purpose of factoring out possible twists from the 
so-called \lq elliptic elements', set then
 $\Gamma_\zeta':=\Gamma\cap \Big(A\cdot N(\mathbb C)\cdot A^{-1}\Big)
=\Gamma_\zeta\cap\Big(A\cdot N(\mathbb C)\cdot A^{-1}\Big)$, which consists
of  parabolic elements in $\Gamma_\zeta$ together with the identity
$I$. Then it is known that there are following 3 possibilities
for $\Gamma_\zeta'$ (see [EGM] for details):

\noindent
(1) $\Gamma_\infty'=\{I\}$; or 

\noindent
(2) $\Gamma_\infty'$ is isomorphic to $\mathbb Z$; or

\noindent 
(3) $\Gamma_\infty'$ is a lattice in $N(\mathbb C)\simeq \mathbb C$.

It is this final case 3) that we want to pursue.  
In this case, further, there are 3 subcases:

\noindent
(i) $\Gamma_\infty=\Gamma_\infty'$;

\noindent
(ii) $\Gamma_\infty$ is conjugate in $B(\mathbb C)$ to a group of the form
$$\bigg\{\left(\begin{matrix} \varepsilon&\varepsilon b\\ 
0&\varepsilon^{-1}\end{matrix}\right):\ b\in\Lambda,\ 
\varepsilon =\{1,i\}\bigg\},$$
where $\Lambda\subset \mathbb C$ is an arbitrary lattice. 
(The abstract group $\Gamma_\infty$ is isomorphic 
to $\mathbb Z^2\times \mathbb Z/2\mathbb Z$ where the 
nontrivial element of $\mathbb Z/2\mathbb Z$ acts by 
multiplication by $-1$);

\noindent
(iii)  $\Gamma_\infty$ is conjugate in $B(\mathbb C)$ to a group of the form
$$\Gamma(n,t):=\bigg\{\left(\begin{matrix} \varepsilon&\varepsilon b\\ 
0&\varepsilon^{-1}\end{matrix}\right):\ b\in
\mathcal O_K,\ \varepsilon=\exp\big(\frac{\pi ivt}{n}\big),\  
1\leq v\leq 2n\bigg\},$$
where $n=4$ or 6 and $t|n, \mathcal O_n$ is the ring of integers 
in the quadratic number field $\mathbb Q\Big(\exp
\big(\frac{2\pi i}{n}\big)\Big).$ (Hence as an abstract group 
$\Gamma_\infty$ is isomorphic to $\mathbb Z^2\times
 \mathbb Z/n\mathbb Z$ for some $n\in\{1,2,3,4,6\}$. An
 element $\varepsilon\in \mathbb Z/m\mathbb Z\simeq\Big\{
\exp\big(\frac{\pi iv}{m}\big):v \in \mathbb Z\Big\}$ acts on 
$\mathbb Z^2\simeq \mathcal O_{m'}^+$ by multiplication 
with $\varepsilon^2$ where $m'=4$ in case $m=1,2,4$ and 
$m'=6$ otherwise.)
\vskip 0.20cm
By definition, an element $\zeta\in \mathbb P^1(\mathbb C)$ 
is called a {\it cusp} of a discrete group 
$\Gamma\subset SL(2,\mathbb C)$ if $\Gamma_\zeta$ contains a 
free abelian group of rank 2. We write 
$\mathcal C_\Gamma$ for the set of cusps of $\Gamma$. Clearly, 
the group $\Gamma$ leaves $\mathcal C_\Gamma$ 
invariant, and breaks $\mathcal C_\Gamma$ into $\Gamma$-classes.
 Moreover, it is known that if $\Gamma\subset SL(2,\mathbb C)$ 
is a discrete group of finite covolume, then $\Gamma$ has only
finitely many $\Gamma$-classes of cusps.  Write 
them as $\eta_1,\ldots,\eta_h$ and fix $A_1,\ldots,A_h\in 
SL(2,\mathbb C)$ such that $\eta_1=A_1\cdot\infty,\ldots,
\eta_h=A_h\cdot \infty.$
(For detailed argument, please see [EGM], which we unconditionally 
follow here.)

\subsection{Rank Two $\mathcal O_K$-Lattices}

We further divide this subsection into two: A) for totally real fields
and B) for general number fields, for the purpose to indicate how a
general theory is built up on the classics. As such, if the reader
feels that there is a part in A (resp. in B) which is a bit clumsy, then 
she or he is suggested to refer the corresponding part in B (resp. in A)
for a clearer explanation. 

\subsubsection{ A. Totally Real Fields}

We unconditionally follow Siegel's presentation [S] in part A) 
(here and in the sequel). As such, there is a discrepency in terms of 
notations. But this is not serious as the content can be understood
without too much difficulty. So we will leave them as they are.
\vskip 0.20cm
Let $K$ be a totally real algebraic number field of degree $n$ over 
$\mathbb Q$,  and let $K^{(1)}=K,K^{(2)},\ldots, 
K^{(n)}$ be the conjugates of $K$. Naturally, this gives then an embedding
$$\mathbb P^1(K)\hookrightarrow \Big(\mathbb P^1(\mathbb R)\Big)^n\subset
\Big(\widehat{\mathcal H}\Big)^n,\qquad \lambda\mapsto(\lambda^{(1)},\cdots, 
\lambda^{(n)})=:\lambda.$$ For a given 
$z=(z_1,\ldots,z_n)\in\mathcal H^n$, 
the norm $N(z)$ shall stand for $\prod_{i=1}^nz_i$ 
and the trace $\mathrm{Tr}(z)$ for $\sum_{i=1}^n z_i$. 
If $\lambda\in K$, $N(\lambda)$ and $\mathrm{Tr}(\lambda)$ coincide with 
the usual norm and trace in $K_1$ respectively.

Let $G=SL(2,K)$ with $Z=\{\pm I\}$. Then for the factor group $G/Z$, 
we have a faithfulful representation as the group of mapping
$\Big(\widehat {\mathcal H}\Big)^n\to \Big(\widehat{\mathcal H}\Big)^n$ 
defined 
by  $z=(z_1,\ldots,z_n)\mapsto z_M:=(z_1^*,\ldots,z_n^*)$ with 
$z_j^*:=\frac{\alpha^{(j)} z_j+\beta^{(j)}}
{\gamma^{(j)} z_j+\delta^{(j)}}$ corresponding to each 
$M=\begin{pmatrix} \alpha&\beta\\ \gamma&\delta\end{pmatrix}\in G.$ We 
shall denote $z_M$ as $\frac{\alpha z+\beta}{\gamma z+\delta}$ symbolically.
In particular, for $\infty=\frac{1}{0}$, 
$\infty_M=\Big(\frac{\alpha^{(1)}}{\gamma^{(1)}},\ldots, 
\frac{\alpha^{(n)}}{\gamma^{(n)}}\Big).$
 Let further $z_j=x_j+iy_j$,  
$z_j^*=x_j^*+iy_j^*$. Writing $z=x+iy$ where $x=(x_1,\ldots,x_n)$ and 
$y=(y_1,\ldots,y_n)$, we shall denote $(x_1^*,\ldots,x_n^*)$ by $x_M$ and 
$(y_1^*,\ldots,y_n^*)$ by $y_M$ so that $z_M=x_M+iy_M$. It is easy to see 
that for $M_1,M_2\in G$, $z_{M_1M2}=(z_{M_2})_{M_1}$.

Let $\Gamma:=SL(2,\mathcal O_K)$ be the subgroup of $G$ consisting of 
$\begin{pmatrix} \alpha&\beta\\ \gamma&\delta\end{pmatrix}$ with 
$\alpha,\beta, \gamma, \delta\in \mathcal O_K.$ The factor group 
$\Gamma/Z$ is precisely the {\it inhomogeneous Hilbert modular group}, which 
we shall denote  by $\Gamma_K$.

One can consider more general group than $\Gamma$, say, the group 
$\Gamma_0=GL^+(2,\mathcal O_K),$ consisting of
$\begin{pmatrix} \alpha&\beta\\ \gamma&\delta
\end{pmatrix}$ with $\alpha,\beta, \gamma, \delta\in \mathcal O_K$ and 
$\alpha\cdot\delta-\beta\gamma=\varepsilon$ with $\varepsilon$ being a 
totally positive unit in $K$. Let $Z_0$ be the subgroup of $\Gamma_0$ 
consisting of matrices of the  form $\begin{pmatrix} \varepsilon&0\\ 
0&\varepsilon\end{pmatrix}$ with $\varepsilon>0$ a unit in $K$. Further, 
let $\rho_1=1,\,\rho_2,\,\ldots,\,\rho_{\mu(K,2)}$ be a complete set of 
representatives of the group $U_K^+$ of units $\varepsilon>0$ in $K$ 
modulo the subgroup $U_K^2$of squares of units in $K$. Then it is clear 
that $\Gamma_0/Z_0$ is isomorphic to the
 group of substitutions $z\mapsto \rho\cdot z_M$, where $z\mapsto z_M$ is a 
Hilbert modular substitution and $\rho=\rho_i$ for some $i$. Thus the study 
of $\Gamma_0/Z_0$ can be reduced to that of $\Gamma_K.$ The group 
$\Gamma_K$ is in general smaller than $\Gamma_0/Z_0$ and is called 
therefore the {\it narrow Hilbert modular group} usually.
\vskip 0.20cm
Two elements $\lambda,\ \mu$ in $\mathbb P^1(K)$ are called
{\it equivalent} (in symbols, $\lambda\sim\mu$), if for some 
$M=\begin{pmatrix} \alpha&\beta\\ 
\gamma&\delta\end{pmatrix}\in\Gamma$, $\mu=\lambda_M=
\frac{\alpha\lambda+\beta}{\gamma\lambda+\delta}$. This is a genuine 
equivalence relation. We shall presently show that $\mathbb P^1(K)$ falls into 
$h$ equiavlence classes, where $h=h_K$ is the class number of $K$.
\vskip 0.20cm
\noindent
{\bf{\Large Proposition}}. {\it There exist $\lambda_1,\ldots,\lambda_h
\in\mathbb P^1(K)$ such that for any $\lambda\in\mathbb P^1(K)$, we have 
$\lambda\sim \lambda_i$ for some $i$ uniquely determined by $\lambda$.}

\noindent
Proof. Any $\lambda\in\mathbb P^1(K)$ is of the form $\frac{\rho}{\sigma}$ 
with $\rho,\sigma\in \mathcal O_K$. (If $\lambda=\infty$, we may take 
$\rho=1$ and $\sigma=0$.) With $\lambda$, introduce the associate the 
integral ideal 
$\frak a:=(\rho,\lambda)$. We first see that 
 if $\frak a_1=(\theta)\frak a$ 
is any integral ideal in the class of $\frak a$, then $\frak a_1$ is of the 
form $(\rho_1,\sigma_1)$ where $\rho_1$ and $\sigma_1$ are in $\mathcal O_K$ 
such that $\lambda=\frac{\rho_1}{\sigma_1}$.
This is quite obvious, since $\frak a_1=(\rho\theta,\sigma\theta)$ and taking 
$\rho_1=\rho\theta, \sigma_1=\sigma\theta$, our assertion is proved.
Conversely,  any integral ideal $\frak a_1$ associated with 
$\lambda\in\mathbb P^1(K)$ in this way is necessarily in the same ideal-class 
as $\frak a$. In fact let $\lambda$ be written in the form 
$\frac{\rho_1}{\sigma_1}$ with $\rho_1$ and $\sigma_1$  in $\mathcal O_K$
 and let $\frak a_1=(\rho_1,\sigma_1).$ Then we claim that $\frak a_1$ 
is in the same class as $\frak a$. For since $\lambda=\frac{\rho}
{\sigma}=\frac{\rho_1}{\sigma_1}$, we have $\rho\sigma_1=\rho_1\sigma$, 
and hence $\frac{(\rho)}{\frak a}\cdot \frac{\sigma_1}{\frak a_1}=
\frac{(\rho_1)}{\frak a_1}\cdot \frac{\sigma}{\frak a}.$ Now 
$\frac{(\rho)}{\frak a}$ and  $\frac{(\sigma)}{\frak a}$ as also 
$\frac{(\rho_1)}{\frak a_1}$ and $\frac{(\sigma_1)}{\frak a_1}$
 are mutually coprime, as we may assume. Hence we have $\frac{(\rho)}{\frak a}=
\frac{(\rho_1)}{\frak a_1}$ and similarly $\frac{(\sigma)}{\frak a}=
\frac{(\sigma_1)}{\frak a_1}$. This means that  
$\frak a_1=(\theta)\frak a$ for a $\theta\in K$. (Consequently,
$\rho_1=\rho\theta$ and  $\sigma_1=\sigma\theta$.)

We choose in the $h$ ideal classes, fixed integral ideals 
$\frak a_1,\ldots,\frak a_h$ such that $\frak a_i$ is of minimum norm 
among all the integral ideals of its class. (Perhaps, the ideal $\frak a_i$ 
is not uniquely fixed in its class by this condition, but there are at most 
finitely many possibilities for $\frak a_i$ and we choose from these, a fixed 
$\frak a_i$.) It follows from the above that to a given 
$\lambda\in\mathbb P^1(K)$, we can make correspondingly an ideal 
$\frak a_i$ such that 
$\lambda=\frac{\rho}{\sigma}$ and $\frak a_i=(\rho,\sigma)$ for suitable 
$\rho,\sigma\in\mathcal O_K$. Now if $\mu=\frac{\alpha\lambda+\beta}
{\gamma\lambda+\delta}\sim \lambda$, then to $\mu$ again corresponds the 
ideal $(\alpha\lambda+\beta,\gamma\lambda+\delta)$ which is just $\frak a_i$
 since $\begin{pmatrix} \alpha&\beta\\ \gamma&\delta\end{pmatrix}$ is 
unimodular. Thus all elements of an equivalence class in $\mathbb P^1(K)$ 
correspond to the same ideal.

We shall now show that if the same ideal $\frak a_i$ corresponds to 
$\lambda,\ \lambda^*\in\mathbb P^1(K)$, then necessarily 
$\lambda\sim\lambda^*$. Set, for example, $\lambda=\frac{\rho}{\sigma},\
\lambda^*=\frac{\rho^*}{\sigma^*}$ and $\frak a_i=(\rho,
\sigma)=({\rho^*},{\sigma^*})$. It is well-known that there exist elements 
$\xi,\,\eta,\,\xi^*,\,\eta^*$ in $\frak a_i^{-1}$ such that 
$\rho\eta-\sigma\xi=1$ and 
$\rho^*\eta^*-\sigma^*\xi^*=1.$ Consequently, if we set
 $A:=\begin{pmatrix} \rho&\xi\\ 
\sigma&\eta\end{pmatrix}$ and $A^*:=\begin{pmatrix} \rho^*&\xi^*\\ 
\sigma^*&\eta^*\end{pmatrix}$; then $A,\,A^*\in G=SL(2,\mathbb R)$.  
Moreover, 
$A^*A^{-1}=\begin{pmatrix} \rho^*&\xi^*\\ \sigma^*&\eta^*\end{pmatrix}\cdot 
\begin{pmatrix} \eta&-\xi\\ -\sigma&\rho\end{pmatrix}=\begin{pmatrix} 
\alpha&\beta\\ \gamma&\delta\end{pmatrix},$ where 
 $\alpha=\rho^*\eta-\xi^*\sigma,\ \beta=-\rho^*\xi+\xi^*\rho,\ \gamma=
\sigma^*\eta-\eta^*\sigma,\ \delta=-\sigma^*\xi+\eta^*\rho$ clearly are 
integers in $K$ satisfying $\alpha\delta-\beta\gamma=1$. Hence, 
$A^*A^{-1}\in\Gamma$ and thus
$$\lambda^*=\frac{\rho^*}{\sigma^*}=\frac{\alpha\rho+\beta\sigma}
{\gamma\rho+\delta\sigma}=\frac{\alpha\lambda+\beta}{\gamma\lambda+\delta}
=\lambda_{A^*A^{-1}}\,\sim\,\lambda.$$
 
Thus we see that to different equivalence classes in $\mathbb P^1(K)$ 
correspond different ideals $\frak a_i$. It is almost trivial to verify that 
to each ideal $\frak a_i$, there corresponds an equivalence class in 
$\mathbb P^1(K)$. As a direct consequence, there are exactly $h$ equivalence 
classes in $\mathbb P^1(K)$ and our proposition is proved.
\vskip 0.30cm
We now make a convention to be followed (for A) in the sequel. We shall assume 
hat $\lambda_1=\frac{1}{0}=\infty$, without loss of generality. Moreover, 
with $\lambda_i=\frac{\rho_i}{\sigma_i}$, we associate a fixed matrix 
$A_i=\begin{pmatrix} \rho_i&\xi_i\\ 
\sigma_i&\eta_i\end{pmatrix}\in G.$ Let us remark that it is always 
possible to find, though not uniquely numbers $\xi_i,\,\eta_i$ in 
$\frak a_i^{-1}$ such that $\rho_i\eta_i-\xi_i\sigma_i=1$. Further, 
we shall suppose that $A_1=\begin{pmatrix} 1&0\\ 
0&1\end{pmatrix}.$

\subsubsection{B. Genaral Number Fields}

Now the working site is the space 
$SL(\mathcal O_K\oplus \frak a)\Big\backslash\Big({\mathcal H}^{r_1}
\times{\mathbb H}^{r_2}\Big)$ with ${\mathcal H}$ the upper half plane, 
${\mathbb H}$ the upper half space, 
and $SL(\mathcal O_K\oplus \frak a)$ the special automorphism
group defined by
$\bigg\{A\in\left(\begin{matrix} \mathcal O_K&\frak a\\ 
\frak a^{-1}&\mathcal O_K\end{matrix}\right):\det A=1\bigg\}.$
Here the action of $SL(\mathcal O_K\oplus \frak a)$ is 
via the action of $SL(2,K)$ on ${\mathcal H}^{r_1}
\times{\mathbb H}^{r_2}$. More precisely, $K^2$ admits 
natural embeddings $K^2\hookrightarrow 
\Big(\mathbb R^{r_1}\times\mathbb C^{r_2}\Big)^2\simeq 
\big(\mathbb R^2\big)^{r_1}\times\big(\mathbb C^2\big)^{r_2}$ so that
 $\mathcal O_K\oplus \frak a$ naturally embeds into 
$\big(\mathbb R^2\big)^{r_1}\times\big(\mathbb C^2\big)^{r_2}$ 
as a rank two $\mathcal O_K$-lattice. As such, 
$SL(\mathcal O_K\oplus \frak a)$ acts on the image of 
  $\mathcal O_K\oplus \frak a$ in $\big(\mathbb R^2\big)^{r_1}
\times\big(\mathbb C^2\big)^{r_2}$ as automorphisms. Our task 
here is to understand the cusps of this action of 
$SL(\mathcal O_K\oplus \frak a)$ on 
${\mathcal H}^{r_1}\times{\mathbb H}^{r_2}$. For this, we go as folllows.

First, the space ${\mathcal H}^{r_1}\times{\mathbb H}^{r_2}$ 
admits a natural boundary 
$\mathbb R^{r_1}\times\mathbb C^{r_2}$, in which the field $K$ 
is imbedded via Archmidean places in $S_\infty$: 
$K\hookrightarrow \mathbb R^{r_1}
\times\mathbb C^{r_2}$. Consequently, $\mathbb P^1(K)
\hookrightarrow \mathbb P^1(\mathbb R)^{r_1}\times
\mathbb P^1(\mathbb C)^{r_2}$ with $\left[\begin{matrix} 1\\ 
0\end{matrix}\right]:=\infty\mapsto (\infty^{(r_1)},\infty^{(r_2)})$. 
As usual, via fractional linear transformations, $SL(2,\mathbb R)$ 
acts on  $\mathbb P^1(\mathbb R)$, and $SL(2,\mathbb C)$ acts on 
 $\mathbb P^1(\mathbb C)$, hence so does $SL(2,K)$ on 
$$\mathbb P^1(K)\hookrightarrow \mathbb P^1(\mathbb R)^{r_1}\times 
\mathbb P^1(\mathbb C)^{r_2}.$$ Being a 
discrete subgroup of $SL(2,\mathbb R)^{r_1}\times SL(2,\mathbb C)^{r_2}$,
for the action of $SL(\mathcal O_K\oplus \frak a)$ 
 on $\mathbb P^1(K)$, we call the corresponding 
orbits (of $SL(\mathcal O_K\oplus \frak a)$ on $\mathbb P^1(K)$) 
the {\it cusps} (say, due to the fact that if we look at each 
local component ${\mathcal H}_\sigma$ and $\mathbb H_\sigma$, 
the induced orbits
 corresponding to the cusps in the sense of subsections 2.2.1 
and 2.2.2). Very often we also call representatives cusps.
 
As before, we would like to study cusps by transforming 
$\zeta\in \mathbb P^1(K)$ to $\infty$, and hence want to 
assume without loss of generality that $\zeta=\infty$ 
in our discussion. For this becoming posible, we are then
supposed to be able to find, for $\zeta:=\left[\begin{matrix} 
\alpha\\ \beta\end{matrix}\right]\in \mathbb P^1(K)$
 an element $M:=M_\zeta:=\left(\begin{matrix} \alpha&\alpha^*\\ \beta&\beta^*
\end{matrix}\right)\in SL(2,K)$, since then 
 it is clear that $\left(\begin{matrix} \alpha&\alpha^*\\ \beta&\beta^*
\end{matrix}\right)\cdot
 \left[\begin{matrix} 1\\ 0\end{matrix}\right]=\left[\begin{matrix} 
\alpha\\ \beta\end{matrix}\right]$, that is,
$M\cdot\infty=\zeta$. This is clearly possible, because if we set 
$\frak c:=\mathcal O_K\cdot \alpha+\mathcal O_K\cdot\beta$ to be the
 fractional ideal generated by $\alpha$ and $\beta$, then 
$1\in\mathcal O_K=\frak c\cdot\frak c^{-1}=\alpha
 \frak c^{-1}+\beta\frak c^{-1}$. Therefore, there exist 
$\alpha^*,\ \beta^*\in\frak c^{-1}\subset K$ such that 
$\alpha\beta^*-\alpha^*\beta=1$.
\vskip 0.30cm
\noindent
{\bf{\Large Theorem}.} ({\bf Cusp and Ideal Class Correspondence}) 
{\it There is a natural bijection between the ideal class group 
$CL(K)$ of $K$ and the cusps $\mathcal C_\Gamma$ of 
$\Gamma=SL(\mathcal O_K\oplus \frak a)$ 
acting on ${\mathcal H}^{r_1}\times
{\mathbb H}^{r_2}$ given by
$$\mathcal C_\Gamma\to CL(K),\qquad \left[\begin{matrix} \alpha\\ 
\beta\end{matrix}\right]\mapsto 
\Big[\mathcal O_K\,\alpha+\frak a\,\beta\Big].$$}

This type of results are rooted back to Maa$\beta$. 
But we here give a proof following Siegel as presented in 
A) above, while we reminder the reader that 
our case at hand is  more complicated.

Let $h=h_K$ denote the class number of $K$. Choose fixed integral 
$\mathcal O_K$-ideals $\frak a_1,\cdots,
\frak a_h$ representing the ideal class group $CL(K)$. We want 
to show that the elements of $\mathbb P^1(K)$
are divided into $h$ equivalence classes by the action of 
$\gamma=\left(\begin{matrix} a&b\\ c&d\end{matrix}
 \right)\in SL(\mathcal O_K\oplus \frak a)$ on 
$P=\left[\begin{matrix} p\\ s\end{matrix}\right]\in 
 \mathbb P^1(K)$ defined by $\gamma\cdot P=
\left[\begin{matrix} ap+bs\\ cp+ds\end{matrix}\right].$

Let then $P=\left[\begin{matrix} p\\ s\end{matrix}\right]$ be a fixed point
in $\mathbb P^1(K)$ with $p,s\in K$. Define $\pi(P)$ to be the ideal 
class associated to the fractiona ideal $\mathcal O_K\cdot p+\frak a\cdot s$.
\vskip 0.20cm
\noindent
{\bf {\large Claim}}. (1) {\it $\pi: \mathbb P^1(K)\to CL(K)$ is well-defined.}

\noindent
(2) {\it $\pi$ factors through the orbit space 
$SL(\mathcal O_K\oplus \frak a)\backslash \mathbb P^1(K)$.}

\noindent
Proof. (1) Indeed, if $P=\left[\begin{matrix} p_1\\ s_1
\end{matrix}\right]=\left[\begin{matrix} p_2\\ s_2
\end{matrix}\right]$, then, as ideal classes,
$$\Big[\mathcal O_K\cdot p_1+\frak a\cdot s_1\Big]
=\Big[\frac{s_2}{s_1}(\mathcal O_K\cdot p_1+
\frak a\cdot s_1)\Big]
=\Big[\mathcal O_K\cdot s_2\cdot\frac{p_1}{s_1}+\frak a\cdot s_2\Big]=
\Big[\mathcal O_K\cdot p_2+\frak a\cdot s_2\Big].$$ Here, we use 
$[\ ]$ to denote an ideal class.

\noindent
(2) For $\gamma=\left(\begin{matrix} a&b\\ c&d\end{matrix}\right)\in 
SL(\mathcal O_K\oplus \frak a)$,
$\pi\big(\gamma\cdot P\big)=\pi\bigg(\left[\begin{matrix} ap+bs\\ cp+
ds\end{matrix}\right]\bigg)=\Big[\mathcal O_K\cdot (ap+bs) 
+\frak a\cdot(cp+ds)\Big].$ But, by 
definition, $a,\,d\in \mathcal O_K,\ b\in\frak a,\ c\in \frak a^{-1}$. 
Hence,
$$\begin{aligned}
~&\mathcal O_K\cdot(ap+bs)+\frak a\cdot(cp+ds)=
(ap)\cdot\mathcal O_K+(bs)\cdot\mathcal O_K
+(cp)\cdot\frak a+(ds)\cdot \frak a\\
&\subset p\cdot\mathcal O_K+s\cdot\frak a
+p\cdot(\frak a^{-1}\cdot\frak a)+s\cdot\frak a
=p\cdot\mathcal O_K+s\cdot\frak a.\end{aligned}$$ On the other hand,
the inverse inclusion must hold as well because the determinant 
of $\gamma$ is one. Therefore, we have 
$$\pi\big(\gamma\cdot P\big)=\Big[\mathcal O_K\cdot (ap+bs) 
+\frak a\cdot(cp+ds)\Big]=\Big[\mathcal O_K\cdot p +\frak a\cdot s\Big].$$ 
This completes the proof of the claim.
\vskip 0.20cm
Consequently, we get a well-defined map 
$$\Pi:\ SL(\mathcal O_K\oplus \frak a)\Big\backslash \mathbb P^1(K)\to 
CL(K),\qquad \left[\begin{matrix} p\\ s\end{matrix}\right]\mapsto 
[\mathcal O_K\cdot p+\frak a\cdot s].$$
We want to show that $\Pi$ is a bijection. 

To start with, let us first show that $\Pi$ is surjective. 
This then is a direct consequence of the following
\vskip 0.20cm
\noindent
{\bf{\large Lemma}}. {\it For any two fractional $\mathcal O_K$-ideals 
$\frak a,\frak b$, there exist elements $\alpha,\beta\in
 K$ such that $\mathcal O_K\cdot \alpha+\frak a\cdot \beta=\frak b$.}

\noindent
Proof. Recall the following generalization of the 
classical Chinese Reminder Theorem for $\mathbb Z$ to 
all Dedekind domains.
\vskip 0.20cm
\noindent
{\bf {\Large Chinese Reminder Theorem}.} 
{\it Let $\frak p_j$ for $j=1,\cdots, s$ 
denote distinct prime ideals of $\mathcal O_K$, 
and let $e_j$  for $j=1,\cdots, s$ be positive integers. Then 
the map given by the product of the quotient 
maps $f:\mathcal O_K\to\prod_{j=1}^s\mathcal O_K\big/\frak p_j^{e_j}$ 
yields an isomorphism of rings $
\mathcal O_K\Big/\prod_{j=1}^s \frak p_j^{e_j}\simeq
\prod_{j=1}^s\mathcal O_K\Big/\frak p_j^{e_j}.$}

In terms of congruence, this means that given $x_j\in \mathcal O_K$ 
for $j=1,\cdots,s$, there exists
 $x\in\mathcal O_K$ such that $x\equiv x_j$ mod $\frak p_j^{e_j}$; 
and moreover, this uniquely determines the class of $x$ mod
$\prod_{j=1}^s \frak p_j^{e_j}$.
\vskip 0.20cm
With this in mind, let us go back to the proof of the lemma.
We can and hence now assume that both $\frak a$ and $\frak b$ are 
integral. 

\noindent
(First we may assume that $\frak b$ is an integral 
$\mathcal O_K$-ideal. Indeed, there exist $b\in K$ and an
 integral $\mathcal O_K$-ideal $\frak b'$ such that 
$\frak b=b\cdot\frak b'$. Therefore, if there exist $\alpha',\,\beta'$ 
 such that $\mathcal O_K\cdot \alpha'+\frak a\cdot 
\beta'=\frak b'$, then
$\frak b=b\cdot\frak b'=b\cdot(\mathcal O_K\cdot \alpha'+
\frak a\cdot \beta')=\mathcal O_K\cdot (b\alpha')+\frak a
\cdot (b\beta').$ That is to say, $\alpha=b\alpha'$ and 
$\beta=b\beta'$ will do the job.
Then we may further assume that $\frak a$ is integral. 
Indeed, there exists an $a\in K^*$ such that $a\cdot\frak a
=\frak a'$ is integral. Thus if there exist 
$\alpha',\,\beta'\in K$ such that $\mathcal O_K\cdot \alpha'+\frak a'
\cdot \beta'=\frak b$. Then 
$$\begin{aligned}\frak b=&\mathcal O_K\cdot \alpha'+
\frak a'\cdot \beta'=\mathcal O_K\cdot \alpha'+
\frak a'\cdot (a\cdot a^{-1}) \beta'\\
=&\mathcal O_K\cdot \alpha'+
(\frak a'\cdot a)\cdot (a^{-1} \beta')=\mathcal O_K\cdot 
\alpha'+\frak a\cdot (a^{-1} \beta').\end{aligned}$$ 
That is to say, this time, $\alpha=\alpha'$ and 
$\beta=a^{-1}\beta'$ do the job.)

We want to find $\alpha,\,
\beta\in K$ such that $\mathcal O_K\cdot \alpha+\frak a\cdot 
\beta=\frak b$. (Clearly, if done,
then $\alpha,\,\beta$ cannot be both zero at the same time, 
hence define a point 
$\left[\begin{matrix} \alpha\\ \beta\end{matrix}\right]\in\mathbb P^1(K)$.
 Furthermore, we get $\Pi\bigg(\left[\begin{matrix} \alpha\\ 
\beta\end{matrix}\right]\bigg)=\frak b$ as desired.)

Choose now $\beta\in\frak a^{-1}\frak b\backslash\{0\}$ so 
that $\frak a\cdot\beta\subset \frak a\cdot\frak a^{-1}
\frak b\subset\mathcal O_K\cdot \frak b=\frak b$.
Therefore, by the unique factorization theorem of integral 
$\mathcal O_F$-ideals into product of prime ideals,
 we can assume that $$\frak b=\prod_{i=1}^l\frak p_i^{n_i}
\supset\frak a\cdot\beta=\prod_{i=1}^l\frak p_i^{m_i}$$ 
 with $m_i\geq n_i\geq 0,\ 1\leq i\leq l$. Now choose 
$b_i\in \frak p_i^{n_i}\Big/\frak p_i^{n_i+1}$ for all $i=1,
 \cdots, l$. By the Chinese Reminder Theorem just cited above, 
there exists an element $\alpha\in\mathcal O_K$
 such that $\alpha\equiv b_i$ mod $\frak p_i^{n_i+1}$. Since, in terms of
local orders at $\frak p_i$,
$\nu_{\frak p_i}(\alpha)=\nu_{\frak p_i}(\frak b)$ for
  each $i$, we know that $\alpha\in\frak b$. Thus
 $\mathcal O_K\cdot\alpha+\frak a\cdot \beta\subset\frak b$.

On the other hand, if $\frak p$ is a prime ideal of 
$\mathcal O_K$ which does not lie in the set $\big\{\frak p_1,
\cdots,\frak p_l\big\}$, then $0=\nu_{\frak p}(\beta)=\nu_{\frak p}(\frak b)$. 
Thus we have shown that for all primes 
$\frak p$ of $\mathcal O_K$, $\ \nu_\frak p(\mathcal O_K\cdot\alpha
+\frak a\cdot \beta)=\inf\Big\{\nu_\frak p
(\mathcal O_K\cdot\alpha),\nu_\frak p(\frak a\cdot \beta)\Big\}=
\nu_\frak p(\frak b).$ Therefore, 
$\mathcal O_K\cdot\alpha+\frak a\cdot \beta=\frak b$. This completes the proof
of the lemma and hence the theorem.
\vskip 0.20cm
\noindent
{\bf Reamrks}. (1) Taking $\frak a=\mathcal O_K$, we in particular 
see that any fractional $\mathcal O_K$-ideal 
is generated by at most two elements, a simple beautiful fact, used
many times in Siegel's arguments copied in A), that should be included in all 
standard textbook in Algebraic Number Theory. The reader may find it in [FT], 
whose proof we followed in our discussion above.

\noindent
(2) We would like to reminder the reader that during this 
process of studying non-abelian zeta functions for 
number fields, all basic facts, not only the
finiteness results on ideal class group and units, but the Chinese Reminder 
Theorem are used. Is not it beautiful and wonderful?!
\vskip 0.20cm
With this being done, we are left with the injectivity of $\Pi$. 
For this, we use the trick of Siegel in A), following 
the presentation of Terras [Te].

In order to establish the injectivity, one may probably first think 
of the following argument. Suppose that 
$\Pi\bigg(\left[\begin{matrix} p_1\\ s_1\end{matrix}\right]\bigg)=k
\cdot \Pi\bigg(\left[\begin{matrix} p_2\\ s_2\end{matrix}\right]\bigg)$ for
some $k\in K$. Clearly,  we may and hence assume that $p_i, s_i$ 
are all in $\mathcal O_K$ and that $\frak a$
 is integral. Then if writing $k=\frac {\omega}{\tau}$ with 
$\omega,\tau\in\mathcal O_K$,  we see that
$\tau(ap_1+bs_1)=\omega p_2$ and $\tau(cp_1+ds_1)=\omega s_2$ 
for some $\gamma=\left(\begin{matrix} a&b\\ 
c&d\end{matrix}\right)\in GL(\mathcal O_K\oplus\frak a).$ It follows 
that $\gamma \cdot \left[\begin{matrix} p_1\\ 
s_1\end{matrix}\right]=\left[\begin{matrix} ap_1+bs_1\\ cp_1+ds_1
\end{matrix}\right]=\left[\begin{matrix} p_2\\ 
s_2\end{matrix}\right]$.
 
Surely, this says that $\left[\begin{matrix} p_1\\ s_1\end{matrix}\right]$ and
 $\left[\begin{matrix} p_2\\ s_2
\end{matrix}\right]$ are indeed equivalent modulo 
$GL(\mathcal O_K\oplus\frak a)$. But  we 
need to know that they are equivalent modulo 
$SL(\mathcal O_K\oplus\frak a)$. Unfortunately, the difference between special
 and general linear group over 
$\mathcal O_K\oplus\frak a$ has been seen to be possibly very large due to 
 the presence of units.

Since our Mitsubishi Fuso is running with fire and we 
like the brand Mitsubishi, so we make the following 
changes  using the Siegel Model.

Take $\zeta_1:=\left[\begin{matrix} p_1\\ s_1\end{matrix}\right]$ 
and $\zeta_2:=\left[\begin{matrix} p_2\\ s_2
\end{matrix}\right]$ in $\mathbb P^1(K)$ with $p_i,s_i$ in 
$\mathcal O_K$. Then by the discussion just above 
 the theorem, there exist 
$M_1:=\left(\begin{matrix} p_1&p_1^*\\ s_1&s_1^*\end{matrix}\right)
$ and $M_2:=\left(\begin{matrix} p_2&p_2^*\\ 
s_2&s_2^*\end{matrix}\right)$ in $SL(2,K)$ such that
$M_1\cdot\infty=\zeta_1,\ M_2\cdot \infty=\zeta_2$. 
Consequently, $(M_1\cdot M_2^{-1})\zeta_2=\zeta_1$. In 
other words,
$\left[\begin{matrix} p_1\\ s_1\end{matrix}\right]=
\Big(M_1\cdot M_2^{-1}\Big)\cdot\left[\begin{matrix} p_2\\ s_2\end{matrix}
\right]$. Thus by the fact that $p_i,s_i$ are all 
$\mathcal O_K$-integers, easily from the  discussion in A),
we have $M_1\cdot M_2^{-1}\in GL(2,\mathcal O_K\oplus a)$ by writing down
all the entries explicitly. 
Clearly, by definition, $M_1\cdot M_2^{-1}\in SL(2,K)$ 
as well. Hence $M_1\cdot M_2^{-1}\in 
GL(2,\mathcal O_K\oplus a)\cap SL(2,K)=SL(2,\mathcal O_K\oplus a)$. This 
completes the proof.
\vskip 0.30cm
In summary, what we have just established is the following bijection
$$\Pi:SL(2,\mathcal O_K\oplus a)\Big\backslash \mathbb P^1(K)
\simeq CL(K),\qquad \left[\begin{matrix} \alpha\\ \beta
\end{matrix}\right]\mapsto 
\Big[\mathcal O_K\alpha+\frak a\beta:=\frak b\Big].$$
Easily, one checks that the inverse map $\Pi^{-1}$ is given as follows: 
For $\frak b$, choose $\alpha_\frak b,
\beta_\frak b\in K$ such that $\mathcal O_K\cdot\alpha_\frak b+
\frak a\cdot\beta_\frak b=\frak b$; With this, then $\Pi^{-1}([\frak b])$ is 
simply the class of the point 
$\left[\begin{matrix} \alpha_\frak b\\ 
\beta_\frak b\end{matrix}\right]$ in 
$SL(2,\mathcal O_K\oplus a)\Big\backslash \mathbb P^1(K)$. Moreover, there 
always exists $M_{ \left[\begin{matrix} \alpha\\ \beta\end{matrix}\right]}
:= \left(\begin{matrix} 
\alpha&\alpha^*\\ \beta&\beta^*\end{matrix}\right)\in SL(2,K)$ such that 
$M_{\left[\begin{matrix} \alpha\\ 
\beta\end{matrix}\right]}\cdot\infty=
\left[\begin{matrix} \alpha\\ \beta\end{matrix}\right].$

We end this discussion on cusps by mentioning that

\noindent
(1) There is a much finer choice for the matrix 
$M_{\eta}$ for $\eta= \left[\begin{matrix} \alpha\\ \beta\end{matrix}\right]$, 
which will be used in the discussion of fundamental domain and Fourier 
expansion of Eisenstein series later.

\noindent
(2) The reason why we call points in $\mathbb P^1(K)$ cusps will become clear 
after we discuss the stablizer groups of 
them: Just as in the case for $\mathbb H$, 
there is a rank one $\mathcal O_K$ lattice involved.

\section{Stablizer Groups of Cusps}

\subsection{Upper Half Plane}
As said in the previous section,   if  
$\zeta\in \mathbb P^1(\mathbb R)$ is a cusp of a Fuchsian group 
$\Gamma\subset SL(2,\mathbb R)$ of first kind, then there exists
a matrix $A:=A_\zeta\in SL(2,\mathbb R)$ such that 
$A\cdot\infty=\zeta$. Moreover, we have 
$$\Gamma_\zeta=\Gamma\cap \Big(A\cdot  B(\mathbb R)\cdot A^{-1}\Big)
=\Bigg\{A_\zeta\cdot\left(\begin{matrix} 
1&m\cdot d\\ 0&1\end{matrix}\right)\cdot A_\zeta^{-1}:
m\in\mathbb Z\Bigg\}$$ for a certain $d>0$. In particular, a real rank 1
lattice $d\mathbb Z$ is naturally associated.
Consequently, the {\it fundamental domain} of the action of $\Gamma_\zeta$
on $\mathcal H$ is given by $$A_\zeta\cdot\Big\{z=x+iy\in\mathcal H:
-\frac{d}{2}\leq x\leq \frac{d}{2}\Big\}=A_\zeta\cdot
\Big\{z=x+iy\in\mathbb C:y>0,\ 
-\frac{d}{2}\leq x\leq \frac{d}{2}\Big\}.$$
Due to this, without loss of generality, usually we assume in addition   
that $\Gamma$ is reduced at all cusps, that is, for all cusps, the above
constant $d=1$.

\subsection{Upper Half Space}
Recall that from the discussion in the previous section, 
$\zeta\in\mathbb P^1(\mathbb C)$ is a cusp for a discrete subgroup 
$\Gamma\subset SL(2,\mathbb C)$ if and only if the (special) 
stablizer group $\Gamma_\zeta'$ is conjugate to a lattice $D$ 
in $N(\mathbb C)$. (Hence $D$ is isomorphic to a full
lattice in $\mathbb C$.)
As such, a fundamental domain for the action of $\Gamma_\zeta'$ on $\mathbb H$
is given by the above lattice up to the action of units. 
That is, $A_\zeta\cdot \{P=z+rj\in\mathbb H:z\in \mathbb C/D\}$ with 
$A_\zeta\in SL(2,\mathbb C)$ such that $A\cdot\infty=\zeta$.

\subsection{Rank Two $\mathcal O_K$-Lattices}

\subsubsection{A.Totally Real Fields}

\subsubsection{A.1. Stablizer Groups for Cusps}

As above, we follow Siegel [S] to give the presentation. So the same
notations are used as in 2.2.3.A.

As usual, let for $\lambda\in\mathbb P^1(K),\ \Gamma_\lambda$ denote 
the stablizer group of $\lambda$, i.e., the one consisting of Hilbert 
modular substitutions $z\mapsto\frac{\alpha z+\beta}{\gamma z+\delta}$ such 
that $\frac{\alpha \lambda+\beta}{\gamma \lambda+\delta}=\lambda.$
For our use later, we need to determine $\Gamma_\lambda$ explcitly. It 
clearly 
suffices to find $\Gamma_{\lambda_i}$ for $i=1,2,\ldots,h$ since 
$\lambda=(\lambda_i)_M$ for some $M=M_i\in\Gamma$ and $\Gamma_\lambda=
M_i\Gamma_{\lambda_i}M_i^{-1}$. Write $M=\begin{pmatrix}\alpha &\beta\\
\gamma &\delta\end{pmatrix}\in \Gamma$, then 
$\frac{\alpha \lambda_i+\beta}{\gamma \lambda_i+\delta}=\lambda_i.$
We have $\Big(\alpha\rho_i+\beta\sigma_i,\gamma\rho_i+\delta\sigma_i\Big)=
\frak a_i=\Big(\rho_i,\sigma_i\Big)$. Here as above, $\lambda_i=\frac{\rho_i}
{\sigma_i},\ \lambda=\frac{\rho}{\sigma}.$ And since 
$(\alpha\rho_i+\beta\sigma_i)
\sigma_i=(\gamma\rho_i+\delta\sigma_i)\rho_i$, we have
$$\frac{(\alpha\rho_i+\beta\sigma_i)}{\frak a_i}\cdot \frac{(\sigma_i)}
{\frak a_i}=\frac{(\gamma\rho_i+\delta\sigma_i)}{\frak a_i}\cdot
\frac{(\rho_i)}{\frak a_i}.$$ As before, $$\frac{(\alpha\rho_i+
\beta\sigma_i)}{\frak a_i}=\frac{(\rho_i)}{\frak a_i},
\qquad\mathrm{and}\qquad \frac{(\gamma\rho_i+\delta\sigma_i)}{\frak a_i}
= \frac{(\sigma_i)}{\frak a_i}.$$ Hence for a unit $\varepsilon$ in $K$, 
we have $$\alpha\rho_i+\beta\sigma_i=\varepsilon \rho_i,\qquad \gamma\rho_i+
\delta\sigma_i=\varepsilon\sigma_i.$$ This means that, with 
$A_i\cdot \infty=\lambda_i,$ 
$$MA_i=\begin{pmatrix}\rho_i&\xi_i^*\\ \sigma_i&\eta_i^*\end{pmatrix}
\begin{pmatrix} \varepsilon&0\\ 0&\varepsilon\end{pmatrix}$$ where 
$\xi_i^*=(\alpha\xi_i+\beta\eta_i)\varepsilon$ and $\eta_i^*=
(\gamma\xi_i+\delta\eta_i)\varepsilon$ lie in $\frak a_i^{-1}$.
 Further since $\rho_i\eta_i-\sigma_i\xi_i=\rho_i\eta_i^*-
\sigma_i\xi_i^*=1$, we have $$\rho_i(\eta_i^*-\eta_i)=
\sigma_i(\xi^*-\xi)$$ i.e.,
$$\frac{(\rho_i)}{\frak a_i}\cdot\frac{(\eta_i^*-\eta_i)}{\frak a_i^{-1}}
=\frac{(\sigma_i)}{\frak a_i}\cdot\frac{(\xi^*-\xi)}{\frak a_i^{-1}}.$$
Again, since $\frac{(\rho_i)}{\frak a_i}$ is coprime to 
$\frac{(\sigma_i)}{\frak a_i}$, we see that $\frac{(\rho_i)}
{\frak a_i}$ divides $\frac{(\xi^*-\xi)}{\frak a_i^{-1}}$, i.e., 
$(\xi^*-\xi)=\frak a_i^{-2}\frak b(\rho_i)$ for an integral ideal 
$\frak b$. In other words, $$\xi_i^*=\xi_i+\rho_i\zeta,\qquad\zeta\in
\frak a_i^{-2}.$$ As  a result, we obtain also $\eta_i^*=
\eta_i+\sigma_i\zeta.$ We now observe that $$MA_i=A_i
\begin{pmatrix} \varepsilon&\zeta\varepsilon^{-1}\\ 0&\varepsilon^{-1}
\end{pmatrix}.$$ Therefore we have the following
\vskip 0.20cm
\noindent
{\bf{\large Lemma}.}  {\it The stablizer group $\Gamma_{\lambda_i}$ consists 
precisley of the modular 
substitutions $z\mapsto z_M$, where $M=A_i\begin{pmatrix} 
\varepsilon&\zeta\\ 0&\varepsilon^{-1}\end{pmatrix}A_i^{-1}\in M$ 
with $\zeta\in\frak a_i^{-2}$ and $\varepsilon$ being any unit in $K$.}
 
\subsubsection{A.2. Actions of $\Gamma$}

Let $\mathcal H_n:=\mathcal H^n$ denote the product of $n$ copies of the
upper half-plane, namely the set of $z=(z_1,\ldots,z_n)$ with 
$z_j=x_j+iy_j,\, y_j>0$. 
The Hilbert modular group $\Gamma$ has a representation as a group of 
analytic homeomorphisms $z\mapsto \frac{\alpha z+\beta}{\gamma z+\delta}$ 
of $\mathcal H_n$ onto itself. In the following, we shall freely identify, say,
 $M=\begin{pmatrix}\alpha&\beta\\ \gamma&\delta\end{pmatrix}\in \Gamma$ with 
the modular substitution $z\mapsto z_M$ in $\Gamma$ and speak of $M$ 
belonging to $\Gamma$, without risk of confusion.

The points $\lambda=(\lambda^{(1)},\ldots, \lambda^{(n)})$ for 
$\lambda\in\mathbb P^1(K)$ lie on the boundary of $\mathcal H_n$ and are 
called (parabolic) {\it cusps} of $\mathcal H_n$ (since
their stablizer groups, up to the units, are isomorphic to a rank one 
$\mathcal O_K$-lattice).  By the discussion in 2.2.3.A, there 
exist $h$ cusps, $$\lambda_1=(\infty,\ldots,\infty),\ \lambda_2=
(\lambda_2^{(1)},\ \ldots,\ \lambda_2^{(n)}),\ldots, \lambda_h=
(\lambda_h^{(1)},\ldots, \lambda_h^{(n)})$$ which are not equivalent with 
respect to $\Gamma$, and any other (parabolic) cusp of $\mathcal H_n$ is 
equivalent to exactly one of them. For our own convenience,
 $\lambda_1,\ldots,\lambda_h$ are called the {\it base cusps}.

If $V\subset \mathbb C^n$, then for given $M\in\Gamma$, $V_M$ shall 
denote the set of all $z_M$ for $z\in V$. For any $z\in\mathcal H_n$, 
$\Gamma_z$ shall stand for the isotropy group, 
(or the same the stablizer group,) of $z$ in $\Gamma$, namely the group 
of $M\in \Gamma$ for which $z_M=z$. 
\vskip 0.20cm
\noindent
{\bf{\large Lemma}.} {\it For any two compact sets $B,B'$ in $\mathcal H_n$, 
the number of $M\in\Gamma$ for which $B_M$ intersects $B'$ is finite.
In particular, $\Gamma_z$ is finite for all $z\in \mathcal H_n$.}

\noindent
Proof. Let $\Lambda$ denote the set of $M=\begin{pmatrix}\alpha&\beta\\ 
\gamma&\delta\end{pmatrix}\in \Gamma$ such that $B_M$ intersects $B'$, i.e., 
given $M\in\Lambda$, there exists $w=(w_1,\ldots,w_n)$ in $B$ such that 
$w_M=(w_1',\ldots, w_n')\in B'$. Let $w_j=u_j+iv_j$ and $w_j'=u_j'+iv_j'$;
 then $v_j'=\frac{v_j}{|\gamma^{(j)}w_j+\delta^{(j)}|^{2}}$. 
Since $B$ and $B'$ are compact, we see that $|\gamma^{(j)}w_j+\delta^{(j)}|<c,
\ \, j=1,2,\ldots, n$ for a constant $c$ depending only on $B$ and 
$B'$. But then it follows immediately that $\gamma,\delta$ belong to a 
finite set of integers in $K$. Let $I=\begin{pmatrix} 0&1\\ -1&0
\end{pmatrix}$ and let $B^*$ and ${B'}^*$ denote the images of $B$ 
and $B'$ respectively under the modular substitution $z\mapsto z_I=-z^{-1}$. 
The images $B^*$ and ${B'}^*$ are again compact. Further noting that if
 $M=\begin{pmatrix}\alpha&\beta\\ \gamma&\delta\end{pmatrix}\in \Gamma$, 
then $I\cdot M\cdot I^t=M^{-t}=\begin{pmatrix}\delta&-\gamma\\ -\beta&\alpha
\end{pmatrix},$ we can show easily that $B_M\cap B'\not=\emptyset$ if 
and only if $B_{M^{-t}}^*\cap {B'}^*\not=\emptyset$. Now applying the same 
argument as above to the compact sets $B^*$ and ${B'}^*$, we may conclude 
that if for $M=\begin{pmatrix}\alpha&\beta\\ \gamma&\delta\end{pmatrix}\in
 \Gamma$, 
$B_{M^{-t}}^*\cap {B'}^*\not=\emptyset$, then $\alpha,\beta$ belong to a 
finite set of integers in $K$. As a result, we obtain finally that 
$\Lambda$ is finite.

In particular, taking a point $z\in\mathcal H_n$ instead of $B$ and $B'$, 
we deduce that 
$\Gamma_z$ is finite. This completes the proof of the lemma.
\vskip 0.20cm
We are now in a position to prove the following

\noindent
{\bf{\Large Proposition}.} {\it The Hilber modular group $\Gamma$ acts 
properly and  discontinuously on $\mathcal H_n$; in other words, for any 
$z\in\mathcal H_n$, there exists 
a neighborhood $V$ of $z$ such that only for finitely many $M\in \Gamma$, 
$V_M$ intersects $V$ and when $V_M\cap V\not=\emptyset$, then $M\in\Gamma_z$.
Consequently, if $z$ is not a fixed point of any $M$ in $\Gamma$ except 
the identity or in other words, if $\Gamma_z$ consists only of 
$\begin{pmatrix}1&0\\ 0&1\end{pmatrix}$, then there exists a neighborhood 
$V$ of $z$ such that for 
$M\not=\begin{pmatrix}1&0\\ 0&1\end{pmatrix}, V_M\cap V=\emptyset$.}

\noindent
Proof. Let $z\in\mathcal H_n$ and $V$ a neighborhood of $z$ such that the
 closure $\overline V$ of $V$ in $\mathcal H_n$ is compact. 
Taking $\overline V$ for $B$ and $B'$ in the above lemma, 
we note that only for finitely many 
$M\in\Gamma$, say $M_1,\ldots, M_r$, $\,V_M$ intersects $V$. Among these
 $M_i$, let $M_1,\ldots, M_s$ be exactly those which do not belong to 
$\Gamma_z$. Then we can find a neighborhood $W$ of $z$ such that 
$W_{M_i}\cap W=\emptyset,\ i=1,2,\ldots,s$. Now let $U=W\cap V$; then 
$U$ has already the proeprty that for $M\not= M_i,\ i=1,2,\ldots, r$, 
$\ U_M\cap U=\emptyset$. Further, $U_{M_i}\cap U=\emptyset$ for $i=1,2,
\ldots,s$. Thus $U$ satisfies the requirement of the Proposition.

\subsubsection{A.3. Fundamental Domain of $\Gamma_\lambda$ in $\mathcal H_n$}

Let $\lambda=\frac{\rho}{\sigma}$ be a cusp of $\mathcal H_n$ and 
$\frak a=(\rho,\sigma)$ be the integral ideal among $\frak a_1,\ldots,
\frak a_h$ which is associated with $\lambda$. Let $\alpha_1,\ldots,\alpha_n$ 
be a $\mathbb Z$-basis for $\frak a^{-2}$, and further associated with 
$\lambda$, let us choose a fixed $A=
\begin{pmatrix}\rho&\xi\\ \sigma&\eta\end{pmatrix}\in G$ with 
$\xi,\eta$ lying in $\frak a^{-1}$. Moreover, let $\varepsilon_1,\cdots,
\varepsilon_{n-1}$ be $n-1$ independent generators of the group of units 
(up to torsion) in $K$. As a first step towards constructing a fundamental 
domain for $\Gamma$ in $\mathcal H_n$, we shall introduce \lq local 
coordinates' relative to $\lambda$, for every point $z$ in $\mathcal H_n$ 
and then construct a fundamental domain $\mathcal D_\lambda$ for 
$\Gamma_\lambda$ in $\mathcal H_n$.

Let $z=(z_1,\ldots,z_n)$ be any point of $\mathcal H_n$ and 
$z_{A^{-1}}=(z_1^*,\ldots,z_n^*)$ with $z_j^*=x_j^*+iy_j^*$. Denoting
 $y_{A^{-1}}=(y^*_1,\ldots,y^*_n)$ by $y^*$, we define the local coordinates 
of $z$ relative to $\lambda$ by the $2n$ quantities 
$$\frac{1}{\sqrt{N(y^*)}},\ Y_1,\, \ldots,\ Y_{n-1},
\ \ X_1,\,\ldots,\ X_n$$
where $Y_1,\,\ldots,\, Y_{n-1},\ X_1,\,\ldots,\, X_n$  are uniquely determined
 by the linear equations 
$$\begin{aligned}Y_1\log\Big|\varepsilon_1^{(k)}\Big|+\ldots+Y_{n-1}
\log\Big|\varepsilon_{n-1}^{(k)}\Big|=&
\frac{1}{2}\log\Big(\frac{y_k^*}{\sqrt{N(y^*)}}\Big),
\qquad k=1,2,\ldots, n-1,\\
X_1\alpha_1^{(l)}+\ldots+X_n\alpha_n^{(l)}=&x_l^*,\hskip 3.0cm l=1,2,\ldots,n.\end{aligned}$$

The group $\Gamma_\lambda$ consists of all the modular substitutions
 of the form $z\mapsto z_M$ where $M=A\cdot T\cdot A^{-1}$ and 
$T=\begin{pmatrix}
\varepsilon&\zeta\varepsilon^{-1}\\ 0&\varepsilon^{-1}\end{pmatrix}$ 
with $\varepsilon$ a unit in $K$ and $\zeta\in\frak a^{-2}$. The 
transformation $z\mapsto z_M$ is equivalent to the transformation 
$z_{A^{-1}}\mapsto z_{TA^{-1}}=\varepsilon^2 z_{A^{-1}}+\zeta.$ 
It is easily verified that $\Gamma_\lambda$ is generated by the 
dilations $z_{A^{-1}}\mapsto \varepsilon^2_i z_{A^{-1}},\ i=1,2,\ldots, n-1$ 
and the translations $z_{A^{-1}}\mapsto z_{A^{-1}}+\alpha_j,\ i=1,2,\ldots, n$.

Let now $z\mapsto z_M$ with $M=A\cdot T\cdot A^{-1}$, 
$T=\begin{pmatrix}\varepsilon
&\zeta\varepsilon^{-1}\\ 0&\varepsilon^{-1}\end{pmatrix}$ be a modular 
substituation in $\Gamma_\lambda$ and let $\varepsilon=\pm 
\varepsilon_1^{k_1}\cdots\varepsilon_{n-1}^{k_{n-1}}$ and $\zeta=
m_1\alpha_1+\ldots+m_n\alpha_n$, where $k_1,\ldots, k_{n-1},\ m_1,
\ldots, m_n$ are rational integers. It is obvious that 

\noindent
a) The first 
coordinate $\frac{1}{\sqrt{N(y^*)}}$ of $z$ is preserved by the 
modular substitutation $z\mapsto z_M$, since, by definition,
$N(y_{A^{-1}M})=N(y_{TA^{-1}})=N(\varepsilon^2)N(y_{A^{-1}})=N(y_{A^{-1}})
=N(y^*).$ 

\noindent
b) For $Y_i$, the effect of the substitution 
$z\mapsto z_M$ on $Y_1,\ldots, Y_{n-1}$ 
is given by $$\Big(Y_1,Y_2,\ldots, Y_{n-1}\Big)\ \mapsto\ \Big(Y_1+k_1,Y_2+k_2,
\ldots, Y_{n-1}+k_{n-1}\Big)$$ as can be verified easily as well. 

\noindent
c) For $X_j$, there are two cases:

\noindent
i) If $\varepsilon^2=1$, then the effect of the substitution $z\mapsto z_M$ 
on $X_1,\ldots, X_n$ is again just a translation, $$\Big(X_1,\ldots, X_n\Big)
\ \mapsto\ \Big(X_1+m_1,\ldots, X_n+m_n\Big).$$

\noindent
ii) If $\varepsilon^2\not=1$, then the effect of the substitution 
$z\mapsto z_M$ on $X_1,\ldots, X_n$ is not merely a translation but
 an affine transformation given by, $$\Big(X_1,\ldots, X_n\Big)\ \mapsto\ 
\Big(X_1^*+m_1,\ldots, X_n^*+m_n\Big).$$
where $\Big(X_1^*,X_2^*,\ldots, X_n^*\Big)=\Big(X_1,X_2,\ldots, X_n\Big)
\cdot\Big(UVU^{-1}\Big),$
$U$ denotes the $n$-rowed square matrix $(\alpha_i^{(j)})$ and $V$ 
denotes the diagonal matrix $\Big((\varepsilon^{(1)})^2,\ldots,
(\varepsilon^{(n)})^2\Big)$.

We define a point $z\in\mathcal H_n$ to be {\it reduced} with respect to 
$\Gamma_\lambda$, if $$\begin{aligned}-\frac{1}{2}\leq& Y_i<\frac{1}{2},
\qquad i=1,2,\ldots, n-1,\\
-\frac{1}{2}\leq& X_j<\frac{1}{2},\qquad j=1,2,\ldots, n.\end{aligned}
\eqno(*)$$ 
It is first of all 
clear that for any $z\in\mathcal H_n$, there exists an $M\in\Gamma_\lambda$ 
such that the equivalent point $z_M$ is reduced with respect to 
$\Gamma_\lambda$. In fact, for $\varepsilon=\pm \varepsilon_1^{k_1}
\cdots\varepsilon_{n-1}^{k_{n-1}}$ and
$M_1=A\begin{pmatrix}\varepsilon&0\\ 0&\varepsilon^{-1}\end{pmatrix}A^{-1}$,
 the effect of the substitution $z\mapsto z_{M_1}$ is given by 
$$\Big(Y_1,Y_2,\ldots, Y_{n-1}\Big)\mapsto \Big(Y_1+k_1,Y_2+k_2,\ldots,
 Y_{n-1}+k_{n-1}\Big),$$  and hence, by choosing $k_1,\ldots, k_{n-1}$ 
properly, we could suppose that the coordinates $Y_1,\ldots, Y_{n-1}$ 
of $z_{M_1}$ satisfy $(*)$. Again, since, 
for $\zeta=m_1\alpha_1+\ldots+m_n\alpha_n$ and 
$M_2=A\begin{pmatrix}1&\zeta\\ 0&1\end{pmatrix}A^{-1}$, the effect of the 
substitution $z\mapsto z_{M_2}$  on the coordinates $X_1,\ldots,X_n$ of
 $z_{M_1}$ is given by
$$\Big(X_1,X_2,\ldots, X_{n}\Big)\mapsto \Big(X_1+m_1,X_2+m_2,\ldots,
 X_{n}+m_{n}\Big),$$
 we could suppsoe that for suitable $m_1,\ldots,m_n$, the 
coordinates $X_1,\ldots,X_n$ of $z_{M_2M_1}$ is reduced with erspect to 
$\Gamma_\lambda$. 

On the other hand, let $z=x+iy,\ w=u+iv\in\mathcal H_n$ 
be reduced and equivalent with respect to $\Gamma_\lambda$ and let 
the local coordinates of $z$ and $w$ relative to $\lambda$ be respectively
$$\frac{1}{\sqrt{N(y_{A^{-1}})}},\ Y_1,\,\ldots,
\ Y_{n-1},\ \ X_1,\,\ldots,\ X_n$$ and 
$$\frac{1}{\sqrt{N(v_{A^{-1}})}},\ Y_1^*,\,\ldots,\ Y_{n-1}^*,\ \ X_1^*,
\,\ldots,\ X_n^*.$$
Further let $z_{A^{-1}}=\varepsilon^2 w_{A^{-1}}+\zeta$. Then in 
view of the fact that $$Y_i^*\equiv Y_i\pmod {1},\qquad\mathrm{and}\quad 
-\frac{1}{2}\leq Y_i, Y_i^*<\frac{1}{2},\qquad i=1,2,\ldots, n-1,$$ 
we have first $Y_i^*=Y_i,\ i=1,2,
\ldots,n-1$ and hence $\varepsilon^2=1$. Again since we have 
$$X_j^*\equiv X_j\pmod {1},\qquad\mathrm{and}\quad 
-\frac{1}{2}\leq X_j, X_j^*<\frac{1}{2}, 
\qquad i=1,2,\ldots, n$$ we see that $X_j^*=X_j,\ j=1,2,\ldots, n,$ 
and hence $\zeta=0$. Thus $z_{A^{-1}}=w_{A^{-1}}$, i.e., $z=w$.

Denote by $\mathcal D_\lambda$ the set of 
$z\in\mathcal H_n$ whose 
local coordinates $Y_1,\,\ldots,\,Y_{n-1},\ X_1,\,\ldots,\,X_n$ satisfy (*).
Thenwe have the following
\vskip 0.20cm
\noindent
{\bf{\Large Proposition}.} {\it Any $z\in\mathcal H_n$ is equivalent with respect
 to $\Gamma_\lambda$ to a point of $\mathcal D_\lambda$, and no two 
distinct points of $\mathcal D_\lambda$ are equivalent with respect to 
$\Gamma_\lambda$. Consequently, $\mathcal D_\lambda$ is a fundamental domain 
for $\Gamma_\lambda$ in $\mathcal H_n$.}

For $n=1$, the fundamental domain for $\Gamma_\lambda$ in $\mathcal H$
 is just the vertical strip $-\frac{1}{2}\leq x^*<\frac{1}{2},\ y^*>0$ 
with reference to the coordinate $z_{A^{-1}}=x^*+iy^*$, as we know. 
Going back to 
the coordinate $z$, the vertical strips $x^*=\pm\frac{1}{2},\ y^*>0$ 
are mapped into semi-circles passing through $\lambda$ and orthogonal 
to the real axis.
\vskip 0.20cm
For later use, we end this subsection with the following

\noindent
{\bf{\large Lemma}.} {\it All points $z=x+iy\in\mathcal D_\lambda$ 
satisfying $c_1\leq N(y_{A^{-1}})\leq c_2$ lie in a compact set 
in $\mathcal H_n$, depending only on $c_1,c_2$ and on the choices of 
$\varepsilon_1,\ldots,\varepsilon_{n-1}$ and $\alpha_1,\ldots,\alpha_n$ 
in $K$.}

\noindent
Proof. As a matter of fact, from $(*)$, we 
know that  $\ \Big|\log\big(\frac{y_i^*}
{\sqrt{N(y^*)}}\big)\Big|\leq c_3$ for $i=1,2,\ldots, n-1$,
and $\big|x_j^*\big|\leq c_4$ for $j=1,2,\ldots,n$. 
Hence we have $c_5\leq \frac{y_i^*}{\sqrt{N(y^*)}}\leq c_6$ for 
$j=1,2,\ldots, n-1$. Since $c_1\leq N(y^*)\leq c_2$, we have then 
$$c_7\leq\frac{y_i^*}{\sqrt{N(y^*)}}\leq c_8,\qquad i=1,2,\ldots,n.$$
Consequently, we obtain  $$c_9\leq y_i^*\leq c_{10},\qquad 
 |x_j^*|\leq c_4\qquad i,j=1,2,\ldots, n$$ where $c_9,\,c_{10},\,c_4$ 
depend only on $c_1,\,c_2$ and the 
choices of $\varepsilon_1,\ldots,\varepsilon_{n-1}$ and $\alpha_1,
\ldots,\alpha_n$ in $K$. Thus $z_{A^{-1}}$ and therefore $z$ lies in a 
 compact set in $\mathcal H_n$, depending on $c_1,\ c_2$ and $K$.

\subsubsection{B. Rank Two $\mathcal O_K$-Lattices: General Number Fields}

We continue our study of rank two  $\mathcal O_K$-lattices here
aiming at a natural construction of a fundamental domain for 
$SL(\mathcal O_K\oplus\frak a)\Big\backslash 
\Big({\mathcal H}^{r_1}\times{\mathbb H}^{r_2}\Big)$. 
Following Siegel, we start with a  precise construction
of fundamental neighborhoods of cusps here. (The reader will see  the  
reason why we use Siegel's text above: Essentially,
Siegel's treatment works for general number fields as well. 
However, what we are dealing  is much  general and hence more 
complicated. By recalling 
Siegel's text, not only we show the reader how classics
are used in our study and what a kind of refinements or better 
further developments are needed, but how things are naturally 
arranged together in a beautiful way under non-abelian zeta functions.)
\vskip 0.30cm
Recall that under the Cusp-Ideal Class Correspondence, there are exactly $h$
inequivalence cusps $\eta_i,\, i=1,2,\ldots,h$. Moreover, if we write the cusp 
$\eta_i=\left[\begin{matrix} \alpha_i\\ \beta_i\end{matrix}\right]$ for 
suitable $\alpha_i,\,\beta_i\in K$, then the associated ideal class 
is exactly the one for the fractional ideal 
 $\mathcal O_K\alpha_i+\frak a\beta_i=:\frak b_i$. 
Denote the stablizer group of $\eta_i$ by 
$$\Gamma_{\eta_i}:=\Big\{\gamma\in 
 SL(\mathcal O_K\oplus\frak a):\gamma\eta_i=\eta_i\Big\},\qquad 
i=1,2,\ldots,h.$$ 
Quite often, we use $\eta$ as a running symbol for $\eta_i$. 

We want to see the structure of $\Gamma_\eta$. 
As usual, we first shift $\eta$ to $\infty$.
So choose $A=\left(\begin{matrix} \alpha&\alpha^*\\ 
\beta&\beta^*\end{matrix}\right)\in SL(2,K)$. Clearly 
$$A\cdot\infty=A\left[\begin{matrix} 1\\ 0\end{matrix}\right]
=\left[\begin{matrix} \alpha\\ \beta\end{matrix}\right].
$$ Consequently, $\Gamma_\eta=A\cdot\Gamma_\infty \cdot A^{-1}.$
 
Then, we further pin down the choice of $\alpha^*$ and $\beta^*$
appeared in $A$.
For this, we use a trick which according to Elstrodt roots backed to Hurwitz.
\vskip 0.20cm
\noindent
{\bf{\large Lemma}}. {\it Let $\alpha,\,\beta\in K$ such that $\mathcal O_K\alpha+
\frak a\beta=\frak b\not=\{0\}.$
Then there exist $\alpha^*,\beta^*\in K$ such that}

\noindent
(1) $\left(\begin{matrix} \alpha&\alpha^*\\ \beta&\beta^*
\end{matrix}\right)\in SL(2,K)$; {\it and}

\noindent
(2) $\mathcal O_K\beta^*+\frak a^{-1}\alpha^*=\frak b^{-1}.$

\noindent 
Proof. Note that  
 $$1\in \mathcal O_K=\frak b\cdot\frak b^{-1}=(\mathcal O_K\alpha+
\frak a\beta)\cdot \frak b^{-1}=\frak b^{-1}\cdot\alpha+
(\frak a\frak b^{-1})\cdot\beta.$$ As such, we can choose 
$\beta^*\in \frak b^{-1},\,\alpha^*\in \frak a
\frak b^{-1}$ such that $\alpha\beta^*-\beta\alpha^*=1$.
This gives (1). As for (2),  it suffices to show that
$$(\mathcal O_K\beta^*+\frak a^{-1}\alpha^*)\cdot
(\mathcal O_K\alpha+\frak a\beta)=\mathcal O_K.$$
One inclusion is clear. Indeed, by our construction, 
$1\in (\mathcal O_K\beta^*+\frak a^{-1}\alpha^*)
\cdot(\mathcal O_K\alpha+\frak a\beta),$
so $$(\mathcal O_K\beta^*+\frak a^{-1}\alpha^*)
\cdot(\mathcal O_K\alpha+\frak a\beta)\supset\mathcal O_K.$$
As for the inclusion in the other direction, we go as follows: 
Clearly, 
$$(\mathcal O_K\beta^*+\frak a^{-1}\alpha^*)\cdot(\mathcal O_K\alpha+
\frak a\beta)
=\mathcal O_K\cdot\big(\beta^*\alpha\big)+\frak a^{-1}\cdot
\big(\alpha^*\alpha\big)+\big(\frak a\frak a^{-1}\big)\cdot
\big(\alpha^*\beta\big)+\frak a\cdot\big(\beta\beta^*\big).$$
But, by definition, $\frak b=\mathcal O_K\alpha+\frak a\beta$, 
so $\alpha\in\frak b,\, \beta\in\frak a^{-1}\frak b$. This, 
together with $\beta^*\in \frak b^{-1},
\,\alpha^*\in \frak a\frak b^{-1}$, then gives
$$\begin{aligned}(\mathcal O_K\beta^*&+\frak a^{-1}\alpha^*)\cdot
(\mathcal O_K\alpha+\frak a\beta)\\
\subset &\mathcal O_K\cdot \big(\frak b^{-1}\cdot\frak b\big)
+\frak a^{-1}\cdot\big((\frak a\frak b^{-1})\cdot\frak b\big)+
\big(\frak a\frak a^{-1}\big)\cdot\big((\frak a\frak b^{-1})\cdot
(\frak a^{-1}\frak b)\big)+\frak a\cdot\big((\frak a^{-1}\frak b)
\frak b^{-1}\big)\\
=&\mathcal O_K.\end{aligned}$$ This completes the proof.

As a direct consequence, we have the following generalization for the 
structure of the stablizer $\Gamma_\eta$;
\vskip 0.20cm
\noindent 
{\bf \large Corollary}.  {\it With the same notation as above, $$A^{-1}\Gamma_\eta A=\Bigg\{\left(\begin{matrix} u&z\\ 
0&u^{-1}\end{matrix}\right)
:u\in U_K, z\in\frak a\frak b^{-2}\Bigg\}.$$ In particular, the associated 
\lq lattice' for the cusp $\eta$ is given by the fractional ideal $\frak a
\frak b^{-2}$.}
 
\noindent
Proof.  All elements in $A^{-1}\cdot\Gamma_\eta\cdot A$ fix $\infty$, 
so are given by upper triangle matrices.
With this observation, let us now show that $z\in\frak a
\frak b^{-2}$. This is easy. Indeed, by definition, 
$A^{-1}\cdot \Gamma_\eta\cdot A$ 
consists of elements in the form
$\left(\begin{matrix}\beta^*&-\alpha^*\\ -\beta&\alpha\end{matrix}\right)
\cdot\left(\begin{matrix} a&b\\ c&d\end{matrix}\right) 
\cdot \left(\begin{matrix}\alpha&\alpha^*\\ \beta&\beta^*
\end{matrix}\right)=:\left(\begin{matrix} a_{11}&a_{12}\\ 
a_{21}&a_{22}\end{matrix}\right)$ with $a_{21}=0$ 
and 
$$a_{12}=(a-d)\alpha^*\beta^*-c(\alpha^*)^2+b(\beta^*)^2.$$ 
Recall that
$\alpha\in\frak b,\, \beta\in\frak a^{-1}\frak b$ 
and $\beta^*\in \frak b^{-1},\, \alpha^*\in
 \frak a\frak b^{-1}$, and that for $\left(\begin{matrix} a&b\\ 
c&d\end{matrix}\right)\in SL(\mathcal O_K\oplus\frak a)$, 
 $a,\,d\in\mathcal O_K,\ b\in\frak a,\ 
 c\in\frak a^{-1}$, easily we have $$z=a_{12}\subset 
\mathcal O_K\cdot \big((\frak a\frak b^{-1})\cdot\frak b^{-1}\big)
 +\frak a^{-1}\cdot(\frak a\frak b^{-1})^2+\frak a\cdot 
(\frak b^{-1})^2=\frak a\frak b^{-2}.$$ as desired.
 
To complete the proof, we still need to show that $u$ is a unit. 
This may be done as follows. So assume, as we can, that 
$\eta=\left[\begin{matrix} \alpha\\ \beta\end{matrix}\right]$ with $\alpha,\,
\beta\in\mathcal O_K$.
Note that for $\gamma=\left(\begin{matrix} a&b\\ c&d\end{matrix}\right)\in 
SL(\mathcal O_K\oplus\frak a)$ such 
that $\gamma \cdot \left[\begin{matrix} \alpha\\ \beta\end{matrix}\right]
=\left[\begin{matrix} \alpha\\ \beta
\end{matrix}\right],$ we have $$(a\alpha+b\beta)\cdot \beta=
(c\alpha+d\beta)\cdot \alpha.$$ Here $a,\,d\in\mathcal O_K,\ b
\in\frak a,\ c\in\frak a^{-1},$ and  $\alpha\in\frak b,\, \beta\in
\frak a^{-1}\frak b$ with 
$\frak b=\mathcal O_K\alpha+\frak a\beta$.
Thus note that now the ideal generated by $(a\alpha+b\beta)
\beta=(c\alpha+d\beta)\alpha$ is included in
 $\frak a^{-1}\frak b^2$. So dividing by it, we have
$$\frac{(a\alpha+b\beta)}{\frak b}=\frac{(\alpha)}{\frak b}
\qquad\mathrm{and}\quad
\frac{(c\alpha+d\beta)}{\frak a^{-1}\frak b}=
\frac{(\beta)}{\frak a^{-1}\frak b}.\eqno(*)$$ 
On the other hand,
$$\left(\begin{matrix}\alpha&\alpha^*\\ 
\beta&\beta^*\end{matrix}\right)\left(\begin{matrix} u&z\\ 0&u^{-1}
\end{matrix}\right)=\left(\begin{matrix} a&b\\ 
c&d\end{matrix}\right)
\left(\begin{matrix}\alpha&\alpha^*\\ 
\beta&\beta^*\end{matrix}\right),$$ so
$$\left(\begin{matrix} u\alpha&*\\ 
u\beta&*\end{matrix}\right)=\left(\begin{matrix} a\alpha+b\beta&*\\ 
c\alpha+d\beta&*\end{matrix}\right).$$ 
Therefore, $$(u\alpha)=(a\alpha+b\beta),\qquad\mathrm{and}\quad (u\beta)=(c\alpha+d\beta).\eqno(**)$$ 

Clearly, now from the equalities (*) and (**),
as integral ideals $(u\alpha)=(\alpha),
 (u\beta)=(\beta).$ So $u\in U_K$ as desired. This completes the proof.
\vskip 0.30cm
Set now $\Gamma_\eta':=\bigg\{A\left(\begin{matrix} 1&z\\ 
0&1\end{matrix}\right)A^{-1}:z\in\frak a\frak b^{-2}\bigg\},$ 
Then $$\Gamma_\eta=\Gamma_\eta'\times \bigg\{A\left(\begin{matrix} u&0\\ 
0&u^{-1}\end{matrix}\right)A^{-1}:u\in U_K\bigg\}.$$
Note that also componentwisely, $\left(\begin{matrix} u&0\\ 
0&u^{-1}\end{matrix}\right)z=\frac{uz} {u^{-1}}=u^2z$. So, in 
practice, what we really get is the following decomposition
$$\Gamma_\eta=\Gamma_\eta'\times U_K^2$$ with 
$$U_K^2\simeq \bigg\{A\cdot\left(\begin{matrix} u&0\\ 0&u^{-1}\end{matrix}
\right)\cdot A^{-1}:u\in U_K\bigg\}\,\simeq\, 
\bigg\{A\left(\begin{matrix} 1&0\\ 
0&u^2\end{matrix}\right)A^{-1}:u\in U_K\bigg\}.$$ 

Now we are ready to proceed  a construction of a fundamental domain for 
the action of $\Gamma_\eta\subset SL(\mathcal O_K\oplus\frak a)$ on 
${\mathcal H}^{r_1}\times{\mathbb H}^{r_2}$. We follow 
[Ge], in which the field involved is assumed to be totally real, 
to proceed our presentation here. This is based on a construction of 
a fundamental domain for the action of $\Gamma_\infty$ on
 ${\mathcal H}^{r_1}\times{\mathbb H}^{r_2}$. More precisely, 
with an element 
$A=\left(\begin{matrix}\alpha&\alpha^*\\ \beta&\beta^*\end{matrix}\right)
\in SL(2,K)$ used above,

\noindent
i) $A\cdot\infty=\left[\begin{matrix}\alpha
\\ \beta\end{matrix}\right];$ and

\noindent
ii) The isotropy group of $\eta$ in $A^{-1}SL(\mathcal O_K\oplus\frak a)A$ is 
generated by translations $\bold\tau\mapsto\bold\tau+z$ 
with $z\in \frak a\frak b^{-2}$ and by dilations 
$\bold\tau\mapsto u\bold\tau$ where $u$ runs through the group
 $U_K^2$. 

\noindent
(Recall that here, as above, we use $A,\,\alpha,\,\beta,\,\frak b$ as
  running symbols for $A_i,\, \alpha_i,\,\beta_i,\,\frak b_i:=
\mathcal O_K\alpha_i+\frak a\beta_i,\ i=1,\cdots,h$.)

Consider then the map $$\begin{matrix}
\mathrm {ImJ}:&{\mathcal H}^{r_1}\times
{\mathbb H}^{r_2}&\to& \mathbb R_{>0}^{r_1+r_2},\\
&\bold\tau:=(z_1,\cdots,z_{r_1};P_1,\cdots,P_{r_2})&\mapsto&
 (\Im(z_1),\cdots,\Im(z_{r_1});J(P_1),\cdots,J(P_{r_2})),\end{matrix}$$ 
where if $z=x+iy\in{\mathcal H}$ resp. $P=z+rj\in{\mathbb H}$,
 we set $\Im(z)=y$ resp. $J(P)=r$. It induces a map 
$$\Big(A^{-1}\cdot\Gamma_\eta\cdot A\Big)\Big\backslash \Big({\mathcal H}^{r_1}
\times{\mathbb H}^{r_2}\Big)\to U_K^2\Big\backslash 
\mathbb R_{>0}^{r_1+r_2},$$ which exhibits $\Big(A^{-1}\cdot\Gamma_\eta\cdot A\Big)\Big\backslash \Big({\mathcal H}^{r_1}
\times{\mathbb H}^{r_2}\Big)$ as a torus bundle over $U_K^2\Big\backslash 
\mathbb R_{>0}^{r_1+r_2}$ with fiber the $n=r_1+2r_2$ 
dimensional torus $\Big(\mathbb R^{r_1}\times\mathbb C^{r_2}\Big)\Big/
\frak a\frak b^{-2}$. Having factored out the action of the 
translations, we only have to construct a fundamental domain 
for the action of $U_K^2$ on $\mathbb R_{>0}^{r_1+r_2}$. This is 
essentially the same as in I.1.7. We 
look first at the action of $U_K^2$ on the norm-one hypersurface
$\bold S:=\Big\{y\in \mathbb R_{>0}^{r_1+r_2}:N(y)=1\Big\}$. By taking 
logarithms, it is transformed bijectively into a trace-zero hyperplane
 which is isomorphic to the space $\mathbb R^{r_1+r_2-1}$
$$\begin{aligned}\bold S&\buildrel\log\over\to  \mathbb R^{r_1+r_2-1}:=
\Big\{(a_1,\cdots a_{r_1+r_2})\in \mathbb R^{r_1+r_2}:\sum a_i=0\Big\},\\
y&\mapsto\qquad \Big(\log y_1,\cdots,\log y_{r_1+r_2}\Big),\end{aligned}$$ 
where the action of $U_K^2$ on $\bold S$ is carried out 
over an action on $\mathbb R^{r_1+r_2-1}$ by translations:
 $a_i\mapsto a_i+\log\varepsilon^{(i)}$. By Dirichlet's 
Unit Theorem, the logarithm transforms $U_K^2$ into a 
lattice in $\mathbb R^{r_1+r_2-1}$. The exponential map transforms 
a fundamental domain, e.g., a fundamental parallelopiped, for
 this action back into a fundamental domain $\bold S_{U_K^2}$
 for the action of $U_K^2$ on $\bold S$. The cone over 
$\bold S_{U_K^2}\ $, that is,  $\ \mathbb R_{>0}\cdot
\bold S_{U_K^2}\subset \mathbb R_{>0}^{r_1+r_2}$, is a 
fundamental domain for the action of $U_K^2$ on   
$\mathbb R_{>0}^{r_1+r_2}$. If we denote by $\mathcal T$ a fundamental 
domain for the action of the translations by elements of
 $\frak a\frak b^{-2}$ on $\mathbb R^{r_1}\times\mathbb C^{r_2}$,
and set $$\mathrm{ReZ}\,\Big(z_1,\cdots, z_{r_1};
P_1,\cdots,P_{r_2}\Big):=\Big(\Re(z_1),\cdots,\Re(z_{r_1});Z(P_1),
\cdots,Z(P_{r_2})\Big)$$ with
$\Re(z):=x$ resp. $Z(P):=z$ if $z=x+iy\in{\mathcal H}$ resp. 
$P=z+rj\in{\mathbb H}$, 
then what we have just said proves the following
\vskip 0.30cm
\noindent
{\bf{\bf {\Large Theorem.}}} {\it With the same notation as above,}
$$\bold E:=\Big\{\bold\tau\in {\mathcal H}^{r_1}\times
{\mathbb H}^{r_2}:
\mathrm{ReZ}\,(\bold\tau)\in {\mathcal T},\  
\mathrm{ImJ}\,(\bold\tau)\in  \mathbb R_{>0}\cdot
\bold S_{U_K^2}\Big\}$$ 
{\it is a fundamental domain for the action of 
$A^{-1}\Gamma_\eta A$ on ${\mathcal H}^{r_1}\times
{\mathbb H}^{r_2}$.}

For later use, we also set $\mathcal F_\eta:=A_\eta^{-1}\cdot \bold E$.
\vskip 0.20cm
Surely, as in Siegel's discussion, we may introduce $Y_1,\,\cdots,\, Y_{n-1},
\ X_1,\,\cdots,\, X_n$ together with a \lq reduced norm' of $\tau$ 
to precisely written done this fundamental domain in a simple form. 
We leave this for the time being to the reader, while pointing out that 
such a description will play a key role in our application of 
Rankin-Selberg \& Zagier method.

\section{Fundamental Domain}

\subsection{Upper Half Plane}

Let $\Gamma$ be a Fuchsian group and $F$ a connected domain of ${\mathcal H}$. 
W call $F$ a {\it fundamental domain} of $\Gamma$ if 

\noindent
(1) ${\mathcal H}=\cup_{\gamma\in\Gamma}\gamma F$;

\noindent
(2) $F=\overline U$ with an open set $U$ consisting of all 
the interior points of $F$;

\noindent
(3) $\gamma U\cap U=\emptyset$ for any 
$\gamma\in\Gamma\big/ Z(\Gamma)$.

An Fuchsian group $\Gamma$ admits a fundamental domain. A well-known method
for such a construction is that of Fricke, which goes as follows (We here
follow the presentation of [Mi]): Fix a point
 $z_0\in {\mathcal H}$ which is not an elliptic point of 
$\Gamma$, i.e., is not a point fixed by any element 
 $\alpha$ such that $\mathrm{Tr}(\alpha)^2<4\det(\alpha).$ For 
$\gamma\in\Gamma\big/ Z(\Gamma)$, put 
 $$\begin{aligned}F_\gamma:=&\Big\{z\in{\mathcal H}:d(z,z_0)
\leq d(z,\gamma z_0)\Big\};\\
 U_\gamma:=&\Big\{z\in{\mathcal H}:d(z,z_0)<d(z,\gamma z_0)\Big\};\\
C_\gamma:=&\Big\{z\in{\mathcal H}:d(z,z_0)=
d(z,\gamma z_0)\Big\}.\end{aligned}$$ 
Here $d$ denotes the hyperbolic distance on ${\mathcal H}$. Clearly, 
then $C_\gamma$ is a geodesic.

Now define the subset $F$ and $U$ on ${\mathcal H}$ by 
$$\begin{aligned}F:=&\cap_{\gamma\in\Gamma}F_\gamma:=\Big\{z\in{\mathcal H}:
d(z,z_0)\leq d(z,\gamma z_0) \forall \gamma\in\Gamma\Big\},\\
U:=&\cap_{\gamma\in\Gamma\backslash Z(\Gamma)}U_\gamma,\qquad\mathrm{where}\\
U_\gamma:=&\Big\{z\in{\mathcal H}:d(z,z_0)< d(z,\gamma z_0) \forall 
\gamma\in\Gamma\backslash Z(\Gamma)\Big\}.\end{aligned}$$
\vskip 0.20cm 
\noindent
{\Large Theorem}. {\it The subset $F$ of ${\mathcal H}$ is a fundamntal 
domain of $\Gamma$. Moreover}

\noindent 
(1) {\it any geodesic joining two points of $F$ is contained in $F$};

\noindent
(2) {\it Put $L_\gamma:=F\cap \gamma F$ for 
$\gamma\in\Gamma\big/ Z(\Gamma)$. Then $L_\gamma\in C_\gamma$. 
If $L_\gamma\not=\emptyset$, then $L_\gamma$ is only one point or a geodesic;}

\noindent
(3) {\it For any compact subset $M$ of ${\mathcal H}$, 
$\Big\{\gamma\in\Gamma:M\cap \gamma F\not=\emptyset\Big\}$ is finite.}
\vskip 0.20cm 
For each $\gamma\in\Gamma\big/ Z(\Gamma)$, put
 $L_\gamma=F\cap\gamma F$. We call $L_\gamma$ a {\it side} of 
$F$ if $L_\gamma$ is neither a null set nor a point. 
The boundary of $F$ consists of sides of $F$. For two 
distinct sides $L$ and $L'$ of $F$, $L\cap L'$ is either 
null or a point $z$ in which case we call $z$ a 
vertex of $F$ in ${\mathcal H}$. For two sides $L,\,L'$ of $F$, 
we say that $L$ and $L'$ are linked and write 
$L\sim L'$ if either $L=L'$ or there exist distinct sides 
$L_1,\cdots, L_n$ of $F$ with $L=L_1,\, L'=L_n$ and 
$L_\nu\cap L_{\nu+1}\not=\emptyset$, for $1\leq\nu\leq n-1$. 
For a side $L$ of $F$, the connected component of the
 boundary containing $L$ is a union of all sides $L'$ which 
are linked to $L$. When a side $L$ of $F$ has no
  end, we call the intersection points of the extension of $L$ 
with $\mathbb P^1(\mathbb R):=
  \mathbb R\cup\{\infty\}$ the vertices of $F$ on  
$\mathbb P^1(\mathbb R)$. (Here extension is taken by 
  considering $L$ a part of a circle or a line orthogonal 
to the real axis.) One checks that if a vertex $x$ 
  of $F$ on  $\mathbb P^1(\mathbb R)$ is an end of two 
sides and $x$ is fixed by a non-scalar element $\gamma$ 
  of $\Gamma$, then $x$ is a  cusp of $\Gamma$. Moreover,
 if $\Gamma$ is a Fuchsian group of the first kind, 
 any vertex of $F$ on  $\mathbb P^1(\mathbb R)$ is 
a cusp of $\Gamma$ and any cusp of $\Gamma$ is 
  equivalent to a vertex of $F$ on  $\mathbb P^1(\mathbb R)$.

Let $\mathcal C_\Gamma$ be the set of all cusps of $\Gamma$ 
and ${\mathcal H}^*:={\mathcal H}^*_\Gamma:={\mathcal H}
\cup\mathcal C_\Gamma$. Put 
$$U_l:=\Big\{z\in{\mathcal H}:\Im (z)>l\Big\}
\qquad\mathrm{and}\qquad 
E_l^*:=U_l\cup\{\infty\}\qquad\mathrm{for}\qquad l>0.$$
Then we can introduce a topology structure on ${\mathcal H}^*$ as follows:

\noindent
(i) For $z\in{\mathcal H}$, we take as a fundamental 
neighborhood system at $z$ in ${\mathcal H}^*$ that at $z$ in ${\mathcal H}$;

\noindent
(ii) For $x\in\mathcal C_\Gamma$, we take as a  fundamental 
neighborhood system at $z$ in ${\mathcal H}^*$
the family $\Big\{\sigma U_l^*:l>0\Big\}$, where $\sigma\in SL(2,\mathbb R)$ 
such that $\sigma \infty=x$.

\noindent
Then ${\mathcal H}^*$ becomes a Hausdorff space. In fact, put 
$\sigma=\left(\begin{matrix} a&b\\ c&d\end{matrix}
\right)$ and $x=\left[\begin{matrix} \alpha\\ \beta\end{matrix}\right]
=\frac{\alpha}{\beta}$, then
$$\sigma U_l=\Big\{z\in{\mathcal H}:\frac{\Im (z)}{|cz+d|^2}>l\Big\},$$ 
and this is the inside of a circle with the radius 
$(2l^2)^{-1}$ tangent to the real axis at $x$. For 
$x\in\mathcal C_\Gamma$, we call $\sigma U_l$ a 
neighborhood of $x$ in ${\mathcal H}$.

The action of $\Gamma$ on ${\mathcal H}$ extends naturally to 
${\mathcal H}^*$ as we have already seen. 
In particular, under the quotient topology 
$\Gamma\Big\backslash{\mathcal H}^*$ is Hausdorff. Furthermore, it is 
well-known that one  can put a compact Riemann 
surface structure on $\Gamma\Big\backslash{\mathcal H}^*$ when 
$\Gamma$ is a Fuchsian group of the first kind.

Consequently, if $\Gamma$ is a Fuchsian group of the first 
kind with $\infty$ as one of its cusps, then a systems of 
a  fundamental neighborhood system of 
$\Gamma\backslash{\mathcal H}^*$ at $\infty$ is given by $F\cap U_l^*$. 
In precise term, 
$$F\cap U_l^*=\{z\in {\mathcal H}^*:\Im (z)>l\}\Big/\langle z\mapsto z+m\rangle,$$ 
that is, a width $m$ rectangle starting from $y=l$ towards infinity.

\subsection{Upper Half Space}

Following [EGM], a closed subset $\mathcal F\in{\mathbb H}$ 
is called a {\it fundamental domain} of a 
discontinuous group $\Gamma\subset SL(2,\mathbb C)$ if

\noindent
(1) $\mathcal F$ meets each $\Gamma$-orbit at least once;

\noindent
(2) the interior $\mathcal F^0$ meets each $\Gamma$-orbit at most once;

\noindent
(3) the boundary of $\mathcal F$ has Lebesgue measure zero.

Now let $\Gamma\subset SL(2,\mathbb C)$ be a discrete group 
and suppose that $U^\lambda:=
\left(\begin{matrix} 1&\lambda\\ 0&1\end{matrix}\right)\in
\Gamma,$ for some $0\not=\lambda\in\mathbb C$. Choose $M=
\left(\begin{matrix} a&b\\ c&d\end{matrix}\right)\in\Gamma$ with $c\not=0.$ 
Define ${\mathbb H}_\lambda:=
\Big\{z+rj:z\in\mathbb C, r>|\lambda|\Big\}.$ Then 
${\mathbb H}_\lambda\cap \Big(M\cdot {\mathbb H}_\lambda\Big)=\emptyset$. 
Consequently, if $\infty$ is a fixed point of a parabolic element of 
$\Gamma$ and that $0\not=\lambda\in\mathbb C$ is 
such that 
$\left(\begin{matrix} 1&\lambda\\ 0&1\end{matrix}\right)\in\Gamma$, 
with $|\lambda|$ minimal, then two points 
contained in ${\mathbb H}_\lambda$ are $\Gamma$-equivalent if 
and only if they are 
$\Gamma_\infty$-equivalent. Moreover, if $\zeta_1,\,\zeta_2\in
\mathbb P^1(\mathbb C)$ are $\Gamma$-equivalent 
such that $\zeta_1=A_1\infty,\,\zeta_2=A_2\infty$ for $A_1,\,A_2\in
 SL(2,\mathbb C)$ and that there exist 
$\lambda_1,\,\lambda_2\in\mathbb C\backslash\{0\}$ with 
$$A_1\left(\begin{matrix} 1&\lambda_1\\ 0&1
\end{matrix}\right)A_1^{-1},\qquad A_2\left(\begin{matrix} 1&\lambda_2\\ 
0&1\end{matrix}\right)A_2^{-1}\in\Gamma,$$ 
then $\Big(A_1\cdot {\mathbb H}_{\lambda_1}\Big)\cap 
\Big(A_2\cdot\mathbb H_{\lambda_2}\Big)=\emptyset$ for all $\gamma\in\Gamma$.

The set $A\cdot{\mathbb H}_\lambda$ for $A\in SL(2,\mathbb C)$ 
and $\lambda\in\mathbb C\backslash\{0\}$ is
either an open upper half-space or an open ball in 
${\mathbb H}_\lambda$ touching 
$\mathbb P^1(\mathbb C)$. $A\cdot{\mathbb H}_\lambda$ 
are called {\it horoballs}. Taking the usual topology on 
${\mathbb H}$ and the horoballs touching $\mathbb P^1(\mathbb C)$ 
at $\lambda$ as basis for the 
neighborhoods of $\lambda$, we get a topology on 
${\mathbb H}\cup\mathbb P^1(\mathbb C)$. The group 
$SL(2,\mathbb C)$ acts continuously on this space.

As such, if a discrete group $\Gamma\subset SL(2,\mathbb C)$ 
contains a parabolic element, then $\Gamma$ is 
not cocompact, i.e., $\Gamma\Big\backslash {\mathbb H}$ is 
not compact. Furthermore, if $\Gamma$ is of finite
 covolume, then $\Gamma$ is not cocompact if $\Gamma$ contains 
a parabolic element; and if the stablizer
  $\Gamma_\zeta$ of $\zeta\in \mathbb P^1(\mathbb C)$ contains 
a parabolic element, then $\zeta$ is a cusp of
   $\Gamma$. Finally, as said before, finite covolume discrete 
subgroup $\Gamma\subset SL(2,\mathbb C)$ admits 
   only finitely many cusps.

Now for the given a discrete group $\Gamma$ of finite covolume, choose
 $A_1,\cdots,A_h\in SL(2,\mathbb C)$ so that
$$\eta_1=A_1\infty,\,\ldots,\,\eta_h=A_h\infty\in \mathbb P^1(\mathbb C)$$
 are representatives for the 
$\Gamma$-equivalence classes of cusps. Choose further closed fundamental set 
$\mathcal P_i$ (which is supposed to be 
parallelopiped) for the action of the stablizer 
$A_i^{-1}\Gamma_{\eta_i}A_i$ on $\mathbb P^1(\mathbb C)
\backslash\{\infty\}=\mathbb C$. (Recall here that 
$A_i^{-1}\Gamma_{\eta_i}A_i$ contains a lattice of 
$\mathbb C$, being the one associated
with cusp $\eta_i$, by definition.) Define, for $Y>0$, 
$$\widetilde{\mathcal F}_i(Y):=\{z+rj:z\in\mathcal P_i, 
r\geq Y\}.$$ Let further $Y_1,\,\ldots,\,Y_h\in \mathbb R_+^*$ be 
large enough so that the $\mathcal F_i(Y_i):=
A_i\cdot\widetilde {\mathcal F}_i(Y)$ are contained in the horospheres 
$A_i\cdot{\mathbb H}_i$. The $\mathcal F_i(Y_i)$
are called cusp sections. Then we have the following
\vskip 0.20cm
\noindent
{\bf{\bf {\Large Theorem.}}} (See e.g. [EGM]) {\it There exists a compact set 
$\mathcal F_0\subset{\mathbb H}$ so that 
$$\mathcal F=\mathcal F_0\cup \mathcal F_1(Y_1)\cup\cdots\cup 
\mathcal F_h(Y_h)$$ is a fundamental domain for
 $\Gamma$. Furthermore, the compact set $\mathcal F_0$ can be 
chosen so that the intersections 
 $\mathcal F_0\cap \mathcal F_i(Y_i)$ are all contained in 
the boundary of $\mathcal F_0$ and hence have
  Lebesgue measure 0 and the intersections 
$\mathcal F_i(Y_i)\cap \mathcal F_j(Y_j)$ will be empty if $i\not=j$.}

\subsection{Rank Two $\mathcal O_K$-Lattices}

\subsubsection{A. Totally Real Fields}

\subsubsection{A.1. Distance to Cusps}

Back to Siegel's presentation [Sie] for totally real fields again.

Now we introduce the notion of \lq distance of a point $z\in\mathcal H_n$ from
a cusp $\lambda$ of $\mathcal H_n$. We have alreday in $\mathcal H_n$ a 
metric given by $ds^2=\sum_{i=1}^n\frac{dx_i^2+dy_i^2}{y_i^2}$ which is
 non-euclidean in the case $n=1$ and has an invariance property with
 respect to $\Gamma$. But since the cusps lie on the boundary of 
$\mathcal H_n$, the distance relative to this metric of an inner point of 
$\mathcal H_n$ from a cusp is infinite. Hence, this metric is not useful 
for our purposes.

For $z\in\mathcal H_n$ and a cusp $\lambda=\frac{\rho}{\sigma}$ with 
associated $A=\begin{pmatrix}\rho&\xi\\ \sigma&\eta\end{pmatrix}\in G$, 
we define the {\it distance} $\Delta(z,\lambda)$ of $z$ from $\lambda$ by 
$$\Delta(z,\lambda):=\frac{1}{N(y_{A^{-1}}))^{\frac{1}{2}}}
=N\Big(\frac{|-\sigma z+\rho|^2}{y}\Big)^{\frac{1}{2}}=
N\Big(\frac{-(\sigma x+\rho)^2+\sigma^2y^2}{y}\Big)^{\frac{1}{2}}.$$ 
For example, if $\lambda=\infty$, then 
$\Delta(z,\infty)=\frac{1}{\sqrt {N(y)}};$ hence the larger the $N(y)$, the 
closer is $z$ to $\infty$.
\vskip 0.20cm
\noindent
{\bf{\large Lemma}.} (1) {\it $\Delta(z,\lambda)$ has an important invariance property with respect to 
$\Gamma$. Namely, for $M\in\Gamma$, we have $$\Delta(z_M,\lambda_M)=
\Delta(z,\lambda).$$}

\noindent
(2) {\it $\Delta(z,\lambda)$ does not depend on the special choice of $A$ 
associated with $\lambda$.}

\noindent
Proof. These are very easy to verify. Indeed, for (1), since by definition $$
\Delta(z_M,\lambda_M)=\frac{1}{N\Big(\Im(z_M)_{A^{-1}M^{-1}}
\Big)^{\frac{1}{2}}}=
\frac{1}{N\Big(y_{A^{-1}}\Big)^{\frac{1}{2}}}=\Delta(z,\lambda).$$
 
To prove (2), if $A_1=\begin{pmatrix}\rho_1&\xi_1\\ 
\sigma_1&\eta_1\end{pmatrix}\in G$ is associated with $\lambda=
\frac{\rho_1}{\sigma_1}$, then $A^{-1}A_1=\begin{pmatrix}\varepsilon&\zeta\\ 
0&\varepsilon^{-1}\end{pmatrix}$ where $\varepsilon$ is a unit in $K$. So 
$N(y_{A_1^{-1}})=N(y_{A^{-1}})N(\varepsilon^{-2})=N(y_{A^{-1}})$ and hence 
 our assertion is proved.
\vskip 0.20cm
Let, for a given cusp $\lambda$ and $r>0$, $U_{\lambda,r}$ denote the set of 
$z\in\mathcal H_n$ such that $\Delta(z,\lambda)<r$. This defines a \lq
 neighbodhood' of $\lambda$, and all points $z\in\mathcal H_n$ which belong 
to $U_{\lambda,r}$ are inner points of the same. Then

\noindent
(a) {\it The neighborhoods $U_{\lambda,r}$ for
$0<r<\infty$ cover the entire  $\mathcal H_n$}.

\noindent
(b) {\it Each neighborhood $U_{\lambda,r}$ is left invariant by a modular 
substitution in $\Gamma_\lambda$.}

Indeed,  by (1) above, if $M\in\Gamma_\lambda$ then 
$$\Delta(z_M,\lambda)=\Delta(z_M,\lambda_M)=\Delta(z,\lambda)$$ and so 
if $z\in U_{\lambda,r}$, then again $
\Delta(z_M,\lambda)<r$, i..e, $z_M\in U_{\lambda,r}$. Consequently,

\noindent
(c) {\it A fundamental domain for $\Gamma_\lambda$ in 
$U_{\lambda,r}$ is given by $\mathcal D_\lambda\cap U_{\lambda,r}$.}
\vskip 0.30cm
We shall now prove some interesting facts concerning $\Delta(z,\lambda)$ 
which will be useful in constructing a fundamental domain for $\Gamma$ 
in $\mathcal H_n$.

\noindent
{\bf i)} {\it For $z=x+iy\in\mathcal H_n$, there exists a cusp $\lambda_0$ of 
$\mathcal H_n$ such that for all cusps $\mu$ of $\mathcal H_n$, we have 
$$\Delta(z,\lambda_0)\leq\Delta(z,\mu).$$}

\noindent
Proof. If $\lambda$ is a cusp, then $\lambda=\frac{\rho}{\sigma}$ for some
${\rho},{\sigma}\in\mathcal O_K$ such that $({\rho},{\sigma})$ is one 
of the $h$ ideals $\frak a_1,\ldots,\frak a_h$. Then $$\Delta(z,\lambda)=
N\bigg(\frac{-(\sigma x+\rho)^2+\sigma^2y^2}{y}\bigg)^{\frac{1}{2}}.$$
Let us consider the expression $N\bigg(\frac{-(\sigma x+\rho)^2+
\sigma^2y^2}{y}\bigg)^{\frac{1}{2}}$ as a function of the pair of integers 
$({\rho},{\sigma})$. It remains unchanged if ${\rho},{\sigma}$ are
 replaced by ${\rho}\varepsilon,{\sigma}\varepsilon$ for any unit 
$\varepsilon$ in $K$. We shall now show that there exists a pair of 
integers ${\rho}_1,{\sigma}_1$ in $K$ such that
$$N\bigg(\frac{-(\sigma_1 x+\rho_1)^2+
\sigma_1^2y^2}{y}\bigg)^{\frac{1}{2}}\,\leq\, N\bigg(\frac{-(\sigma x+\rho)^2+
\sigma^2y^2}{y}\bigg)^{\frac{1}{2}}\eqno(1)$$ for all pairs 
of integers $({\rho},{\sigma})$. In order to prove it, obviously it 
suffices to show that for given $c_{11}>0$, there are only finitely 
many non-associated pairs of integers $({\rho},{\sigma})$ such that 
$$N\bigg(\frac{-(\sigma x+\rho)^2+
\sigma^2y^2}{y}\bigg)^{\frac{1}{2}}\leq c_{11}.\eqno(2)$$

Now it is known from the theory of algebraic number fields that if 
$\alpha=(\alpha_1,\ldots,\alpha_n)$ is an $n$-tuple of real numbers 
with $N(\alpha)\not=0$, then we can find a unit $\varepsilon$ in $K$ 
such that $$\Big|\alpha_i\cdot\varepsilon^{(i)}\Big|\,\leq\, 
c_{12}\root {n}\of{|N(\alpha)|}\eqno(3)$$ 
for a constant $c_{12}$ depending only on $K$. In view of inequality 
(2), we can 
suppose, after multiplying $\rho$ and $\sigma$ by a suitable unit 
$\varepsilon$, that already we have
$$\frac{(-\sigma^{(i)} x_i+\rho^{(i)})^2+{\sigma^{(i)}}^2y_i^2}{y_i}
\leq c_{13},\qquad i=1,2,\ldots,n$$ for a constant $c_{13}$ depending only 
on $c_{11}, c_{12}$. This
implies that $(-\sigma^{(i)} x_i+\rho^{(i)})$ and ${\sigma^{(i)}}$ and, 
 consequently, $\rho^{(i)}$ and $\sigma^{(i)}$ are bounded for 
$i=1,2,\ldots,n$. Again, in this case, we know from the theory of 
algebraic number theory that there are only finitely many poissibilities 
for $\rho$ and $\sigma$. Hence inequality (2) is true only for finitely many 
non-associated pairs of integers $(\rho,\sigma)$. From these pairs, 
we choose a pair $(\rho_1,\sigma_1)$ such that the valus of
$N\bigg(\frac{-(\sigma x+\rho)^2+
\sigma^2y^2}{y}\bigg)^{\frac{1}{2}}$ is minimum. 
This pair $(\rho_1,\sigma_1)$ now obviously satisfies inequality (1).

Let $(\rho_1,\sigma_1)=\frak b$ and let $\frak b=\frak a_i\cdot(\theta^{-1})$ 
for some $\frak a_i$ and an element $\theta\in K$. Then 
$\frak a_i=(\rho_1\theta,
\sigma_1\theta).$ Now in inequality (1) if we replace $\rho$ and $\sigma$ by 
$\rho_1\theta,\,\sigma_1\theta$ respectively, we get 
$\Big|N(\theta)\Big|\geq 1$. On the 
other hand, since $\frak a_i$ is of minimum norm among the integral ideals 
of its class, $N(\frak a_i)\leq N(\frak b)$ and therefore $|N(\theta)|\leq 1$. 
Thus $\Big|N(\theta)\Big|=1$. Let now $\rho_0=\rho_1\theta,\,
\sigma_0=\sigma_1\theta$ and 
$\lambda_0=\frac{\rho_0}{\sigma_0}$. Then
$\frak a_i=(\rho_0,\sigma_0)$ and by definition, 
$\Delta(z,\lambda_0)=N\bigg(\frac{-(\sigma_0 x+\rho_0)^2+
\sigma_0^2y^2}{y}\bigg)^{\frac{1}{2}}=N\bigg(\frac{-(\sigma_1 x+\rho_1)^2+
\sigma_1^2y^2}{y}\bigg)^{\frac{1}{2}}.$ If we use inequality (1), then we
 see at once that
$\Delta(z,\lambda)\leq \Delta(z,\mu)$ for all cusps $\mu$. This 
completes the proof.
\vskip 0.20cm
For given $z\in\mathcal H_n$, define $\Delta(z)=\inf_\lambda\Delta
(z,\lambda)$. By i), there exists a cusp $\lambda$ such that 
$\Delta(z)=\Delta(z,\lambda)$. In general, $\lambda$ is unique, but
 there are exceptional cases when the minimum is attained for more 
than one $\lambda$. We shall see presently that there exists a constant
 $d>0$, depending only on $K$ such that if $\Delta(z)<d$, then the cusp
 $\lambda$ for which $\Delta(z)=\Delta(z,\lambda)$ is unique. 

\noindent
{\bf ii)} {\it There exists $d>0$ depending only on $K$ such that if for 
$z=x+iy\in\mathcal H_n$, $\Delta(z,\lambda)<d$ and $\Delta(z,\mu)<d$,
 then necessarily $\lambda=\mu$.}

\noindent
Proof. Let $\lambda=\frac{\rho}{\sigma}$ and $\mu=\frac{\rho_1}{\sigma_1}$ 
and let for a real number $d>0$,
$$\begin{aligned}\Delta(z,\lambda)=&N\Big(-(\sigma x+\rho)^2y^{-1}+\sigma^2y
\Big)^{\frac{1}{2}}<d,\\
\Delta(z,\mu)=&N\Big(-(\sigma_1 x+\rho_1)^2y^{-1}+\sigma_1^2y\Big)^{\frac{1}{2}}<d.\end{aligned}$$
After multiplying $\rho$ and $\sigma$ by a suitable unit 
$\varepsilon$ in $K$, we might assume in view of inequality (3) that
$$(\sigma^{(i)} x_i-\rho^{(i)})^2y_i^{-1}+{\sigma^{(i)}}^2y_i
<c_{12}d^{2/n},\qquad i=1,2,\ldots, n.$$ Hence 
$$\Big|-\sigma^{(i)} x_i+\rho^{(i)}\Big|
y_i^{-\frac{1}{2}}<\sqrt{c_{12}}d^{1/n},\qquad
\Big|\sigma^{(i)}\Big|y_i^{\frac{1}{2}}<\sqrt{c_{12}}d^{1/n}.$$
Similarly we have 
$$\Big|-\sigma_1^{(i)} x_i+\rho_1^{(i)}\Big|y_i^{-\frac{1}{2}}
<\sqrt{c_{12}}d^{1/n},\qquad
\Big|\sigma_1^{(i)}\Big|y_i^{\frac{1}{2}}<\sqrt{c_{12}}d^{1/n}.$$
Now $\rho^{(i)}\sigma_1^{(i)}-\rho_1^{(i)}\sigma^{(i)}=
\Big(-\sigma^{(i)} x_i+\rho^{(i)}\Big)y_i^{-\frac{1}{2}}\cdot{\sigma^{(i)}}
y_i^{\frac{1}{2}}-\Big(-\sigma_1^{(i)} x_i+\rho_1^{(i)}\Big)y_i^{-\frac{1}{2}}
\cdot{\sigma^{(i)}}y_i^{\frac{1}{2}}$ and hence 
$$\Big|N(\rho\sigma_1-\rho_1\sigma)\Big|
\,<\,\Big(2c_{12}d^{2/n}\Big)^n.$$ If we set $d=\Big(2c_{12}\Big)^{-n/2}$, 
then $\Big|N(\rho\sigma_1-\rho_1\sigma)\Big|<1$. Since 
$\rho\sigma_1-\rho_1\sigma$ is an integer, 
it follows that $\rho\sigma_1-\rho_1\sigma=0$, i.e., $\lambda=\mu$.

Thus for $d=\Big(2c_{12}\Big)^{-n/2}$, the conditions $\Delta(z,\lambda)<d$, 
$\Delta(z,\mu)<d$ for a $z\in\mathcal H_n$ imply that $\lambda=\mu$. 
Therefore, the neighborhoods $U_{\lambda,d}$ for the various cusps 
$\lambda$ mutually disjoint.
\vskip 0.20cm
We shall now prove that for $z\in\mathcal H_n$, $\Delta(z)$ is 
uniformly bounded in $\mathcal H_n$. To this end, it suffices to prove

\noindent
{\bf iii)} {\it There exists $c>0$ depending only on $K$ such that for any 
$z=x+iy\in\mathcal H_n$, there exists a cusp $\lambda$ with the property 
that $\Delta(z,\lambda)<c$. In particular, 
$\mathcal H_n=\cup_\lambda U_{\lambda,c}$.}

\noindent
Proof. We shall prove the existence of a constant $c>0$ depending only on 
$K$ and a pair of integers $(\rho,\sigma)$ not both zero such that 
$$N\Big(-(\sigma x+\rho)^2y^{-1}+\sigma^2y\Big)^{\frac{1}{2}}<c.$$

Let $\omega_1,\ldots,\omega_n$ be a $\mathbb Z$-basis of $\mathcal O_K$. 
Consider now 
the following system of $2n$ linear inequalities in the $2n$ variables 
$a_1,\,\ldots,\, a_n,\ b_1,\,\ldots,\,b_n$, viz.
$$\begin{aligned}&\Big|y_1^{-\frac{1}{2}}\Big(\omega_1^{(1)}a_1+\ldots+
\omega_n^{(1)}a_n\Big)-x_1y_1^{-\frac{1}{2}}\Big(\omega_1^{(1)}b_1+\ldots+
\omega_n^{(1)}b_n\Big)\Big|\leq\alpha_1\\
&\cdots\\
&\Big|y_n^{-\frac{1}{2}}\Big(\omega_1^{(n)}a_1+\ldots+
\omega_n^{(n)}a_n\Big)-x_ny_n^{-\frac{1}{2}}\Big(\omega_1^{(n)}b_1+\ldots+
\omega_n^{(n)}b_n\Big)\Big|\leq\alpha_n\end{aligned}$$
$$\begin{aligned}&\Big|y_1^{\frac{1}{2}}\Big(\omega_1^{(1)}b_1+\ldots+
\omega_n^{(1)}b_n\Big)\Big|\leq\beta_1\\
&\cdots\\
&\Big|y_n^{\frac{1}{2}}\Big(\omega_1^{(n)}b_1+\ldots+
\omega_n^{(n)}b_n\Big)\Big|
\leq\beta_n.\end{aligned}$$
The determinant of this system of linear forms is 
$(\omega_i^{(j)})^2=\Delta_K$ where $\Delta_K$ is the absolute value of
the discriminant of $K$. 
By Minkowski's theorem on linear forms, this system of linear 
inequalities has a non-trivial solution in rational integers 
$a_1,\,\ldots,\, a_n,\ b_1,\,\ldots,\,b_n$ if $\alpha_1,\,\ldots,\, \alpha_n, 
\ \beta_1,\,\ldots,\,\beta_n\geq \root{2n}\of {\Delta_K}$. 
Taking $\alpha_i=\beta_j=
\root{2n}\of {\Delta_K},\ \, i,j=1,2,\ldots,n$, in particular, this system 
of inequalities has a non-trivial solution in rational integers,
 say, $\,a_1',\,\ldots,\, a_n',\ b_1',\,\ldots,\,b_n'$. Let us take 
$\rho=a_1'\omega_1+\ldots+a_n'\omega_n$ and $\sigma=b_1'\omega_1+
\ldots+b_n'\omega_n$. Then we obtain $$\Big|(-\sigma^{(i)}x_i+
{\rho^{(i)}})^2y_i^{-1}+{\sigma^{(i)}}^2y_i\Big|\leq 2\root{n}\of{\Delta_K}, 
\qquad i=1,2,\ldots, n.$$
Hence $N\Big(-(\sigma x+\rho)^2y^{-1}+\sigma^2y\Big)^{\frac{1}{2}}\leq 
2^{n/2}\Delta_K^{\frac{1}{2}}=:c$ say.

Now let $(\rho,\sigma)=\frak b;\ \frak b=\frak a_i(\theta)^{-1}$ 
for some $\frak a_i$ and a $\theta\in K$. Since $\frak a_i$ 
is of minimum norm among the integral ideals of its class, 
$\big|N(\theta)\big|\leq 1$. Further $\frak a_i=(\rho\theta,\sigma\theta)$. 
Now, for the cusp $\lambda=\frac{\rho\theta}{\sigma\theta}=
\frac{\rho}{\sigma},$ we have
$$\Delta(z,\lambda)=\Big|N(\theta)\Big|\cdot 
N\Big(-(\sigma x+\rho)^2y^{-1}+\sigma^2y\Big)^{\frac{1}{2}}\leq c,$$ which was what we wished to prove. Consequently, 
we deduce that $\mathcal H_n=\cup_\lambda U_{\lambda,c}.$
\vskip 0.20cm
\noindent
{\bf iv)} {\it For $z\in\mathcal H_n$ and $M\in\Gamma$, 
$\Delta(z_M)=\Delta(z)$.}

\noindent
Proof. In fact, $$\Delta(z_M)=\inf_\lambda\Delta(z_M,\lambda)=
\inf_\lambda\Delta(z,\lambda_{M^{-1}})=\inf_\lambda\Delta(z,\lambda)=
\Delta(z).$$

\subsubsection{A.2. Fundamental Domain for $\Gamma$ in $\mathcal H_n$}

We now have all the necessary material for the construction of a fundamental 
domain for $\Gamma$ in $\mathcal H_n$.

A point $z\in\mathcal H_n$ is {\it semi-reduced} (with respect to a cusp 
$\lambda$), if $\Delta(z)=\Delta(z,\lambda)$. If $z$ is semi-reduced with
 respect to $\lambda$, then for all cusps $\mu$, we have 
$\Delta(z,\mu)\geq\Delta(z,\lambda).$

Let $\lambda_1(=(\infty,\ldots,\infty)),\lambda_2,\ldots,\lambda_h$ 
be the $h$ inequivalent base cusps of $\mathcal H_n$. We denote by 
$\mathcal F_{\lambda_i}$, the set of all $z\in\mathcal H_n$ which are 
semi-reduced with respect to $\lambda_i$. Clearly, 
$\mathcal F_{\lambda_i}\subset U_{\lambda_i,r}$, in view of iii) above. 
The set $\mathcal F_{\lambda_i}$ is invariant under the modular substitution
 $z\mapsto z_M$ for $M\in\Gamma_{\lambda_i}$. For by iv) above, 
$\Delta(z_M)=\Delta(z)$ and further
$\Delta(z)=\Delta(z,\lambda_i)=\Delta\big(z_M,(\lambda_i)_M\big)=
\Delta(z_M,\lambda_i).$ Thus $\Delta(z_M)=\Delta(z_M,\lambda_i)$ 
and hence $z_M\in \mathcal F_{\lambda_i}$ for $M\in\Gamma_{\lambda_i}$.
\vskip 0.30cm
Let $\overline{\mathcal D}_{\lambda_i}$ denote the closure in $\mathcal H_n$ 
of the set ${\mathcal D}_{\lambda_i}$ and $\mathcal D_i=
\mathcal F_{\lambda_i}\cap \overline{\mathcal D}_{\lambda_i}$. Then 
$\mathcal D_i$ is explicitly defined as the set of $z\in\mathcal H_n$ 
whose local coordinates $X_1,\ldots,X_n, Y_1,\ldots,Y_{n-1}$ relative
 to the cusp $\lambda_i$ satisfy the conditions 
$$-\frac{1}{2}\leq X_k,\ \ Y_l\leq\frac{1}{2},\qquad k=1,2,\ldots,n,
\ l=1,2,\ldots,n-1,$$ and further
for all cusps $\mu$, $$\Delta(z,\mu)\geq\Delta(z,\lambda_i).$$ Let 
$\mathcal F:=\cup_{i=1}^h\mathcal D_i$. We see that the $\mathcal D_i$ 
as also $\mathcal F$ are closed in $\mathcal H_n$.

A point $z\in\mathcal D_i$ is an {\it inner point} of $\mathcal D_i$ if in all
 the conditions above, strictly inequality holds, viz
$$-\frac{1}{2}< X_k,\ \ Y_l<\frac{1}{2},\qquad k=1,2,\ldots,n,
\ l=1,2,\ldots,n-1,$$ and
$$\Delta(z,\mu)>\Delta(z,\lambda_i)\qquad\forall\mu\not=\lambda_i.$$
If equality holds even in one of these conditions, then $z$ is said to 
be a {\it boundary point} of $\mathcal D_i$. The set of boundary points of 
$\mathcal D_i$ constitute the boundary of $\mathcal D_i$, which may be 
denoted by $\partial \mathcal D_i$. We denote by ${\mathcal D_i^0}$ the set 
of inner points of $\mathcal D_i$.

It is clear that the $\mathcal D_i$'s do not overlap and intersect at most 
on their boundary.

A point $z\in\mathcal F$ may now be called an {\it inner point} of 
$\mathcal F$, if $z$ is an inner point of some $\mathcal D_i$; similarly 
we may define a {\it boundary point} of $\mathcal F$ and denote the set of 
boundary points of $\mathcal F$ by $\partial\mathcal F$.

We say that a point $z\in\mathcal H_n$ is {\it reduced} (with respect to 
$\Gamma$) if, in the first place, $z$ is semi-reduced with respect to 
some one of the $h$ cusps $\lambda_1,\ldots,\lambda_h$, say $\lambda_i$ 
and then further $z\in\overline {\mathcal D}_{\lambda_i}$. 
Clearly $\mathcal F$ 
is just the set of all $z\in\mathcal H_n$ reduced with respect to $\Gamma$.
\vskip 0.20cm
Before we proceed to show that $\mathcal F$ is a fundamental domain for 
$\Gamma$ in $\mathcal H_n$ we shall prove the following result concerning 
the inner points of $\mathcal F$, namely,

\noindent
{\bf{\large Lemma}}. {\it The set of inner points of $\mathcal F$ is open 
in $\mathcal H_n$.}

\noindent
Proof. It is enough to show that each ${\mathcal D_i^0}$ is open in 
$\mathcal H_n$. Let then $z_0=x_0+iy_0\in\mathcal D_j^0$; we have to prove 
that there exists
a neighborhood $V$ of $z_0$ in $\mathcal H_n$ which is wholly contained in 
${\mathcal D_j^0}$.

Recall that for each $z=x+iy\in\mathcal H_n$ and a  cusp 
$\mu=\frac{\rho}{\sigma}$, by defintion,
 we have $\Delta(z,\lambda)=N\Big((-\sigma x+\rho)^2y^{-1}+
\sigma^2y\Big)^{\frac{1}{2}}.$ Using the fact that for each $j=1,2,\ldots,n$,
$\ (-\sigma^{(j)} x_j+\rho^{(j)})^2y_j^{-1}+
{\sigma^{(j)}}^2y_j$ is a positive-definite binary quadratic form in 
$\sigma^{(i)}$ and $\rho^{(j)}$, we can easily show that 
$$\Delta(z,\mu)\geq\alpha\cdot N(\sigma^2+\rho^2)^{\frac{1}{2}}$$ where 
$\alpha=\alpha(z)$ depends continuously on $z$ and does not depend on 
$\mu$. We can find a sufficiently small neighborhood $W$ of $z_0$ such 
that for all $z\in W$,
$$\Delta(z,\mu)\geq\frac{1}{2}\alpha_0\cdot N(\sigma^2+\rho^2)^{\frac{1}{2}}$$ 
where $\alpha_0=\alpha(z_0)$. Moreover, we can assume $W$ so chosen that 
for all $z\in W$, $$\Delta(z,\lambda_j)\leq 2\Delta(z_0,\lambda_j).$$ 
Thus for all cusps $\mu=\frac{\rho}{\sigma}$ and $z\in W$,
$$\Delta(z,\mu)-\Delta(z,\lambda_j)\geq \frac{1}{2}\alpha_0\cdot
N(\sigma^2+\rho^2)^{\frac{1}{2}}-2\Delta(z_0,\lambda_j).$$

Now employing an argument used already in the previous subsection, 
we can show that 
there are only finitely many non-associated pairs of integers 
$({\rho},{\sigma})$ such that 
$$\frac{1}{2}\alpha_0\cdot N(\sigma^2+\rho^2)^{\frac{1}{2}}\leq 2
\Delta(z_0,\lambda_j).$$
 It is an immediate consequence that, except for finitely many
 cusps $\mu_1,\ldots,\mu_r$, we have  
$$\Delta(z,\mu)>\Delta(z,\lambda_j)\qquad\forall\mu\not=\lambda_j.$$ 
Now since $\Delta(z_0,\mu)>
\Delta(z_0,\lambda_j)$ for all cusps $\mu\not=\lambda_j$, we can, 
in view of the continuity in $z$ of $\Delta(z,\mu_i)$ for $i=1,2,
\ldots,r,$ find a neighborhood $U$ of $z_0$ such that 
$$\Delta(z,\mu_k)>\Delta(z,\lambda_j),\qquad k=1,2,\ldots,r,
\ \forall z\in U\ \&\ \mu_k\not=\lambda_j.$$ We then have finally
$$\Delta(z,\mu)>\Delta(z,\lambda_j),\qquad\forall z\in U\cap W\ \ \mathrm{provided}\ \ \mu\not=\lambda_j.$$ 
Further, we could have chosen $U$ such that for all $z\in U$, the inequalities 
$(*)$ in 2.3.3.A.3 is satisfied, in addition. Thus the neighborhood $V=U\cap W$ of $z_0$
 satisfies our requirements and so $\mathcal D_j^0$ is open. This completes
 the proof.
\vskip 0.20cm
It may now be seen that the closure $\overline{\mathcal D_i^0}$ of 
${\mathcal D_i^0}$ is just $\mathcal D_i$. In fact, 
let $z\in\mathcal D_i$ and let 
$\lambda_i=\infty_{A_i}, \ z^*=z_{A_i^{-1}},\ \mu_{A_i^{-1}}=\nu=
\frac{\rho}{\sigma}\not=\infty$. Then we have $$\sigma\not=0,\qquad
\frac{\Delta(z,\mu)}{\Delta(z,\lambda_i)}=\frac{\Delta(z^*,\nu)}
{\Delta(z^*,\infty)}=\bigg(N\Big((-\sigma x^*+\rho)^2+
(\sigma y^*)^2\Big)\bigg)^{\frac{1}{2}}.$$
If $y^*$ is replaced by $ty^*$ where $t$ is a positive scalar factor, 
then the expression above is a strictly monotonic increasing function of
 $t$, whereas the coordinates $X_1,\ldots,X_n, Y_1,\ldots,Y_{n-1}$ remain 
unchanged. Thus, if $z\in \mathcal D_i$, and $z_{A_i^{-1}}=z^*=x^*+iy^*,$ 
then for $z^{(t)}=(x^*+ity^*)_{A_i},\ t>1$, the inequalities 
$\Delta(z,\mu)>\Delta(z,\lambda_i),\ \mu\not=\lambda_i$ 
are satisfied and moreover for a 
small change in the coordinates $X_1,\ldots,X_n, Y_1,\ldots,Y_{n-1}$, the 
inequalities $(*)$ are also satisfied. As a consequence, if 
$z\in\mathcal D_i$ then every neighborhood of $z$ intersects 
${\mathcal D_i^0}$; in other words, $\overline{\mathcal D_i^0}=\mathcal D_i$.

From above, clearly, we have that if $z\in \mathcal D_i$, then the entire 
curve defined by $z^{(t)}=(x^*+ity^*)_{A_i},\ t\geq 1$ lies in $\mathcal D_i$ 
and hence, in particular, $z^{(t)}$ for large $t>0$ belongs to 
${\mathcal D_i^0}$. Using this it may be shown that $\mathcal D_i$ and 
similarly ${\mathcal D_i^0}$ are connected.

The existence of inner points of $\mathcal D_i$ is an immediate consequence 
of result ii). For $j=1$, it may be verified that $z=(it,\ldots, it),\ t>1$ 
is an inner point of $\mathcal D_i$.
\vskip 0.20cm
Now we proceed to prove the following
\vskip 0.30cm
\noindent
{\bf {\Large Theorem.}} {\it  $\mathcal F$ is a fundamantal domain for 
$\Gamma$ in $\mathcal H_n$.}
\vskip 0.30cm
\noindent
We have first to show that

\noindent
(a) {\it The images $\mathcal F_M$ of $\mathcal F$ for $M\in\Gamma$ cover 
$\mathcal H_n$ without gaps.}

\noindent
Proof. This is obvious from the very method of construction of $\mathcal F$. 
First, for any $z\in\mathcal H_n$, there exists a cusp $\lambda$ such that 
$\Delta(z)=\Delta(z,\lambda)$. Let $\lambda=(\lambda_i)_M$ for some
 $\lambda_i$ and $M\in\Gamma$. We have then $$
\Delta(z_{M^{-1}})=\Delta(z)=\Delta(z,\lambda)=\Delta(z_{M^{-1}},\lambda_i)$$ 
and hence $z_{M^{-1}}\in \mathcal F_{\lambda_i}$. Now we can find 
$N\in\Gamma_{\lambda_i}$ such that $(z_{M^{-1}})_N$ is reduced with 
respect to $\Gamma_{\lambda_i}$ and thus $z_{NM^{-1}}\in\mathcal F$.
\vskip 0.20cm
\noindent
Next we need to show that

\noindent
(b) {\it The images $\mathcal F_M$ of $\mathcal F$ for $M\in\Gamma$ cover 
$\mathcal H_n$ without overlaps.}

\noindent
Proof. Let $z_1,z_2\in\mathcal F$ such that $z_1=(z_2)_M$ for an 
$M\not=\pm I_2$ in $\Gamma$ and let $z_1\in \mathcal D_i$ and $z_2\in 
\mathcal D_j.$
Now since $z_1\in \mathcal F_{\lambda_i}$, we have 
$$\Delta(z_1,\lambda_i)\leq\Delta\big(z_1,(\lambda_j)_{M^{-1}}\big)
=\Delta(z_2,\lambda_j).$$
 Similarly
$$\Delta(z_2,\lambda_j)\leq\Delta\big(z_2,(\lambda_i)_M\big)
=\Delta(z_1,\lambda_i).$$
Therefore, we obtain $$\Delta(z_1,\lambda_i)=\Delta(z_2,\lambda_j)=
\Delta\big(z_1,(\lambda_j)_M\big).$$
Two cases have to be discussed.

\noindent
($\alpha$) Let us first suppose that 
$$\Delta(z_1,\lambda_i)=\Delta(z_2,\lambda_j)<d.$$
Then since $\Delta(z_1,\lambda_i)<d$ as also
$\Delta(z_1,(\lambda_j)_M)<d$, we infer from the result ii) that 
$\lambda_i=(\lambda_j)_M$. But $\lambda_i$ and $\lambda_j$ for $i\not=j$ 
are not equivalent with respect to $\Gamma$. Therefore $i=j$ and 
$M\in\Gamma_{\lambda_i}$. Again, since both $z_1$ and $z_2$ are in
 $\overline{\mathcal D}_{\lambda_i}$ and further $z_1=(z_2)_M$ with 
$M\in\Gamma_{\lambda_i}$, we conclude that necessarily $z_1$ and $z_2$
 belong to $\partial \mathcal F$ and indeed their local coordinates $
X_1,\,\ldots,\,X_n,\ Y_1,\,\ldots,\,Y_{n-1}$ relative to 
$\lambda_i$ satisfy at 
least one of the conditions $\,X_1=\pm\frac{1}{2},\,\ldots,
\,X_n=\pm\frac{1}{2},\ Y_1=\pm\frac{1}{2},\,\ldots,\,Y_{n-1}=\pm\frac{1}{2}.$ 
Further $M$ clearly belongs to a finite set $M_1,\ldots, M_r$ of elements 
in $\cup_{i=1}^h\Gamma_{\lambda_i}$.

\noindent
($\beta$) We have now to deal with the case
$$d\leq\Delta(z_1,\lambda_i),\qquad \Delta(z_2,\lambda_j)\leq c$$
where we may suppose that $\lambda_i\not=(\lambda_j)_M$. Let for 
$i=1,2,\ldots,h$, $B_i$ denote the set of $z\in\mathcal H_n$ for which $
d\leq\Delta(z_1,\lambda_i)\leq c$ and $z\in\overline{\mathcal D}_{\lambda_i}$.
Then from Lemma 2.3.3.A.3, $B_i$ is compact and so is $B=\cup_{i=1}^h B_i$.
 Now both $z_2$ and $z_1=(z_2)_M$ belong to the compact set $B$. We may 
then deduce from Lemma 2.3.3.A.2 that $M$ belongs to a finite set of elements 
$M_{r+1},\ldots, M_s$ in $\Gamma$, depending only on $B$ and hence only on
 $K$. Moreover $z_1$ satisfies $$\Delta(z_1,\lambda_i)=
\Delta\big(z_1,(\lambda_j)_M\big)$$ with $(\lambda_j)_M\not=\lambda_i$. 
Hence $z_1$ 
and similarly $z_2$ belongs to $\partial \mathcal F$. As a result, 
arbitrarily near $z_1$ and $z_2$, there exist points $z$
such that $\Delta(z_1,\lambda_i)\not=\Delta\big(z,(\lambda_j)_M\big)$ for 
$M=M_{r+1},\ldots, M_s$.

Thus, finally, no two inner points of $\mathcal F$ can be equivalent 
with respect to $\Gamma$. Further $\mathcal F$ intersects only finitely
 many of its neighbours $\mathcal F_{M_1},\ldots, \mathcal F_{M_s}$ and 
indeed only on its boundary. We have therefore established (b).
\vskip 0.30cm
From (a) and (b) above, it follows that $\mathcal F$ is a fundamental 
domain for $\Gamma$ in $\mathcal H_n$. It consists of $h$ connected \lq 
pieces' corresponding to the $h$ inequivalent base cusps $\lambda_1,
\ldots,\lambda_h$ and is bounded by a finite number of manfolds of the form
$$\Delta(z,\lambda_i)=\Delta\big(z,(\lambda_j)_M\big),\qquad i,j=1,2,\ldots,h, 
\ \ M=M_{r+1},\ldots, M_s,$$
and hypersurfaces defined by
$$X_i^{(k)}=\pm\frac{1}{2},\qquad Y_j^{(k)}=\pm\frac{1}{2},\qquad i=1,2,\ldots,n;\ \  j=1,2,\ldots,n-1$$ where $X_1^{(k)},\ldots,X_n^{(k)}, Y_1^{(k)},
\ldots,Y_{n-1}^{(k)}$ are local coordinates relative to the base cusp
 $\lambda_k$. This completes the proof of the Theorem.
\vskip 0.20cm
The manifolds defined by the $\Delta$ equation above are seen to be 
generalizations of the isometric circles in the sense of Ford, 
for a fuchsian group. In fact 
if $\lambda_i=\frac{\rho_i}{\sigma_i}$ and $(\lambda_j)_M=\frac{\rho}
{\sigma}$ then the condition becomes $$N\Big(\big|-\sigma_i z+\rho_i\big|\Big)
=N\Big(\big|-\sigma z+\rho\big|\Big).$$
If we set $n=1$ and $\lambda_i=\infty$ or equivalently $\rho_i=1$ and
 $\sigma_i=0$, then the condition reads as $$\Big|-\sigma z+\rho\Big|=1$$
which is the familiar \lq isometric circle' corresponding to the 
transformation $z\mapsto \frac{\eta z-\xi}{-\sigma z+\rho}$ of 
$\mathcal H_1$ on itself.

The conditions $\Delta(z,\lambda)\geq \Delta(z,\lambda_i)$ by which 
$\mathcal F_{\lambda_i}$ was defined, simply mean for $n=1$ and 
$\lambda_i=\infty$ that $\Big|\gamma z+\delta\Big|\geq 1$ for all pairs of 
coprime rational integers $(\gamma,\delta)$. Thus, just as the 
points of the well-known fundamental domain in $\mathcal H=\mathcal  H_1$ for 
the elliptic modular group lie in the exterior of the isometric circles 
$\Big|\gamma z+\delta\Big|=1$ corresponding to the same group, 
$\mathcal F_{\lambda_i}$ lies in the \lq exterior' of the generalized 
isometric circles $\Delta(z,\lambda)= \Delta(z,\lambda_i)$.

Consequently, we have the following important result on truncations.
\vskip 0.20cm
\noindent
{\bf{\Large Proposition}}. (1) {\it Let $\mathcal F^*$ denote the set of $z\in\mathcal F$ for which
$\Delta(z,\lambda_i)\geq e_i>0,\ \ i=1,2,\ldots,h.$ Then
$\mathcal F^*$ is compact in $\mathcal H_n$.}

\noindent
(2) {\it For any compact set $C$ in $\mathcal H_n$, there exists
 a constant $b=b(C)>0$ such that $C\cap U_{\mu,b}=\emptyset$ for all 
cusps $\mu$.}

\noindent
Proof. (1) The proof is almost trivial in the light of Lemma 2.3.3.A.3. 
Indeed, let 
$B_i$ denote the set of $z\in\mathcal G_{\lambda_i}$ for which $e_i\leq 
 \Delta(z,\lambda_i)\leq c$. Then $B_i$ as also $B=\cup_{i=1}^hB_i$ is 
compact in $\mathcal H_n$. Hence $\mathcal F^*$ which is closed and 
contained in $B$ is again compact.

\noindent
(2) Since $C$ is compact, we can find $\alpha,\beta>0$ depending 
only on $C$ such that, for $z=x+iy\in C$, we have $\beta\leq N(y)\leq\alpha$.
 Now for any cusp $\lambda=\frac{\rho}{\sigma},\ {\rho},{\sigma}\in
\mathcal O_K$ and $z=x+iy\in C$, it is clear that $\Delta(z,\lambda)=
\Big(\big(N((-\sigma z+\rho)^2y^{-1}+\sigma^2 y\big)\Big)^{1/2}$ satisfies 
$$\Delta(z,\lambda)\geq\begin{cases}|N(\sigma)|N(y)^{\frac{1}{2}}
\geq\beta^{\frac{1}{2}}&\sigma\not=0\\
|N(\rho)|N(y)^{-\frac{1}{2}}\geq\alpha^{-\frac{1}{2}}&
\sigma=0\end{cases}$$ If we choose $b$ for which 
$0<b<\min(\alpha^{-\frac{1}{2}},\beta^{\frac{1}{2}})$, then it is 
obvious that for all cusps $\lambda$, $\ U_{\lambda,b}\cap C=\emptyset$.
This completes the proof.
\vskip 0.20cm
\noindent
{\bf Remarks.} (1) It was Blumenthal who first gave a method of 
constructing a fundamental 
domain for $\Gamma$ in $\mathcal H_n$, but his proof contained an error 
since he obtained a fundamental domain with just one cusp and not $h$ cusps. 
This error was set right by Maass.

\noindent
(2) Siegel's method of constructing the fundamental domain 
$\mathcal F$ is essentially different from the well-known method of Fricke
 for constructing a normal polygon for a discontinuous group of analytic
 automorphisms of a bounded domain in the complex plane. This method uses 
only the notion of distance of a point of $\mathcal H_n$ from a cusp, 
whereas we require a metric invariant under the group, for Fricke's method
cited in 2.4 based on the fact 
that $\mathcal H_n$ carries a Riemannian metric which is 
invariant under $\Gamma$.  For our later 
purposes, we require a fundamental domain whose nature near the cusps
 should be well known. Therefore we see in the first place that the 
adaptation of Fricke's method to our case is not practical in view of
 the fact that the distance of a point of $\mathcal H_n$ from the 
cusps relative to the Riemannian metric is infinite. In the second 
place, it is advantageous to adapt Fricke's method only when the 
fundamental domain is compact, whereas we know that the fundamental 
domain is not compact in our case. Moreover, Siegel's method of construction
of the fundamental domain uses the deep and intrinsic properties of 
algebraic number fields.

\subsubsection{B. General Number Fields}
Guided by Siegel's discussion on totally real fields,  
we are now ready to construct  fundamental domains for general number fields. 
We largely follow [Ge] for the presentation even though our field may not 
be totally real.

So we are dealing with  rank two 
$\mathcal O_K$-lattices whose underlying projective modules $P$ are all
given by the same $P=P_\frak a:=\mathcal O_K\oplus \frak a$ for a fixed 
fractional $\mathcal O_K$-ideal $\frak a$.
This then leads to  the space $SL(\mathcal O_K\oplus\frak a)
\Big\backslash \Big(\mathcal H^{r_1}\times\mathbb H^{r_2}\Big).$ 

To facilitate ensuring discussion, recall that
for $\bold\tau=(z_1,\ldots,z_{r_1};P_1,\cdots,P_{r_2})
\in {\mathcal H}^{r_1}\times{\mathbb H}^{r_2}$, we set 
$$\mathrm{ImJ}(\tau):=\Big(\Im(z_1),\ldots,\Im(z_{r_1}),J(P_1),\ldots, 
J(P_{r_2})\Big)\in\mathbb R^{r_1+r_2}$$ where
$\Im(z)=y$ resp. $J(P)=v$ for $z=x+iy\in\mathcal H$ resp. $P=z+vj\in\mathbb H$.
For our own convenience, we now set 
$$N(\tau):=N\Big(\mathrm{ImJ}(\tau)\Big)=\prod_{i=1}^{r_1}
\Im (z_i)\cdot\prod_{j=1}^{r_2} J(P_j)^2=\Big(y_1\cdot\ldots\cdot y_{r_1}\Big)
\cdot\Big(v_1\cdot\ldots\cdot v_{r_2}\Big)^2.$$
Then by an obvious computation, (see e.g. 2.1.1 and 2.1.2,) we have, for all 
$\gamma=\begin{pmatrix} a&b\\ c&d
\end{pmatrix}\in SL(2,K)$,  $$N\Big(\mathrm{ImJ}(\gamma\cdot\tau)\Big)
=\frac{N(\mathrm{ImJ}(\bold\tau))}{\|N(c\bold \tau+d)\|^2}.\eqno(*)$$
In particular, only the second row of $\gamma$ appears.

As the first step to construct a fundamental domain, we need to have 
a generalization of Siegel's \lq distance to cusps'. 
For this, recall that for a cusp $\eta=\left[\begin{matrix}\alpha\\ 
\beta\end{matrix}\right]\in \mathbb P^1(K)$, by the Cusp-Ideal Class
Correspondence, we have a natural corresponding  ideal class associated to
the fractional ideal  
$\frak b:=\mathcal O_K\cdot\alpha+\frak a\cdot\beta$. Moreover, by assuming 
that $\alpha,\beta$ appeared above are all contained in $\mathcal O_K$, 
as we may, we know that the corresponding stablzier group $\Gamma_\eta$ 
can be described by
$$A^{-1}\cdot \Gamma_\eta\cdot A=\bigg\{\gamma=\begin{pmatrix} u&z\\ 0&u^{-1}
\end{pmatrix}\in\Gamma: u\in U_K, z\in\frak a\frak b^{-2}\bigg\},$$
where $A\in SL(2,K)$ satisfying $A\infty=\eta$ which may be further chosen in 
the form
$A=\begin{pmatrix} \alpha&\alpha^*\\ \beta&\beta^*\end{pmatrix}\in SL(2,K)$
so that $\mathcal O_K\beta^*+\frak a^{-1}\alpha^*=\frak b^{-1}$.

Now we  define the {\it reciprocal distance 
$\mu(\eta,\bold\tau)$  from a point
$\bold\tau=\Big(z_1,\ldots,z_{r_1};P_1,\cdots,P_{r_2}\Big)$
 in ${\mathcal H}^{r_1}\times{\mathbb H}^{r_2}$  to the cusp 
$\eta=\left[\begin{matrix}\alpha\\ \beta\end{matrix}\right]$  in 
$\mathbb P^1(K)$} by 
$$\begin{aligned}\mu(\eta,\bold\tau):=&
N\Big(\frak a^{-1}\cdot (\mathcal O_K\alpha+\frak a\beta)^2\Big)
\cdot\frac{\Im(z_1)\cdots \Im(z_{r_1})\cdot J(P_1)^2\cdots 
J(P_{r_2})^2}{\prod_{i=1}^{r_1}|(-\beta^{(i)}z_i+
\alpha^{(i)})|^2\prod_{j=1}^{r_2}
\|(-\beta^{(j)}P_j+\alpha^{(j)})\|^2}\\
=&\frac{1}{N(\frak a\frak b^{-2})}\cdot\frac{N(\mathrm{ImJ}(\bold\tau))}
{\|N(-\beta\bold \tau+\alpha)\|^2}.\end{aligned}$$ 
This is well-defined. Indeed, if
$\eta=\left[\begin{matrix}\alpha\\ \beta\end{matrix}\right]=
\left[\begin{matrix}\alpha'\\ \beta'\end{matrix}\right]$ in $\mathbb P^1(K)$, 
then, there exists $\lambda\in K^*$ such that $\alpha'=\lambda\cdot\alpha,\
\beta'=\lambda\cdot\beta$. Therefore, $\mu(\eta,\bold\tau)$ in terms of 
$\left[\begin{matrix}\alpha'\\ \beta'\end{matrix}\right]$ is given by
$$\frac{1}{N(\frak a\frak {b'}^{-2})}\cdot\frac{N(\mathrm{ImJ}(\bold\tau))}
{\|N(-\beta'\bold \tau+\alpha')\|^2}$$ where $\frak b'=\mathcal O_K
\alpha'+\frak a\beta'=(\lambda)\cdot\frak b$.
Hence, $\mu(\eta,\bold\tau)$ in terms of 
$\left[\begin{matrix}\alpha'\\ \beta'\end{matrix}\right]$ becomes 
$$\begin{aligned}&\frac{N(\lambda)^2}{N(\frak a\frak {b}^{-2})}\cdot
\frac{N(\mathrm{ImJ}(\bold\tau))}
{N(\lambda)^2\cdot \|N(-\beta\bold \tau+\alpha)\|^2}\cr
=&\frac{1}{N(\frak a\frak {b}^{-2})}\cdot\frac{N(\mathrm{ImJ}(\bold\tau))}
{\|N(-\beta\bold \tau+\alpha)\|^2},\end{aligned}$$ 
which is nothing but $\mu(\eta,\bold\tau)$ in terms of 
$\left[\begin{matrix}\alpha\\ \beta\end{matrix}\right]$. We are done.

As such, our definition is clearly a generalization and more importantly a 
normalization of Siegel's distance to cusps. In particular, this definition is
enviromentally free. Say no assumption such as $\alpha,\beta$ are 
$\mathcal O_K$-integers is needed.

Just as for the case of totally 
real fields, this distance plays also a key role in the sequel.
Before go further, let us show how basic properties work here.
\vskip 0.20cm
\noindent
{\bf {\large Lamma 1.}} {\it $\mu$ is invariant under the action of 
$SL(\mathcal O_K\oplus\frak a)$. That is to say,
$$\mu(\gamma\eta,\gamma\bold\tau)=
\mu(\eta,\bold\tau),\qquad\forall \gamma\in SL(\mathcal O_K\oplus\frak a).$$}
Proof. By the well-defined argument above, we may simply assume that
for a cusp $\eta$, $\alpha, \beta$ are fixed.
Then the proof is based on the following observation.
For the cusp $\eta=\left[\begin{matrix}\alpha\\ \beta\end{matrix}\right]\in
\mathbb P^1(K)$, we may choose $A_\eta=\left(\begin{matrix}\alpha&\alpha^*\\ 
\beta&\beta^*\end{matrix}\right)\in SL(2,K)$ such that $A\infty=\eta$.
(Surely, $A_\eta$ is not unique, however this does not matter.) 
Clearly, $A_\eta^{-1}=\left(\begin{matrix}\beta^*&-\alpha^*\\ 
-\beta&\alpha\end{matrix}\right).$ Therefore, by defintion,
$$\mu(\eta,\bold\tau)
=\frac{1}{N(\frak a\frak b^{-2})}\cdot
N\Big(\mathrm{ImJ}\big(A_\eta^{-1}(\bold\tau)\big)\Big).\eqno(**)$$
Note that now even $A_\eta$ is not unique, as said above, 
with a fixed $\tau$, from (*),
$N\Big(\mathrm{ImJ}\big(A_\eta^{-1}(\bold\tau)\big)\Big)$ depends only on
the second row of $A_\eta^{-1}$, which is simply $(-\beta,\alpha)$, uniquely
determined by the cusp $\eta$.

With (**), the proof may be completed easily as follows. First, let us 
consider 
the factor $N(\frak a\frak b^{-2})$. Clearly, with the change from $\eta$ to
$\gamma\eta$ for $\gamma\in SL(\mathcal O_K\oplus\frak a)$, the fractional 
ideal $\frak a\frak b^{-2}$ does not really change, so this factor remains 
unchanged. Therefore, it suffices to consider the second factor
$N\Big(\mathrm{ImJ}\big(A_\eta^{-1}(\bold\tau)\big)\Big)$. 
By an easy calculation, $A_{\gamma\eta}=\gamma A_\eta$. Consequently,
$$\begin{aligned}A_{\gamma \eta}^{-1}(\gamma \bold\tau)
=&\Big(\gamma A_\eta\Big)^{-1}(\gamma\bold\tau)
=A_\eta^{-1}\gamma^{-1}(\gamma\eta)=A_\eta^{-1}
\Big(\gamma^{-1}\gamma\eta\Big)\\
=&A_\eta^{-1}(\eta).\end{aligned}$$ Done.
\vskip 0.20cm
\noindent
{\bf{\large Lemma 2}.} {\it There exists a positive constant $C$ 
depending only on $K$ and $\frak a$ 
such that if $\mu(\eta,\bold\tau)>C$ and 
$\mu(\eta',\bold\tau)>C$ for $\bold\tau\in  
{\mathcal H}^{r_1}\times{\mathbb H}^{r_2}$ and 
$\eta,\,\eta'\in\mathbb P^1(K)$, then $\eta=\eta'$.}
\vskip 0.20cm
\noindent
{\bf Remark.} An effective version of this lemma will be given in 2.5.4 below.

\noindent
Proof. Set $\mu(\eta,\bold\tau)=\frac{1}{N(\frak a\frak b^{-2})}\cdot\frac{1}
{\Delta(\eta,\tau)}$. 
Since $N(\frak a^{-1}\frak b)\geq N(\frak a^{-1})$, it suffices
to show that {\it there exists a positive constant $c$ 
depending only on $K$ 
such that if $\Delta(\eta,\bold\tau)<c$ and 
$\Delta(\eta',\bold\tau)<c$ for $\bold\tau\in  
{\mathcal H}^{r_1}\times{\mathbb H}^{r_2}$ and 
$\eta,\,\eta'\in\mathbb P^1(K)$, then $\eta=\eta'$.}

By the Cusp-Ideal Class correspondence and the 
invariance property just proved, we can write 
$\eta=\left[\begin{matrix}\alpha\\ \beta\end{matrix}\right],$ 
$\ \eta'=\left[\begin{matrix}\alpha'\\ \beta'\end{matrix}\right]$ 
with $\mathcal O_K$-integers $\alpha,\beta,\alpha',\beta'$ 
such that $\frak b:=\mathcal O_K\alpha+\frak a\beta$ and 
 $\frak b':=\mathcal O_K\alpha'+\frak a\beta'$ have norm 
less than a constant $C$ depending only on $K$. For every
 $(r_1+r_2)$-tuple $(t_1,\cdots,t_{r_1+r_2})$ of non-zero
 real numbers, by Dirichlet's Unit Theorem, there exists a 
unit $\varepsilon\in K$ such that 
$$\Big|t_i\varepsilon^{(i)}\Big|\leq c\cdot|N(t)|^{\frac{1}{r_1+r_2}}$$ 
where $N(t):=\prod_{i=1}^{r_1}t_i\cdot\prod_{j=r_1+1}^{r_1+r_2}t_j^2$ 
with  $c$ a constant depending only on $K$.
 Hence, after multiplying $\alpha$ and $\beta$ by a suitable uint, 
we have $$\begin{aligned}\max\Big\{&\,\Im(z_i)^{-1}
\Big|-\beta^{(i)}z_i+\alpha^{(i)}\Big|,\ \,
J(P_j)^{-2}\Big\|-\beta^{(j)}P_j+\alpha^{(j)}\Big\|^2\Big\}\\
&\leq c\cdot\Delta(\eta,\bold\tau)^{-\frac{1}{r_1+r_2}}\cdot
 C^{\frac{2}{r_1+r_2}}\leq c\cdot T^{-\frac{1}{r_1+r_2}}\cdot 
C^{\frac{2}{r_1+r_2}}.\end{aligned}$$ This gives
$$\begin{aligned}\max\Big\{&\,\Big|-\beta^{(i)}\Re(z_i)+\alpha^{(i)}\Big|
\cdot\Im(z_i)^{-1/2},\ \,
\Big\|-\beta^{(j)}Z(P_j)+\alpha^{(j)}\Big\|\cdot J(P_j)^{-1}\Big\}\\
&\qquad\leq c^{1/2}\cdot T^{-\frac{1}{2(r_1+r_2)}}\cdot 
C^{\frac{1}{r_1+r_2}}\end{aligned}$$
and $$\max\Big\{\,\Big|\beta^{(i)}\Big|\cdot\Im(z_i)^{1/2},\ \,
\Big\|\beta^{(j)}\Big\|\cdot J(P_j)\Big\}\leq 
c^{1/2}\cdot T^{-\frac{1}{2(r_1+r_2)}}\cdot C^{\frac{1}{r_1+r_2}}.$$
For $\alpha'$ and $\beta'$, we obtain similar inequalities. But now, 
for real places
$$\begin{aligned}\alpha^{(i)}(\beta')^{(i)}-\beta^{(i)}(\alpha')^{(i)}
=&\Big(-\beta^{(i)}\Re(z_i)+\alpha^{(i)}\Big)\Im(z_i)^{-1/2}
\cdot (\beta')^{(i)}\Im(z_i)^{1/2}\\
&\qquad-\Big(-(\beta')^{(i)}\Re(z_i)+(\alpha)^{(i)}\Big)\Im(z_i)^{-1/2}
\cdot \beta^{(i)}\Im(z_i)^{1/2},\end{aligned}$$ while for complex places,
$$\begin{aligned}\alpha^{(j)}(\beta')^{(j)}-\beta^{(j)}(\alpha')^{(j)}
=&\Big(-\beta^{(j)}Z(P_j)+\alpha^{(j)}\Big)J(P_j)\cdot
(\beta')^{(j)}J(P_j)\\
&\qquad-\Big(-(\beta')^{(j)}Z(P_j)+(\alpha')^{(j)}\Big)J(P_j)\cdot
\beta^{(j)}J(P_j).\end{aligned}$$ Consequently
$$N\Big(\alpha\beta'-\beta\alpha'\Big)\leq (2c)^{r_1+r_2}
\cdot T^{-1}\cdot C^2.$$ So if $T>(2c)^{r_1+r_2}\cdot C^2$, 
the norm of the algebraic integr $\alpha\beta'-\beta\alpha'$ 
has absolute value less than 1, that is, 
$\alpha\beta'-\beta\alpha'=0$. This implies that  $\eta=\eta'$ as desired.
\vskip 0.30cm
More correctly, we should consider $\frac{1}{\mu(\eta,\bold\tau)^{1/2}}$ as 
the \lq distance' of $\bold\tau$ to the cusp $\eta$. For example, if
$\eta=\infty$, the distance is just $\frac{1}{N(\bold\tau)^{1/2}}
\cdot \frac{1}{N(\frak a)^{1/2}}$, since by definition, 
$\mu(\infty,\bold\tau)=\frac{N(\mathcal O_K\cdot 1+\frak a\cdot 0)^2N
(\bold\tau)}{|N(-0\bold\tau+1)|^2}=N(\bold\tau)$. As also for totally real 
fields, this distance is universally bounded as well.
\vskip 0.20cm
\noindent
{\bf{\large Lemma 3}.} {\it There exists a positive real number $T:=T(K)$ 
depending only on $K$ such that for 
$\bold\tau\in {\mathcal H}^{r_1}\times {\mathbb H}^{r_2}$, 
there exists a cusp $\eta$ such that $\mu(\eta,\bold\tau)>T$.}

\noindent
Proof. Since $N(\frak a^{-1}\frak b^2)\geq N(\frak a^{-1})$, 
and there are finitely many
inequivalent cusps,
it is sufficient to find a solution of $\alpha,\beta$ in $\mathcal O_K$
satisfying the inequality 
$$\Big|N(-\beta\bold\tau+\alpha)\Big|^2\cdot 
N\Big(\mathrm{ImJ}(\bold\tau)\Big)^{-1}\leq T^{-1}.$$ 
Consider the inequalities
$$\begin{aligned}\Big|-\beta^{(i)}\Re(z_i)+\alpha^{(i)}\Big|\cdot \Im(z_i)^{-1/2}\leq& c_i,\\
 \Big|\beta^{(i)}\Big|\cdot \Im(z_i)^{1/2}\leq& d_i,\qquad i=1,\cdots, r_1\\
\Big\|-\beta^{(j)}Z(P_j)+\alpha^{(j)}\Big\|\cdot J(P_j)^{-1}
\leq&  c_j,\\
\Big\|\beta^{(j)}\Big\|\cdot J(P_j)\leq&  d_j,\qquad j=1,\cdots,r_2,
\end{aligned}$$
which we may write, using a $\mathbb Z$-basis 
$\omega_1,\cdots, \omega_{r_1+r_2}$ of $\mathcal O_K$ as a system of 
$r_1+2r_2$ linear inequalities 
(by changing the last $r_2$ to the $2r_2$ inequalities 
involving only real numbers with respect to complex conjugations).
According to a theorem of Minkowski, we can find a solution
 $\alpha=\sum a_i \omega_i,\ \beta=\sum b_i\omega_i$ with $a_i,\, b_i\in 
\mathbb Z$ provided that $\Big(\prod c_i\cdot \prod d_j^2\Big)$ is no less 
than the absolute of the determinant of this system. Clearly, this absolute
value is simply $\Big|\omega_i^{(k)}\big|^2=\Delta_K$, the discriminant
 of $K$. So we can take $c_i=d_j=\Delta_K^{\frac{1}{r_1+2r_2}}$, 
and hence $T=2^{r_2}\cdot \Delta_K.$ This completes the proof.
\vskip 0.30cm
Now for the cusp $\eta=\left[\begin{matrix}\alpha\\ 
\beta\end{matrix}\right]\in\mathbb P^1(K)$, we define the 
\lq sphere of influence' of $\eta$ by 
$$F_\eta:=\Big\{\bold\tau\in {\mathcal H}^{r_1}\times
{\mathbb H}^{r_2}:\mu(\eta,\bold\tau)\geq 
\mu(\eta',\bold\tau),\forall\eta'\in \mathbb P^1(K)\Big\}.$$

\noindent
{\bf{\large Lemma 4}.} {\it The action of $SL(\mathcal O_K\oplus\frak a)$ 
in the interior $F_\eta^0$ of $F_\eta$ reduces to that 
of the isotropy group $\Gamma_\eta$ of $\eta$, i.e., if 
$\bold\tau$ and $\gamma\bold\tau$ both belong to $F_\eta^0$, 
then $\gamma\bold\tau=\bold\tau.$}

\noindent
Proof. We have $$\begin{matrix} \mu(\gamma^{-1}\eta,
\bold\tau)&\leq&\mu(\eta,\bold\tau)\\
\|&&\|\\
\mu(\eta,\gamma\bold\tau)&\geq&\mu(\gamma \eta,\gamma\bold\tau)
\end{matrix}$$ for $\bold\tau,\, \gamma\bold\tau\in F_\eta^o$, 
and the inequalities are strict if $\gamma\eta\not=\eta$.

Consequently, the boundary of $F_\eta$ consists of pieces of \lq generalized 
isometric circles' given by equalities 
$\mu(\eta,\bold\tau)=\mu(\eta',\bold\tau)$ with $\eta'\not=\eta$.
\vskip 0.30cm
Using above discussion, we arrive at the following way to 
decompose the orbit space $SL(\mathcal O_K\oplus\frak a)
\Big\backslash \Big({\mathcal H}^{r_1}\times{\mathbb H}^{r_2}\Big)$ into
$h$ pieces glued in some way along pants of their boundary. 
\vskip 0.20cm
\noindent
{\bf {\Large Theorem.}} {\it With the same notation as above,
let $$i_\eta:\Gamma_\eta\Big\backslash F_\eta\hookrightarrow 
SL(\mathcal O_K\oplus\frak a)\Big\backslash\Big({\mathcal H}^{r_1}
\times{\mathbb H}^{r_2}\Big)$$ be the natural map. Then 
$$SL(\mathcal O_K\oplus\frak a)\Big\backslash 
\Big({\mathcal H}^{r_1}\times{\mathbb H}^{r_2}\Big)=\cup_\eta 
i_\eta\Big(\Gamma_\eta\Big\backslash F_\eta\Big),$$
where the union is taken over a set of $h$ cusps 
representing the ideal classes of $K$. Each piece corresponds
to an ideal class of $K$.}
\vskip 0.30cm
Note that the action of $\Gamma_\eta$ on ${\mathcal H}^{r_1}
\times{\mathbb H}^{r_2}$ is free. Consequently, all fixed 
points of $SL(\mathcal O_K\oplus\frak a)$ on ${\mathcal H}^{r_1}
\times{\mathbb H}^{r_2}$ lie on the boundaries of $F_\eta$.
 
Further, we may give a more precise description of the fundamental domain, 
based on our understanding of the fundamental domains for stablizer groups of
cusps. To state it, denote by $\eta_1,\,\ldots,\,\eta_h$ 
inequivalent cusps for the action of $SL(\mathcal O_K\oplus\frak a)$ on 
${\mathcal H}^{r_1}\times{\mathbb H}^{r_2}$. Choose $A_{\eta_i}\in SL(2,K)$
such that $A_{\eta_i}\infty=\eta_i,\ i=1,2,\ldots,h$. Write $\bold S$ for
the norm-one hypersurface
$\bold S:=\Big\{y\in \mathbb R_{>0}^{r_1+r_2}:N(y)=1\Big\}$, and 
$\bold S_{U_K^2}$ for
 the action  of $U_K^2$ on $\bold S$. Denote by $\mathcal T$ a fundamental 
domain for the action of the translations by elements of
 $\frak a\frak b^{-2}$ on $\mathbb R^{r_1}\times\mathbb C^{r_2}$, 
and $$\bold E:=\Big\{\bold\tau\in {\mathcal H}^{r_1}\times
{\mathbb H}^{r_2}:
\mathrm{ReZ}\,(\bold\tau)\in {\mathcal T},\  
\mathrm{ImJ}\,(\bold\tau)\in  \mathbb R_{>0}\cdot
\bold S_{U_K^2}\Big\}$$ 
for a fundamental domain for the action of 
$A_\eta^{-1}\Gamma_\eta A_\eta$ on ${\mathcal H}^{r_1}\times
{\mathbb H}^{r_2}$. Easily, we know that 
the intersections of $\bold E$ 
with $i_\eta(F_\eta)$ are connected. Consequently,
we have the following
\vskip 0.20cm
\noindent
{\bf{\Large Theorem.$'$}} (1) {\it $A_\eta^{-1} {\bold E}\cap F_\eta$ is a fundamental domain for the 
action of $\Gamma_\eta$ on $F_\eta$ which we call $D_\eta$;}

\noindent
(2) {\it There exist $\alpha_1,\,\cdots,\,\alpha_h\in 
SL(\mathcal O_K\oplus\frak a)$ such that 
$\cup_{i=1}^h\alpha(D_{\eta_i})$ is connected and hence a
 fundamental domain for $SL(\mathcal O_K\oplus\frak a)$.}
\vskip 0.20cm
We may put this concrete discussion on fundamental domains 
in a more theoretical way. For this, we first introduce
a natural geometric truncation for the fundamental domain.
So define a compact manifold with boundary 
$$S_T:=SL(\mathcal O_K\oplus\frak a)\Big\backslash
\Big\{\bold\tau\in {\mathcal H}^{r_1}\times{\mathbb H}^{r_2}:
\mu(\eta,\bold\tau)\leq T\ \forall \eta\in \mathcal 
C_{SL(\mathcal O_K\oplus\frak a)}\Big\},$$ where $C_{SL(\mathcal O_K\oplus\frak a)}$ denotes the collections of cusps, and $T$ is so 
large that for all cusps $\eta$, $\ W(\eta,T):=\Big\{\bold\tau\in 
{\mathcal H}^{r_1}\times{\mathbb H}^{r_2}:
\mu(\eta,\bold\tau)\leq T\Big\}$ is contained in $F_\eta$, 
so disjoint for different classes $\eta$ and $\eta'$. 
Clearly, then the boundary $\partial S_T$ consists of $h$ component 
manifold $i_\eta\Big(\Gamma_\eta\Big\backslash \partial W(\eta,T)\Big)$ 
of dimension $2r_1+3r_2-1$.
Moreover, let
$\Sigma:=\Big\{(t_1,\cdots,t_{r_1};s_1,\cdots, s_{r_2})\in 
\mathbb R_{>0}^{r_1+r_2}:
\,\prod_{i=1}^{r_1} t_i\prod_{i=1}^{r_2} s_j^2=1\Big\}$ act on 
$\mathbb R^{r_1}\times\mathbb C^{r_2}$ by component-wise 
multiplication. The semi-direct product 
$\mathcal E=\Big(\mathbb R^{r_1}\times\mathbb C^{r_2}\Big)\times \Sigma$ 
acts on $ {\mathcal H}^{r_1}\times{\mathbb H}^{r_2}$ by
$$\Big((u_i,v_j),(t_i,s_j)\Big)\cdot \big(\bold \tau=(z_i;P_j)\big)
:=\Big(\lambda_iz_i+u_i;s_jP_j+v_j\Big).$$ The boundary 
$\partial W(\infty,T)$ is a partial homogeneous space
 for this semi-direct product. We view 
$A_\eta^{-1}\Gamma_\eta A_\eta\Big\backslash\partial W(\infty,Y)$ 
as the quotient of $\mathcal E$ by the discrete subgroup 
$A_\eta^{-1}\Gamma_\eta A_\eta$. It is a $r_1+2r_2$-torus bundle 
over $U_K^2\Big\backslash \Sigma$ with fiber 
$\mathbb R^{r_1}\times\mathbb C^{r_2}$ modulo the translations 
in $A_\eta^{-1}\Gamma_\eta A_\eta$. The manifold with boundary
 $S_T$ is homotopically equivalent to 
$SL(\mathcal O_K\oplus\frak a)\Big\backslash\Big( 
{\mathcal H}^{r_1}\times{\mathbb H}^{r_2}\Big)$. (See e.g. [Ga].) 
Consequently, we have 
$$SL(\mathcal O_K\oplus\frak a)\Big\backslash\Big({\mathcal H}^{r_1}
\times{\mathbb H}^{r_2}\Big)=S_T\cup_{\partial S_T}
\Big(\partial S_T\times[0,\infty)\Big),$$ i.e., 
$SL(\mathcal O_K\oplus\frak a)\Big\backslash\Big(\mathcal
 H^{r_1}\times{\mathbb H}^{r_2}\Big)$ is topologically 
a manifold with $h$ \lq ends' of the form 
$T^{r_1+2r_2}$-bundle over $T^{r_1+r_2-a}\times [0,\infty).$ 
\vskip 0.30cm
With all this, we may end our long discussion on the 
fundamental domain for the action of 
$SL(\mathcal O_K\oplus\frak a)$ on ${\mathcal H}^{r_1}
\times{\mathbb H}^{r_2}$. The essentials are, of 
course, that a fundamental domain may be given as 
$S_Y\cup \mathcal F_1(Y_1)\cup\cdots\cup \mathcal F_h(Y_h)$ 
with $\mathcal F_i(Y_i)=A_i\cdot \widetilde  {\mathcal F}_i(Y_i)$ and
 $$\widetilde{\mathcal F}_i(Y_i):=\Big\{\bold\tau\in {\mathcal H}^{r_1}
\times{\mathbb H}^{r_2}:\mathrm{ReZ}(\bold\tau)\in\Sigma,\,
\mathrm{ImJ}(\bold\tau)
\in\mathbb R_{>T}\cdot \bold S_{U_K^2}\Big\}.$$
Moreover, all $\mathcal F_i(Y_i)$'s are  
disjoint from each other when $Y_i$ are sufficiently large.

\section{Stability}

\subsection{Upper Half Plane}

So we are working with rank two $\mathbb Z$-lattice of volume 1. The space, 
i.e., the moduli space of all such 
lattices, is simply $SL(2,\mathbb Z)\Big\backslash SL(2,\mathbb R)\Big/SO(2)$,
or better, $SL(2,\mathbb Z)\Big\backslash\mathcal H$. For it,
we have a well-known  fundamental domain $\mathcal D$ 
whose closure is given by 
$\overline{\mathcal D}:=\Big\{z\in\mathcal H:|z|\geq 1, 
 |x|\leq\frac {1}{2}\Big\}$. Our question then is:
\vskip 0.20cm
\noindent
{\it What are the points in $\mathcal D$ corresponding to isometric 
classes of rank 2 semi-stable lattices of volume 1}? 

The answer is given by  classical reduction theory.
For any rank two $\mathbb Z$-lattice $\Lambda$ of volume 1 in $\mathbb R^2$ 
(equipped with the standard 
Euclideal metric), fix $\bold x\in\Lambda\backslash\{0\}$ such that its length 
gives the first Minkowski successive minimum
$\lambda_1=\lambda_1(\Lambda)$ of $\Lambda$. Then via rotation when 
necessary, we may assume that $\bold x=(\lambda_1,0)$. 
Furthermore, classical reduction theory tells us that 
$\frac{1}{\lambda_1}\Lambda$ 
is simply the lattice of the volume $\lambda_1^{-2}=:y_0$ generated by the 
vectors $(1,0)$ and $\omega:=x_0+iy_0\in\overline{\mathcal D}$. 
In particular,
with one generator $(1,0)$ being fixed, all lattices are parametrized
by only one vector, i.e., the (other) generator 
$\omega=x_0+iy_0\in \overline{\mathcal D}$. 
Consequently, our problem now becomes: 
\vskip 0.20cm
\noindent
{\it What are the points $\omega\in\overline{\mathcal D}$ whose
corresponding lattices, i.e.,
those generated by $(1,0)$ and $\omega$, are semi-stable}?

To answer this, set $\mathcal D_T:=\Big\{z\in\overline{\mathcal D}:
y=\Im(z)\leq T\Big\}$. 
Then by the above discussion,  up to points on the boundary,
the points in $\mathcal D_T$ are in one-to-one corresponding with  rank two 
$\mathbb Z$-lattices (in $\mathbb R^2$) of 
volume one whose first Minkowski successive minimums $\lambda_1$ satisfying 
$\lambda_1^{-2}\leq T$, since $\lambda_1^{-2}=y_0\leq T$. 
Write this condition in a better form, we have
$\lambda_1(\Lambda)\geq T^{-1/2}$, 
or equivalently, $\mathrm{deg}(\Lambda)\leq \frac{1}{2}\log T$. 
Then what we have just siad may be restated in a more 
theoretical form as the following
\vskip 0.20cm
\noindent
{\bf {\Large Fact}} (VI$)_{\mathbb Q}$  
({\bf Grometric Truncation=Algebraic Truncation}) 
{\it Up to a subset of measure zero, there is a natural one-to-one 
and onto morphism 
$$\mathcal M_{\mathbb Q,2}^{\leq \frac{1}{2}\log T}\Big[1\Big]
\simeq \mathcal D_T,$$ 
where $\mathcal M_{\mathbb Q,2}^{\leq \frac{1}{2}\log T}\Big[1\Big]$ 
denotes the moduli space of rank two 
$\mathbb Z$-lattices $\Lambda$ of volume 1 (over $\mathbb Q$) whose 
sublattices of rank one all have degrees $\leq\frac {1}{2}\log T$.
In particular, 
$$\mathcal M_{\mathbb Q,2}^{\leq 0}\Big[1\Big]=\mathcal M_{\mathbb Q,2}
\Big[1\Big]\simeq \mathcal D_1.$$}

That is to say, the moduli space of rank 2 semi-stable lattices of volume 1 
corresponds to the part $\mathcal D_1$ of $\mathcal D$ bounded under the 
line $y=1$.

\subsection{Upper Half Space Model}

Here we are supposed to work with imaginary quadratic fields. Our question is:

\noindent
{\it In $SL(\mathcal O_K\oplus\frak a) \backslash \bold\mathbb H$, where are 
rank two semi-stable $\mathcal O_K$-lattices of volume  $N(\frak a)\Delta_K$}?
However, as we cannot really gain anything by assuming that the fields involved
are imaginary quadratic fields, so we omit the entire discussion here by going 
directly to the most general case.

\subsection{Rank Two $\mathcal O_K$-Lattices: Level Two}

We start with our discussion by citing a result of Tsukasa Hayashi [Ha]. 

Let $\Lambda$ be a rank two $\mathcal O_K$-lattice of volume 
$N(\frak a)\cdot\Delta_K$ with underlying projective module 
$\mathcal O_K\oplus \frak a$.
Recall that, by definition, $\Lambda$  is semi-stable if for any rank one 
$\mathcal O_K$-sublattice $\Lambda_1$ of $\Lambda$,  equipped with the induced 
metric,  $$\mathrm{Vol}(\Lambda_1)^2\geq N(\frak a)\Delta_K.$$
To understand this condition, let us first understand the structure of 
rank one $\mathcal O_K$-sublattices $\Lambda_1$ of $\Lambda$.

By the Lemma in \S I.1, any rank one $\mathcal O_K$-submodule of $\Lambda$ has 
the form $\frak c\cdot\left(\begin{matrix} x\\ y\end{matrix}\right)$ where 
$\frak c$ is a fractional $\mathcal O_K$-ideal and 
$\frak c\cdot \left(\begin{matrix} x\\ y\end{matrix}\right)\in 
\mathcal O_K\oplus\frak a.$ 

Set now $\frak b=\mathcal O_K x+\frak a^{-1}y$. Since
$\frak c\cdot x\in \mathcal O_K,\ \frak c\,y\in\frak a$, we have 
$$\frak b\cdot\frak c\subset 
\Big(\mathcal O_K x+\frak a^{-1} y\Big)\cdot\frak c
=\frak c\cdot x+\frak a^{-1}(\frak c\cdot 
y)\subset\mathcal O_K+\frak a^{-1}\cdot\frak a=\mathcal O_K.$$ 
Therefore, $$\frak c\subset\frak b^{-1}.$$
This then proves (1) of the following
\vskip 0.20cm
\noindent
{\bf{\large Lemma}.}([Ha]) (1) {\it Any rank one sublattice of 
$\Lambda=\Big(\mathcal O_K\oplus\frak a,\rho_\Lambda\Big)$ is contained in 
$\frak b^{-1}\left(\begin{matrix} x\\ y\end{matrix}\right)\cap\Lambda$ where 
$\left(\begin{matrix} x\\ y\end{matrix}\right) \in K^2
\Big\backslash\bigg\{\left(\begin{matrix} 0\\ 0
\end{matrix}\right)\bigg\}$ and $\frak b=\mathcal O_K x+\frak a^{-1}y$;}

\noindent
(2) {\it $\Lambda$ is semi-stable if and only if 
$$\prod_{\sigma\in S_\infty}\bigg\|\left(\begin{matrix} x_\sigma\\ 
y_\sigma\end{matrix}\right)\bigg\|_{\Lambda_\sigma}^2\geq 
N\Big(\frak a\frak b^{2}\Big) 
=N\Big(\mathcal O_Kx+\frak a^{-1}y\Big)\cdot N\Big(\mathcal O_Ky+\frak ax\Big),
\qquad\forall \left(\begin{matrix} x\\ y\end{matrix}\right)\in K^2
\Big\backslash\bigg\{\left(\begin{matrix} 0\\ 0
\end{matrix}\right)\bigg\}.$$}
Proof. From (1), it suffices to check the 
semi-stable condition for all rank one sublattices $\Lambda_1$
induced from the submodules $\frak b^{-1}\left(\begin{matrix} x\\ y
\end{matrix}\right)$, where 
$\left(\begin{matrix} x\\ y\end{matrix}\right)\in K^2
\Big\backslash\bigg\{\left(\begin{matrix} 0\\ 0
\end{matrix}\right)$ with  $\frak b:=\mathcal O_Kx+\frak a^{-1} y$. 
Now, by the Arakelov-Riemann-Roch formula,
$$\mathrm{Vol}(\Lambda_1)=N(\frak c)\cdot
\Delta_K^{1/2}\cdot\prod_\sigma\bigg\|\left(\begin{matrix} x_\sigma\\ 
y_\sigma\end{matrix}\right)\bigg\|_{\Lambda_\sigma}.$$
Therefore, the semi-stable condition becomes
$$\bigg(\Big(N\big(\frak b^{-1}\big)\Delta_K^{1/2}\Big)\cdot 
\prod_{\sigma\in 
S_\infty}\bigg\|\left(\begin{matrix} x_\sigma\\ 
y_\sigma\end{matrix}\right)\bigg\|\bigg)^2\geq N(\frak a)\cdot\Delta_K.$$ 
That is to say, 
$$\begin{aligned}
&\prod_{\sigma\in S_\infty}\bigg\|\left(\begin{matrix} x_\sigma\\ 
y_\sigma\end{matrix}\right)\bigg\|_{\Lambda_\sigma}^2
\geq N(\frak a\frak b^2)=N\Big(\frak a(\mathcal O_Kx+\frak a^{-1} y)\cdot
\frak b\Big)\\
=&N\Big((\mathcal O_Ky+\frak a x)\frak b\Big)\\
=&N\Big(\frak ax+\mathcal O_Ky\Big)\cdot 
N\Big(\mathcal O_Kx+\frak a^{-1}y\Big).\end{aligned}$$ 
This completes the proof.

\subsection{Rank Two $\mathcal O_K$-Lattices: Normalization or Convention}

In the discussion on semi-stable lattices above, for a vector
$\left(\begin{matrix} x\\ y\end{matrix}\right)\in K^2
\Big\backslash\bigg\{\left(\begin{matrix} 0\\ 0
\end{matrix}\right)$, we introduced the fractional ideal $\frak b$ 
to be $$\frak b:=\mathcal O_K x+\frak a^{-1} y,$$ 
while in the discussion on cusps and ideal class correspondence,
for a cusp 
$\eta:=\left[\begin{matrix} \alpha\\ \beta\end{matrix}\right]
\in\mathbb P^1(K)$, the ideal class is defined to be the one associated
to the fractional idea $\frak b$  defined to be 
$$\frak b:=\mathcal O_K\alpha+\frak a\beta.$$ 
There is a discrepency among these two definitions of $\frak b$.
On the other hand, the similarity among the above definitions of $\frak b$
suggests that there may be some intrinsic relations between stability
and cusps. This is indeed the case. But before we expose this, let us make sure
that in our discussion, right normalizations and compactible notations
are used. It is for this reason, we add this rather elementary subsection.
(Technically, the point is, as it will become clear later, about how left 
multiplication and right multiplication are co-operated with each other.)
So rather then give all details here, we often use examples to indicate 
what are points we want the reader to see.
The experienced reader may skip this subsection entirely.

\subsubsection{A. $SL(2,\mathbb Z)$ Acts on the Upper Half Plane}

For $z\in \mathcal H,\gamma=\left(\begin{matrix} a&b\\ c&d\end{matrix}\right)
\in SL(2,\mathbb Z)$, set $\gamma\cdot z:=\frac{az+b}{cz+d}.$ This left action of 
$SL(2,\mathbb Z)$ is compactible with the following action of 
$SL(2,\mathbb Z)$ on $\mathbb Z^2$ whose elements are written as column 
vectors. That is to say, we have, for 
$\left(\begin{matrix} \alpha\\ \beta\end{matrix}\right)\in\mathbb Z^2$,
$$\gamma\,\left(\begin{matrix} \alpha\\ \beta\end{matrix}\right)=
\left(\begin{matrix} a&b\\ c&d\end{matrix}\right)\cdot 
\left(\begin{matrix} \alpha\\ \beta\end{matrix}\right)=
\frac{a\alpha+b\beta}{c\alpha+d\beta}=\frac{a\frac{\alpha}{\beta}+b}
{c\frac{\alpha}{\beta}+d}.$$
As such, under the natural quitient map $\mathbb Z^2\Big\backslash
\bigg\{\left(\begin{matrix} 0\\ 0\end{matrix}\right)\bigg\}\to 
\mathbb P^1(\mathbb Q)$
sending $\left(\begin{matrix} \alpha\\ \beta\end{matrix}\right)$ to 
$\left[\begin{matrix} \alpha\\ \beta\end{matrix}\right]$, the above action of
 $SL(2,\mathbb Z)$ on $\mathbb Z^2$ decends to the following action of 
$SL(2,\mathbb Z)$ on $\mathbb P^1(\mathbb Z)$;
$$\gamma\cdot \left[\begin{matrix} \alpha\\ \beta\end{matrix}\right]
=\left[\begin{matrix} a\alpha+b\beta\\ c\alpha+d\beta\end{matrix}\right].$$
Here as usual, with the following identifications, we view 
$\mathbb Q$ as a subset of $\mathbb P^1(\mathbb Q)$:
$$\left[\begin{matrix} \alpha\\ \beta\end{matrix}\right]=\left[\begin{matrix} 
\frac{\alpha}{\beta}\\ 1\end{matrix}\right]\mapsto\frac{\alpha}{\beta}.$$
That is to say, $$\mathbb P^1(\mathbb Q)=\mathbb Q\cup\{\infty\}
=\bigg\{\left[\begin{matrix} \frac{\alpha}{\beta}\\ 1\end{matrix}\right]:
z\in\mathbb Q\bigg\}\cup\bigg\{\left[\begin{matrix} 1\\ 
0\end{matrix}\right]\bigg\}$$ with 
$z=\left[\begin{matrix} z\\ 1\end{matrix}\right]$ for $z\in\mathbb Q$ and 
$\infty=\frac{1}{0}=\left[\begin{matrix} 1\\ 0\end{matrix}\right]$.
 
\subsubsection{B. Identification Between 
$SL(2,\mathbb R)\Big/SO(2)$ and $\mathcal H$}

Similarly, the operation $\gamma\cdot z=\frac{az+b}{cz+d}$
for $\gamma=\left(\begin{matrix} a&b\\ c&d\end{matrix}\right)\in 
SL(2,\mathbb R)$ and $z=x+iy\in\mathcal H$ defines an action of 
$SL(2,\mathbb R)$
on the upper half plane $\mathcal H:=\Big\{z=x+iy\in\mathbb C:x, y\in\mathbb R,
 y>0\Big\}$ since $$\begin{aligned}\frac{az+b}{cz+d}=&\frac{(ax+b)(cx+d)+acy^2}
{(cx+d)^2+c^2y^2}+i\frac{ay(cx+d)-(ac+b)cy}
{(cx+d)^2+c^2y^2}\\
=&\frac{(ax+b)(cx+d)+acy^2}
{(cx+d)^2+c^2y^2}+i\frac{(ad-bc)y}
{(cx+d)^2+c^2y^2}\\
=&\frac{(ax+b)(cx+d)+acy^2}
{(cx+d)^2+c^2y^2}+i\frac{y}
{(cx+d)^2+c^2y^2}\in\mathcal H\end{aligned}$$ and 
$$\begin{aligned}
&\bigg(\left(\begin{matrix} a_2&b_2\\ 
c_2&d_2\end{matrix}\right)\left(\begin{matrix} a_1&b_1\\ c_1&d_1
\end{matrix}\right)\bigg)z\\
=&\left(\begin{matrix} a_2a_1+b_2c_1&a_2b_1+b_2d_1\\ 
c_2a_1+d_2c_1&c_2b_1+d_2d_1\end{matrix}\right)z\\
=&\frac{(a_2a_1+b_2c_1)z+
(a_2b_1+b_2d_1)}{(c_2a_1+d_2c_1)z+(c_2b_1+d_2d_1)}\\
=&\frac{a_2(a_1z+b_1)
+b_2(c_1z+d_1)}{c_2(a_1z+b_1)+d_2(c_1z+d_1)}\\
=&\left(\begin{matrix} a_2&b_2\\ 
c_2&d_2\end{matrix}\right)\Big(\frac{a_1z+b_1}{c_1z+d_1}\Big)\\
=&
\left(\begin{matrix} a_2&b_2\\ c_2&d_2\end{matrix}\right)
\bigg(\left(\begin{matrix} a_1&b_1\\ c_1&d_1\end{matrix}\right) z\bigg).\end{aligned}$$

Note that $\gamma\cdot\sqrt{-1}=\sqrt{-1}$ if and only if 
$\frac{a\sqrt{-1}+b}{c\sqrt{-1}+d}=\sqrt{-1}$, that is, $a=d,b=-c$. 
Since $a^2+b^2=\det\left(\begin{matrix} a&b\\ -b&a\end{matrix}\right)=1$
and 
$\left(\begin{matrix} a&b\\ -b&a\end{matrix}\right)^t
\cdot\left(\begin{matrix} a&b\\ -b&a\end{matrix}\right)
=\left(\begin{matrix} a^2+b^2&0\\ 0&a^2+b^2\end{matrix}\right)
=\left(\begin{matrix} 1&0\\ 0&1\end{matrix}\right),$
so, $\gamma=\left(\begin{matrix} a&b\\ -b&a\end{matrix}
\right)\in SO(2)$. Consequently,
 the stablizer group of $\sqrt{-1}$ in $SL(2,\mathbb R)$ 
is exactly $SO(2)$. Write it formally,
 $$\mathrm{Stab}_{\sqrt{-1}}SL(2,\mathbb R)
=SL(2,\mathbb R)_{\sqrt{-1}}=\bigg\{\left(\begin{matrix} a&b\\ -b&a\end{matrix}
\right)\in SL(2,\mathbb R)\bigg\}=SO(2).$$ 
Moreover, for any $z=x+iy\in\mathcal H$, there exists 
$\gamma=\left(\begin{matrix} a&b\\ c&d\end{matrix}\right)\in SL(2,\mathbb R)$ 
such that $\gamma \sqrt{-1}=z$ since if we choose a point $(c,d)$ on the 
circle $c^2+d^2=\frac{1}{y}$, then
from the equations $ca+db=\frac{x}{y}$ and $da-cb=1$, the determinant of whose
 coefficient matrix is simply $-c^2-d^2=-\frac{1}{y}$, we can
construct a real solution for $a,b$ as well. Obviously, with such choices of 
$a,\,b,\,c,\,d$ $$\gamma \sqrt{-1}=\frac{bd+ac}
{c^2+d^2}+i\frac{ad-bc}{c^2+d^2}=x+iy=z$$ and $ad-bc=1$. Consequently, 
$\gamma \sqrt{-1}=z$ has a solution in $SL(2,\mathbb R)$. (Clearly,
the choice of the matrix $\gamma$ is essential unique modulo 
a choice of points on the circle, whose group structure is exactly that of 
$SO(2)$ viewed as a group consisting of rotations.)

That is to say, the 
action of $SL(2,\mathbb R)$ on $\mathcal H$ is also transitive, i.e., 
for any two $z_1,z_2\in\mathcal H$, there exists $\gamma_{12}\in 
SL(2,\mathbb R)$ such that $\gamma_{12}z_1=z_2$. Therefore, there is  
a natural identification $SL(2,\mathbb R)\Big/SO(2)\to
 \mathcal H$  sending $[\gamma]$ to $\gamma(\sqrt{-1})$.
 
\subsubsection{C. Metrized Structures on $\mathbb R^2$}
 
This is the first place we have to be a bit careful: The point is how to view 
a vector in $\mathbb R^2$ so that the consideration for matrices is 
compactible with that of
the action of $SL(2,\mathbb R)$ on $\mathcal H$. (By saying this, i.e., by 
assuming that the determinant is 1, we in practical term assume that the 
volumes of the lattices are fixed to be 1.)

Let us start with the column consideration. Say, write an element 
$\left(\begin{matrix} x\\ y\end{matrix}\right)$ for a 2 dimensional real 
vector space $\mathbb R^2$, and denote its norm is by 
$\bigg\|\left(\begin{matrix} x\\ y\end{matrix}\right)\bigg\|_{\rho}.$  Then,
there is a natural map from $SL(2,\mathbb R)$ to the 
space of metric structures on $\mathbb R^2$ defined by 
$g\mapsto g^t\cdot g$. Indeed, then,  
$$\bigg\|\left(\begin{matrix} x\\ y\end{matrix}\right)\bigg\|_{\rho(g)}=
(x,y)g^t\cdot g\left(\begin{matrix} x\\ y\end{matrix}\right)=
\bigg(g\left(\begin{matrix} x\\ y\end{matrix}\right)\bigg)^t\cdot 
\bigg(g\left(\begin{matrix} x\\ y\end{matrix}\right)\bigg).$$ 
Moreover, suppose $g_1,g_2\in SL(2,\mathbb R)$ give the same 
metric structure so that $$(x,y)g_1^t\cdot g_1\left(\begin{matrix} x\\ 
y\end{matrix}\right)=(x,y) g_2^t\cdot g_2\left(\begin{matrix} x\\ 
y\end{matrix}\right),$$  then there exists an orthogonal matrix $u$ such that
$ug_1=g_2$. In particular, being multiplicated from the left,
we have the identification $SO(2)\backslash SL(2,\mathbb R)$ with the 
metric structures on $\mathbb R^2$.
\vskip 0.20cm
However, such an identification is not compactible with the identification
$SL(2,\mathbb R)\Big/SO(2)\simeq \mathcal H$ we used. 
Consequently, when we consider real vector spaces and their associated 
metrized structures, for compactibility, 
it is better to view vectors in a vector space as {\it row vectors}. 
For example, with $(x,y)\in \mathbb R_2$, we have 
$$\|(x,y)\|_{\rho(g)}^2=(x,y)g\cdot g^t\left(\begin{matrix} x\\ 
y\end{matrix}\right)=(x,y)g\cdot(uu^t)\cdot g^t\left(\begin{matrix} x\\ 
y\end{matrix}\right)=\Big((x,y)gu\Big)\cdot \Big((x,y)(gu)\Big)^t,$$ 
since then a metrized structure corresponds exactly to the class $gu$ 
with $u\in SO(2)$. 
 
Accordingly, when we talk about lattices, all vectors are  understood 
to be row vectors in order to make everything compactible.
In particular, for our rank two $\mathcal O_K$-lattices 
$\Lambda=(\mathcal O_K\oplus\frak a,\rho_\Lambda)$, 
we have to use this convention and hence write 
$(x,y)\in \mathcal O_K\oplus\frak a$ with
$$\Big\|(x,y)\Big\|_\Lambda^2=\Big\|(x,y)\Big\|_{\rho_\Lambda}^2
=\prod_{\sigma\in S_\infty}\Bigg(\Big((x_\sigma,
y_\sigma)g_\sigma\Big)\cdot\overline{\Big((x_\sigma,y_\sigma)
g_\sigma\Big)^t}\Bigg)^{N_\sigma},$$ where $N_\sigma=1$ resp. 2 if
$\sigma$ is real (resp. complex).

\subsubsection{D. Automorphism Group $\mathrm{Aut}_{\mathcal O_K}(\mathcal O_K
\oplus\frak a)$}

For $A=\left(\begin{matrix} a&b\\ c&d\end{matrix}\right)\in 
\mathrm{Aut}_{\mathcal O_K}(\mathcal O_K\oplus\frak a)$, 
we should have $$(x,y)\left(\begin{matrix} a&b\\ 
c&d\end{matrix}\right)=\Big(ax+cy,bx+dy\Big)\in \mathcal O_K\oplus\frak a$$ 
for all $(x,y)\in \mathcal O_K\oplus\frak a$. Since $x\in\mathcal O_K$ and 
$y\in \frak a$, we have then $c\in\frak a^{-1}$ and $d\in \mathcal O_K$, 
say by taking $x=0$.
Similarly, one sees that $a\in\mathcal O_K$ while $b\in\frak a$. Therefore,
we conclude that
$\mathrm{Aut}_{\mathcal O_K}(\mathcal O_K\oplus\frak a)=
GL(\mathcal O_K\oplus\frak a)$ is the subgroup contained in 
$GL(2,K)\cap \left(\begin{matrix} \mathcal O_K&\frak a\\ 
\frak a^{-1}&\mathcal O_K\end{matrix}\right)$ consisting of elements whose 
determinents are units of $K$. Therefore, 
$$SL\Big(\mathcal O_K\oplus\frak a\Big)=SL(2,K)\cap
\left(\begin{matrix} \mathcal O_K&\frak a\\ 
\frak a^{-1}&\mathcal O_K\end{matrix}\right).$$
(If the vectors involved were column vectors, we should use 
$$\left(\begin{matrix} \mathcal O_K&\frak a^{-1}\\ 
\frak a&\mathcal O_K\end{matrix}\right)$$ instead.)

\subsubsection{E. Cusp-Ideal Class Correspondence}

So we are working over $\mathcal H^{r_1}\times{\mathbb H}^{r_2}$.
With  boundaries $\partial \mathcal H=\mathbb R,\,\partial{\mathbb H}
=\mathbb C$, via natural imbedding, we have the following commutative diagram
 $$\begin{matrix}
K&\hookrightarrow& \mathbb R^{r_1}\times\mathbb C^{r_2}&=
(\partial \mathcal H)^{r_1}\times (\partial {\mathbb H})^{r_2}\\
\downarrow&&\downarrow&\\
\mathbb P^1(K)&\hookrightarrow&\mathbb P^1(\mathbb R)^{r_1}\times
\mathbb P^1(\mathbb C)^{r_2}.&\end{matrix}$$ Hence,
we may view $K$ and hence 
$\mathbb P^1(K)$ as subsets of the boundary of $\mathcal H^{r_1}
\times {\mathbb H}^{r_2}$. As before, we have used the 
convention that $\frac{\alpha}{\beta}\mapsto \left[\begin{matrix} 
\frac{\alpha}{\beta}\\ 1\end{matrix}\right]=\left[\begin{matrix} \alpha\\ 
\beta\end{matrix}\right]$ and $\infty=\left[\begin{matrix} 1\\ 
0\end{matrix}\right]$.

Now motivated by the discussion on stability,
we assume that the Cusp-Ideal Class Correspondence were given by the map 
$\left[\begin{matrix} \alpha\\ \beta\end{matrix}\right]\mapsto [\frak b]$ 
with $\frak b:=\mathcal O_K\alpha+\frak a^{-1}\beta$, instead of
 $\frak b:=\mathcal O_K\alpha+\frak a\beta$,
let us see what would happen. 

\noindent
(i) {\it It is well-defined.} 
Indeed, suppose that $\left[\begin{matrix} \alpha\\ 
\beta\end{matrix}\right]=\left[\begin{matrix} \alpha_1\\ \beta_1
\end{matrix}\right]$ with $\alpha,\beta\in K$, then there exists 
$\lambda\in K$ such that $\alpha_1=\lambda\alpha,\ \,
\beta_1=\lambda\beta,$ so that $$\frak b'=\mathcal O_K\alpha_1+
\frak a^{-1}\beta_1=\mathcal O_K(\lambda\alpha)+\frak a^{-1}(\lambda\beta)
=\lambda(\mathcal O_K\alpha+\frak a^{-1}\beta)=\lambda\frak b.$$
Hence, $[\frak b]=[\frak b']$. Done;

\noindent
(ii) {\it Surjectivity}. As said before in the main text, this is a direct 
consequence of the Chinese Reminder Theorem. No problem;

\noindent 
(iii) {\it Injectivity.} This in fact will cause a certain difficulty 
due to our 
definition of $\frak b$ above. Indeed, suppose that both elements
$\left[\begin{matrix} \alpha\\ \beta\end{matrix}\right]$ and 
$\left[\begin{matrix} \alpha_1\\ \beta_1\end{matrix}\right]$ are sending 
to the same ideal class $[\frak b]$ with $\frak b:=\mathcal O_K\alpha+
\frak a^{-1}\beta$. Then by the main text, there exist matrices 
$\left(\begin{matrix} \alpha&*\\ \beta&*\end{matrix}\right)$ and 
$\left(\begin{matrix} \alpha_1&*\\ \beta_1&*\end{matrix}\right)$ in 
$SL(2,K)$ such that
$$\left(\begin{matrix} \alpha&*\\ \beta&*\end{matrix}\right)
\left[\begin{matrix} 1\\ 0\end{matrix}\right]=\left[\begin{matrix} \alpha\\ 
\beta\end{matrix}\right],\qquad
\left(\begin{matrix} \alpha_1&*\\ \beta_1&*\end{matrix}\right)
\left[\begin{matrix} 1\\ 0\end{matrix}\right]=\left[\begin{matrix} 
\alpha_1\\ \beta_1\end{matrix}\right].$$ Then one checks that
by an obvious calculation (see also below) that
$$A=\left(\begin{matrix} \alpha_1&*\\ 
\beta_1&*\end{matrix}\right)\left(\begin{matrix} \alpha&*\\ 
\beta&*\end{matrix}\right)^{-1}\in
\left(\begin{matrix} \mathcal O_K&\frak a\\ \frak a^{-1}&\mathcal O_K
\end{matrix}\right)$$ which in general is not in 
$SL(\mathcal O_K\oplus\frak a)$ as indicated in D).

Therefore,  the above definition for $\frak b$ as 
$\mathcal O_K\alpha+\frak a^{-1}\beta$ does not compactible with our
discussion. Surely, as used in the main text, the right one is the following:
\vskip 0.20cm
\noindent
{\bf {\Large Cusp-Ideal Class Correspondence}:} {\it There is a  natural bijection 
between cusps and ideal classes. More precisely, we have $$\begin{matrix}
SL(\mathcal O_K\oplus\frak a)\Big\backslash \mathbb P^1(K)&\to& CL(K)\\
\left[\begin{matrix} \alpha\\ \beta\end{matrix}\right]&\mapsto&
\mathcal O_K\alpha+\frak a\beta.\end{matrix}$$}

In particular, take the matrices $\left(\begin{matrix} \alpha_1&\alpha_1^*\\ 
\beta_1&\beta^*\end{matrix}\right)$ and $\left(\begin{matrix} 
\alpha&\alpha^*\\ \beta&\beta^*\end{matrix}\right)$ with $\alpha_1,\,
\alpha\in \frak b,\ \beta_1,\,\beta\in\frak a^{-1}\frak b,$ and
$\alpha_1^*,\,\alpha^*\in \frak a\frak b^{-1},
\ \beta_1^*,\,\beta^*\in\frak b^{-1}.$
We have, $$\begin{aligned}
A=&\left(\begin{matrix} \alpha_1&\alpha_1^*\\ 
\beta_1&\beta_1^*\end{matrix}\right)\left(\begin{matrix} \alpha&\alpha^*\\ 
\beta&\beta^*\end{matrix}\right)^{-1}\\
=&\left(\begin{matrix} \alpha_1&
\alpha_1^*\\ \beta_1&\beta_1^*\end{matrix}\right)\left(\begin{matrix} 
\beta^*&-\alpha^*\\ -\beta&\alpha\end{matrix}\right)\\
\subset&
\left(\begin{matrix} \frak b&\frak a\frak b^{-1}\\ \frak a^{-1}\frak b
&\frak b^{-1}\end{matrix}\right)\left(\begin{matrix} \frak b^{-1}&
\frak a\frak b^{-1}\\ \frak a^{-1}\frak b&\frak b\end{matrix}\right)\\
\subset& 
\left(\begin{matrix} \mathcal O_K&\frak a\\ \frak a^{-1}&\mathcal O_K
\end{matrix}\right)\end{aligned},$$ which is 
compactible with our convention on 
the group $SL(\mathcal O_K\oplus\frak a)$.

\subsubsection{F. Actions of $SL(\mathcal O_K\oplus\frak a)$ on Lattices and 
on $\mathcal H^{r_1}\times{\mathbb H}^{r_2}$}

With the above discussion, the actions of $SL(\mathcal O_K\oplus\frak a)$ 
on both lattices and on points in $\mathcal H^{r_1}
\times{\mathbb H}^{r_2}$
are given by the matrix multiplication from the left. In other words, we have
$$\begin{matrix} \Big(SL(2,\mathbb R)/SO(2)\Big)^{r_1}\times 
\Big(SL(2,\mathbb C)/SU(2)\Big)^{r_2}&\to&\Big\{(\mathcal O_K\oplus\frak a),
\rho_\Lambda=g\cdot \overline{g^t})\Big\}\Big/\equiv\\
\downarrow&&\downarrow\\
\Big(SL(2,\mathbb R)/SO(2)\Big)^{r_1}\times \Big(SL(2,\mathbb C)/SU(2)
\Big)^{r_2}&\to&\Big\{(\mathcal O_K\oplus\frak a),\rho_\Lambda=g\cdot 
\overline{g^t})\Big\}\Big/\equiv\end{matrix}$$
given by
$$\begin{matrix}[g]&\mapsto&g\cdot g^t\\
\downarrow&&\downarrow\\
[Ag]&\mapsto&Ag\cdot\overline{(Ag)}^t.\end{matrix}$$

\subsubsection{G. Stablizer Groups of Cusps}

Clearly $\left(\begin{matrix} \alpha&\alpha^*\\ \beta&\beta^*\end{matrix}
\right)\left[\begin{matrix} 1\\ 0\end{matrix}\right]=\left[\begin{matrix} 
\alpha\\ \beta\end{matrix}\right]$ which we call $\eta$ for simplicity.  
Then the stablizer group $\Gamma_\eta=\left(\begin{matrix} \alpha&\alpha^*\\ 
\beta&\beta^*\end{matrix}\right)\Gamma_\infty \left(\begin{matrix} 
\alpha&\alpha^*\\ \beta&\beta^*\end{matrix}\right)^{-1}.$
Note that when we set $\frak b=\mathcal O_K\alpha+\frak a\beta$, then 
$\alpha\in\frak b,\ \beta\in\frak a^{-1}\frak b$ while 
$\alpha^*\in\frak a\frak b^{-1},\ \beta^*\in\frak b^{-1}$. Thus 
$$\Gamma_\infty\subset\left(\begin{matrix} \frak b&\frak a
\frak b^{-1}\\ \frak a^{-1}\frak b&\frak b^{-1}\end{matrix}\right)^{-1}
\left(\begin{matrix} \mathcal O_K&\frak a\\ \frak a^{-1}&\mathcal O_K
\end{matrix}\right)\left(\begin{matrix} \frak b&\frak a\frak b^{-1}\\ 
\frak a^{-1}\frak b&\frak b^{-1}\end{matrix}\right).$$
Consequently,
$$\Gamma_\infty\subset\left(\begin{matrix} \frak b^{-1}&\frak a\frak b^{-1}\\ 
\frak a^{-1}\frak b&\frak b\end{matrix}\right)
\left(\begin{matrix} \mathcal O_K&\frak a\\ \frak a^{-1}&\mathcal O_K
\end{matrix}\right)\left(\begin{matrix} \frak b&\frak a\frak b^{-1}\\ 
\frak a^{-1}\frak b&\frak b^{-1}\end{matrix}\right)\subset 
\left(\begin{matrix} \mathcal O_K&\frak a\frak b^{-2}\\ 
\frak a^{-1}\frak b^2&\mathcal O_K\end{matrix}\right).$$
Hence, 
$$\Gamma_\infty\subset \left(\begin{matrix} \mathcal O_K&
\frak a\frak b^{-2}\\ \frak a^{-1}\frak b^2&\mathcal O_K\end{matrix}\right)\cap
 \left(\begin{matrix} \mathcal O_K&\frak a\\ 0&\mathcal O_K
\end{matrix}\right)$$since $$\left(\begin{matrix} a&b\\ c&d
\end{matrix}\right)\left[\begin{matrix} 1\\ 0\end{matrix}\right]
=\left[\begin{matrix} 1\\ 0\end{matrix}\right]\Leftrightarrow c=0.$$ 
Therefore,
$$\Gamma_{\left[\begin{matrix} \alpha\\ \beta\end{matrix}\right]}=\bigg\{
\left(\begin{matrix} \alpha&\alpha^*\\ \beta&\beta^*\end{matrix}\right)
\left(\begin{matrix} u&\omega\\ 0&u^{-1}\end{matrix}\right) 
\left(\begin{matrix} \alpha&\alpha^*\\ \beta&\beta^*
\end{matrix}\right)^{-1}:u\in U_K,\ \omega\in\frak a\frak b^{-2}\bigg\}.$$ 
That is to say, we have shown the following

\noindent
{\bf{\large Lemma}.} {\it The corresponding \lq lattice' for the cusp
$\left[\begin{matrix} \alpha\\ \beta\end{matrix}\right]$ with 
$\alpha, \beta\in\mathcal O_K$ (as we may asusme) is given by 
$\frak a\frak b^{-2}$ with $\frak b=\mathcal O_K\alpha+\frak a\beta$.}

\subsection{Rank Two $\mathcal O_K$-Lattices: Stability and Distance to Cusps}

In this subsection, we expose an intrinsic relation in Geometric Arithmetic, 
which connects stability and distance to cusps in a very beautiful way.

Assume that $\Lambda=(\mathcal O_K\oplus\frak a,\rho_\Lambda)$ is semi-stable.
 Then for any non-zero element $(x,y)\in K\oplus K$, set 
$\frak b_0:=\mathcal O_K x+
\frak a^{-1}y$ so that $x\in\frak b_0,\ y=\frak a\frak b_0$. 
Thus
$\frak b_0^{-1}x\subset\mathcal O_K$ and $\frak b_0^{-1}y\subset\frak a$ and 
$$\frak b_0^{-1}\Big(x,y\Big)\subset \Big(\frak b_0^{-1}x,\frak b_0^{-1}y\Big)
\subset \mathcal O_K\oplus\frak a.$$ Moreover, if $P_1$ is a projective 
$\mathcal O_K$-submodule of rank 1 in 
$\mathcal O_K\oplus\frak a$, then $P_1=\frak c(x,y)$ with $\frak c$ a 
fractional ideal and $(x,y)\in K\oplus K\backslash\{(0,0)\}.$ Since 
$\frak c x\subset\mathcal O_K$ and 
$\frak c y\subset\frak a$, we have $$\frak c\cdot \frak b_0=\mathcal O_K \cdot 
\frak c x+\frak a^{-1}\cdot \frak c y\subset
\mathcal O_K \cdot \mathcal O_K+\frak a^{-1}\cdot \frak a=\mathcal O_K.$$ 
Hence $\frak c\subset\frak b_0^{-1}$. Consequently, $P_1=\frak c\Big(x,y\Big)
\subset \frak b_0^{-1}\Big(x,y\Big).$
Therefore,

\noindent
(i) {\it $\frak b_0^{-1}(x,y)$ is a projective 
$\mathcal O_K$-submodule of rank 1 
in $\mathcal O_K\oplus\frak a;$} and

\noindent
(ii) {\it Any projective $\mathcal O_K$-submodule of rank 1 in 
$\mathcal O_K\oplus\frak a$ is contained in $\frak b_0^{-1}(x,y)$.}
\vskip 0.30cm
Consequently, the semi-stability condition becomes
$$\Big(\mathrm{Vol}(\frak b_0^{-1}(x,y),\rho_\Lambda)\Big)^2\geq 
\mathrm{Vol}(\mathcal O_K\oplus\frak a,\rho_\Lambda).$$ That is,
$$\Big(N(\frak b_0)^{-2}\cdot(\Delta_K^{\frac{1}{2}})^2\Big)\cdot
\Big\|(x,y)\Big\|_{\rho_\Lambda}^2\geq N(\frak a)\cdot 
\Delta_K^{2\times \frac{1}{2}}$$ or better
$$\Big\|(x,y)\Big\|_{\rho_\Lambda}^2\geq N(\frak a\frak b_0^2).\eqno(*)$$

On the other hand, 
for $g_\Lambda=\left(\begin{matrix} a&b\\ c&d\end{matrix}\right)$
such that $\rho_\Lambda=\rho(g_\Lambda)$,
$$\Big\|(x,y)\Big\|_{\rho_\Lambda}^2=\Big\|(x,y)g_\Lambda\Big\|^2=
\prod_{\sigma\in S_\infty}\Big||a_\sigma x_\sigma+c_\sigma y_\sigma|^2
+|b_\sigma x_\sigma+d_\sigma y_\sigma|^2\Big|^{N_\sigma}=N
\bigg(\left(\begin{matrix} *&*\\ x&y\end{matrix}\right)\left(\begin{matrix} 
a&b\\ c&d\end{matrix}\right)\mathrm{ImJ}\bigg)^{-1},\eqno(**)$$
where $N_\sigma=1$ resp. 2 if $\sigma$ is real resp. complex,
and $\mathrm{ImJ}:=(i,\cdots,i,\, j,\cdots,j)\in\mathcal H^{r_1}\times
\mathbb H^{r_2}$ with $i=\sqrt{-1}\in\mathcal H$ and $j=(0,0,1)\in\mathbb H$.
(Recall that we have set $N(\tau):=N(\mathrm{ImJ}(\tau)).$

Indeed, 
$$\begin{aligned}&\left(\begin{matrix} *&*\\ x&y\end{matrix}\right)\left(\begin{matrix} 
a&b\\ c&d\end{matrix}\right)\mathrm{ImJ}
=\left(\begin{matrix} *&*\\ ax+cy&bx+dy\end{matrix}\right)\mathrm{ImJ}\\
=&\bigg(
\left(\begin{matrix} *&*\\ a_{\sigma_1}x_{\sigma_1}+c_{\sigma_1}y_{\sigma_1}
&b_{\sigma_1}x_{\sigma_1}+d_{\sigma_1}y_{\sigma_1}\end{matrix}\right)\,i,
\cdots,
\left(\begin{matrix} *&*\\ a_{\sigma_{r_1}}x_{\sigma_{r_1}}+c_{\sigma_{r_1}}
y_{\sigma_{r_1}}
&b_{\sigma_{r_1}}x_{\sigma_{r_1}}+d_{\sigma_{r_1}}y_{\sigma_{r_1}}
\end{matrix}\right)\,i,\\
&\qquad \left(\begin{matrix} *&*\\ a_{\tau_1}x_{\tau_1}+c_{\tau_1}y_{\tau_1}&b_{\tau_1}
x_{\tau_1}+d_{\tau_1}y_{\tau_1}\end{matrix}\right)j,\cdots, \left(\begin{matrix} *&*\\ a_{\tau_{r_2}}x_{\tau_{r_2}}+c_{\tau_{r_2}}y_{\tau_{r_2}}&b_{\tau_{r_2}}
x_{\tau_{r_2}}+d_{\tau_{r_2}}y_{\tau_{r_2}}\end{matrix}\right)j\bigg)\\
&\hskip 5.0cm\in\mathcal H^{r_1}\times
\mathbb H^{r_2},\end{aligned}$$ where $\sigma_{1},\cdots,\sigma_{r_1}$ 
(resp. $\tau_1,\cdots,\tau_{r_2}$) denote real places (resp. complex places) 
in $S_\infty$.
From here, to get (**), we use the following obvious calculations:

\noindent
(a) For reals,
if $M=\left(\begin{matrix} A&B\\ C&D\end{matrix}\right)
\in SL(2,\mathbb R)$, 
for $z=X+Yi\in \mathcal H$ with $X, Y\in\mathbb R$, set 
$M(X+iY)=X^*+Y^*i$ with $X^*,\, Y^*\in\mathbb R$. 
Then   $$Y^*:=\frac{Y}{(CX+D)^2+C^2Y^2}.$$ In particular, 
when applied to the local factor
for real $\sigma$ in (**), we have $C=a_\sigma x_\sigma+c_\sigma y_\sigma,\, 
D=b_\sigma x_\sigma+d_\sigma y_\sigma$ and $X=0, Y=1$. Therefore,
the corresponding $Y^*$ is simply 
$$\begin{aligned}&\frac{1}{\Big((a_\sigma x_\sigma+b_\sigma y_\sigma)\cdot 0+
(b_\sigma x_\sigma+d_\sigma y_\sigma)\Big)^2+(a_\sigma x_\sigma+c_\sigma 
y_\sigma)^2\cdot 1^2}\\
=&\frac{1}{(b_\sigma x_\sigma+d_\sigma y_\sigma)^2+(a_\sigma x_\sigma
+c_\sigma y_\sigma)^2},\end{aligned}$$ as desired;

\noindent
(b) For complexes,
if $M=\left(\begin{matrix} A&B\\ C&D\end{matrix}\right)
\in SL(2,\mathbb C)$, 
for $P=Z+Vj\in \mathbb H$ with $Z\in\mathbb C, V\in\mathbb R$, set 
$M(Z+Vj)=Z^*+V^*i$ with $Z^*\in\mathbb C,\, V^*\in\mathbb R$. 
Then   $$V^*:=\frac{V}{|CZ+D|^2+|C|^2V^2}=\frac{V}{\|CP+D\|^2}.$$
In particular, when applied to the local factor for complex $\tau$ in (**),
we have $C=a_\tau x_\tau+c_\tau y_\tau,\, 
D=b_\tau x_\tau+d_\tau y_\tau$ and $Z=0,\, V=1$. Therefore,
the corresponding $(V^*)^2$ is simply 
$$\begin{aligned}&\Big(\frac{1}
{|(a_\tau x_\tau+c_\tau y_\tau)\cdot 0+(b_\tau x_\tau+d_\tau y_\tau)|^2
+|a_\tau x_\tau+c_\tau y_\tau|^2\cdot 1^2}
\Big)^2\\
=&
\Big(\frac{1}{|b_\tau x_\tau+d_\tau y_\tau|^2+|a_\tau x_\tau+c_\tau y_\tau|^2}
\Big)^2\end{aligned}$$ as desired.

Consequently, the relation (**), together with (*), implies  

\noindent
(iii) {\it The lattice 
$\Lambda=\Big(\mathcal O_K\oplus\frak a,\rho_\Lambda(g)\Big)$ with 
$g:=\left(\begin{matrix} a&b\\ c&d\end{matrix}\right)$ 
is semi-stable if and only if for any non-zero $(x,y)\in K\oplus K$,
$$N\bigg(\left(\begin{matrix} *&*\\ x&y\end{matrix}\right)
\left(\begin{matrix} a&b\\ c&d\end{matrix}\right)\mathrm{ImJ}\bigg)\cdot N
(\frak a\frak b_0^2)\leq 1,$$
where $\frak b_0:=\mathcal O_Kx+\frak a^{-1}y$.}
\vskip 0.30cm
But, by definition, for the lattice $\Lambda=\Big(\mathcal O_K\oplus\frak a,
\rho(g_\Lambda)\Big)$, the corresponding point 
$\tau_\Lambda\in\mathcal H^{r_1}\times\mathbb H^{r_2}$ is given by
$g_\Lambda \Big(\mathrm{ImJ}\Big)$. Hence, we have the following equivalent

\noindent
(iii$'$) {\it The lattice 
$\Lambda=\Big(\mathcal O_K\oplus\frak a,\rho_\Lambda(g)\Big)$
is semi-stable if and only if for any non-zero $(x,y)\in K\oplus K$,
$$N\bigg(\left(\begin{matrix} *&*\\ x&y\end{matrix}\right)
\tau_\Lambda \bigg)\cdot N
(\frak a\frak b_0^2)\leq 1,$$
where $\frak b_0:=\mathcal O_Kx+\frak a^{-1}y$.}

Set now $x=-\beta$ and $y=\alpha$. Then $\frak b_0=\mathcal O_K\beta+
\frak a^{-1}\alpha.$ 
In particular, $\beta\in\frak b_0$ and $\alpha\in\frak a\frak b_0$. 
So if we define 
$$\frak b:=\frak a\frak b_0.$$ 
Then $\alpha\in\frak b,\ \beta\in\frak a^{-1}\frak b$, and
$$\mathcal O_K\alpha+\frak a\beta\subset\frak b=\frak a\frak b_0=\frak a
\cdot\Big(\mathcal O_K\beta+\frak a^{-1}\alpha\Big)
\subset \mathcal O_K\cdot\frak a
\beta+\frak a^{-1}\cdot\frak a\alpha=\frak a\beta+\mathcal O_K\alpha.$$ 
Therefore, $\frak b=\mathcal O_K\alpha+\frak a\beta$, and
$$\frak a\frak b_0^2=\frak a\cdot(\frak a^{-1}\frak b)^2=
\frak a^{-1}\frak b^2.$$
Consequently, the semi-stability condition (iii$'$) becomes
for any cusp $\eta=\left[\begin{matrix}\alpha\\ 
\beta\end{matrix}\right]\in\mathbb P^1(K)$,
$$\mu(\eta,\tau_\Lambda)=
N\bigg(\left(\begin{matrix} *&*\\ -\beta&\alpha\end{matrix}\right)
\tau_\Lambda \bigg)\cdot 
N(\frak a^{-1}\frak b^2)\leq 1.$$ Or better,
in terms of {\it distance to cusp},
$$d(\eta,\tau_\Lambda):=\frac{1}{\mu(\eta,\tau_\Lambda)}\geq 1.$$
In this way, we arrive at the following  fundamental result, which
exposes a beautiful intrinsic relation between stability and the 
distance to cusps.
\vskip 0.30cm
\noindent
{\bf {\Large Fact}} (VII) 
{\it The lattice $\Lambda$ is semi-stable if and only if the distances  of
corresponding point $\tau_\Lambda\in \mathcal H^{r_1}\times{\mathbb H}^{r_2}$ 
to all cusps are all bigger or equal to 1.}
\vskip 0.30cm
\noindent
{\bf Remark.} One can never overestimate the importance of this relation.
Being stable, lattices should be away from cusps. More generally,
while the stability condition is defined in terms of sublattices, the 
relation above transforms these volumes inequalities 
in terms of  distances to cusps. In a more theoretical term 
for higher rank lattices, the essence of this fact is that,
sublattices and cusps, as two different aspects of parabolic subgroups,  
are naturally corresponding to each other: the stability conditions for 
various sublattices 
are naturally related with generalized distances to all types of cusps.

\subsection{Moduli Space of Rank Two Semi-Stable $\mathcal O_K$-Lattices}

For a rank two 
$\mathcal O_K$-lattice $\Lambda$, denote by $\tau_\Lambda\in
 \mathcal H^{r_1}\times\mathbb H^{r_2}$ the corresponding point. 
Then, from the Fact in the 
previous subsection,  $\Lambda$ is semi-stable if and only if 
for all cusps $\eta$,  
$d(\eta, \tau_\Lambda):=\frac{1}{\mu(\eta, \tau_\Lambda)}$ 
are bigger than or equal to 1.
This then leads to the consideration of the following
truncation of the fundamental domain $\mathcal D$ of $SL(\mathcal O_K\oplus
\frak a)\Big\backslash \Big(\mathcal H^{r_1}
\times\mathbb H^{r_2}\Big)$:
For $T\geq 1$, denote by $$\mathcal D_T:=\Big\{\tau\in\mathcal D:
d(\eta,\tau_\Lambda)\geq T^{-1},\ \forall\mathrm{cusp}\ \eta\Big\}.$$
 
The space $\mathcal D_T$ may be precisely described in terms of $\mathcal D$ 
and certain neighborhood of cusps. To explain this, we first establish the 
following

\noindent
{\bf{\large Lemma}.} {\it For a cusp $\eta$, denote by $$X_\eta(T):=\Big\{
\tau\in \mathcal H^{r_1}\times\mathbb H^{r_2}: d(\eta,\tau)<T^{-1}\Big\}.$$
Then for $T\geq 1$, $$X_{\eta_1}(T)\cap X_{\eta_2}(T)\not=\emptyset\qquad
\Leftrightarrow \eta_1=\eta_2.$$}

\noindent
{\bf Remark.} This result is an effective version of ii) of 2.4.3.A) and
Lemma 2 of 2.4.3.B.

\noindent
Proof. One direction is clear. Hence, it suffices to show that 
if $\tau\in \mathcal H^{r_1}
\times\mathbb H^{r_2}$ satisfies $d(\tau,\eta_1)<1$ and $d(\tau,\eta_2)<1,$
then $\eta_1=\eta_2$.

For this, let  
$\eta_1=\left[\begin{matrix}\alpha_1\\ \beta_1\end{matrix}\right],
\ \eta_2=\left[\begin{matrix}\alpha_2\\ \beta_2\end{matrix}\right]$ and $\frak b_1=\mathcal O_K\alpha_1+\frak a \beta_1,\ \frak b_2=\mathcal O_K\alpha_2+\frak a \beta_2.$
Clearly,  $$\begin{aligned}&N\bigg(\begin{pmatrix}*&*\\ -\beta_1&\alpha_1\end{pmatrix}\tau
\bigg)\cdot N(\frak a^{-1}\frak b_1^2)\\
&=N\bigg(\begin{pmatrix}*&*\\ 
-\beta_1&\alpha_1\end{pmatrix}\cdot\bigg(
\begin{pmatrix}\alpha_2&\alpha_2^*\\ \beta_2&\beta_2^*\end{pmatrix}\cdot
\begin{pmatrix}\beta_2^*&-\alpha_2^*\\ -\beta_2&\alpha_2
\end{pmatrix}\bigg)\cdot
\tau\bigg)\cdot N(\frak a^{-1}\frak b_1^2)\\
&=N\bigg(\bigg(\begin{pmatrix}*&*\\ -\beta_1&\alpha_1\end{pmatrix}\cdot
\begin{pmatrix}\alpha_2&\alpha_2^*\\ \beta_2&\beta_2^*\end{pmatrix}\bigg)\cdot
\begin{pmatrix}*&*\\ -\beta_2&\alpha_2\end{pmatrix}\tau\bigg)\cdot 
N(\frak a^{-1}\frak b_1^2)\\
&=N\bigg(\begin{pmatrix}*&*\\ c&d\end{pmatrix}\cdot
\begin{pmatrix}*&*\\ -\beta_2&\alpha_2\end{pmatrix}\tau\bigg)\cdot 
N(\frak a^{-1}\frak b_1^2)\\
&=\frac{N\bigg(\begin{pmatrix}*&*\\ -\beta_2&\alpha_2\end{pmatrix}\tau\bigg)}
{\bigg\|c\cdot \begin{pmatrix}*&*\\ 
-\beta_2&\alpha_2\end{pmatrix}\tau+d\bigg\|^2}\cdot 
N(\frak a^{-1}\frak b_1^2),\end{aligned}$$ where 
$c=\alpha_1\beta_2-\beta_1\alpha_2$.

We want to show that $c=0$, since then $\eta_1=\eta_2$. Thus to continue, let us recall that we have the following conditions ready to use:

$$\begin{aligned}&N\bigg(\begin{pmatrix}*&*\\ -\beta_1&\alpha_1\end{pmatrix}\tau\bigg)\cdot N(a^{-1}\frak b_1^2)>1,\\
&N\bigg(\begin{pmatrix}*&*\\ -\beta_2&\alpha_2\end{pmatrix}\tau\bigg)\cdot N(a^{-1}\frak b_2^2)>1.\end{aligned}$$

As such, set  $\tau'=\begin{pmatrix}*&*\\ -\beta_2&\alpha_2\end{pmatrix}\tau$,
then what we need to show becomes the following
\vskip 0.20cm
\noindent
{\bf{\large Lemma}$'$}. {\it With the same notaion as above, if}

\noindent
(i) $N(\tau')\cdot N(\frak a^{-1}\frak b_2^2)>1,$

\noindent 
(ii) $N(\tau')\cdot N(\frak a^{-1}\frak b_1^2)>\|c\tau'+d\|^2,$ and

\noindent
(iii) $c=\alpha_1\beta_2-\beta_1\alpha_2$ {\it with} $\alpha_1,\,\beta_1,\,\alpha_2,
\,\beta_2\in\mathcal O_K$, 

{\it Then $c=0$.}
\vskip 0.20cm
\noindent
Proof. First note that $\alpha_1\in\frak b_1,\ 
\beta_1\in\frak a^{-1}\frak b_1$ and $\alpha_2\in\frak b_2,\ 
\beta_2\in\frak a^{-1}\frak b_2$, we have $c\in\frak a^{-1}\frak b_1\frak b_2$.
Thus $$N(c)\geq N(\frak a^{-1}\frak b_1\frak b_2).\eqno(*)$$

Then we use the following
\vskip 0.30cm
\noindent
{\bf {\large Sublemma}.} {\it $\|c\tau'+d\|^2\geq N(c)^2\cdot N(\tau')^2$.}
\vskip 0.20cm
\noindent
Proof. Indeed, if suffices to prove this inequality locally. This is however
an obvious calculation. Say for real $\sigma$, by definition,
$$\Big\|c_\sigma z_\sigma+d_\sigma\Big\|^2
=(c_\sigma x_\sigma+d_\sigma)^2+c_\sigma^2
y_\sigma^2\geq c_\sigma^2y_\sigma^2,$$ done. 
(We leave the complex case to the reader.)

Thus by (ii), we have $$N(\tau')\cdot N(\frak a^{-1}\frak b_1^2)>
N(c)^2\cdot N(\tau')^2.$$ That is to say, 
$$N(\frak a^{-1}\frak b_1^2)>
N(c)^2\cdot N(\tau').$$ Consequently, by (i), we have
$$N(\frak a^{-1}\frak b_1^2)\cdot N(\frak a^{-1}\frak b_2^2) >
N(c)^2,$$ or better $N(\frak a^{-1}\frak b_1\frak b_2)>N(c),$ contrads 
with (*). This completes the proof of the Lemma.
\vskip 0.30cm
All in all, then we have exposed the following
\vskip 0.20cm
\noindent
{\bf {\Large Fact}} (VI$)_K$. {\it There is a natural identification between

\noindent
(a) the moduli space of rank two semi-stable 
$\mathcal O_K$-lattices of volume $N(\frak a)\Delta_K$
with underlying projective module 
$\mathcal O_K\oplus\frak a$ and 

\noindent
(b) the truncated compact domain 
$\mathcal D_1$ consisting of points in the fundamental domain $\mathcal D$
whose distances to all cusps are bigger than 1.}

In other words, the truncated compact domain  $\mathcal D_1$
is obtained from the fundamental domain $\mathcal D$
of $SL(\mathcal O_K\oplus
\frak a)\Big\backslash \Big(\mathcal H^{r_1}
\times\mathbb H^{r_2}\Big)$ by delecting the disjoint open neighborhoods
$\cup\cup_{i=1}^h\mathcal F_i(1)$ associated to inequivalent cusps 
$\eta_1,\eta_2,\ldots,\eta_h$,
where $\mathcal F_i(T)$ denotes the neighborhood of $\eta_i$ consisting
of $\tau\in\mathcal D$ whose distance to $\eta_i$ is strictly 
less than $T^{-1}$.

For later use, we set also 
$$\mathcal D_T:=\mathcal D\backslash \cup\cup_{i=1}^h\mathcal F_i(T),
\qquad T\geq 1.$$

\chapter{Epstein Zeta Functions and Their Fourier Expansions}

\section{Upper Half Plane}

We will follow Kubota [Kub] for the presentation. (So
we claim no originality in any sense here --  We add this elementary
section for the purpose of indicating how general theory is built from 
the classics.)

Let $\Gamma$ be a discontinuous group reduced at infinity. 
Recall that a function $f(z)$ is called an automorphic function with respect
 to $\Gamma$ if $f(\gamma z)=f(z)$ for all $\gamma\in\Gamma$.
If $\kappa$ is a cusp of $\Gamma$, then there exists an
$A\in G$ such that $A\infty=\kappa$ and such that 
$A^{-1}\Gamma_\kappa A=\Gamma_\infty$. Thus for automorphic
 $f(z)$, $\ f(Az)$ is a periodic function with period 1, i.e., 
$f(A(z+1))=f(Az)$ due to the fact that
 $\Gamma_\infty=\Big\{\left(\begin{matrix} 1&n\\ 0&1\end{matrix}\right):
n\in\mathbb Z\Big\}$, a typical element
 $\left(\begin{matrix} 1&n\\ 0&1\end{matrix}\right)$ of which acts on 
$z$ by a shift $z\mapsto z+n$, in particular, 
 $\left(\begin{matrix} 1&1\\ 0&1\end{matrix}\right)z=z+1$. Therefore,
 putting $e(x)=\exp(2\pi \sqrt{-1}x)$ and 
 $a_m(y)=\int_0^1 f(Az)e(-mx)dx,\ z=x+\sqrt{-1}y$, we have a 
Fourier expansion $$f(Az)=\sum_{m\in\mathbb Z}a_m(y) e(mx),$$
whenever $f$ satisfies some natural analaytic conditions. This 
will be called the {\it Fourier expansion} of $f$ at the cusp $\kappa$.

Let us now find an explicit formula for the Fourier expansion of 
an Eisenstein series at a cusp. Denote by
 $\Gamma$ a discontiuous group of finite  type, and by 
$\kappa_1,\cdots,\kappa_h$ a complete set of 
 inequivalent cusps of $\Gamma$. Choose $\sigma_i\in SL(2,\mathbb R)$
such that $\sigma_i\infty=\kappa_i,\ i=1,2,\ldots,h$. 
Define the {\it Eisenstein series} $E_i(z,s)$ of $\Gamma$ at $\kappa_i$  by 
$$E_i(z,s):=\sum_{\gamma\in\Gamma_i\backslash \Gamma} 
y^s(\sigma_i^{-1}\gamma z),$$ where $y(x+iy):=y$. Then 
the Fourier expansion of 
$E_i(z,s)$ at $\kappa_j$ is given in the form
$$E_i(\sigma_jz,s)=\sum_{m\in\mathbb Z} a_{ij,m}(y,s)\,e(mx),$$
with $a_{ij,m}(y,s):=\int_0^1E_i(A_j z,s)e(-mx)dx.$

To give $a_{ij,m}$ an explicit form, recall that $\Gamma$ is 
reduced at $\infty$, that is to say, 
$\Gamma_\infty$ is generated by $\left(\begin{matrix} 1&n\\ 
0&1\end{matrix}\right)$. We will calculate 
$a_{11,m}$. However, for the sake of convenience, we omit the 
index \lq 1' in the notation relative to 
$\infty$; for example, we write $E$ for $E_1$, and $a$ for $a_{11}$, and 
$$a_m(y,s)=\int_0^\infty\sum_{\gamma\in\Gamma_\infty\backslash\Gamma} 
y(\gamma z)^s e(-mx)dx.$$

To give an explicit expression for $a_m(y,s)$,
let us make use of the double coset decomposition 
$\Gamma_\infty\Big\backslash\Gamma\Big/\Gamma_\infty$. Since
 $\Gamma$, as a discontinuous group reduced at $\infty$, cannot 
contain any hyperbolic transformation of the
  form $\left(\begin{matrix} t&b\\ 0&t^{-1}\end{matrix}\right),\ t>0$, 
we see that the double coset 
  $\Gamma_\infty\gamma\Gamma_\infty$ with 
$\gamma=\left(\begin{matrix} a&b\\ c&d\end{matrix}\right)\in\Gamma$
is equal to $\Gamma_\infty \left(\begin{matrix} 1&b\\ 
0&1\end{matrix}\right) \Gamma_\infty$ whenever $c=0$.
 Furthermore, if $c\not=0$, or more specifically if $c>0$, 
which may be assumed without loss of generality,
  then $\Gamma_\infty\gamma\Gamma_\infty=\Gamma_\infty\gamma'\Gamma_\infty$ 
for $\gamma=\left(\begin{matrix} a&b\\ c&d\end{matrix}\right),\ \, 
\gamma'=\left(\begin{matrix} a'&b'\\ 
c'&d'\end{matrix}\right)$ if and only if $c=c'$ and 
$d\equiv  d' \pmod {c}$. As easily to be seen, for a 
given $c>0$, there exist only a finite number of $d$ incongruent 
$\pmod {c}$ such that $(c,d)$ is the second row of some $\gamma\in \Gamma$.
Thus we get a decomposition $$\Gamma_\infty\Big\backslash\Gamma
\Big/\Gamma_\infty=\Gamma_\infty\cup\Big(\cup_{c,d}
(\Gamma_\infty\gamma\Gamma_\infty)\Big),\qquad c>0,\ d\pmod{c},\ \, 
\gamma=\left(\begin{matrix} a&b\\
 c&d\end{matrix}\right)\in\Gamma.$$
This is a kind of {\it Bruhat decomposition}. It also follows from elementary 
arguments that $\Gamma_\infty
\gamma\gamma_0$ and $\Gamma_\infty\gamma\gamma_0'$ 
$(\gamma\not\in\Gamma_\infty)$ are different cosets 
in $\Gamma_\infty\backslash \Gamma$, whenever $\gamma$ and $\gamma'$ 
are different elements of
 $\Gamma_\infty$. Hence, $\Gamma_\infty\cap \gamma\Gamma_\infty 
\gamma^{-1}=\{1\}$ for $c\not=0$. So
$$\begin{aligned}a_m(y,s)=&\int_0^1\sum_{\gamma\in \Gamma_\infty\backslash \Gamma}
y(\gamma z)^s e(-mx)dx\\
=&\delta_{0m}y^s+
\int_{-\infty}^\infty \sum_{\gamma=\left(\begin{matrix} a&b\\ 
c&d\end{matrix}\right)\in \Gamma_\infty\backslash
 \Gamma/\Gamma_\infty, c\not=0}y(\gamma z)^s e(-mx)dx\end{aligned}$$ with 
Kronecker's $\delta$, and this latter quantity 
 is further equal to 
$$\begin{aligned}\delta_{0m}y^s&+\sum_{c,d}\int_{-\infty}^\infty 
\frac{y^s}{|cz+d|^{2s}} e(-mx)dx\\
=&\delta_{0m}y^s+\sum_{c}\frac{1}{|c|^{2s}}\Big(\sum_d e(\frac{md}{c})\Big)
\cdot\int_{-\infty}^\infty 
\frac{y^s}{|x^2+y^2|^{2s}} e(-mx)dx\\
=&\delta_{0m}y^s+y^{1-s}
\sum_{c}\frac{1}{|c|^{2s}}\cdot\sum_d 
e(\frac{md}{c})\cdot\int_{-\infty}^\infty \frac{1}{(1+t^2)^{s}} 
e(-mt)dt,\\
&\hskip 3.0cm c>0,\ d\pmod {c},\  
\left(\begin{matrix} *&*\\ c&d\end{matrix}\right)\in\Gamma.\end{aligned}$$ 
We now set $$\phi_m(s)=\sum_{c}\frac{1}
{|c|^{2s}}\cdot\sum_d e(\frac{md}{c}),\qquad
 c>0,\ d\ \pmod {c},\ \, \left(\begin{matrix} *&*\\ c&d\end{matrix}\right)
\in\Gamma$$ and
recall $$\int_{-\infty}^\infty \frac{1}{(1+t^2)^{s}} e(-nt)dt=2\pi^s
|n|^{s-\frac{1}{2}}\Gamma(s)^{-1}
K_{s-\frac{1}{2}}(2\pi|n|),\qquad n\in \mathbb R\backslash\{0\}$$ and 
$$\int_{-\infty}^\infty \frac{1}{(1+t^2)^{s}}dt=
\pi^{\frac{1}{2}}\frac{\Gamma(s-\frac{1}{2})}{\Gamma(s)}.$$ Here 
$K_s$ denotes the so-called modified
 Bessel function defined by $$K_s(z)=\frac{\pi}{2}
\frac{I_{-s}(z)-I_s(z)}{\sin\,s\pi}\quad\mathrm{with}\quad
I_s(z):=\sum_{m=0}^\infty\frac{(\frac{1}{2}z)^{s+2m}}{m!\Gamma(s+m+2)}.$$ 
We have thus obtained
$$a_m(y,s)=2\pi^s\Big|m\Big|^{s-\frac{1}{2}}\Gamma(s)^{-1}\cdot
y^{\frac{1}{2}}K_{s-\frac{1}{2}}(2\pi|m|y)\phi_m(s),\qquad m\not=0$$ 
and $$a_0(y,s)=y^s+
\phi(s)y^{1-s}$$ with
$$\phi(s):=\pi^{\frac{1}{2}}\frac{\Gamma(s-\frac{1}{2})}
{\Gamma(s)}\,\phi_0(s).$$ 
The general case for other cusps can be treated in almost the same 
way using the double coset decomposition 
$\Gamma_\infty\Big\backslash\Big(A_i^{-1}\Gamma A_i\Big)\Big/\Gamma_\infty$. 
The only thing which we need to note is that the coset
 $\Gamma_\infty\left(\begin{matrix} 1&0\\ 0&1\end{matrix}\right)\Gamma_\infty$ 
does not appear in the above 
 decomposition unless $i=j$. For, otherwise, $\kappa_i$ and $\kappa_j$, 
$i\not=j$, would be equivalent. 
 The result of the calculation in general case is
$$a_{ij,m}(y,s)=2\pi^s|m|^{s-\frac{1}{2}}\Gamma(s)^{-1}\cdot y^{\frac{1}{2}}
K_{s-\frac{1}{2}}(2\pi|m|y)\phi_{ij,m}(s),\qquad m\not=0$$ and 
$$a_{ij,0}(y,s)=\delta_{ij}y^s+\phi_{ij}(s)y^{1-s}$$
 with 
$$\phi_{ij,m}(s)=\sum_{c}\frac{1}{|c|^{2s}}\cdot\sum_d e(\frac{md}{c}),\qquad
 c>0,\ d\ \pmod {c},\ \,\left(\begin{matrix} *&*\\ c&d\end{matrix}\right)
\in A_i^{-1}\Gamma A_i,$$
and $$\phi_{ij}(s):=\pi^{\frac{1}{2}}\frac{\Gamma(s-\frac{1}{2})}{\Gamma(s)}\phi_{ij,0}(s).$$ 
The functions $\phi_{ij,m}(s)$ are all Dirichlet series in a rather general sense.

The matrix $\Phi(s):=(\phi_{ij}(s))$ of functions appearing in the 
constant term of the Fourier expansion of
 Eisenstein series has a very important manning. Due to the 
 involution map $x\mapsto x^{-1}$ of $G$, it is clear that $\Phi(s)$ is 
a symmetric matrix.
\vskip 0.30cm
In fact, we may also use the fact that $E$ is a solution of the partial 
differential equation 
$\Delta E=\lambda E$ with $\lambda=s(s-1)$ to see that the coefficients 
$a_{ij,m}(y,s)$ satisfy
$$\frac{d^2 a_{ij,m}}{dy^2}-\Big(4\pi^2m^2+\frac{\lambda}{y^2}\Big)
a_{ij,m}=0.$$ Such a second order ordinary
differential equation may be solved by using standard solutions. 
With these solutions, by looking at the
growth conditions, we will arrive  also at the conclusion that the 
constant terms are certain combinations of 
$y^s$ and $y^{1-s}$ while the rest are coming from Bessel-K functions. 
For details, please see the next section.

\section{Upper Half Space}

We will reverse the ordering from the previous section 
to presume the case at hand. The presentation follows [EGM]. (We make this 
decision  based on the same reason as the one when we use Siegel's original 
text for fundamental domains: Classics are already parts of our culture. So no
new writting is needed. On the other hand, by adding them here,
we present the reader in a {\it single volume} on how classics work; and more
importantly, based on it how new theory should be developed. 
We claim no originality 
about the classics, rather we claim the responsibility for our particular 
choices from the classics.)

Suppose $\Lambda$ is a lattice in $\mathbb C$.  Let $f:{\mathbb H}
\to\mathbb C$ be a $\Lambda$-invariant 
$C^2$-function, that is $f(P+z)=f(P)$ for all $z\in\Lambda$ and $f$ is
differentiable up to the second order, 
satisfying the differential equation $-\Delta f=\lambda f$. 
Choose $s\in\mathbb C$ with $\lambda=2s(2-2s).$ 
Assume further that $f(z+rj)$ is of polynomial growth as 
$r\to\infty$, that is $f(z+rj)=O(r^k)$ as $r\to\infty$
for some constant $k$ uniformly with respect to $z\in\mathbb C$. 
Then in case $s\not=0$, $f$ possesses a
 Fourier expansion $$f(z+rj)=a_0r^{2s}+b_0r^{2-2s}+
\sum_{0\not=\mu\in\Lambda^\vee}a_\mu r\cdot K_{2s-1}\Big(2\pi|\mu|r\Big)
\, e\Big(\langle\mu,z\rangle\Big)$$ wheras in case $s=0$, 
 $$f(z+rj)=a_0r^2+b_0r^2\log r^2+\sum_{0\not=\mu\in\Lambda^\vee}
 a_\mu r\cdot K_{2s-1}\Big(2\pi|\mu|r\Big)\cdot 
e\Big(\langle\mu,z\rangle\Big).$$ Here 
$\langle\ , \rangle$ denotes the usual scalar 
 product on $\mathbb R^2=\mathbb C$ and $\Lambda^\vee$ denotes 
the dual lattice of $\Lambda$, $$\Lambda^\vee:=
\Big\{\mu\in\mathbb C: \langle\mu,z\rangle\in\mathbb Z,\forall 
z\in\Lambda\Big\}.$$

Indeed, since the function $z\mapsto f(z+rj)$ is real analytic and 
$\Lambda$-invariant, it has a Fourier
 expansion $$f(z+rj)=\sum_{\mu\in\Lambda^\vee}q_\mu(r)\cdot 
e\Big(\langle\mu,z\rangle\Big),$$
which may be differentiated termwise up to the second order. 
Hence our formula for $\Delta$ in the 
coordinaters $z,\,r$ yields the ordinary differential equation 
$$\Big(r^2\frac{d^2}{dr^2}-r\frac{d}{dr}+
\lambda-4\pi^2|\mu|^2r^2\Big)\cdot q_\mu(r)=0\eqno(*).$$
For $\mu=0,\, s\not=0$, the function $r^{2s}, r^{2-2s}$ form a  
fundamental system of solutions of $(*)$. For
 $\mu=0,\, s=0$, we take the functions $r^2$, $r^2\log r^2$ as fundamental 
system. This gives the constant terms of
  the Fourier expansion above.

We now study what happens for $\mu\not=0$.  If $Z_n(u)$ is an 
arbitrary solution of Bessel's differential 
equation of order $n$, the function $w=u^\alpha\cdot Z_n(\beta u)$ 
satisfies the differential equation 
$$u^2\frac{d^2w}{du^2}+\Big(1-2\alpha\Big)u\frac{dw}{du}+\Big((\beta u)^2+
\alpha^2-r^2\Big)w=0.$$ As such, in our case 
we may choose $\alpha=1,\ \beta=2\pi i|\mu|,\ n=s$. Hence we arrive at 
the solution $g_\mu(r)=r Z_{2s-1}(2\pi i|\mu|r)$. The function $Z_s$ can 
be written as a linear combination of the fundamental system 
$K_s$ and $I_s$. So the general solution for us is $$g_\mu(r)
=a_\mu r\cdot K_{2s-1}(2\pi|\mu|r^2)
+b_\mu r\cdot I_{2s-1}(2\pi|\mu|r^2),$$ where $a_\mu, b_\mu$ are constants.

But $$g_\mu(r)=\frac{1}{\mathrm{Vol}(\mathcal P)}\int_{\mathcal P}
f(z+rj)\cdot e\Big(-\langle \mu,z\rangle\Big) dx\,dy$$ 
where $\mathrm{Vol}(\mathcal P)$ denotes the Euclidean area of a fundamental 
parallelopiped $\mathcal P$ of $\Lambda$. By our 
assumption that $f$ is of polynomial growth, $q_\mu(r)$ is of  
polynomial growth. The function $K_s(x)$ 
decreases exponentially, and $I_s(x)$ increases exponentially as 
$x\to\infty$. Hence $b_\mu=0$, and we are done.
\vskip 0.20cm
Next, we proceed to determine an explicit Fourier expansion for 
Eisenstein series.

By definition, if $\Gamma\subset SL(2,\mathbb C)$ is a discrete 
group and $\eta=A\infty\in\mathbb 
P^1(\mathbb C)$ is one of its cusps, the {\it Eisenstein series} of 
$\Gamma$ at $\eta$ is defined as
$$E_A(P,s):=\sum_{M\in\Gamma_\eta'\backslash\Gamma}r(AMP)^{2s}.$$ 
This series converges for $\Re(s)>1$. 
If $\eta=B^{-1}\infty\in\mathbb P^1(\mathbb C)$ is another cusp of 
$\Gamma$, the function $P\mapsto 
E_A(B^{-1}P,s)$
is invariant under the action of the lattice $\Lambda$ corresponding 
to $(B\Gamma B^{-1})_\infty'
=B\Gamma_\eta' B^{-1}$, that is, $$B\Gamma_\eta' B^{-1}
=\Bigg\{\left(\begin{matrix} 1&\omega\\ 0&1\end{matrix}\right):
\omega\in\Lambda\Bigg\}.$$ As before, we write $\langle\ ,\ \rangle$ for 
the Euclidean inner product on $\mathbb
 R^2=\mathbb C$ and $\Lambda^\vee$ for the lattice dual to $\Lambda$ 
with respect to this inner product. 
 Writing $P=z+rj\in{\mathbb H}$, being of slow growth and 
satisfying the PDE $\Delta E=-\lambda E$ as
  can be checked by standard method, the above discussion of 
Fourier expansion ensures the existence of an 
  expansion $$E_A(B^{-1}P,s)=
a_0r^{2s}+b_0r^{2-2s}+\sum_{0\not=\mu\in\Lambda^\vee}a_\mu r
\cdot K_{2s-1}(2\pi|\mu|r^2)\cdot e\Big(\langle\mu,z\rangle\Big).$$
We shall give here an explicit formula for the coefficients $a_\mu$.
 In the formulation we shall use the 
Kronecker symbol $\delta_{\eta,\zeta}=1$ if 
$\eta\equiv \zeta \pmod {\Gamma}$, 0 otherwise for cusps 
$\eta,\, \zeta$ of $\Gamma$.
\vskip 0.20cm
\noindent
{\bf {\Large Theorem 1.}} ([EGM]) {\it For $\Re(s)>1$, the 
Eisenstein series $E_A(B^{-1}P,s)$ has the Fourier expansion
$$\begin{aligned}&E_A(B^{-1}P,S)=\\
=&\delta_{\eta,\zeta}
[\Gamma_\zeta:\Gamma_\zeta']
|d_0|^{-4s}r^{2s}+\frac{\pi}
{\mathrm{Vol}(\mathcal P)^s}\Big(\sum_{\left(\begin{matrix} *&*\\ 
c&d\end{matrix}\right)\in\mathcal R}
|c|^{-4s}\Big)r^{2-2s}\\
&+\frac{2\pi^{2s}}{\mathrm{Vol}(\mathcal P)
\Gamma(2s)}\sum_{0\not=\mu\in\Lambda^\vee}
|\mu|^{2s-1}\cdot\Big(\sum_{\left(\begin{matrix} *&*\\ c&d\end{matrix}\right)
\in\mathcal R}\frac{e(-\langle\mu,
\frac{d}{c}\rangle)}{|c|^{4s}}\Big)
\cdot r\cdot 
K_{2s-1}\Big(2\pi|\mu|r^2\Big)\cdot
e\Big(\langle\mu,z\rangle\Big),\end{aligned}$$ where
$\mathcal R$ denotes a system of representatives 
$\left(\begin{matrix} *&*\\ c&d\end{matrix}\right)$ of the 
double cosets in $$\Big(A\Gamma_\zeta'A^{-1}\Big)
\Big\backslash \Big(A\Gamma B^{-1}\Big)\Big/\Big(B\Gamma_\eta'B^{-1}\Big),$$ 
such that $c\not=0$, 
$\mathcal P$ is a fundamental parallelopiped for $\Lambda$ with Euclidean 
area $\mathrm{Vol}(\mathcal P)$. If
$\eta$ and $\zeta$ are $\Gamma$-equivalent, $d_0$ is defined by 
$\left(\begin{matrix} *&*\\ 
0&d_0\end{matrix}\right)=AL_0B^{-1}$ for a certain $L_0\in\Gamma$ 
satisfying that $L_0\eta=\zeta$.}

\noindent
Proof. This will be done by a direct calculation. So being 
$\Lambda$-invariant, $E_A(B^{-1}P,S),$ $\,\Re(s)>1$ 
admits a Fourier expansion of the form $$E_A(B^{-1}P,S)
=\sum_{\mu\in\Lambda^\vee}a_\mu(r,s)\,
e\Big(\langle\mu,z\rangle\Big),\quad P=z+rj.$$
Clearly, by definition, $$a_\mu(r,s)=\frac{1}{\mathrm{Vol}(\mathcal P)}
\sum_{M\in\Gamma_\zeta'
\backslash\Gamma}\int_{\mathcal P}r\Big(AMB^{-1}(z+rj)\Big)^{2s}\cdot 
e\Big(-\langle\mu,z\rangle\Big) dx\,dy.$$

 First, we want to reduce our system of elements $AMB^{-1}$ modulo 
$B\Gamma_\eta'B^{-1}$ from the right. 
 Consider $AMB^{-1}=\left(\begin{matrix} *&*\\ c&d\end{matrix}\right)$. 
Then there exists an $M\in\Gamma$ in the
  same class such that $c=0$ if and only if $AMB^{-1}$ fixes $\infty$, 
that is, if and only if   $\eta=B^{-1}\infty$ and 
 $\zeta=A^{-1}\infty$ are equiavalent modulo $\Gamma$. 
 
 If $\eta$ and $\zeta$ are $\Gamma$-equivalent, let $L_0\in\Gamma$ 
be so chosen that $L_0\eta=\zeta$ and put
$AL_0B^{-1}= \left(\begin{matrix} *&*\\ 0&d_0\end{matrix}\right).$ Then 
for all $M\in\Gamma$ such that $AMB^{-1}$
has the form $\left(\begin{matrix} *&*\\ c&d\end{matrix}\right)$ with $c=0$, 
we have $|d|=|d_0|$, and there 
are exactly $[\Gamma_\zeta:\Gamma_\zeta']$ different elements in 
$\Gamma_\zeta'\backslash \Gamma$ with this
 property.

Now if $\mu\not=0$, the contribution of these terms to the 
integration we are computing equals to zero, 
since exp is periodic and its average over $[0,1]$ is 0; whereas 
for $\mu=0$, the contribution is equal 
to $\delta_{\eta\zeta}[\Gamma_\zeta:\Gamma_\zeta']|d_0|^{-4s}r^{2s}.$ 
Hence we are left with the computation of the sum
$$\frac{1}{\mathrm{Vol}(\mathcal P)}\sum_{AMB^{-1}=\left(\begin{matrix} *&*\\ 
c&d\end{matrix}\right),\, c\not=0
}\int_{\mathcal P}\Big(\frac{r}{\|c(z+rj)+d\|^2}\Big)^{2s}\cdot 
e\Big(\langle\mu,z\rangle\Big) dx\,dy,$$
where the summation extends over $M\in \Gamma_\zeta'\backslash \Gamma$ 
such that $c\not=0$, that is, over 
all cosets in $A\Gamma_\zeta'A^{-1}\backslash A\Gamma B^{-1}.$ Hence 
we obtain that this latest summation 
is equal to
$$\begin{aligned}\frac{1}{\mathrm{Vol}(\mathcal P)}&\sum_{\left(\begin{matrix} *&*\\ 
c&d\end{matrix}\right)\in\mathcal R}
\sum_{\omega\in\Lambda}\int_{\mathcal P}\Big(\frac{r}{\|c(z+\omega)
+d\|^2+|c|^2r^2}\Big)^{2s}\cdot 
e\Big(\langle\mu,z\rangle\Big) dx\,dy\\
=&\frac{1}{\mathrm{Vol}(\mathcal P)}\sum_{\left(\begin{matrix} *&*\\ 
c&d\end{matrix}\right)\in\mathcal R}
\int_{\mathbb C}\Big(\frac{r}{|c z+d|^2+|c|^2r^2}\Big)^{2s}\cdot 
e\Big(\langle\mu,z\rangle\Big) dx\,dy\\
=&\frac{1}{\mathrm{Vol}(\mathcal P)}\sum_{\left(\begin{matrix} *&*\\ 
c&d\end{matrix}\right)\in\mathcal R}
\frac{e(\langle\mu,\frac{d}{c}\rangle)}{|c|^{4s}}
\int_{\mathbb C}\Big(\frac{r}{|z|^2+r^2}\Big)^{2s}\cdot 
e\Big(|\mu|\cdot x\Big) dx\,dy,\end{aligned}$$ where we have applied 
an orthogonal linear transformation of $\mathbb R^2$ sending $\mu$ to 
$(|\mu|,r)$.

For $\mu=0$, we obtain then $$\begin{aligned}\int_{\mathbb C}&
\Big(\frac{r}{|z|^2+r^2}\Big)^{2s}dx\,dy\\
=r^{2-2s}&
\int_{\mathbb C}(|z|^2+1)^{-2s}dx\,dy\\
=&2\pi r^{2-2s}\int_0^\infty\frac{\rho d\rho}{(\rho^2+1)^{2s}}\\
=&\frac{\pi}{2s-1}r^{2-2s}.\end{aligned}$$ 

For $\mu\not=0$, we obtain then
$$\begin{aligned}\int_{\mathbb C}&\Big(\frac{r}{|z|^2+r^2}\Big)^{2s} \cdot
e\Big(|\mu|x\Big)dx\,dy\\
=&r^{2-2s}\int_{-\infty}^\infty\int_{-\infty}^\infty\frac{dy}
{\Big(y^2+(1+x^2)\Big)^{2s}}\cdot e\Big(|\mu|x\Big)dx\\
=&r^{2-2s}\int_{-\infty}^\infty\frac{dt}{(1+t^2)^{2s}}\int_{-\infty}^\infty 
\frac{e(r|\mu|x)}{(1+x^2)^{2s-1/2}}dx\\
=&\frac{2\pi^{2s}|\mu|^{2s-1}}{\Gamma(2s)}r\cdot K_{2s-1}(2\pi|\mu|r^2).
\end{aligned}$$
This then completes the proof.

\noindent
{\bf Remark.} As to be point out clearly, the appearence of 
$[\Gamma_\zeta:\Gamma_\zeta']$ is due to the fact that
we want to factor out the twists from elliptic points. Indeed, in the 
case of $SL(\mathcal O_K\oplus\frak a)$, these elliptic points
are induced from units, since the fields involved now are imaginary 
quadratic fields (hence by Dirichlet's Unit 
Theorem, there are only finitely many of them.) So the picture is 
coherent. We thank Elstrodt for explaining 
elliptic point is the motivation to the introduction of the group 
$\Gamma_\zeta'$. Hence, in general, a 
modified approach has to be taken since the unit group 
is really quite large.
\vskip 0.30cm
We continue this discussion by looking at closely the Fourier expansion 
of Eisenstein series for $SL(2)$ 
over imaginary quadratic integers. Write $\mathcal M=\mathcal M_K$ 
for the group of fractional ideals of 
$\mathcal O_K$. For $\frak m\in\mathcal M$, $s\in\mathbb C$ with 
$\Re(s)>1$, define
$$E_\frak m(P,s):=N\frak m^{2s}
\sum_{c,d\in K, \langle c,d\rangle=\frak m}
\Big(\frac{r}{\|cP+d\|^2}\Big)^{2s},
$$ where the summation extends over all pairs $(c,d)$ of 
generators of $\frak M$ as an $\mathcal O_K$-module,
 and $$\widehat  E_\frak m(P,s):=N\frak m^{2s}\sum_{c,d\in\frak m}'
\Big(\frac{r}{\|cP+d\|^2}\Big)^{2s},$$ where 
 the prime indicates that the summation extends over all pairs 
$c,d$ in $\frak M$ with $(c,d)\not=(0,0)$.
$E_\frak m$ and $\widehat  E_\frak m$ are called {\it Eisenstein series} for 
$SL(2,\mathcal O_K)$ associated with 
$\frak m$. One checks easily they satisfy the following properties:

\noindent
{\bf{\Large Proposition.}}
(1) {\it The Eisenstein series $E_\frak m$ and $\widehat  E_\frak m$ depend 
only on the class of $\frak m$ in $CL(K)=\mathcal M/K^*;$}

\noindent
(2) {\it They converge uniformly on compact sets for $\Re(s)>1$;}

\noindent
(3) {\it For $\Re(s)>1$, they are $SL(2,\mathcal O_K)$-invariant 
functions that satisfy the differential equation
 $\Delta E=2s(2s-2)E$.}
\vskip 0.30cm
To go further, we first express $\widehat  E_\frak m$ as a linear 
combination of the $E_\frak n$
for $\frak n\in CL(K)$. Here certain zeta functions come up that we define now.

For $\frak m,\frak n\in\mathcal M$, let $\zeta(\frak m,\frak n;s)
:=N(\frak m\frak n^{-1})^s\sum_{\lambda\in 
\frak m\frak n^{-1}}'N(\lambda)^{-s}$ and for 
$\frak m^{\#}\in CL(K)$, set $\zeta(\frak m^{\#};s):=
\sum_{\frak a\in\frak m^{\#},\frak a\subset \mathcal O_K}
N(\frak a)^{-s}.$ 

One checks easily that $\zeta(\frak m,\frak n;s)$ and 
$\zeta(\frak m^{\#};s)$ are well-defined for $\Re(s)>1$ 
and $\zeta(\frak m,\frak n;s)$ depends only on the classes of 
$\frak m$ and $\frak n$ in $CL(K)$. Note also that
$$\zeta(\frak m,\frak n;1)=\zeta(\frak m\cdot \frak n^{-1},
\mathcal O_K;s)\qquad \mathrm{and}\qquad
\zeta(\frak a,\mathcal O_K;s)=\#U_K\cdot \zeta([\frak a^{-1}],s).$$ 
Thus we have
$\zeta(\frak m,\frak n;s)=\#U_K\cdot \zeta([\frak m\frak n^{-1}],s).$
\vskip 0.20cm
\noindent
{\bf{\large Lemma 1}.} {\it  For $\frak m\in\mathcal M,\ P\in{\mathbb H}$, 
we have $$\# U_K\cdot \widehat  E_\frak m(P,s)=\sum_{[\frak n]\in CL(K)}
\zeta(\frak m,\frak n;2s)E_\frak n(P,s)\qquad\Re(s)>1.$$}

\noindent
Proof. Let $\frak n$ run through a representative system $V$ of $CL(K)$. 
Consider a pair $(\gamma,\delta)$ of
 generators of an arbitrary element $\frak n\in V$ and an arbitrary 
$\lambda\in\frak m\frak n^{-1}$, and
  consider the map
$$\Big(\lambda,(\lambda,\delta)\Big)\mapsto 
(c,d):=(\lambda \gamma,\lambda\delta)\in\frak m\oplus\frak m\backslash
 \{(0,0)\}.$$
This map is surjective, and every $(c,d)$ has precisely $\#U_K$ diferent 
inverse images. This yields the assertion.
\vskip 0.20cm
We now show that the function $E_\frak m(P,s)$ agree with the 
Eisenstein series $E_A(P,S)$ up to elementary 
factors.

First we introduce the following notation. For 
$A=\left(\begin{matrix} \alpha&\beta\\ \gamma&\delta\end{matrix}
\right)\in SL(2,K)$, let $\frak u_A=\langle\gamma,
\delta\rangle\in\mathcal M,\ \frak v_A=\langle\alpha,\beta
\rangle\in\mathcal M$. The maps $SL(2,K)\to\mathcal M$ given by 
$A\mapsto \frak u_A$ and $A\mapsto\frak v_A$ 
are surjective.
\vskip 0.20cm
\noindent
{\bf{\large Lemma 2}.} {\it If $\zeta\in\mathbb P^1(K)$ is a cusp of 
$\Gamma=SL(2,\mathcal O_K)$ and $A=\left(\begin{matrix} 
\alpha&\beta\\ \gamma&\delta\end{matrix}\right)\in SL(2,K)$ with 
$A\zeta=\infty$, then
$$E_A(P,s)=\frac{1}{2}(N\frak u_A)^{-2-2s}E_{\frak u_A}(P,s)$$ for all 
$P\in{\mathbb H},\ s\in\mathbb C$
 with $\Re(S)>1$.}

\noindent
Proof. Consider the set $L$ of pairs $(c,d)\in K^2$ which generate the 
$\mathcal O_K$-module $\frak u_A$. 
For every $(c,d)\in L$, there exists an $M\in SL(2,\mathcal O_K)$ 
such that $M\left(\begin{matrix} d\\ 
-c\end{matrix}\right)=\left(\begin{matrix} \delta\\ 
-\gamma\end{matrix}\right).$ We use this fact to construct the map
$$\phi:L\to \Gamma_\zeta'\backslash\Gamma,\qquad (c,d)\mapsto \Gamma_\zeta'M.$$ 
Note that $\phi$ is well-defined. 
If $(c,d)\in L$ and $\phi\Big((c,d)\Big)=\Gamma_\zeta' M$ with $M\in\Gamma$, 
then $AM=\left(\begin{matrix} *&*\\
c&d\end{matrix}\right)$. Conversely, if
$\Gamma_\zeta'M\in\Gamma_\zeta'\backslash \Gamma$ with  
$AM=\left(\begin{matrix} *&*\\ c&d\end{matrix}\right)$,
 then $\phi^{-1}\Big(
\Gamma_\zeta'M\Big)=\Big\{(c,d),(-c,-d)\Big\}$. This completes the proof.
\vskip 0.30cm
The explicit computation of the Fourier expansion of $E_\frak m$ 
turns out to be rather clumsy, however the
Fourier  expansion of $\widehat  E_\frak m$ can be determined much more easily.
\vskip 0.30cm
\noindent
{\bf {\Large Theorem 2.}} ([EGM]) {\it Suppose that $\frak m\in\mathcal M$ and 
$\eta=B^{-1}\infty$ with $B\in SL(2,K)$  a cusp
of $SL(2,\mathcal O_K)$. Let $\Lambda$ be the lattice in $\mathbb C$ 
corresponding to the unipotent stablizer 
$SL(2,\mathcal O_K)_\eta'$ of the cusp, that is 
$B\cdot SL(2,\mathcal O_K)_\eta'\cdot B^{-1}=\Bigg\{\left(\begin{matrix} 
1&\omega\\ 0&1\end{matrix}\right):\omega \in\Lambda\Bigg\}.$
Let $\Lambda^\vee$ be the dual lattice of $\Lambda$. Then the
 $\Lambda$-invariant function $\widehat  E_\frak m(B^{-1}P,s)$ has the 
Fourier expansion
$$\begin{aligned}&\widehat  E_\frak m(B^{-1}P,s)=\\
=&N(\frak u_B)^{2s}\zeta(\frak m,\frak u_B, 2s)\cdot r^{2s}
+\frac{\pi N(\frak m)^{2s}}{|\Lambda|^{2s-1}}\Big(\sum_{(c,d)
\in\mathcal R_0}|c|^{-4s}\Big)r^{2-2s}\\
+&\frac{2\pi^{2s}N(\frak m)^{2s}}{|\Lambda|\cdot\Gamma(2s)}
\sum_{0\not=\omega'\in\Lambda^\vee}|\omega'|^{2s-1}\cdot
\sum_{(c,d)\in \mathcal R_0}\frac{e\Big(-\langle\omega',
\frac{d}{c}\rangle\Big)}{|c|^{4s}}\cdot r\cdot 
K_{2s-1}\Big(2\pi|\omega'|r^2\Big)\cdot 
e\Big(\langle\omega',z\rangle\Big),\end{aligned}$$ 
where $\mathcal R_0$ is a maximal system of 
representatives of $(c,d)$ of $\Big(\frak m\oplus\frak m\Big)
B^{-1}/B\Gamma_\eta' B^{-1}$ with $c\not=0$.}

\noindent
Proof.  
There is first of all a Fourier 
expansion of the form $$\widehat  E_\frak m(B^{-1}P,s)=
\sum_{\omega'\in\Lambda^\vee}a_{\omega'}(r,s)\cdot 
e\Big(\langle\omega',z\rangle\Big).$$ The Fourier coefficients are computed 
as $$a_{\omega'}(r,s):=
\frac{N(\frak m)^{2s}}{|\Lambda|}\sum_{(c,d)\in 
\Big(\frak m\oplus\frak m\Big)B^{-1}, (c,d)\not=(0,0)}
\int_Q\Big(\frac{r}{\|cP+d\|^2}\Big)^{2s}\cdot 
e\Big(-\langle\omega',z\rangle\Big)\, dx\,dy$$ with $Q$ the fundamental 
parallelpoid of $\Lambda$ in $\mathbb C^2$.
Accodring to whether $c=0$ or not, we further get
$$\begin{aligned}\frac{N(\frak m)^{2s}}{|\Lambda|}&\sum_{(0,d)\in 
\Big(\frak m\oplus\frak m\Big)B^{-1}, d\not=0}
\int_Q\Big(\frac{r}{\|d\|^2}\Big)^{2s}\cdot 
e\Big(\langle\omega',x\rangle\Big)\, dx\,dy\\
&+\frac{N(\frak m)^{2s}}{|\Lambda|}\sum_{(c,d)\in \Big(\frak m\oplus\frak m\Big)
B^{-1}, c\not=0}\int_Q\Big(\frac{r}
{\|cP+d\|^2}\Big)^{2s}\cdot
e\Big(-\langle\omega',z\rangle\Big)\, dx\,dy.\end{aligned}\eqno(**)$$ The first 
term on the right hand side vanishes 
termwise for $\omega'\not=0$ and for $\omega'=0$ it is equal to 
$$N(\frak m)^{2s}\Big(\sum_{(0,d)\in 
\Big(\frak m\oplus\frak m\Big)B^{-1}, d\not=0}|d|^{-4s}\Big)r^{2s}.
\eqno(***)$$ 
If $d\not=0$, then $(0,d)\in 
\Big(\frak m\oplus\frak m\Big)B^{-1}$ if and only if $(0,d)B\in 
\Big(\frak m\oplus\frak m\Big)$, that is if and only if 
$d\frak u_B\subset \frak m$. Hence, the sum (***) is equal to 
$\frak u_B^{2s}\zeta(\frak m,
\frak u_B; 2s)\cdot r^{2s}.$ This is the first term on the 
RHS of the theorem. 

The second sum above is treated in the same way as in the proof 
of the previous Theorem.
\vskip 0.30cm
Even though higher order terms are not needed for our present discussion
on rank two zeta functions. But for
 completeness, we decide to include them. In fact, even to give a 
number theoretic interpretation for 
 the constant terms above, more works has to be done. 

Let us then fix some notation. Suppose that $\frak m\in\mathcal M$, 
and let $\eta\in\mathbb P^1(K)$ be 
a cusp of $\Gamma$, $\eta=B^{-1}\infty$ with $B=\left(\begin{matrix} 
\alpha&\beta\\ \gamma&\delta\end{matrix}
\right)\in SL(2,K)$. Let $\Lambda\subset \mathbb C$ be the lattice 
such that $$B\cdot SL(2,\mathcal O_K)_\eta'\cdot 
B^{-1}=\Bigg\{\left(\begin{matrix} 1&\omega\\ 0&1\end{matrix}\right):
\omega\in\Lambda\Bigg\}.$$
As above, denote the $\mathcal O_K$-modules generated by the row vectors of 
$B$ by $$\frak u:=\frak u_B:=\langle 
\gamma,\delta\rangle,\qquad \frak v:=\frak v_B:=\langle \alpha,\beta\rangle.$$ 
The set of row vectors $(c,d)$ in 
the sum (**) is contained in the $\mathcal O_K$-module $\mathcal L:=
\Big(\frak m\oplus\frak m\Big)B^{-1}.$ Note that
 trivially $\mathcal L\subset (\frak m\frak u\oplus\frak m\frak v).$
For $0\not= c_0\in \frak m\frak u$, let $$\mathcal L(c_0):=\Big(\{c_0\}
\times K\Big)\cap\mathcal L=\Big\{(c,d)\in \mathcal L:c=c_0\Big\}.$$

\noindent
{\bf{\large Sublemma 1}.} {\it With the preceding notations we have}

\noindent
(1) $\Lambda=\frak u^{-2}$;

\noindent
(2) {\it if $(c,d)\in\mathcal L$ and $\omega\in\Lambda$, then 
$(c,c\omega+d)\in\mathcal L$;}

\noindent
(3) {\it if $0\not=c_0\in\frak m\frak u$, then 
$\mathcal L(c_0)\not=\emptyset$;}

\noindent
(4) $\frak m\frak u^{-1}\subset \frak m\frak v$.

\noindent
Proof. (1): Clearly, $\omega\in\Lambda$ if and only if 
$$B^{-1}\left(\begin{matrix} 1&\omega\\ 0&1\end{matrix}
\right) B=\left(\begin{matrix} 1+\gamma\delta\omega&\delta^2\omega\\ 
-\gamma^2\omega&1-\gamma\delta\omega
\end{matrix}\right)\in SL(2,\mathcal O_K).$$ This holds if and only if 
$\delta^2\omega,\, \gamma^2\omega,\, \gamma\delta\omega
\in \mathcal O_K$. Since $\delta^2,\, \gamma^2,\, \gamma\delta$ generate the 
$\mathcal O_K$-module $\frak u^2$, we obtain $\lambda=\frak u^{-2}$.

\noindent
(2): We have from definition
$$(c,c\omega+d)=(c,d) B\cdot B^{-1} \left(\begin{matrix} 1&\omega\\ 
0&1\end{matrix}\right) B\cdot B^{-1}\in 
\Big(\frak m\oplus\frak m\Big)B^{-1}=\mathcal L.$$

\noindent
(3): We have $c_0=x\delta-y\gamma$ for some $x,y\in\frak m$. Defining 
$d:=-x\beta+y\alpha$, then $(c_0,d)=(x,y)B^{-1}\in\mathcal L.$

\noindent
(4): By definition, $\mathcal O\subset \frak u\frak v$. 
This completes the proof.

From (2), the group $\Lambda=\frak u^2$ acts on $\mathcal L$ by 
$(c,d)\mapsto (c,c\omega+d)$ where
$(c,d)\in\mathcal L$, $\omega\in\Lambda$.
We compute the number of orbits of the restriction of this group 
action to $\mathcal L(c_0)$.

\noindent
{\bf {\large Sublemma 2.}} {\it If $0\not=c_0\in\frak m\frak u$, then 
$$\#\Big(\mathcal L(c_0)/\Lambda\Big)
=\frac{N(c_0)}{N(\frak m) \,N(\frak u)}.$$}

\noindent
Proof. Consider the homomorphism of $\mathcal O_K$-modules 
$$\phi:\frak m\frak v\to(\frak m
 \frak u\frak v\oplus\frak m \frak u\frak v)/(\frak m\oplus \frak m), 
\qquad x\mapsto (\gamma x,\delta x)+
 (\frak m\oplus \frak m).$$ The range of $\phi$ is well-defined by (4), 
and the same result  implies also that 
 $\mathrm{Ker}\phi=\frak m\frak v\cap \frak m\frak u^{-1}=\frak m
\frak u^{-1}$ and hence 
 $c_0\frak u^{-2}\subset \mathrm {Ker}\phi$. Thus $\phi$ induces a 
homomorphism $\overline{\phi}:
 \frak m\frak v/c_0 \frak u^{-1}\to (\frak m \frak u\frak v\oplus
\frak m \frak u
 \frak v)/(\frak m\oplus \frak m)$. The element $\lambda_0:=(c_0\alpha,
 c_0\beta)+
 (\frak m\oplus \frak m)\in (\frak m \frak u\frak v\oplus\frak m 
\frak u\frak v)/(\frak m\oplus \frak m)$ 
 belongs to the image of $\overline{\phi}$ by (3) and the map 
$$\mathcal L(c_0)/\frak u^{-2}\to\overline{\phi}^{-1}
 (\lambda_0),\qquad \Big\{(c_0,c_0\omega+d):\omega\in\frak u^{-2}\Big\}
\mapsto -d+c_0\frak u^{-2}$$ is a bijection.
  Hence $$\#\Big(\mathcal L(c_0)/\frak u^{-2}\Big)
=\#\Big(\overline{\phi}^{-1}(\lambda_0)\Big)
=\#\Big(\mathrm{Ker}\overline{\phi}\Big)=\Big[\frak m 
  \frak u^{-1}:c_0\frak u^{-2}\Big]=\frac{N(c_0)}{N(\frak m) \,N(\frak u)}.$$

\noindent
{\bf {\large Sublemma 3.}} {\it We have} $$\frac{\pi N(\frak m)^{2s}}
{|\Lambda| (2s-1)}\sum_{(c,d)\in\mathcal R_0}|c|^{-4s}=
\frac{2\pi}{\sqrt{\Delta_K}}N(\frak u_B)^{2-2s}\zeta (\frak m,
\frak u_B^{-1},2s-1).$$

\noindent
Proof. The set $\mathcal R_0$ is a maximal set of representatives 
$(c,d)\in\mathcal L,\, c\not=0$ for 
the action above, and for a fixed entry $c_0$ of some 
element of $\mathcal R_0$, the number of 
different $d$ with $(c_0,d)\in\mathcal R_0$ is given by the 
Sublemma 2 above. Hence we obtain
$$\begin{aligned} N(\frak m)^{2s}&\sum_{(c,d)\in\mathcal R_0}|c|^{-4s}= 
N\frak m^{2s}\sum_{c\in \frak m\frak u}'
\frac{\#\Big(\mathcal L(c)/\Lambda\Big)}{N(c)^{2s}}\\
=& \frac{N(\frak m)^{2s-1}}{N(\frak u)}\sum_{c\in \frak m\frak u}'
Nc^{1-2s}=N\frak u_B^{-2s}\zeta (\frak m,\frak u,2s-1).\end{aligned}$$ Since 
$\Big\{1,\frac{\Delta_K+\sqrt{\Delta_k}}{2}\Big\}$ is a 
$\mathbb Z$-basis of $\mathcal O_K$, we have $|\mathcal O_K|
=\frac{1}{2}\sqrt {\Delta_K}$ and hence $|\Lambda|=
\frac{1}{2}\sqrt{\Delta_K}N(\frak u)^{-2}$. This completes the proof.
\vskip 0.20cm
This the end our discussion here for the constant terms. Our 
next aim is the explicit computation of the
 higher Fourier coefficients. This computation is more 
complicatd for the following reason. Note that the 
 choice of $B$ is quite arbitrary; for instance, $B$ may be 
multiplied from the left by any translation 
 $\left(\begin{matrix} 1&\lambda\\ 0&1\end{matrix}\right), \lambda\in K$. 
Such a change of $B$ leaves 
 $\Lambda$ unchanged and means that the higher Fourier coefficients 
are multiplied by $e\Big(\langle\omega', \lambda\rangle\Big)$. 
We shall circumvent this technical inconvenience 
later by a suitable chioce of $B$. As 
 the moment, a normalizaton of $B$ is not yet necessary.

For a fixed $0\not=c_0\in\frak m\frak u$, consider the sum with 
respect to $(c_0,d)\in\mathcal R_0$
in the third term of the right hand side of the Fourier expansion 
and define the {\it Kloosterman-like sum}
$$S(\omega',c_0):=\sum_{(c_0,d)\in\mathcal R_0}e\Big(-\langle \omega',
\frac{d}{c_0}\rangle\Big),\qquad 0\not=\omega'\in
\Lambda^\vee,\ 0\not=c_0\in\frak m\frak u.$$

\noindent
{\bf{\large Sublemma 4}.} {\it For $0\not=\omega'\in\Lambda^\vee,\
 0\not=c_0\in\frak m\frak u$, $$S(\omega',c_0)=0\qquad\mathrm{unless}\qquad 
\frac{\omega'}{\bar c_0}\in(\frak m\frak u^{-1})^\vee.$$}
Proof. Let $(c_0,d)\in\mathcal L$ and $x\in \frak m\frak u^{-1}$. 
Then $(c_0,d+x)B=(c_0,d)B+(0,x)B\in 
\frak m\oplus \frak m$ because $x\frak u\subset \frak m$. Hence if 
$(c_0,d)$ runs through a system of 
representatives for $\mathcal L(c_0)/\Lambda$, then $(c_0,d+x)$ does 
the same for every fixed $x\in 
\frak m\frak u^{-1}$. This implies 
$S(\omega',c_0)=e\Big(\langle \omega',\frac{d}{c_0}\rangle\Big)
\cdot S(\omega',c_0)$ 
for all $x\in \frak m\frak u^{-1}$. We conclude that $S(\omega',c_0)=0$ 
unless the condition 
$e(-\langle \omega',\frac{d}{c_0}\rangle)=1$ for all
$x\in \frak m\frak u^{-1}$ holds. The latter is equivalent to  
$\frac{\omega'}{\bar c_0}
\in(\frak m\frak u^{-1})^\vee.$
\vskip 0.20cm
\noindent
{\bf{\large Sublemma 5}.} 
{\it Suppose that $0\not=\omega'\in\Lambda^\vee,\ 0\not=
c_0\in\frak m\frak u,\ \frac{\omega'}{\bar c_0}
\in(\frak m\frak u^{-1})^\vee$ and $(c_0,d_0)\in \Big(\frak m\oplus\frak m\Big)
B^{-1}$. Then  
 $$S(\omega',c_0)=\frac{Nc_0}{N(\frak m)\, N(\frak u)}\cdot
e\Big(-\langle \omega',\frac{d_0}{c_0}\rangle\Big)$$ where 
 $e\Big(-\langle \omega',\frac{d}{c_0}\rangle\Big)$ is a  root of unity.}
 
\noindent
Proof. If $(c_0,d)\in\mathcal R_0$, we have $\Big((c_0,d_0)-(c_0,d)\Big)B
=\Big(0,d_0-d\Big)B\in \Big(\frak m\oplus\frak m\Big)$, 
 i.e., $d_0-d\in \frak m\frak u^{-1}$. Hence all the terms in the 
sum defining $S$ are equal, and the  number is given 
 by Sublemma 2.
\vskip 0.20cm
\noindent
{\bf{\large  Sublemma 6.}} {\it If $n\in\mathcal M$, then 
$\frak n^\vee=\frac {2}{\sqrt{\Delta_K}}\bar{\frak n}^{-1}.$ In particular, 
 $\Lambda^\vee=\frac {2}{\sqrt{\Delta_K}}\bar{\frak u}^{2}.$}

\noindent
Proof. By definition, the dual $\mathbb Z$-lattice $\frak n^\vee$ 
is the set of all $\lambda\in K$ such that
 $$\langle\lambda, x\rangle=\frac{\lambda\bar x+\bar\lambda x}{2}=
\frac{Tr(\bar\lambda x)}{2}$$ is a rational
  integer for all $x\in \frak n$. Hence $\frak n^\vee=2\bar{\frak n}^*$ 
where $\frak n^*$ is the complementary
   module with respect to the trace form. It is known that 
$\frak n^*=\frak D^{-1}\frak n^{-1}$, where
    $\frak D=\sqrt\Delta_K\mathcal O_K$ is the different of $K$. 
This completes the proof.
\vskip 0.20cm
We now define a normalization condition on $B$ that will enable us 
to compute the higher Fourier coefficients
 explicitly. Remember that $\frak u\frak v\subset\mathcal O_K$ by
Sublemma 1.4). By definition, a matrix $B\in SL(2,K)$ is called 
 quasi-integral if $\frak u_B\frak v_B=\mathcal O_K$.
Maintaining our notation  $B=\left(\begin{matrix} \alpha&\beta\\ 
\gamma&\delta\end{matrix}\right)$. We see that 
$B$ is quasi-integral if and only if $\alpha\gamma,\,\alpha\delta,
\,\beta\delta\in\mathcal O_K$.
\vskip 0.20cm
\noindent
{\bf {\large Sublemma 7}.} (1) {\it For $\gamma^*,\,
\delta^*\in K,\ (\gamma^*,\delta^*)\not=(0,0)$,
 there exists a quasi-integral matrix 
$B^*\in SL(2,K)$ such that $B^*=\left(\begin{matrix} *&*\\ 
\gamma^*&\delta^*\end{matrix}\right)$;}

\noindent
(2) {\it For every $\eta\in\mathbb P^1(K)$ there exists a quasi-integral 
matrix $B^*\in SL(2,K)$ such that $B\eta=\infty$;}

\noindent
(3) {\it For every $\frak n\in\mathcal M$, there exists a quasi-integral 
matrix $B^*\in SL(2,K)$ such that
 $\frak n=\frak u_B.$}

\noindent
Proof. A much more general version was proved before in our discussion about
fundamental domains.
\vskip 0.20cm
The following type of divisor sum will come up in the final 
formula for the higher Fourier coefficients; 
For $\frak a,\,\frak b\in\mathcal M$, $\,s\in\mathbb C$ and 
$\omega\in K^*$, let $$\sigma_s(\frak a,\frak b,
\omega)=N\frak a^{-s}\sum_{\lambda\in\frak a\frak b,\,
\omega\in\lambda \frak a^{-1}\frak b}N\lambda^{s}.$$ 
This sum is a finite one. It is empty unless $\omega\in\frak b^2$. 
For $\frak a=\mathcal O_K,\, 
\frak b\subset\mathcal O_K$ an ideal, and $\omega\in\frak b^2,\,
\omega\not=0$, the sum extends over all 
divisors $\lambda$ of $\omega$ such that $\lambda\in\frak b$ and
 $\frac{\omega}{\lambda}\in \frak b$. If
$\mu\in K^*$, then $\sigma_s(\mu\frak a,\frak b,\omega)=
\sigma_s(\frak a,\frak b,\omega).$ Moreover we
have the followng reciprocity formula
$|\omega|^{-s}\sigma_s(\frak a,\frak b,\omega)=|\omega|^s
\sigma_{-s}(\frak a^{-1},\frak b,\omega).$
\vskip 0.30cm
\noindent
{\bf {\Large Theorem 3.}} ([EGM]) {\it Suppose that $\frak m\in\mathcal M$ and 
$\eta\in\mathbb P^1(K)$ is a cusp of $SL(2,\mathcal O_K)$.
Choose a quasi-integral matrix $B=\left(\begin{matrix} \alpha&\beta\\ 
\gamma&\delta\end{matrix}\right)\in 
SL(2,K)$ such that $\eta=B^{-1}\infty$ and let $\frak u_B:=\langle\gamma,
\delta\rangle$. Then for $\Re(s)>1$,
$\ \widehat  E_\frak m(B^{-1}P,s),\ P=z+rj\in{\mathbb H}$ has the Fourier expansion
$$\begin{aligned}&\widehat  E_\frak m(B^{-1}P,s)=\\
&N(\frak u_B)^{2s}\zeta(\frak m,\frak u_B, 2s)r^{2s}
+\frac{2\pi}{\sqrt{\Delta_K}\cdot (2s-1)}
 N(\frak u_B)^{2-2s}\zeta(\frak m,\frak u_B^{-1},2s-1)r^{2-2s}\\
+&\frac{2^{1+2s}\pi^{2s}N(\frak u_B)}{\Delta_K^{s}\Gamma(2s)}
\sum_{0\not=\omega\in\frak u^2}|\omega|^{2s-1}\cdot\sigma_{1-2s}(\frak m,
\frak u_B,\omega)\cdot r\cdot 
K_{2s-1}\Big(\frac{4\pi|\omega|r^2}{\sqrt{\Delta_K}}\Big)
\cdot e\Big(\langle\frac{2\omega}
{\sqrt{\Delta_K}},z\rangle\Big).\end{aligned}$$}

\noindent
Proof. The coefficients of $r^{2s}$ and $r^{2-2s}$ are given already. 
We compute the higher Fourier 
coefficients. Let  $0\not=\omega'\in\Lambda^\vee,\ 0\not=c_0\in\frak m\frak u$.
 If  
$\frac{\omega'}{\bar c_0}\not\in(\frak m\frak u^{-1})^\vee$, then 
$S(\omega',c)=0$. Assume now
$\frac{\omega'}{\bar c_0}\in(\frak m\frak u^{-1})^\vee$, and let 
$(c,d)\in\mathcal R_0$. Then 
$$\langle\omega',\frac{d}{c}\rangle=\langle\frac{\omega'}{\bar c},
d\rangle\in\mathbb Z$$ because 
$d\in\frak m\frak v=\frak m\frak u^{-1}$ since $B$ is quasi-integral. 
This means that all terms in the 
sum defining $S$ are equal to 1 and hence $
S(\omega',c)=\frac{N(c)}{N(\frak m)\,N(\frak u_B)}.$ Also from above, 
the map from $\frak u_B^2$ to $\Lambda^\vee$
 defined by $\omega\mapsto \omega':=\frac{2}{\sqrt \Delta_K}\bar\omega
\in\Lambda^\vee$ is bijective. Replacing 
 $\omega'$ by $\frac{2}{\sqrt {\Delta_K}}\bar\omega$, we obtain
$$\begin{aligned}N\frak m^{2s}
\sum_{(c,d)\in\mathcal R_0}&\frac{e\Big(-\langle\omega',
\frac{d}{c}\rangle\Big)}{|c|^{4s}}\\
=&
N\frak m^{2s}\sum_{c\in \frak m\frak u_B,c\not=0}S(\omega',c)N(c)^{-2s}\\
=&\frac{N(\frak m)^{2s-1}}{N(\frak u_B)}
\sum_{c\in \frak m\frak u_B, \frac{\omega'}{\bar c}
\in(\frak m\frak u_B^{-1})^\vee}Nc^{1-2s}\\
=&\frac{N(\frak m)^{2s-1}}
{N(\frak u_B)}\sum_{c\in \frak m\frak u_B, 
\omega\in c\frak m^{-1}\frak u_B^{-1}}N(c)^{1-2s}\\
=&\frac{1}{N(\frak u_B)}\sigma_{1-2s}(\frak m,\frak u,\omega).\end{aligned}$$
This completes the proof.

\section{Rank Two $\mathcal O_K$-Lattices}

\subsection{Epstein Zeta Function and Eisenstein Series}

We start with a relation between Epstein zeta function and Eisenstein
series on $\mathcal H^{r_1}\times {\mathbb H}^{r_2}$.
 
Motivated by our study on non-abelian zeta functions for number fields in 
Chapter 1, for a fixed integer $r\geq 1$ and a fractional ideal $\frak a$
of a number field $K$, let us define {\it the Epstein type zeta function} 
$\widehat E_{r,\frak a;\Lambda}(s)$
associated to an $\mathcal O_K$-lattice $\Lambda$  with underlying 
projective module $P_\frak a=\mathcal O_K^{(r-1)}\oplus\frak a$ to be
$$\begin{aligned}\widehat E_{r,\frak a;\Lambda}(s):=
&\Big(\pi^{-\frac{rs}{2}}\Gamma
(\frac{rs}{2})\Big)^{r_1}\Big(2\pi^{-rs}\Gamma(rs)\Big)^{r_2}\cdot
\Big(N(\frak a)\Delta_K^{\frac{r}{2}}\Big)^s
\cdot\sum_{\bold x\in \mathcal O_K^{(r-1)}\oplus\frak a/U_{K,r}^+,
\bold x\not=(0,\cdots,0)}
\frac{1}{\|\bold x\|_\Lambda^{rs}}\end{aligned}$$
where $U_{K,r}^+:=\Big\{\varepsilon^r:
\varepsilon\in U_K,\,\varepsilon^r\in U_K^+\Big\}
=U_K^+\cap U_K^r$.
For example, note that in the case $r=2$, $U_{K,2}^+=U_K^2$, we have
$$\begin{aligned}\widehat E_{2,\frak a;\Lambda}(s):=
&\Big(\pi^{-s}\Gamma(s)\Big)^{r_1}
\Big(2\pi^{-2s}\Gamma(2s)\Big)^{r_2}\cdot \Big(N(\frak a)\Delta_K\Big)^s
\cdot\sum_{\bold x\in \mathcal O_K\oplus\frak a/U_K^2,\bold x\not=(0,0)}
\frac{1}{\|\bold x\|_\Lambda^{2s}}.\end{aligned}$$ 
From now on, we will concentrate on  this rank 2 case.
 
We want to relate the rank 2 Epstein zeta function defined in terms 
of lattices to an Eisenstein series defined over $\mathcal H^{r_1}\times 
{\mathbb H}^{r_2}$. This is based on the following simple but key 
observation, which serves as a bridge between lattices model and the upper 
half space model. (See also our discussion on stability and distance to cusps.)

Recall that, for any non-zero vector $(x,y)\in\mathcal O_K\oplus\frak a$, 
the lattice norm of $(x,y)$ associated with the lattice
$\Lambda=(\mathcal O_K\oplus\frak a,\rho_\Lambda(g))$ where 
$g=\left(\begin{matrix} a&b\\ c&d
\end{matrix}\right)$ is given by
$$\begin{aligned}&\Big\|(x,y)\Big\|_\Lambda^2=\Bigg\|(x,y)\left(\begin{matrix} a&b\\ c&d\end{matrix}\right)\Bigg\|^2\\
=&\prod_{\sigma:\mathbb R}\Big((a_\sigma x_\sigma+c_\sigma 
y_\sigma)^2+(b_\sigma x_\sigma+d_\sigma y_\sigma)^2\Big)\cdot
\prod_{\tau:\mathbb C}\Big(|a_\tau x_\tau+c_\tau y_\tau|^2+|b_\tau x_\tau+
d_\tau y_\tau|^2\Big)^2\\
=&\Big(\frac{N(g_\Lambda(\mathrm{ImJ}))}
{\Big\|x\cdot g_\Lambda (\mathrm{ImJ})+y\Big\|^2}\Big)^{-1},\end{aligned}$$ 
where, by ImJ, we mean the point 
$\mathrm{ImJ}:=(\overbrace{i,\,\ldots,\,i}^{r_1-\mathrm{times}},
\overbrace{j,\,\ldots,\,j}^{r_2-\mathrm{times}})\in
\mathcal H^{r_1}\times {\mathbb H}^{r_2}.$ (Recall that we have set
$N(\tau):=N(\mathrm{ImJ}(\tau)).$) Here for $X\in K$, we set
$\|X\|:=N(X):=\prod_{\sigma:\mathbb R}|X_\sigma|\cdot
\prod_{\tau:\mathbb C}|X_\tau|^2.$
Also change the action of units 
to the one induced from the diagonal action.
Then,
$$\begin{aligned}\widehat E_{2,\frak a;\Lambda}(s):=
&\Big(\pi^{-s}\Gamma(s)\Big)^{r_1}
\Big(2\pi^{-2s}\Gamma(2s)\Big)^{r_2}\cdot
\Big(N(\frak a)\Delta_K\Big)^s
\cdot
\sum_{(x,y)\in \mathcal O_K\oplus\frak a/U_K,(x,y)\not=(0,0)}
\Bigg(\frac{N(\mathrm{ImJ}(\tau_\Lambda))}{\Big\|x\cdot \tau_\Lambda+y\Big\|^2}\Bigg)^s.
\end{aligned}$$
Set then for $\Re(s)>1$,
 $$\begin{aligned}\widehat  E_{2,\frak a}(\tau,s):=&\Big(\pi^{-s}
\Gamma(s)\Big)^{r_1}\Big(2\pi^{-2s}\Gamma(2s)\Big)^{r_2}\cdot
\Big(N(\frak a)\Delta_K\Big)^s\cdot
\sum_{(x,y)\in \mathcal O_K\oplus\frak a/U_K,(x,y)\not=(0,0)}
\Bigg(\frac{N(\mathrm{ImJ}(\tau_\Lambda))}{\Big\|x\cdot \tau+y\Big\|^2}\Bigg)^s.
\end{aligned}$$ 
Then we have just completed the proof of the following

\noindent
{\bf{\large Lemma}.} {\it For a rank two
$\mathcal O_K$-lattice $\Lambda=(\mathcal O_K\oplus a,\rho_\Lambda)$, 
denote by $\tau_\Lambda$ the corresponding point in  the moduli space
$SL(\mathcal O_K\oplus a)\Big\backslash\Big( 
\mathcal H^{r_1}\times {\mathbb H}^{r_2}\Big)$. Then 
$$\widehat E_{2,\frak a;\Lambda}(s)
=\widehat E_{2,\frak a}(\tau_\Lambda,s).$$}

Consequently, to understand rank two non-abelian zeta functions,
we need to study the Eisenstein series 
$\widehat E_{2,\frak a}(\tau,s)$ for $\tau\in 
SL(\mathcal O_K\oplus a)\Big\backslash\Big( 
\mathcal H^{r_1}\times {\mathbb H}^{r_2}\Big)$.

\subsection{Fourier Expansion: Constant Term}

For simplicity, introduce the standard Eisenstein series by setting 
$$E_{2,\frak a}(\tau,s):=
\sum_{(x,y)\in \mathcal O_K\oplus\frak a/U_K,(x,y)\not=(0,0)}
\Bigg(\frac{N(\mathrm{ImJ}(\tau))}{\Big\|x\cdot
 \tau+y\Big\|^2}\Bigg)^s,\qquad\Re(s)>1.$$
Then the completed one becomes $$\widehat E_{2,\frak a}(\tau,s)=
\Big(\pi^{-s}\Gamma(s)\Big)^{r_1}\Big(2\pi^{-2s}\Gamma(2s)\Big)^{r_2}
\cdot \Big(N(\frak a)\Delta_K\Big)^s
\cdot E_{2,\frak a}(\tau,s).\eqno(*)$$ Following the classics, we 
in this subsection give  an explicit expression of 
Fourier expansion for the  Eisenstein series to find. (But, the final result 
shows, as the reader will find, that it is the completed Eisenstein series 
which makes the whole theory more elegent.)

As before, for the cusp  $\eta=\left[\begin{matrix}
 \alpha\\ \beta\end{matrix}\right]$, choose a (normalized) matrix 
$A=\left(\begin{matrix} \alpha&\alpha^*\\ 
\beta&\beta^*\end{matrix}\right)\in SL(2,F)$ such that 
and that if $\frak b=\mathcal O_K\alpha+\frak a\beta$, then 
$\mathcal O_K\beta^*+\frak a\alpha^*=\frak b^{-1}$. Clearly, 
$A\infty=\eta$, and moreover,
$$A^{-1}\Gamma_\eta' A=\Bigg\{\left(\begin{matrix} 
1&\omega\\ 0&1\end{matrix}\right):\omega\in \frak a\frak b^{-2}\Bigg\}.$$
Since $\widehat  E_{2,\frak a}(\tau,s)$, and hence $E_{2,\frak a}(\tau,s)$, is
$SL(\mathcal O_K\oplus \frak a)$-invariant, $E_{2,\frak a}(\tau,s)$ is
$\Gamma_\eta'\subset SL(\mathcal O_K\oplus \frak a)$-invariant. Therfeore 
$E_{2,\frak a}(A\tau,s)$ is $\frak a\frak b^{-2}$-invariant, that is, 
$E_{2,\frak a}(A\tau,s)$ is invariant under parallel transforms by elements
of $\frak a\frak b^{-2}$. As a direct consequence, we have the Fourier 
expansion
$$E_{2,\frak a}(A\tau,s)=\sum_{\omega'\in(\frak a\frak b^{-2})^\vee}
a_{\omega'}\Big(\mathrm{ImJ}(\tau),s\Big)\cdot e^{2\pi i\langle\omega',\mathrm{ReZ}(\tau)\rangle},$$ 
where $(\frak a\frak b^{-2})^\vee$ denotes the dual lattice of 
$\frak a\frak b^{-2}$. Thus, if we use $Q$ to denote a fundamental 
parallolgram of $\frak a\frak b^{-2}$ in 
$\mathbb R^{r_1}\times\mathbb C^{r_2}$, then
$$\begin{aligned}a_{\omega'}(\mathrm{ImJ}(\tau),s):=&\frac{1}{\mathrm{Vol}(\frak a\frak b^{-2})}
\sum_{(c,d)\in (\mathcal O_K\oplus \frak a)A/U_K, (c,d)\not=(0,0)}\\
&\int_Q
\Big(\frac{N(\mathrm{ImJ}(\tau))}{\|c\tau+d\|^2}\Big)^s\cdot
e^{-2\pi i\langle\omega',\mathrm{ReZ}(\tau)\rangle}\prod_{\sigma:\mathbb R}
dx_\sigma\cdot\prod_{\tau:\mathbb C}dx_\tau dy_\tau.\end{aligned}$$
(As such, we are using in fact the standard Lebesgue measure, rather than the
canonical one. So the notation may cause a bit confusion.
However, since the canonical metric and the Lebesgue one
differ by a constant factor depending only on  the field $K$, 
we, up to such a constant factor, may ignore the actual difference:
Even this makes results in this section not as explicit as possible, 
it serves our purpose of understanding  rank two non-abelian zetas 
quite well.)

Now, let us compute the Fourier coefficients in more details. For this, we 
break the summation in $a_{\omega'}$ into two cases according to whether 
$c=0$ or not.

\noindent
1) {\it Case when $c=0$.} Then the contribution becomes
$$\begin{aligned}\frac{1}{\mathrm{Vol}(\frak a\frak b^{-2})}&\sum_{(0,d)\in 
(\mathcal O_K\oplus \frak a)A/U_K, d\not=0}\int_Q
\Big(\frac{N(\mathrm{ImJ}(\tau))}
{\|d\|^2}\Big)^s\cdot
e^{-2\pi i\langle\omega',\mathrm{ReZ}(\tau)\rangle}\prod_{\sigma:\mathbb R}
dx_\sigma\cdot\prod_{\tau:\mathbb C}dx_\tau dy_\tau\\
=&
\frac{1}{\mathrm{Vol}(\frak a\frak b^{-2})}
\sum_{(0,d)\in (\mathcal O_K\oplus \frak a)A/U_K, d\not=0}
\Big(\frac{N(\mathrm{ImJ}(\tau))}{\|d\|^2}\Big)^s\int_Q
e^{-2\pi i\langle\omega',\mathrm{ReZ}(\tau)\rangle}\prod_{\sigma:\mathbb R}
dx_\sigma\cdot\prod_{\tau:\mathbb C}dx_\tau dy_\tau.\end{aligned}$$
So according to whether $\omega'=0$ or not, this case may further be 
classified into two subcases.

\noindent
1.a) {\it Subcase when $\omega'\not=0$}. Then, $$\int_Q
e^{-2\pi i\langle\omega',\mathrm{ReZ}(\tau)\rangle}\prod_{\sigma:\mathbb R}
dx_\sigma\cdot\prod_{\tau:\mathbb C}dx_\tau dy_\tau=0.$$ That is to say, 
the corresponding Fourier coefficient $a_{\omega'}=0$. So in this subcase, 
there is no contribution at all.

\noindent
1.b) {\it Subcase when $\omega'=0$.} Then $$\int_Q
\prod_{\sigma:\mathbb R}dx_\sigma\cdot\prod_{\tau:\mathbb C}dx_\tau 
dy_\tau=\mathrm{Vol}(\frak a\frak b^{-2}).$$ 
(We are illegally using canonical metric then. But then up to a fixed
constant factor argument will save us.) Hence,  in this subcase,
accordingly, $$\begin{aligned}a_{0}\Big(\mathrm{ImJ}(\tau),s\Big)=&
\sum_{(0,d)\in (\mathcal O_K\oplus \frak a)A/U_K, d\not=0}
\Big(\frac{N(\mathrm{ImJ}(\tau))}{\|d\|^2}\Big)^s\\
=&\Big(
\sum_{(0,d)\in (\mathcal O_K\oplus \frak a)A/U_K, d\not=0}N(d)^{-2s}\Big)
\cdot N(\mathrm{ImJ}(\tau))^s.\end{aligned}$$

To go further, let us look at the summation 
$$\sum_{(0,d)\in (\mathcal O_K\oplus \frak a)
A/U_K, d\not=0}N(d)^{-2s}$$ more carefully.

By definition, $$(\mathcal O_K\oplus \frak a)A=(\mathcal O_K\oplus \frak a)
\left(\begin{matrix} \alpha&\alpha^*\\ \beta&\beta^*\end{matrix}\right)$$ with 
$\left(\begin{matrix} \alpha&\alpha^*\\ \beta&\beta^*\end{matrix}\right)
\in SL(2,F)$ such that if $\frak b=\mathcal O_K\alpha+\frak a\beta$, then 
$\mathcal O_K\beta^*+\frak a\alpha^*=\frak b^{-1}$.

\noindent
{\bf \large Claim.}  $\Big(\mathcal O_K\oplus \frak a\Big)A\Big/U_K=
\Big(\mathcal O_K\alpha+
\frak a\beta,\mathcal O_K\alpha^*+\alpha\beta^*\Big)\Big/U_K
=\Big(\frak b\oplus 
\frak a\frak b^{-1}\Big)\Big/U_K.$
 
Indeed, by defintion $\mathcal O_K\alpha+\frak a\beta=\frak b$. So it 
suffices to prove that $$\mathcal O_K\alpha^*+\frak a\beta^*=\frak a
\frak b^{-1}.$$
Clearly $\alpha^*\in \frak a\frak b^{-1},\,\beta^*\in\frak b^{-1}$ 
(as already used several times), so  $\mathcal O_K\alpha^*+\alpha\beta^*
\subset \frak a\frak b^{-1}$. On the other hand, as showed before 
$\frak b^{-1}=\mathcal O_K\beta^*+\frak a^{-1}\beta^*$ so 
$$\frak a\frak b^{-1}=\frak a\cdot (\mathcal O_K\beta^*+\frak a^{-1}\beta^*)
\subset\frak a\beta^*+\mathcal O_K\alpha^*,$$  we are done.

As such, then the corresponding summation in the coeffcient $a_0$ becomes 
the one over $\Big(\frak a\frak b^{-1}\backslash \{0\}\Big)/U_K$. 
Now we use the following

\noindent
{\bf{\large Lemma}.} {\it For a fractional $\mathcal O_K$ ideal $\frak a$, 
denote by $\frak R$ the ideal class associated with $\frak a^{-1}$. Then} 

\noindent
(1) {\it there is a natural bijection}
$$\begin{matrix} \Big(\frak a\backslash\{0\}\Big)/U_K&\to&\Big\{\frak b\in 
[\frak a^{-1}]=\frak R:\frak b\ \textrm{integral}\ \mathcal O_K-
\mathrm{ideal}\Big\}\\
\overline a&\mapsto&\frak b:=a\frak a^{-1}\end{matrix}.$$

\noindent
(2) {\it For $\zeta(\frak R,s):=\sum_{\frak b\in \frak R:\frak b\ 
\textrm{integral}\ \mathcal O_K-\mathrm{ideal}}N(\frak b)^{-s},$ we have
$$\zeta(\frak R,s)=N(\frak a)^s\cdot \Big(\sum_{a\in \Big(\frak 
a\backslash\{0\}\Big)/U_K}N(a)^{-s}\Big).$$}

\noindent
Proof. All are standard. For example, (1) may be found in [Neu], 
while (2) is a direct consequence of (1).
\vskip 0.20cm
Therefore, we arrive at the following

\noindent
{\bf \Large Proposition.} {\it For the subcases at hand, the corresponding Fourier
coefficient is given by}
$$a_0\Big(\mathrm{ImJ}(\tau),s\Big)=\Big(N(\frak a^{-1}\frak b)^{2s}
\cdot\zeta([\frak a^{-1}\frak b],2s)\Big)\cdot N(\mathrm{ImJ}(\tau))^s.$$
\vskip 0.20cm
\noindent
2) {\it Case when $c\not=0$}. In this case,
$$\begin{aligned}a_{\omega'}\Big(\mathrm{ImJ}(\tau),s\Big):=&\frac{1}{\mathrm{Vol}(\frak a\frak b^{-2})}
\sum_{(c,d)\in (\mathcal O_K\oplus \frak a)A/U_K, c\not=0}\\
&\qquad\int_Q
\Big(\frac{N(\mathrm{ImJ}(\tau))}{\|c\tau+d\|^2}\Big)^s\cdot
e^{-2\pi i\langle\omega',\mathrm{ReZ}(\tau)\rangle}\prod_{\sigma:\mathbb R}
dx_\sigma\cdot\prod_{\tau:\mathbb C}dx_\tau dy_\tau.\end{aligned}$$
To compute this, consider the coset of $A^{-1}\Gamma_\eta'A$ among 
$\left(\begin{matrix} *&*\\ c&d\end{matrix}\right)$. 

\noindent
{\bf \large Claim.} {\it For $(c,d)\in (\mathcal O_K\oplus \frak a)A/U_K,\ c\not=0$ and $\omega\in \frak a\frak b^{-2}$, we have $(c,c\omega+d)\in 
(\mathcal O_K\oplus \frak a)A/U_K,\ c\not=0$.}

It suffices to deal the component of $d$ and $c\omega+d$. Note that 
$A=\left(\begin{matrix} \alpha&\alpha^*\\ \beta&\beta^*\end{matrix}\right)
\in SL(2,F)$ with
$\alpha\in\frak b,\,\beta\in\frak a^{-1}\frak b,\,
\alpha^*\in\frak a\frak b^{-1},\,\beta^*\in\frak b^{-1}$, 
we have $$c\in \mathcal O_K\cdot \frak b+\frak a
\cdot \frak a^{-1}\frak b=\frak b\qquad\mathrm{and} 
\qquad d\in \mathcal O_K\cdot \frak a\frak b^{-1}+\frak a\cdot \frak b^{-1}=
\frak a\, \frak b^{-1}.$$ So, we should show that
with $c\in\frak b$, $d\in \frak a\, \frak b^{-1}$ and $\omega\in \frak a
\frak b^{-2}$, we have $c\omega+d\in \frak a\, \frak b^{-1}$. 
But this is clear since $c\omega+d\in \frak b\cdot \frak a\frak b^{-2}+\frak a\, \frak b^{-1}=\frak a\, \frak b^{-1}.$ This completes the proof of the Claim.
 
Now since
$$\left(\begin{matrix} *&*\\ c&d\end{matrix}\right)\left(\begin{matrix} 
1&\omega\\ 0&1\end{matrix}\right)=\left(\begin{matrix} *&*\\ c&c\omega+
d\end{matrix}\right)$$
 with $(c,d)\in (\mathcal O_K\oplus \frak a)A/U_K,\ c\not=0$ and $\omega\in 
\frak a\frak b^{-2}$. Consequently, if we let $\mathcal R$ to be a system of 
representatives of $\left(\begin{matrix} *&*\\ c&d\end{matrix}\right)$ 
modulo the right action of $A^{-1}\Gamma_\eta' A$, or better, to be a 
system of representatives of $(c,d)$ modulo the relation 
$$(c,d)\sim (c,c\omega+d)$$ with 
$(c,d)\in (\mathcal O_K\oplus \frak a)A/U_K,\, c\not=0$ and $\omega\in 
\frak a\frak b^{-2}$, then in the case at hand, note that for $\tau\in \mathcal H^{r_1}\times\mathbb H^{r_2}$, 
$\mathrm{ImJ}(\tau+\omega)=\mathrm{ImJ}(\tau)$, we have
$$\begin{aligned}&a_{\omega'}\Big(\mathrm{ImJ}(\tau),s\Big)\\
=&\frac{1}{\mathrm{Vol}
(\frak a\frak b^{-2})}
\sum_{
\left(\begin{matrix} *&*\\ c&d\end{matrix}\right)
\in\mathcal R}
\int_{\omega\in \frak a\frak b^{-2}}\int_Q
\Bigg(\frac{N(\mathrm{ImJ}(\tau))}
{\|c(\tau+\omega)+d\|^2}\Bigg)^s\cdot
e^{-2\pi i\langle\omega',\mathrm{ReZ}(\tau)\rangle}\prod_{\sigma:\mathbb R}
dx_\sigma\cdot\prod_{\tau:\mathbb C}dx_\tau dy_\tau\\
=&\frac{1}{\mathrm{Vol}(\frak a\frak b^{-2})}\sum_{
\left(\begin{matrix} *&*\\ c&d\end{matrix}\right)\in
\mathcal R}\int_{\mathbb R^{r_1}\times\mathbb C^{r_2}}
\Bigg(\frac{N(\mathrm{ImJ}(\tau))}{\|c\tau+d\|^2}\Bigg)^s
\cdot
e^{-2\pi i\langle\omega',\mathrm{ReZ}(\tau)\rangle}\prod_{\sigma:\mathbb R}
dx_\sigma\cdot\prod_{\tau:\mathbb C}dx_\tau dy_\tau\\
=&\frac{1}{\mathrm{Vol}(\frak a\frak b^{-2})}\sum_{
\left(\begin{matrix} *&*\\ c&d\end{matrix}\right)\in
\mathcal R}\frac{1}{N(c)^{2s}}\int_{\mathbb R^{r_1}\times\mathbb C^{r_2}}
\Bigg(\frac{N(\mathrm{ImJ}(\tau))}
{\|\tau+\frac{d}{c}\|^2}\Bigg)^s\cdot
e^{-2\pi i\langle\omega',\mathrm{ReZ}(\tau)+\frac{d}{c}-\frac{d}{c}\rangle}
\prod_{\sigma:\mathbb R}dx_\sigma\cdot\prod_{\tau:\mathbb C}dx_\tau dy_\tau
\\
=&\frac{1}{\mathrm{Vol}(\frak a\frak b^{-2})}\sum_{
\left(\begin{matrix} *&*\\ c&d\end{matrix}\right)\in
\mathcal R}\frac{e^{2\pi i\langle \omega',\frac{d}{c}\rangle}}{N(c)^{2s}}
\int_{\mathbb R^{r_1}\times\mathbb C^{r_2}}
\Bigg(\frac{N(\mathrm{ImJ}(\tau))}
{\|\tau\|^2}\Bigg)^s\cdot
e^{-2\pi i\langle\omega',\mathrm{ReZ}(\tau)\rangle}\prod_{\sigma:\mathbb R}
dx_\sigma\cdot\prod_{\tau:\mathbb C}dx_\tau dy_\tau
\end{aligned}$$

\noindent
2.a) {\it Subcase when $\omega'=0$}.
Then $$\begin{aligned}a_0\Big(\mathrm{ImJ}(\tau),s\Big)=&\frac{1}
{\mathrm{Vol}(\frak a\frak b^{-2})}
\sum_{
\left(\begin{matrix} *&*\\ c&d\end{matrix}\right)\in\mathcal R}
\frac{1}{N(c)^{2s}}\cdot\int_{\mathbb R^{r_1}\times\mathbb C^{r_2}}
\Bigg(\frac{N(\mathrm{ImJ}(\tau))}{\|\tau\|^2
}\Bigg)^s\prod_{\sigma:\mathbb R}dx_\sigma\cdot\prod_{\tau:\mathbb C}dx_\tau 
dy_\tau.\end{aligned}$$
So according to whether $\sigma:\mathbb R$ or $\tau:\mathbb C$, we have to 
compute the following integrations:

\noindent
2.a.i) {\it For reals},
$$\begin{aligned}\int_{\mathbb R}\Big(\frac{y}{x^2+y^2}\Big)^sdx
=&\frac{1}{y^s}\int_{\mathbb R}\Big(\frac{1}{(\frac{x}{y})^2+1}\Big)^sd
\frac{x}{y}\cdot y\\
=&y^{1-s}\int_{\mathbb R}\frac{dt}{(1+t^2)^s}\\
=&y^{1-s}\cdot\pi^{\frac{1}{2}}
\frac{\Gamma(s-\frac{1}{2})}{\Gamma(s)};\end{aligned}$$

\noindent
2.a.ii) {\it For complexes},
$$\begin{aligned}\int_{\mathbb C}&\Big(\frac{r}{|z|^2+r^2}\Big)^{2s}dx\,dy\\
=&\frac{1}{r^{2s}}\int_{\mathbb C}\Big(\frac{1}{(\frac{|z|}{r})^2+1}
\Big)^{2s}d\frac{x}{r}\cdot r d\frac{y}{r}\cdot r\\
=&r^{2-2s}\int_{\mathbb C}\frac{dx\,dy}{(1+|z|^2)^{2s}}\\
=&r^{2-2s}\cdot
\frac{\pi}{2s-1};\end{aligned}$$

\noindent
2.b) {\it Subcase when $\omega'\not=0$.}
Then $$\begin{aligned}a_{\omega'}\Big(\mathrm{ImJ}(\tau),s\Big)
=&\frac{1}{\mathrm{Vol}(\frak a\frak b^{-2})}\sum_{
\left(\begin{matrix} *&*\\ c&d\end{matrix}\right)\in
\mathcal R}\frac{e^{-2\pi i\langle \omega',\frac{d}{c}\rangle}}{N(c)^{2s}}\\
&\times \int_{\mathbb R^{r_1}\times\mathbb C^{r_2}}
\Bigg(\frac{N(\mathrm{ImJ}(\tau))}{\|\tau\|^2}\Bigg)^s\cdot
e^{-2\pi i\langle\omega', \mathrm{ReZ}(\tau)\rangle}
\prod_{\sigma:\mathbb R}dx_\sigma\cdot
\prod_{\tau:\mathbb C}dx_\tau dy_\tau.\end{aligned}$$
So according to whether $\sigma:\mathbb R$ or $\tau:\mathbb C$, we have to 
compute the following integrations:

\noindent
2.b.i) {\it For reals},
$$\begin{aligned}&\int_{\mathbb R}\Big(\frac{y}{x^2+y^2}\Big)^se^{-2\pi i|\omega'|\cdot x}dx\\
=&\frac{1}{y^s}\int_{\mathbb R}\Big(\frac{1}{(\frac{x}{y})^2+1}\Big)^s
e^{-2\pi i|\omega'|\cdot \frac{x}{y}y}
d\frac{x}{y}\cdot y\\
=&y^{1-s}\int_{\mathbb R}\frac{1}{(1+t^2)^s}e^{-2\pi i|\omega'|yt}dt\\
=&y^{1-s}\cdot \Big(2\pi^s|\omega'|^{s-\frac{1}{2}}\cdot y^{s-\frac{1}{2}}
\cdot\frac{1}{\Gamma(s)}\cdot K_{s-\frac{1}{2}}(2\pi|\omega'|y)\Big)\\
=&2\pi^s|\omega'|^{s-\frac{1}{2}}\cdot y^{\frac{1}{2}}\cdot\frac{1}{\Gamma(s)}
\cdot K_{s-\frac{1}{2}}(2\pi|\omega'|y);\end{aligned}$$

\noindent
2.b.ii) {\it For complexes},
$$\begin{aligned}&\int_{\mathbb C}
\Big(\frac{r}{|z|^2+r^2}\Big)^{2s}e^{-2\pi i|\omega'|\cdot x}dx\,dy
=\int_{\mathbb R^2}\Big(\frac{r}{x^2+y^2+r^2}\Big)^{2s}e^{-2\pi i|\omega'|
\cdot x}dx\,dy\\
=&\int_{\mathbb R^2}\Big(\frac{1}{(\frac{x}{r})^2+(\frac{y}{r})^2+1}
\Big)^{2s}r^{-2s}e^{-2\pi i|\omega'|\cdot \frac{x}{r}\cdot r}d
\frac{x}{r}\,d\frac{y}{r}\cdot r^2\\
=&r^{2-2s}\int_{\mathbb R^2}\frac{e^{-2\pi i|\omega'|\cdot x\cdot r}}
{(x^2+y^2+1)^{2s}}dx\,dy\\
=&r^{2-2s}\int_{\mathbb R}\int_{\mathbb R}\frac{dy}{(y^2+x^2+1)^{2s}}\cdot 
e^{-2\pi i|\omega'| r x}dy\,dx\\
=&r^{2-2s}\int_{\mathbb R}\Big(\int_{\mathbb R}\frac{1}{(\Big(\frac{y}
{\sqrt{x^2+1}}\Big)^2+1)^{2s}}d\frac{y}{\sqrt{x^2+1}}\cdot \sqrt{x^2+1}
\cdot\frac{1}{(x^2+1)^{2s}}\Big)\cdot e^{-2\pi i|\omega'| r x}dx\\
=&r^{2-2s}\int_{\mathbb R}\Big(\int_{\mathbb R}\frac{dt}{(1+t^2)^{2s}}\Big)
\cdot \frac{e^{-2\pi i|\omega'| r x}}{(x^2+1)^{2s-\frac{1}{2}}}dx\\
=&r^{2-2s}\cdot\Big(\pi^{\frac{1}{2}}\cdot\frac{\Gamma(2s-\frac{1}{2})}
{\Gamma(2s)}\Big)\cdot
\int_{\mathbb R}\frac{e^{-2\pi i|\omega'| r x}}{(x^2+1)^{2s-\frac{1}{2}}}dx\\
=&r^{2-2s}\cdot\Big(\pi^{\frac{1}{2}}\cdot\frac{\Gamma(2s-\frac{1}{2})}
{\Gamma(2s)}\Big)\cdot\Big(2\pi^{2s-\frac{1}{2}}|\omega'|^{2s-1}r^{(2s-1)}
\frac{1}{\Gamma(2s-\frac{1}{2})}K_{2s-1}(2\pi|\omega'|r)\Big)\\
=&\frac{2\pi^{2s}|\omega'|^{2s-1}}{\Gamma(2s)}rK_{2s-1}(2\pi|\omega'|r)
\end{aligned}$$ 
by using the calculation for reals. Or more directly, 
$$\begin{aligned}&\int_{\mathbb C}\Big(\frac{r}{|z|^2+r^2}\Big)^{2s}
e^{-2\pi i|\omega'|\cdot x}dx\,dy\\
=&r^{-2s}\int_{\mathbb C}\Big(\frac{1}{(\frac{|z|}{r})^2+1}
\Big)^{2s}de^{-2\pi i|\omega'|\cdot \frac {x}{r}\cdot r}
\frac{x}{r}\cdot r d\frac{y}{r}\cdot r\\
=&r^{2-2s}\int_{\mathbb C}\frac{1}{(1+|z|^2)^{2s}}e^{-2\pi i|\omega'|\cdot x
\cdot r}\cdot dx\,dy\\
=&r^{2-2s}\cdot\frac{2\pi^{2s}|\omega'|^{2s-1}}
{\Gamma(2s)}r^{2s-1}K_{2s-1}(2\pi|\omega'|r)\\
=&\frac{2\pi^{2s}|\omega'|^{2s-1}}{\Gamma(2s)}r K_{2s-1}(2\pi|\omega'|r).
\end{aligned}$$

All in all, we have then obtain the following
\vskip 0.30cm
\noindent
{\bf {\Large Theorem.}} {\it With the same notation as above, 
we have the following Fourier expansion for the 
 Eisenstein series}
$$\begin{aligned}&E_{2,\frak a}(A\tau,s)=
\zeta([\frak a^{-1}\frak b],2s)\cdot N(\frak a\frak b^{-1})^{-2s}
\cdot N(\mathrm{ImJ}(\tau))^s\\
&+\frac{1}{\mathrm{Vol}(\frak a\frak b^{-2})}\sum_{
\left(\begin{matrix} *&*\\ c&c\omega+d\end{matrix}\right)\in
\mathcal R}
\frac{1}{N(c)^{2s}}
\cdot(\pi^{\frac{1}{2}})^{r_1}\cdot\Big(\frac
{\Gamma(s-\frac{1}{2})}{\Gamma(s)}\Big)^{r_1}
\cdot \Big(\frac{\pi}{2s-1}
\Big)^{r_2}\cdot
N(\mathrm{ImJ}(\tau))^{1-s}\\
&+\frac{1}{\mathrm{Vol}(\frak a\frak b^{-2})}
\sum_{
\left(\begin{matrix} *&*\\ c&d\end{matrix}\right)\in
\mathcal R}
\frac{e^{2\pi i\langle\omega',\frac{d}{c}\rangle}}{N(c)^{2s}}
\cdot N(\mathrm{ImJ}(\tau))^{\frac{1}{2}}\cdot N(\omega')^{s-\frac{1}{2}}\\
&\times \Big(\frac{2\pi^s}{\Gamma(s)}\Big)^{r_1}
\prod_{\sigma:\mathbb R}K_{s-\frac{1}{2}}(2\pi |\omega'|_\sigma y_\sigma)
\cdot\Big(\frac{2\pi^{2s}|\omega'|^{2s-1}}{\Gamma(2s)}\Big)^{r_2}\cdot
\prod_{\tau:\mathbb C} K_{2s-1}(2\pi |\omega'|_\tau r_\tau).
\end{aligned}$$

\noindent 
{\bf Warning:} This is not the final version of the Fourier expansion
we expect:
We are supposed to have a formula as in Theorem 3 of 3.2. Furthermore, 
if one wants to have an analogue of Kronecker 
limit formula here,  Kloosterman type sums appeared above have to be studied.
However, this does not matter for our limited purpose here: 
For rank two zetas, only the first 
coefficience, i.e., the coefficience of $N(\mathrm{ImJ}(\tau))^s$, 
palys a key role.
 (I intend to come back to this point later so as to also 
give precise expressions for all the coefficiants in terms of $K$ and 
the associated fractional ideals  in the style of [EGM].)

\chapter{Explicit Formula for Rank Two Zeta Functions: Rankin-Selberg \& Zagier Method}

The original Rankin-Selberg method gives a way to express the 
Mellin transfrom of the constant term in the
Fourier expansion of an automorphic function as the scalar product 
of the automorphic function with an
Eisenstein series when the automorphic function is very small when 
approach to the cusp. In a paper of 
Zagier, this method is extended to a much board type of automorphic 
functions, certain kind of slow growth 
functions. 

In this chapter, we  use a generalization of
the Rankin-Selberg \& Zagier method to give an explicit expression
for rank two zeta functions of number fields in terms of 
Dedekind zeta functions.

\section{Upper Half Plane}

We here will give two different approaches, both are due to Zagier. 
The first uses a particular truncation of the fundamenmtal 
domain,  while the second uses the original Rankin-Selberg method even when
slow growth functions are involved.

\subsection{Geometric Approach}

The original paper of Zagier deals with $SL(2,\mathbb Z)$, but the method 
clearly works for general Fuchsian groups as well.
Here, we use Gupta's beautiful exposition in her Journal of 
Number Theory paper, which follows closely Zagier's paper. 

Let $\Gamma$ be a congruence subgroup of $SL(2,\mathbb Z)$ which is
reduced at infinity so that
 $\Gamma_\infty$ is generated by $\left(\begin{matrix} 1&1\\
 0&1\end{matrix}\right)$. 
 $\Gamma$ acts naturally on $\mathcal H$. Let  $\kappa_1=\infty,
\,\kappa_2,\,\ldots,\,\kappa_h$ be the set of inequivalent cusps.
For $i=1,\,\ldots,\,h$, denote the isotropy groups by 
$\Gamma_i:=\Gamma_{\kappa_i}$, and choose 
$A_i\in GL(2,\mathbb Q)$ such that $A_i\infty=\kappa_i$. Then
$A_i^{-1}\Gamma_i A_i=\Gamma_\infty$.
 
Recall that the Eisenstein series $E_i(z,s)$ of $\Gamma$ at 
$\kappa_i$ is defined by 
$$E_i(z,s):=\sum_{\gamma\in\Gamma_i\backslash \Gamma} y^s(A_i^{-1}\gamma z),$$ 
where $y(x+iy):=y$. Set 
$$\bold E(z,s):=\left[\begin{matrix}
E_1(z,s)\\ E_2(z,s)\\ \cdots\\ E_h(z,s)\end{matrix}\right].$$ Then,  
there is a matrix $\Phi(s)$ of 
functions such that $$\bold E(z,s)=\Phi(s) \bold E(z,s), \qquad
\Phi(s)\Phi(1-s)=I_{h\times h}$$ where 
$$\Phi(s)=(\phi_{ij})_{h\times h}\qquad \mathrm{and}\qquad 
\phi_{ij,0}(s)=\sum_{c>0}
\frac{1}{|c|^{2s}}\sum_{d\equiv \pmod{c}, 
\left(\begin{matrix} *&*\\ c&d\end{matrix}\right)\in A_i^{-1}\Gamma A_j}1.$$
 Note that $E_1(z,s)$ is the Eisenstein series at $\infty$. We will 
write it as $E(z,s)$ for short.

Zagier's idea is as follows. For a continuous function $F(z)$ invariant 
under the action of $\Gamma$
with the Fourier expansion $F(A_iz)=\sum_{m\in\mathbb Z}a_m^i(y)e(mx)$ 
such that it is of slow growth, then we want to know
the Mellin tranforms of the constant term $a_0^i(y)$ along the line of 
Rankin-Selberg. However, usually, slow 
growth condition is not enough. (Recall that $F$ being slow growth means that
$$F(A_iz)=\psi_i(y)+O(y^{-N}), \qquad\mathrm{as}\ 
y=\Im (z)\to\infty, \qquad \forall N.$$
What Zagier added, in addition to the slow 
growth condition, is
that $$\psi_i(y)=\sum_{j=1}^l
 \frac{c_{ij}}{n_{ij}!}y^{\alpha_{ij}}\log^{n_{ij}}y,\qquad\mathrm{where}\  
c_{ij},\ \alpha_{ij}\in\mathbb C,\ n_{ij}\in  \mathbb Z_{\geq 0}.$$
As such, we may then modify the Mellin tranform for the constant term 
to the one for $a_0^i(y)-\psi(y)$ but 
shift by $-1$. More precisely, we define the {\it Zagier tranform} $R_i(F,s)$
of $F$ at the cusp $\kappa_i$ by
$$R_i(F,s):=\int_0^\infty\Big(a_0^i(y)-\psi_i(y)\Big)y^{s-2}dy
=\int_0^\infty\int_0^1\Big[F(A_iz)-\psi_i(y)\Big]y^s\cdot 
\frac{dx\wedge dy}{y^2}.$$ 
Note that since the slowth growth part 
has been truncated, the integration makes perfect sense for say $\Re(s)$ 
is sufficiently large. Now the 
original Rankin-Selberg method can be applied. In fact, much more is true.
To state it,  set $$\bold R(F,s):=\left[\begin{matrix}
R_1(F,s)\\ R_2(F,s)\\ \cdots\\ R_h(F,s)\end{matrix}\right]\qquad\mathrm{and} 
\qquad \bold h(s):=\left(\begin{matrix}-\sum_{j=1}^l
\frac{c_{1j}}{(1-\alpha_{1j}-s)^{n_{1j}}+1}\\
-\sum_{j=1}^l
\frac{c_{2j}}{(1-\alpha_{2j}-s)^{n_{2j}}+1}\\
\cdots\\
-\sum_{j=1}^l
\frac{c_{hj}}{(1-\alpha_{hj}-s)^{n_{hj}}+1}\end{matrix}\right).$$ 
\vskip 0.30cm
\noindent
{\bf{\Large Proposition}.} ([Z] and [Gu]) (1) {\it $\bold R(F,s)$ is well-defined for say $\Re(s)$ 
is sufficiently large;}

\noindent
(2) (Functional equation) $\bold R(F,s)=\Phi(s)\bold R(F,1-s);$

\noindent
(3) $\xi(2s)\bold R(F,s)=\xi(2s)\bold h(s)+\xi(2s)\Phi(s)\bold h(1-s)+
\frac{\mathrm{entire\ function\ of}\ 
s}{s(s-1)}.$ 

\noindent
Proof. Indeed, for $\mathcal D$ the standard fundamental domain 
for the action of
$SL(2,\mathbb Z)$ on $\mathcal H$, 
let $\mathcal D_\Gamma$ be a fundamental domain of $\Gamma$ with 
$|x|\leq \frac{1}{2}$. Let $$\begin{aligned}S_\infty(T):=&
\Big\{z\in\mathcal H:\Im(z)>T,|x|\leq\frac{1}{2}\Big\},\\  
S_{\kappa_i}(T):=&\Big\{z\in\mathcal H:A_i^{-1}z\in S_\infty(T)\Big\}
=A_iS_\infty(T).\end{aligned}$$ Consider then the truncated domain 
$\mathcal D_T=\mathcal D_\Gamma\Big\backslash
\Big(\cup_i S_{\kappa_i}(T)\Big)$, for sufficiently large $T$. Then 
$\mathcal D_T$ is the fundamental domain for the action of 
$\Gamma$ on 
$$\begin{aligned}\mathcal H_T:=&\cup_{\gamma\in\Gamma}\gamma\mathcal D_T\\
=&
\Big\{z\in\mathcal H:\max_{\delta\in\Gamma, i\geq 1}
\Im (A_i^{-1}\delta z)\leq T\Big\}\\
=&\Big\{z\in\mathcal H:\Im (z)\leq T\Big\}\Big\backslash 
\Big(\cup_{c\geq 1}\cup_{a\in\mathbb Z,
(a,c)=1}S_{a/c}\Big),\end{aligned}$$ where 
$S_{a/c}:=\delta^{-1}A_i\,\Big\{\Im(z)>T\Big\}$, 
$\,\delta\in\Gamma$ and $\delta(a/c)=\kappa_i$ 
for some $i\geq 1$. (Easily one sees that $S_{a/c}$ is in fact an open disk 
in the upper half plane tangent to
 $x$-axis at the point $(a/c,0)$.)

Thus $$\Gamma_\infty\Big\backslash \mathcal H_T=
\Big\{x+iy:0\leq x\leq 1, 0\leq y\leq T\Big\}\Big\backslash 
\Big(\cup_{c\geq 1}\cup_{a\,(\mathrm{mod} c),\,(a,c)=1}S_{a/c}\Big).$$ 
Note that $F$ is $\Gamma$-invariant, so is $F(z)\cdot\chi_T$ where
 $\chi_T$ denotes the characteristic 
function of $\mathcal H_T$. As such, using the unfolding trick, we 
have arrived at
$$\begin{aligned}
\int_{\mathcal D_T}&F(z)E(z,s)d\mu=\int_{\Gamma_\infty\backslash\Gamma} 
F(z) y^s d\mu\\
=&\int_0^T\int_0^1 F(z) y^sd\mu-\sum_{c=1}^\infty
\sum_{a(\mathrm{mod} c),(a,c)=1}\iint_{S_{a/c}}F(z) y^sd\mu.\end{aligned}
\eqno(*)$$ To compute the summation on the right hand side of 1), 
we divide it into two cases.

\noindent
(1) Case $a/c\sim \infty$.

Let $\gamma_0=\left(\begin{matrix} a&b\\ c&d\end{matrix}\right)\in\Gamma$. 
Then $\gamma_0^{-1}S_{a/c}=\Big\{z\in
\mathcal H:\Im(z)>T\Big\}.$ We have $$\begin{aligned}\iint_{S_{a/c}}F(z) y^sd\mu=&\int_T^\infty
\int_{-\infty}^\infty F(z)
\Im(\gamma_0z)^s d\mu\\
=&\int_T^\infty\int_{-1/2}^{1/2} F(z)\sum_{n=-\infty}^\infty 
\Im(\gamma_0(z+n))^s d\mu\\
=&\iint_{S_\infty(T)}F(z)\sum_{\gamma=\left(\begin{matrix} a&*\\ 
c&*\end{matrix}\right)\in\Gamma}
\Im(\gamma z)^s d\mu,\end{aligned}$$ where the sum is over all 
$\gamma\in\Gamma$ 
with the form $\gamma=\left(\begin{matrix}
 a&*\\ c&*\end{matrix}\right)$, all of which are in the form 
$\gamma_0\left(\begin{matrix} 1&n\\ 0&1\end{matrix}
 \right)$ for some $n\in\mathbb Z$.
 But 
$$\Bigg\{
\Gamma_\infty\backslash\Gamma/\{\pm 1\}
\Bigg\}=
\Bigg(
\cup_{c>0}
\cup_{a\,\mathrm{mod} z,\, a/c\sim \infty}
 \cup_{\left(\begin{matrix} a&*\\ c&*\end{matrix}\right)\in\Gamma}
\left(\begin{matrix} a&*\\ c&*\end{matrix}\right)
\Bigg)
\cup
 \Bigg\{I_2\Bigg\}.$$ Thus,
 $$\sum_{c=1}^\infty\sum_{a(\mathrm{mod}c), (a,c)=1, a/c\sim \infty}
\iint_{S_{a/c}}F(z) y^sd\mu=
 \iint_{S_\infty(T)}F(z)\cdot \Big(E(z,s)-y^s\Big)d\mu.$$

\noindent
(2) Case $a/c\sim \kappa_{i\geq 2}\not\sim\infty$. 

So $\delta(a/c)=\kappa_i$. Then 
$$\begin{aligned}\sum_{c=1}^\infty&\sum_{a(\mathrm{mod}c), (a,c)=1, a/c\sim \kappa_i}
\iint_{S_{a/c}}F(z) y^sd\mu\\
=&
 \iint_{S_{\kappa_i}(T)}F(z)E(z,s)d\mu\\
 =&\iint_{A_iS_\infty(T)}F(z)E(z,s)d\mu\\
=&\iint_{S_\infty(T)}F(A_iz)E(A_iz,s)d\mu.
 \end{aligned}$$
 Therefore, $$\int_{\mathcal D_T}F(z)E(z,s)d\mu=\int_0^Ta_0^\infty(y)y^{s-2}dy-
\sum_{i=1}^h\iint_{S_\infty(T)} F(A_iz)\Big[E(A_iz,s)-\delta_{i\infty}y^s\Big]d\mu,$$
since $\int_0^T\int_0^1 F(z) y^sd\mu=\int_0^Ta_0^\infty(y)y^{s-2}dy$ 
as to be easily checked.
\vskip 0.20cm
Now let us use the properties of Eisenstein series to simplify the 
right hand side.

Let $e_{j\infty}:=e_{j\infty}(y,s):=\int_0^1E(A_jz,s)d\mu$ be the 
constant term in the Fourier expsnsion of 
$E(z,s)$ at $\kappa_j$. It is well-known that $e_{i\infty}=
\delta_{i\infty}y^s+\phi_{i\infty}y^{1-s}.$ Thus 
for $\Re(s)$ sufficiently large,
$$\begin{aligned}\int_{\mathcal D_T}F(z)&E(z,s)d\mu
=\int_0^Ta_0^\infty(y)y^{s-2}dy\\
&-
\sum_{i=1}^h\iint_{S_\infty(T)}F(A_iz)\Big[E(A_iz,s)-e_{i\infty}\Big]\,d\mu
-\sum_{i=1}^h\iint_{S_\infty(T)}F(A_iz)\phi_{i\infty}y^{1-s}d\mu\\
&\qquad\qquad=\int_0^Ta_0^\infty(y)y^{s-2}dy\\
&-
\sum_{i=1}^h\iint_{S_\infty(T)}F(A_iz)\Big[E(A_iz,s)-e_{i\infty}\Big]\,d\mu
-\sum_{i=1}^h\phi_{i\infty}\int_T^\infty a_0^i(y)y^{-1-s}dy.\end{aligned}$$
The difference $\Big[E(A_iz,s)-e_{i\infty}\Big]$ is an entire function of $s$ 
and is of rapid decay with respect to $y$.
Thus we have obtained the following

\noindent
{\bf \large Equation I.} {\it With the same notation as above,}
$$\begin{aligned}\int_{\mathcal D_T}F(z)E(z,s)d\mu&+
\sum_{i=1}^h\iint_{S_\infty(T)}F(A_iz)\Big[E(A_iz,s)-e_{i\infty}\Big]d\mu\\
=&\int_0^Ta_0^\infty(y)y^{s-2}dy
-\sum_{i=1}^h\phi_{i\infty}\int_T^\infty a_0^i(y)y^{-1-s}dy.\end{aligned}$$

We further evaluate the right hand side. Write 
$$\int_0^Ta_0^i(y)y^{s-2}dy=\int_0^T\Big(a_0^i(y)-\psi_i(y)\Big)
y^{s-2}dy+\int_0^T\psi_i(y)y^{s-2}dy,$$ and set 
$$h_T^i(s):=\int_0^T\psi_i(y)y^{s-2}dy=\sum_{j=1}^l
\frac{c_{ij}}{n_{ij}!}\frac{\partial^{n_{ij}}}
{\partial s^{n_{ij}}}\Big(\frac{T^{s+\alpha_{ij}-1}}{s+
\alpha_{ij}-1}\Big)$$ and $\bold h_T(s)=(h_T^i(s)).$
Then $$\int_0^Ta_0^i(y)y^{s-2}dy=R_i(F,s)-\int_T^\infty
\Big(a_0^i(y)-\psi_i(y)\Big)y^{s-2}dy+h_T^i(s).$$

On the other hand, write $$\int_T^\infty a_0^i(y)y^{-s-1}dy=
\int_T^\infty\Big(a_0^i(y)-\psi_i(y)\Big)y^{-s-1}dy+
\int_T^\infty\psi_i(y)y^{-s-1}dy.$$ Since 
$\int_0^\infty \psi(y)y^{-s-1}dy=0$,  $$\int_T^\infty\psi_i(y)
y^{-s-1}dy=-\int_0^T\psi_i(y)y^{-s-1}dy=-h_T^i(1-s).$$ Thus, 
$$\int_T^\infty a_0^i(y)y^{-s-1}dy=\int_T^\infty
(a_0^i(y)-\psi_i(y))y^{-s-1}dy-h_T^i(1-s).$$
Consequently, we get the following

\noindent
{\bf \large Equation II.} {\it With the same notation as above,}
$$\begin{aligned}\int_{\mathcal D_T}F(z)E(z,s)d\mu&+
\sum_{i=1}^h\iint_{S_\infty(T)}F(A_iz)\Big[E(A_iz,s)-e_{i\infty}\Big]d\mu\\
=&R_\infty(F,s)-\int_T^\infty \Big(a_0^\infty(y)-\psi_\infty(y)\Big)
y^{s-2}dy\\
&+h_T^\infty(s)
-\sum_{i=1}^h\phi_{i\infty}\Big[\int_T^\infty (a_0^i(y)-
\psi_i(y))y^{-1-s}dy-h_T^i(1-s)\Big].\end{aligned}$$

Since $$\iint_{S_\infty(T)}\sum_{i=1}^h 
F(A_iz)e_{i\infty}d\mu=\int_T^\infty a_0^\infty(y)y^{s-2}dy+\sum_i
\phi_{i\infty}\int_T^\infty a_0^i(y)y^{-s-1}dy,$$ we have 
$$\begin{aligned}~&R_\infty(F,s)
+h_T^\infty(s)+\sum_{i=1}^h\phi_{i\infty}h_T^i(1-s)\\
=&
\int_{\mathcal D_T}F(z)E(z,s)d\mu+
\sum_{i=1}^h\iint_{S_\infty(T)}F(A_iz)\Big[E(A_iz,s)-e_{i\infty}\Big]d\mu\\
&+
\int_T^\infty \Big(a_0^\infty(y)-\psi_\infty(y)\Big)
y^{s-2}dy+\sum_{i=1}^h\phi_{i\infty}\int_T^\infty \Big(a_0^i(y)-
\psi_i(y)\Big)y^{-1-s}dy\\
=&\int_{\mathcal D_T}F(z)E(z,s)d\mu+\sum_{i=1}^h\iint_{S_\infty(T)}
\Big(F(A_iz)E(A_iz,s)-\psi_i(y)e_{i\infty}\Big)d\mu.\end{aligned}$$
Thus by similarly working over other cusps, we have obtained the following

\noindent
{\bf \large Equation III.} {\it With the obvious change of notations, we have}
$$\begin{aligned}\int_{\mathcal D_T}&F(z)E_\kappa(z,s)d\mu+
\sum_{i=1}^h\iint_{S_\infty(T)}
\Big(F(A_iz)E_\kappa(A_iz,s)-\psi_i(y)e_{i\kappa}\Big)\,d\mu\\
=&R_\kappa(F,s)
+h_T^\kappa(s)+\sum_{i=1}^h\phi_{i\kappa}h_T^i(1-s).\end{aligned}$$ {\it Or in vector form,}
$$\begin{aligned}\int_{\mathcal D_T}&F(z)\bold E(z,s)d\mu+\sum_{i=1}^h\iint_{S_\infty(T)}
\Big(F(A_iz)\bold E(A_iz,s)-\psi_i(y)\bold e_{i}(y,s)\Big)d\mu\\
=&\bold R(F,s)
+\bold h_T(s)+\Phi(s)\bold h_T(1-s).\end{aligned}$$

From here it is easy to get the conclusions stated above, by noticing that
 $$\bold h_T(s)-\bold h_1(s)=\Bigg(\sum_{j=1}^l\frac{c_{ij}}{n_{ij}!}
\frac{\partial^{n_{ij}}}{\partial 
 s^{n_{ij}}}\Big(\frac{T^{s+\alpha_{ij}-1}-1}{s+\alpha_{ij}-1}\Big)\Bigg)$$ 
is entire in $s$.
This then completes the proof.

As a direct consequence, taking the spacial case with $F=1$ the 
constant function, we have the following

\noindent
{\bf Corollary.} $$\begin{aligned}\int_{\mathcal D_T}E_\kappa(z,s)\,d\mu
=&\int_0^Ta_0^\infty(y)y^{s-2}dy-\sum_{i=1}^h\phi_{i\infty}\int_T^\infty
a_0^i(y)y^{-1-s}dy\\
=&h_T^\kappa(s)
+\sum_{i=1}^h\phi_{i\kappa}h_T^i(1-s).\end{aligned}$$ 

In fact, the first equality is a direct consequence of 
Equation I. As for the second one, with the case in hand,

\noindent
(i) $F(A_iz)E_\kappa(A_iz,s)-\psi_i(y)e_{i\kappa}=
E_\kappa(A_iz,s)-e_{i\kappa}$ so its integration over $S_\infty(T)$ 
is simply zero;

\noindent
(ii) $R_\kappa(F,s)$ is simply zero by definition.

\subsubsection{Appendix: Rankin-Selberg \& Zagier Method (II)}

Even through we may use the above geometrically oriented method to 
study rank two non-abelian 
zeta functions, we yet give another method of Zagier. Simply put, this 
second one, explained to me by Zagier in 2004,
much simpler than the one outlined above, 
uses the classical Rankin-Selberg for rapid decreasing functions
(to deal with slow increasing functions).

To explain the idea, for simplicity, let us here concentrate with the 
simplest case. Thus,
we at the beginning assume that $\Gamma=SL(2,\mathbb Z)$. 
(The method works in general.)
 
Let $F(z)$ be a $\Gamma$-invariant funtion. Assume first that 

\noindent
(i) $F(z)$ is of slow growth near the cusp $\infty$, that is to say, 
$F(x+iy)=O(y^N)$. 

\noindent
Since this is too weak, we need a bit more stronger condition, say

\noindent
(ii) $F(z)=\phi(y)+O(y^{-N})$ for any $N>0$, and $\phi(y)=y^\alpha$.  

\noindent
(Again here, for simplicity, we have taken $\phi$ to be of the simplest form.)

Being $z\mapsto z+1$ invariant, $F(z)$ admits a Fourier expansion 
$$F(z)=\sum_{m=-\infty}^\infty A_m(y)
e^{2\pi im x}.$$ Thus in particular, by our assumption (i) and (ii),
 $A_0(y)=\phi(y)+O(y^{-N})$ while
$A_{n\not=0}(y)=O(y^{-N})$.

Clearly, $F(z)$ may be written as $y^\alpha$ modulo the 
terms which are quite small. Moreover, 
$\Delta F(z)$, still a $\Gamma$-invariant function, naturally decomposite
as $\alpha(1-\alpha)y^\alpha$ plus the remaining terms, 
which are assumed to be small as well. This latest statement is not a 
direct consequence of (i) and (ii), so we
 may well take this as the additional condition (iii).

Set then $\widehat  F(z):=\Big(\Delta-\alpha(1-\alpha)\Big)(z)$ with 
$\Delta=-y^2(\frac{\partial^2}{\partial x^2}
+\frac{\partial^2}{\partial y^2})$. 
(There is a sign difference bwteen here and elsewhere.) Then we see that 
$\widehat F(z)$ is $\Gamma$-invariant and small. In 
particular, we have a small function $$\widehat  A_0(y)=
-y^2A''_0(y)-\alpha(\alpha-1)A_0(y).$$
Therefore, it makes sense to define its Mellin transform 
$$R(\widehat  F,s):=\int_0^\infty\widehat  A_0(y)y^{s-2}dy,$$
in acoordance with the  Zagier transform 
$$R(F,s):=\int_0^\infty (A_0(y)-\phi(y))y^{s-2}dy$$
introduced above. Set then $$R^*(F,s):=\xi(2s)\cdot R(F,s)
\qquad\mathrm{and}\qquad R^*(\widehat  F,s):=
\xi(2s)\cdot R(\widehat  F,s).$$ Then we have the following:

\noindent
{\bf{\large Lemma}.} (Zagier) {\it With the sam notation as above,} 

\noindent
(1) $$R(F,s)=\frac{1}{(s-\alpha)(1-s-\alpha)}R(\widehat  F,s)\ and 
\ R^*(F,s)=\frac{1}{(s-\alpha)(1-s-\alpha)}
R^*(\widehat  F,s);$$

\noindent
(2) $R^*(F,s)$ {\it admits a meromorphic continuation to the whole complex 
$s$-plane whoes only singularities are 
simple poles at $s=0,1,\alpha, 1-\alpha$;}

\noindent
(3) $R^*(F,s)=R^*(F,1-s)$.

\noindent
Proof. Modulo the clasical Rankin-Selberg method, which can be applied to
 $\widehat  F$ directly, it suffices to
 prove (1). But this is  simply a consequence of  integration by parts. 
Indeed, by definition, $$\begin{aligned}R(\widehat  F,s)
=&\int_0^\infty\widehat  A_0(y)y^{s-2}dy\\
=&-\int_0^\infty A''_0(y)y^sdy+\alpha(\alpha-1)\int_0^\infty A_0(y)y^{s-2}dy\\
=&-\int_0^\infty \Big(A''_0(y)-\phi''(y)\Big)y^sdy+\alpha(\alpha-1)
\int_0^\infty \Big(A_0(y)-\phi(y)\Big)y^{s-2}dy\end{aligned}
$$ since $y^2\phi''(y)=\alpha(\alpha-1)\phi(y)$. Thus, integrating by parts,
$$\begin{aligned}R(\widehat  F,s)=&-\int_0^\infty 
\Big(A_0(y)-\phi(y)\Big)(y^s)''dy+
\alpha(\alpha-1)\int_0^\infty \Big(A_0(y)-\phi(y)\Big)y^{s-2}dy\\
=&\Big(s(1-s)-\alpha(1-\alpha)\Big)\int_0^\infty 
\Big(A_0(y)-\phi(y)\Big)y^{s-2}dy\\
=&\Big(s-\alpha\Big)\Big(1-s-\alpha\Big)R(F,s).\end{aligned}$$ 
This completes the proof.
 
For the reader who does not know the classical Rankin-Selberg, 
let us give a few detail. As a by-product, we give the following

\noindent
{\bf \large Corollary.} (Zagier) $\mathrm{Res}_{s=\alpha}R^*(F,s)=
\frac{1}{s-2\alpha}R^*(\widehat  F,\alpha)=\xi(2\alpha-1).$

\noindent
Proof. Apply the unfolding trick to the function $\widehat  F$, we have
$$R^*(\widehat  F,s)=\xi(2s)\int_0^\infty\widehat  A_0(y)y^{s-2}dy
=\int_{\Gamma\backslash \mathcal H}\widehat  F(z) E^*(z,s)d\mu$$ where 
$d\mu=\frac{dx\,dy}{y^2}$ and $E^*(z,s):=
\xi(2s) E(z,s)$.

Thus $$R^*(\widehat  F,\alpha)=\int_{\mathcal D}\widehat  F(z) E^*(z,s)d\mu
=\int_{\mathcal D_T}\Big(\Delta-\alpha(1-\alpha)\Big)F(z)E^*(z,s)d\mu+O(1).$$
Now using the Green's formula, we obtain, for any two reasonably nice 
functions $F$ and $G$,
 $$\int_{\mathcal D}\Delta F\cdot Gd\mu=\int_{\mathcal D} F\cdot\Delta G
d\mu+\int_{\partial \mathcal D}[F,G]$$
where $$[F,G]:=\frac{\partial F}{\partial \bar z}Gd\bar z+\frac{\partial G}
{\partial z}Fd z
=\Big(\frac{\partial G}{\partial y}F-\frac{\partial F}{\partial y}G\Big)dx
+(\cdots)dy.$$
So going back to  $R^*(\widehat  F,\alpha)$, we get
 $$\begin{aligned}R^*(\widehat  F,\alpha)=&\int_{\mathcal D_T}
\Big(\Delta-\alpha(1-\alpha)\Big)E^*(z,s)\cdot F(z)d\mu
 +\int_{\partial \mathcal D_T}[E^*(z,s),F(z)]+ O(1)\\
=&
\int_{\partial \mathcal D_T}[E^*(z,s),F(z)]+ O(1)\end{aligned}$$ 
since $\Delta E^*(z,s)=\alpha(1-\alpha)E^*(z,s).$
Therefore, note that the contributions coming from the vertical 
boundary and that on $|z|=1$ cancal out 
in pairs,
 $$\begin{aligned}R^*(\widehat  F,\alpha)=&\int_{iT-\frac{1}{2}}^{iT+\frac{1}{2}}
[E^*(z,s),F(z)]+ O(1)\\
 =&\int_{iT-\frac{1}{2}}^{iT+\frac{1}{2}}\Big[\xi(2\alpha)y^\alpha+
\xi(2\alpha-1)y^{1-\alpha}+\exp(*),y^\alpha+
 \exp(*)\Big]+O(1)\end{aligned}$$ where $\exp(*)$ denotes a
certain function with exponentially decay. 
This completes the proof.
\vskip 0.20cm
As Zagier explained, all this may be understood in the framework 
of what may be better called generalized 
Mellin tranform, which we recall here following him.

So let $\phi(t)$ be a nice continuous function, say at least it 
is of polynomial growth on $(0,\infty)$.
Set $$\widetilde \phi_+(s;t_0):=\int_{t_0}^\infty \phi(t) t^{s-1}dt,
\qquad \Re(s)<< 0$$ and 
$$\widetilde \phi_-(s;t_0):=\int_0^{t_0} \phi(t) t^{s-1}dt,
\qquad \Re(s)>> 0.$$
Even though in the definition, $\widetilde \phi_+(s;t_0)$ resp. 
$\widetilde \phi_-(s;t_0)$ are defined for the half 
plane on the far left resp. on the far right, if 
both of these two functions admit  
meromorphic continuations to a common strip $S$ somewhere 
in the middle. Then we may define a new function 
$$\widetilde \phi(s):=\widetilde \phi_+(s;t_0)+\widetilde \phi_+(s;t_0)\qquad
s\in S.$$  As the notation suggests, this function is independent 
of $t_0$.
Moreover, if $\phi$ is of rapidly decreasing, then we get 
$$\widetilde \phi(s)=\int_0^\infty \phi(t)t^{s-1}dt$$
which is just the standard Mellin transform.
Let us give two examples.

\noindent
{\bf \large Ex.} (1) Take $\phi(t)=t^\alpha$, or $t^\alpha(\log t)^n,  
\,\alpha\in \mathbb C,\ n\in\mathbb Z_{\geq 0}.$
Then $\widetilde \phi(s)\equiv 0$ and $S=\mathbb C$. Say, 
$$\int_{t_0}^\infty t^{s+\alpha-1}dt=\frac{1}
{s+\alpha}t^{s+\alpha}\Big|_{t_0}^\infty=-\frac{1}{s+\alpha} 
t_0^{s+\alpha},$$ while $$\int_0^{t_0}t^{s+\alpha-1}dt
=\frac{1}{s+\alpha} t_0^{s+\alpha}.$$

\noindent
{\bf \large Ex.} (2) Take $F$ as in the Rankin-Selberg \& Zagier method,
 we have $$R(F,s)=\widetilde A_0(s-1),\qquad s\in\mathbb C.$$

\noindent 
At this point, I recall in one series of lectures, a kind of 
mini course  around 1990-1991 at MPI f\"ur
Mathematik at Bonn on HyperGeometric Functions and 
Differential Equations,
Zagier reminded the audience, as usual in his extremely fast fashion,
that Ex could mean example, and could also well mean 
exercise. So let us follow him and leave the details 
of Ex.2 to the reader. 

\subsection{Rank Two Non-Abelian Zeta Function For $\mathbb Q$}

Recall that if we set 
$\mathcal D_T:=\Big\{x\in\mathcal D:y=\Im(z)\leq T\Big\}$,  
the points in $\mathcal D_T$ are in one-to-one corresponding with  rank two 
$\mathbb Z$-lattices (in $\mathbb R^2$) of 
volume one whose first Minkowski successive minimums $\lambda_1$ satisfying 
$\lambda_1(\Lambda)\geq T^{-1/2}$. 
Thus if we set $\mathcal M_{\mathbb Q,2}^{\leq \frac{1}{2}\log T}\Big[1\Big]$ 
be the moduli space of rank two 
$\mathbb Z$-lattices $\Lambda$ of volume 1 (over $\mathbb Q$) whose 
sublattices of rank one have degree 
$\leq\frac {1}{2}\log T$, then up to a measure zero subset,
{\it there is a natural one-to-one and onto morphism 
$$\mathcal M_{\mathbb Q,2}^{\leq \frac{1}{2}\log T}\Big[1\Big]
\simeq \mathcal D_T.$$} 
In particular, the corresponding
moduli space of semi-stable lattices is given by
$$\mathcal M_{\mathbb Q,2}^{\leq 0}\Big[1\Big]=\mathcal M_{\mathbb Q,2}
\Big[1\Big]\simeq \mathcal D_1.$$
Moreover, motivated by our definition of non-abelian zeta functions, 
we introduce a (generalized) rank two zeta function 
$\xi_{\mathbb Q,2}^T(s)$ by setting
$$\xi_{\mathbb Q,2}^T(s):=\int_{\mathcal D_T}\widehat E(z,s)\,
\frac{dx\wedge dy}{y^2},\qquad 
\Re(s)>1.$$ Then we have the following 
\vskip 0.30cm
\noindent
{\bf \Large Fact} (VIII$)_{\mathbb Q}$ {\it For the generalized 
zeta function $\xi_{\mathbb Q,2}^T(s)$,
$$\xi_{\mathbb Q,2}^T(s)=\frac{\xi(2s)}{s-1}\cdot T^{s-1}-
\frac{\xi(2s-1)}{s}\cdot T^{-s}.$$  In particular, 
the rank two non-abelian zeta function $\xi_{\mathbb Q,2}(s)$ for 
the field of rationals $\mathbb Q$ is given by}
$$\xi_{\mathbb Q,2}(s)=\frac{\xi(2s)}{s-1}-
\frac{\xi(2s-1)}{s},\qquad\Re(s)>1.$$

\noindent
Proof. This is a direct consequence of the 
Corollary in 4.1.1. Indeed,  it is well known that the Fourier 
expansion of 
$\widehat E(z,s)$ is given by
$$\widehat E(z,s)=\xi(2s)y^{s}+\xi(2s-1)y^{1-s}+\mathrm{non-constant\ term}.$$
Hence we have $$\begin{aligned}\xi_{\mathbb Q,2}^T(s)=&\int_0^T
\Big(\xi(2s)y^{s}\Big)\frac{dy}{y^2}-
\int_T^\infty\Big(\xi(2s-1)y^{1-s}\Big)\frac{dy}{y^2}\\
=&\frac{\xi(2s)}{s-1}\cdot T^{s-1}-
\frac{\xi(2s-1)}{s}\cdot T^{-s}.\end{aligned}$$
This completes the proof.

\noindent
{\bf Remarks.} (1) Even though originally $T\geq 1$, we may extend it as a 
function
of complex variable $T$ in terms of the right hand side. Denote this resulting 
function also by $\xi_{\mathbb Q,2}^T(s)$. 
Surely, for $\Re(s)>1$, if $T$ is 
real and $T\geq 1$, then $\xi_{\mathbb Q,2}^T(s)$ is simply the integration of
$\widehat E(z,s)$ over the domain $\mathcal D_T$.
 Based on this, even when
 $T$ is real and  $0< T\leq 1$, we  have a geometric interpretation for 
$\xi_{\mathbb Q,2}^T(s)$: it is simply the combination 
$$\Big(\int_{D_{1,T}}-\int_{D_{-1,T}}\Big) \widehat E(z,s)\cdot
\frac{dx\wedge dy}{y^2},$$ where $D_{1,T}:=\mathcal D\cap 
\Big\{z=x+iy: y\leq T, |x|\leq\frac{1}{2}\Big\}$ and $D_{-1,T}:=
\Big\{z\in\mathcal H: |z|\leq 1\Big\}\cap \Big\{z=x+iy: y\geq T, |x|\leq\frac{1}{2}\Big\};$

\noindent
(2) By taking the residue at $s=1$, we have   
$$\Big(\mathrm{Res}_{s=1}\widehat E(z,s)\Big)\cdot 
\mathrm{Vol}\Big(\mathcal D_1\Big)=\xi(2)-\mathrm{Res}_{s=1}\xi(2s-1);$$

\noindent
(3) We see, in particular, for half positive integers $n\geq \frac {3}{2}$,
$$\Big((n-1)n\Big)\cdot \xi_{\mathbb Q,2}(n)
=n\cdot \xi(2n)-(n-1)\cdot \xi(2n-1).$$ So the special values of the 
Riemann zeta function at two successive integers are related naturally
via the special values of  rank two non-abelian zeta function. 
This clearly is a fact which should be taken very seriously. In particular, in view of Remark (1) above, we suggest the reader to see what happens for small
$n$'s by writting out the non-abelian zeta in terms of the integrations for 
the terms defining Eisenstein series. With this, it is very likely that
the reader will be convinced that, say, when talking about specail values 
of $\xi(s)$ at odd integers, it is better to distinguish the values at 
$4\mathbb Z_{>0}-1$ from these at $4\mathbb Z_{> 0}+1$.

\section{Upper Half Space Model: Rankin-Selberg Method}

In this subsection, we discuss a generalization of the original Rankin-Selberg 
method for upper half space. This should be known to experts. But as we can 
hardly find any details in the literature, so we decide to write all the 
details down.

For the hyperbolic upper half space $\mathbb H$, with respect to the 
hyperbolic metric 
$ds^2:=\frac{dx^2+dy^2+dr^2}{r^2}$, the volume form is 
$d\mu:=\frac{dx\wedge dy\wedge dr}{r^3}$ while the corresponding 
Laplace operator becomes
$$\Delta:=r^2\Big(\frac{\partial^2}{\partial x^2}+\frac{\partial^2}
{\partial y^2}
+\frac{\partial^2}{\partial r^2}\Big)-r\frac{\partial}{\partial r}.$$
Also recall that the standard Stokes' formula over $\mathbb R^3$ is simply
$$\iiint_D\Big(\frac{\partial}{\partial x}A+\frac{\partial}{\partial y}B+
\frac{\partial}{\partial r}C\Big)\,dx\wedge dy\wedge dr=\iint_{\partial D}
\Big(Ady\wedge dr+Bdr\wedge dx+Cdx\wedge dy\Big)$$ for a 3 dimensional domain 
$D$ in $\mathbb R^3$ with boundary $\partial D$. Note in particular 
that for any two nice functions $f$ and $g$ on $D$,
$$\begin{aligned}
\frac{\partial}{\partial x}\Big(\frac{\partial}{\partial x}f\cdot g\Big)
-\frac{\partial}{\partial x}f\cdot \frac{\partial}{\partial x}g=&
\frac{\partial^2}{\partial x^2}f\cdot g,\\ 
\frac{\partial}{\partial x}\Big(f\cdot \frac{\partial}{\partial x} g\Big)
-\frac{\partial}{\partial x}f\cdot \frac{\partial}{\partial x}g=&
f\cdot 
\frac{\partial^2}{\partial x^2}g.\end{aligned}$$ Similar formulas
holds with respect to
$\frac{\partial}{\partial y}$ and $\frac{\partial}{\partial r}.$
Consequently, 
$$\begin{aligned}&\Big(\frac{\partial^2}{\partial x^2}+
\frac{\partial^2}{\partial y^2}
+\frac{\partial^2}{\partial r^2}\Big)f\cdot g-f\cdot \Big(\frac{\partial^2}
{\partial x^2}+\frac{\partial^2}{\partial y^2}+\frac{\partial^2}
{\partial r^2}\Big)g\\
=&
\frac{\partial}{\partial x}\Big(\frac{\partial}{\partial x}f\cdot g
-f\cdot \frac{\partial}{\partial x}g\Big)+\frac{\partial}{\partial y}
\Big(\frac{\partial}{\partial y}f\cdot g-f\cdot \frac{\partial}{\partial y}
g\Big)
+\frac{\partial}{\partial r}\Big(\frac{\partial}{\partial r}f\cdot g
-f\cdot \frac{\partial}{\partial r}g\Big).\end{aligned}$$
Clearly, the left hand side is simply 
$$\begin{aligned}&\Big(\frac{1}{r^2}\Delta+
\frac{1}{r}\frac{\partial}{\partial r}\Big)f\cdot g-f\cdot 
\Big(\frac{1}{r^2}\Delta+\frac{1}{r}\frac{\partial}{\partial r}\Big)g\\
=&\frac{1}{r^2}\Big(\Delta f\cdot g-f\Delta g\Big)
+\frac{1}{r}\Big(\frac{\partial f}
{\partial r}\cdot g-f\cdot \frac{\partial g}{\partial r}\Big).\end{aligned}$$ 
Hence
$$\begin{aligned}&\frac{1}{r^2}\Big(\Delta f\cdot g-f\Delta g\Big)\\
=&
\frac{\partial}{\partial x}\Big(\frac{\partial}{\partial x}f\cdot g-f\cdot 
\frac{\partial}{\partial x}g\Big)+\frac{\partial}{\partial y}
\Big(\frac{\partial}{\partial y}f\cdot g-f\cdot \frac{\partial}
{\partial y}g\Big)\\
&+\frac{\partial}{\partial r}\Big(\frac{\partial}{\partial r}f\cdot g
-f\cdot \frac{\partial}{\partial r}g\Big)
-\frac{1}{r}\Big(\frac{\partial f}{\partial r}\cdot g-f\cdot \frac{\partial g}
{\partial r}\Big).\end{aligned}$$ This implies that
$$\begin{aligned}&\iiint_D\Big(\Delta f\cdot g-f\Delta g\Big)\,\frac{dx\wedge dy\wedge dr}{r^3}\\
=&\iiint\Big[
\frac{\partial}{\partial x}\Big(\frac{\partial}{\partial x}f\cdot g
-f\cdot \frac{\partial}{\partial x}g\Big)+\frac{\partial}{\partial y}
\Big(\frac{\partial}{\partial y}f\cdot g-f\cdot \frac{\partial}
{\partial y}g\Big)\\
&+\frac{\partial}{\partial r}\Big(\frac{\partial}{\partial r}f\cdot g
-f\cdot \frac{\partial}{\partial r}g\Big)
-\frac{1}{r}\Big(\frac{\partial f}{\partial r}\cdot g-f\cdot 
\frac{\partial g}{\partial r}\Big)\Big]\,\frac{dx\wedge dy\wedge dr}{r}\\
=&\iiint\Big[
\frac{\partial}{\partial x}\Big(\frac{1}{r}
\Big(\frac{\partial}{\partial x}f\cdot g
-f\cdot \frac{\partial}{\partial x}g\Big)\Big)
+\frac{\partial}{\partial y}\Big(\frac{1}{r}
\Big(\frac{\partial}{\partial y}f\cdot g-f\cdot \frac{\partial}
{\partial y}g\Big)\Big)\\
&+\frac{\partial}{\partial r}\Big(\frac{1}{r}\big(\frac{\partial}{\partial r}
f\cdot g-f\cdot \frac{\partial}{\partial r}g\big)\Big)
+\frac{1}{r^2}\Big(\frac{\partial f}{\partial r}\cdot g-f\cdot 
\frac{\partial g}{\partial r}\Big)
-\frac{1}{r^2}\Big(\frac{\partial f}{\partial r}\cdot g-f\cdot 
\frac{\partial g}{\partial r}\Big)
\Big] \, dx\wedge dy\wedge dr\\
=&\iiint\Big[
\frac{\partial}{\partial x}\Big(\frac{1}{r}
\Big(\frac{\partial}{\partial x}f\cdot g
-f\cdot \frac{\partial}{\partial x}g\Big)\Big)
+\frac{\partial}{\partial y}\Big(\frac{1}{r}
\Big(\frac{\partial}{\partial y}f\cdot g
-f\cdot \frac{\partial}{\partial y}g\Big)\Big)\\
&+\frac{\partial}{\partial r}\Big(\frac{1}{r}\big(\frac{\partial}{\partial r}
f\cdot g-f\cdot \frac{\partial}{\partial r}g\big)\Big)\Big] \, dx\wedge dy
\wedge dr\\
=&\iint_{\partial D}
\Big(\frac{1}{r}\frac{\partial}{\partial x}f\cdot g-\frac{1}{r}f\cdot
 \frac{\partial}{\partial x}g\Big)dy\wedge dr
+\Big(\frac{1}{r}\frac{\partial}{\partial y}f\cdot g-\frac{1}{r}f\cdot
 \frac{\partial}{\partial y}g\Big)dr\wedge dx\\
&+\Big(\frac{1}{r}\frac{\partial}{\partial r}f\cdot g-\frac{1}{r}f\cdot
 \frac{\partial}{\partial r}g\Big)dx\wedge dy,
\end{aligned}$$ where in the last equality, we used the Stokes formula.
Note that $$\begin{aligned}&\iiint_D\Big((\Delta-2s(2s-2) )f\cdot g
-f(\Delta-2s(2s-2) g\Big)\frac{dx\wedge dy\wedge dr}{r^3}\\
=&\iiint_D\Big(\Delta f\cdot g-
f\Delta g\Big)\,\frac{dx\wedge dy\wedge dr}{r^3},\end{aligned}$$ 
we get the following

\noindent
{\bf \large Basic Formula I.} 
({\bf Stokes' Formula for Hyperbolic Geometry}). 
{\it With the same notationa as above,
$$\begin{aligned}&\iiint_D\Big((\Delta-2s(2s-2) )f\cdot g
-f(\Delta-2s(2s-2) g\Big)
\frac{dx\wedge dy\wedge dr}{r^3}\\
=&\iint_{\partial D}
\Big(\frac{\partial}{\partial x}f\cdot g-f\cdot \frac{\partial}
{\partial x}g\Big)dy\wedge \frac{dr}{r}
+\Big(\frac{\partial}{\partial y}f\cdot g-f\cdot \frac{\partial}
{\partial y}g\Big)\frac{dr}{r}\wedge dx\\
&+\Big((\frac{1}{r}\frac{\partial}{\partial r})f\cdot g-f\cdot 
(\frac{1}{r}\frac{\partial}{\partial r})g\Big)dx\wedge dy.\end{aligned}
$$}

Taking an example with $f=1$ and $g=\widehat E(s)=\widehat E(P,s)$, 
by the fact that 
$\Delta \widehat E(s)=2s(2s-2)\widehat E(s)$, we  get 
$$\begin{aligned}&-2s(2s-2)\iiint_D\widehat E(s)
\frac{dx\wedge dy\wedge dr}{r^3}\\
=&\iiint_D\Big(-2s(2s-2)\cdot \widehat E(s)-1\cdot 0\Big)
\frac{dx\wedge dy\wedge dr}{r^3}\\
=&
\iiint_D\Big((\Delta-2s(2s-2))1\cdot \widehat E(s)-1\cdot (\Delta-2s(2s-2)) 
\widehat E(s)\Big)\frac{dx\wedge dy\wedge dr}{r^3}\\
=&\iint_{\partial D}
\Big[-\frac{\partial}{\partial x}\widehat E(s)\, dy\wedge \frac{dr}{r}
-\frac{\partial}{\partial y}\widehat E(s)\, \frac{dr}{r}\wedge dx
-(\frac{1}{r}\frac{\partial}{\partial r})\widehat E(s)\, dx\wedge dy.
\end{aligned}$$ 
That is to say,
$$\begin{aligned}&\iiint_D\widehat E(s)\frac{dx\wedge dy\wedge dr}{r^3}\\
=&\frac{1}{2s(2s-2)}
\iint_{\partial D}\Big[\frac{\partial}{\partial x}\widehat E(s)\, dy\wedge 
\frac{dr}{r}
+\frac{\partial}{\partial y}\widehat E(s)\, \frac{dr}{r}\wedge dx
+(\frac{1}{r}\frac{\partial}{\partial r})\widehat E(s)\, dx\wedge dy.\end{aligned}$$

Now let us concentrate on the case we are really intertested in, i.e., taking 
$D=\mathcal D_T$, where $\mathcal D_T$ is the compact part of the fundamental 
domain $\mathcal D$ of $\Gamma$ by cutting off the cusp neighborhoods defined 
by the condition $r>T$. That is to say 
$$\mathcal D_T=\mathcal D\,\backslash\,\cup_{i=1}^h
\widetilde{\mathcal F_{\eta_i}}(T).$$ 
Note that the boundary of $\mathcal D_T$ consists of 
surfaces given by $\mu(\sigma,P)=\mu(\tau,P)$ with
 $\sigma\not=\tau$, so in fact we can even pair the boundary surfaces 
together while with opposite norm directions, due to that fact that the 
angles in the hyperbolic geometry are really the same as their 
corresponding angles when measured in terms of Euclidean metric. 
Therefore, after such a cancellation, what left on the right hand side of
the above Stokes' formula for $\widehat E$ is only the part concerning
the integration over the surfaces $P_i(T)$ obtained by intersecting 
the cuspidal neighborhood with the surfaces $r=T$. 

Furthermore, note that in the integration, $$\iint_{P_i(T)}
\frac{\partial}{\partial x}\widehat E(s)\, dy\wedge \frac{dr}{r},$$ 
being taking derivatives along $x$ direction, the constant terms in which only
variable $r$ is involved do not contribute, while for non-constant terms, 
the avarage on 
$y$ for the exponential function $e^{2\pi i\langle \omega',z\rangle}$ 
contributes exactly zero, so $$\iint_{P_i(T)}\frac{\partial}{\partial x}
\widehat E(s)\, dy\wedge \frac{dr}{r}=0.$$ Similarly,
$$\iint_{P_i(T)}\frac{\partial}{\partial y}\widehat E(s)\,  \frac{dr}{r}
\wedge dx=0.$$ Thus we are left with the integration $$
\iint_{P_i(T)}\Big(\frac{1}{r}\frac{\partial}{\partial r}\Big)\widehat E(s)\, dx\wedge dy.$$
That is to say, we have arrived at
the following 

\noindent
{\bf\large  Basic Formula II.} {\it With the same notation as above, we have
$$\iiint_D\widehat E(s)\frac{dx\wedge dy\wedge dr}{r^3}=\frac{1}{2s(2s-2)}
\sum_{i=1}^h\iint_{P_i(T)}(r\frac{\partial}{\partial r})\widehat E(s)\, 
dx\wedge dy.$$}

Clearly, for this last integral, only constant terms contribute since 
the average of the exponential over $dx$ or $dy$ gives exactly zero.
In this way, if we set the constant term of $\widehat E(s)$ at the $i$-th cusp
to be $$a_i\cdot r^{2s}+b_i\cdot r^{2-2s},\qquad a_i,\ b_i\ \mathrm{constants},$$ then
$$\begin{aligned}&\frac{1}{2s(2s-2)}\iint_{P_i(T)}
\Big(\frac{1}{r}\frac{\partial}{\partial r}\Big)
\widehat E(s)\, 
dx\wedge dy\\
=&\Big(\frac{1}{r}\frac{\partial}{\partial r}\Big)\widehat E(s)|_{r=T}
\cdot \iint_{P_i(T)}dx\wedge dy\\
=&\Big(\frac{a_i}{s-1}\cdot T^{2s-2}+\frac{b_i}{-s}\cdot T^{-2s}\Big)
\mathrm{Vol}(Q_i).\end{aligned}$$
Therefore, we have
\vskip 0.20cm
\noindent
{\bf{\Large Proposition}.} {\it With the same notation as above,
$$\iiint_{\mathcal D_T}\widehat E(s)\frac{dx\wedge dy\wedge dr}{r^3}= 
\sum_{i=1}^h\Big(\frac{a_i}{s-1}\cdot T^{2s-2}-\frac{b_i}{s}\cdot T^{-2s}\Big)
\mathrm{Vol}(Q_i).$$}

\section{Rank Two $\mathcal O_K$-Lattices}
\subsection{Rankin-Selberg Method}
In this section, $T$ is assumed to be a positive real number $\geq 1$.

Now let us compute the integration $$\iiint_{\mathcal D_T}
\widehat E_{2,\frak a}(\tau,s)d\mu(\tau).$$
Here $\mathcal D_T$ is the compact part obtained from the fundamental domain
$\mathcal D$ for $SL(\mathcal O_K\oplus\frak a)\Big\backslash
\Big(\mathcal H^{r_1}\times \mathbb H^{r_2}\Big)$ by cutting off the cusp 
neighbouhoods defined by the conditions that the distance to cusps is
less than $T^{-1}$. (Recall that, as such, $\mathcal D_1$ is simply the part 
corresponding to semi-stable lattices.) 
We use the Rankin-Selberg \& Zagier method,
but in its simplest form as a generalization of the one stated in
the previous section (hoping that we will come back later for 
general cases).

For doing so, let us first formulate the integration
$$\iiint_{\mathcal D_T}\Big(\Delta_K\,\widehat E_{2,\frak a}(\tau,s)\Big)\,
d\mu(\tau)$$
where $$\Delta_K:=\sum_{\sigma:\mathbb R}\Delta_\sigma+
\sum_{\tau:\mathbb C}\Delta_\tau$$ with $$\Delta_\sigma:=
y_\sigma^2\Big(\frac{\partial^2}{\partial x_\sigma^2}+
\frac{\partial^2}{\partial y_\sigma^2}\Big)$$ and $$
\Delta_\tau:=r_\tau^2\Big(\frac{\partial^2}{\partial x_\tau^2}+
\frac{\partial^2}{\partial y_\tau^2}+\frac{\partial^2}{\partial r_\tau^2}\Big)
-r\frac{\partial}{\partial r_\tau}.$$ (For the time being, by an abuse
of notation, we use $\Delta_K$
to denote the hyperbolic Laplace operator for the space
$\mathcal H^{r_1}\times \mathbb H^{r_2}$, not the absolute value of the
discriminant of $K$ which accordingly is changed to $D_K$.)

Note that 
$$\Delta_\sigma\Big(y_\sigma^s\Big)=s(s-1)\cdot y_\sigma^s,$$ while
$$\Delta_\tau \Big(r_\tau^{2s}\Big)=2s(2s-2)\cdot r_\tau^{2s}$$ by the 
$SL$-invariance of the metrics, we conclude hence that 
$$\Delta_K\Big(\widehat E_{2,\frak a}(\tau,s)\Big)=
\Big(r_1\cdot\Big(s(s-1)\Big)+r_2\cdot\Big(2s(2s-2)\Big)\Big)\cdot 
\widehat E_{2,\frak a}(\tau,s),\qquad\Re(s)>1.$$ Hence 
$$\iiint_{\mathcal D_T}\widehat E_{2,\frak a}(\tau,s)\,d\mu(tau)=
\frac{r_1+4r_2}{s(s-1)}\iiint_{\mathcal D_T}\Delta_K
\widehat E_{2,\frak a}(\tau,s)\,d\mu(\tau).$$

On the other hand, using Stokes' Formula, we have
$$\begin{aligned}&\iiint_{\mathcal D_T }\Delta_K\widehat E_{2,\frak a}
(\tau,s)d\mu(\tau)\\
=&
\iiint_{\mathcal D_T }\Big(\Delta_K\widehat E_{2,\frak a}(\tau,s)\Big)
\cdot 1d\mu(\tau)-
\iiint_{\mathcal D_T } \widehat E_{2,\frak a}(\tau,s)
\cdot\Big(\Delta_K 1\Big)\,d\mu(\tau)\\
=&
\iiint_{\mathcal D_T }\Bigg(\Big(\Delta_K\widehat E_{2,\frak a}(\tau,s)\Big)
\cdot 1- \widehat E_{2,\frak a}(\tau,s)
\cdot\Big(\Delta_K 1\Big)\Bigg)\,d\mu(\tau)\\
=&\iint_{\partial D(T)}\Big(
\frac{\partial \widehat E_{2,\frak a}(\tau,s)}{\partial \nu}\cdot 1-
\widehat E_{2,\frak a}(\tau,s)\cdot
\frac{\partial 1}{\partial \nu}\Big)d\mu\\
=&\iint_{\partial D(T)}
\frac{\partial \widehat E_{2,\frak a}(\tau,s)}{\partial \nu}
d \mu\end{aligned}$$
where $\frac{\partial}{\partial \nu}$ is the outer normal derivative and 
$d \mu$ is the volume element of the boundary $\partial \mathcal D_T $. 

To calculate this latest integration, we start with a trick initially
 used by Siegel (and hence used in many other places such as [Ge] 
and [Ef]) to make the following change of variables at the cusps. 
As such, our presentation will follow them, in particular, that of [Ef].
(Here the reader can have really a good test on, say at least, 
how our mathematics  is different from others and 
how mathematics as a part of our culture is developed -- 
By choosing the way of 
presenting our work in such a form, we here are trying to indicate that 
what a kind of  mathematics we pursue: elegant, fundamental, 
and deeply rooted in the classics. The researcher may be changed, 
but the beauty and the 
essence of mathematics should not and will not be changed. 
It is not about fashionable or popular, but about whether it becomes 
a permanent part of mathematics.)

As usual, we proceed our discussion by transforming cusps $\eta$ to 
$\infty$ using suitable conjugations.
 
Two directions have to be studied: the ReZ direction for $x_\sigma$ resp.
$z_\tau=x_\tau+iy_\tau$,
and the ImJ directions for $y_\sigma$ resp. for $r_\tau=:v_\tau$ when $\sigma$ 
is real resp. $\tau$ is complex. 
As used in the discussion for fundamental domains, 
the change with respect to the ReZ direction is simpler, while
 the change with respect to the ImJ direction is a bit complicated.
More precisely, to deal with ReZ direction, there are two options, namely

\noindent
(1) use a $\mathbb Z$-basis $\omega_1,\ldots,\omega_{n=r_1+2r_2}$
for the $\mathbb Z$-lattice $\frak a\frak b^{-2}$
associated to the cusp $\eta=\left[\begin{matrix}\alpha\\ \beta\end{matrix}
\right]$,
where $\frak b:=\mathcal O_K\alpha+\frak a\beta$ so that under this
change the fundamental domain for the $\mathbb Z$-lattice 
$\frak a\frak b^{-2}$ is changed to the one given by
$$|x_\sigma|\leq \frac{1}{2},\ \sigma\ \mathrm{real};\qquad
|x_\tau|\leq \frac{1}{2},\ |y_\tau|\leq \frac{1}{2},\ \tau\ \mathrm{complex};$$

\noindent
(2) do not make any change so that the corresponding $\mathbb Z$-lattice
remaining to be $\frak a\cdot\frak b^{-2}$.
 
Clearly, all of these two are just affine transformations 
hence easy to be handled. 
We decide here to adopt the second one, that is, 
we will not change the variables $x_{\sigma_1},\cdots, x_{\sigma_{r_1}}, 
z_{\tau_1},\cdots, z_{\tau_{r_2}}$, while pointing out that (1) is the 
one used by Siegel as cited in 2.4.

Now let us turn to the ImJ direction. Recall that here all
components are positive. In particular,
$(y_{\sigma_1},\ldots, y_{\sigma_{r_1}},v_{\tau_1},\ldots, v_{\tau_{r_2}})
\in\mathbb R^{r_1+r_2}_+$
resulting from the ImJ direction of a point 
$(z_1,\cdots, z_{r_1},\,P_1,\cdots,P_{r_2})$ in 
$\mathcal H^{r_1}\times\mathbb H^{r_2}$ admits a natural
norm 
$$\begin{aligned}N(y_{\sigma_1},&\ldots, y_{\sigma_{r_1}},v_{\tau_1},
\ldots, v_{\tau_{r_2}})=\\
=&
\Big(y_{\sigma_1}\cdot\ldots\cdot y_{\sigma_{r_1}}\Big)
\cdot\Big(v_{\tau_1}\cdot
\ldots\cdot v_{\tau_{r_2}}\Big)^2\\
=&\prod_{\sigma:\,\mathbb R}y_\sigma\cdot
\prod_{\tau:\,\mathbb C}v_\tau^2.\end{aligned}$$ 
The key here is that we need to find a variable change
for the ImJ directions so that

\noindent 
(a) the outer normal direction will be seen more clearly; and

\noindent
(b) the fundamental domain for the stablizer group of cusps can be written 
in a very simple way.
 
The generalized version of Siegel's change of variables in the discussion 
of fundamental domains
does exactly this. It is carried out by replacing the 
original variables  $y_{\sigma_1},\cdots, 
y_{\sigma_{r_1}}, v_{\tau_1},\cdots, v_{\tau_{r_2}}$ with the
new variables $Y_0, Y_1,\cdots, Y_{r_1+r_2-1}$. To give a precise definition,
let $\varepsilon_1,\cdots,\varepsilon_{r_1+r_2-1}$ be a generator of the unit
group $U_K$ (modulo the torsion). Then by Dirichlet's Unit Theorem, the matrix
$$\left(\begin{matrix} 1&\log|\varepsilon_1^{(1)}|&\cdots&
\log|\varepsilon_{r_1+r_2-1}^{(1)}|\\
\cdots&\cdots&\cdots&\cdots\\
1&\log|\varepsilon_1^{(r_1+r_2)}|&\cdots&
\log|\varepsilon_{r_1+r_2-1}^{(r_1+r_2)}|\end{matrix}\right)$$
is invertible. Set $(e_j^{(i)})_{i=0,j=1}^{r_1+r_2-1,r_1+r_2}$ be its
inverse. Then by definition and an obvious calculation, 

\noindent
i) the entries of the
first row is given by $e_j^{(0)}=\frac{1}{r_1+r_2},\ j=1,2,\cdots, r_1+r_2$;

\noindent
ii) $\sum_{j=1}^{r_1+r_2}e_j^{(i)}=0,\qquad i=1,\cdots, r_1+r_2-1$;
and

\noindent
iii) $\sum_{j=1}^{r_1+r_2}e_j^{(i)}\,\log|\varepsilon_k^{(j)}|=\delta_{ik},
\qquad i,k=1,\cdots, r_1+r_2-1.$

\noindent
In particular, $$(e_j^{(i)})=
\left(\begin{matrix}
\frac{1}{r_1+r_2}&\frac{1}{r_1+r_2}&\cdots&\frac{1}{r_1+r_2}\\
e_1^{(1)}&e_2^{(1)}&\cdots&e_{r_1+r_2}^{(1)}\\
\cdots&\cdots&\cdots&\cdots\\
e_1^{(r_1+r_2-1)}&e_2^{(r_1+r_2-1)}&\cdots&e_{r_1+r_2}^{(r_1+r_2-1)}
\end{matrix}\right).$$

With this, make a change of variables by 
$$\begin{aligned}Y_0:=&N(y_{\sigma_1},\ldots, y_{\sigma_{r_1}},
v_{\tau_1},\ldots, v_{\tau_{r_2}})
=\prod_{\sigma:\mathbb R}y_{\sigma}\cdot\prod_{\tau:\mathbb C} v_{\tau}^2,\\
Y_1:=&\frac{1}{2}\Big(\sum_{i=1}^{r_1}e_i^{(1)}\log y_{\sigma_i}+
 \sum_{j=1}^{r_2}e_{r_1+j}^{(1)}\log v_{\tau_j}^2\Big)\\
&\cdots\ \cdots\ \cdots\ \cdots\\
Y_{r_1+r_2-1}:=&\frac{1}{2}\Big(\sum_{i=1}^{r_1}e_i^{(r_1+r_2-1)}
\log y_{\sigma_i}+\sum_{j=1}^{r_2}e_{r_1+j}^{(r_1+r_2-1)}\log v_{\tau_j}^2\Big)
\end{aligned}$$
Consequently, by inverting these relations, we have
$$\begin{aligned}y_{\sigma_i}=& Y_0^{\frac{1}{r_1+r_2}}\prod_{q=1}^{r_1+r_2-1}
\Big|\varepsilon_q^{(i)}\Big|^{2Y_q},\qquad i=1,\cdots, r_1,\\
v_{\tau_j}^2=& Y_0^{\frac{1}{r_1+r_2}}\prod_{q=1}^{r_1+r_2-1}
\Big(\Big|\varepsilon_q^{(r_1+j)}\Big|^{2Y_q}\Big)^2,
\qquad j=1,\cdots, r_2.\end{aligned}$$

Further, by taking the fact that $N_\tau=2$ for complex places $\tau$, 
for later use,
we set $$t_j:=v_j^2=Y_0^{\frac{1}{r_1+r_2}}
\prod_{q=1}^{r_1+r_2-1}
\Big(\Big|\varepsilon_q^{(r_1+j)}\Big|^2\Big)^{2Y_q},\qquad j=1,\cdots, r_2.$$

With this change of variables,
from the precise construction of the fundamental domain 
for the action of $\Gamma_\eta$ in $SL(\mathcal O_K\oplus \frak a)$ on
$\mathcal H^{r_1}\times\mathbb H^{r_2}$ in 2.4, it now becomes clear that
this fundamental domain for the action of $\Gamma_\eta$ on 
$\mathcal H^{r_1}\times\bold{\mathbb H}^{r_2}$ is simply given by 
$$\begin{aligned} &0<Y_0<\infty,\qquad -\frac{1}{2}\leq Y_1,\cdots, 
Y_{r_1+r_2-1}\leq\frac{1}{2},\\
&(x_{\sigma_1},\cdots, x_{\sigma_{r_1}}; z_{\tau_1},\cdots, z_{\tau_{r_2}})\in 
\mathcal F_\eta(\frak a\frak b^{-2}),\end{aligned}$$ where 
$\mathcal F_\eta(\frak a\frak b^{-2})$ denotes a fundamental parallelopiped 
associated with the lattice $\frak a\frak b^{-2}$ in 
$\mathbb R^{r_1}\times \mathbb C^{r_2}$.

To go further, we need to know the precise change of the volume forms 
in accordance with the above change of variables.
So we must compute some of the Riemannian geometric invariants in terms 
of these coordinates at the cusps. Clearly the hyperbolic metric on 
$\mathcal H^{r_1}\times\mathbb H^{r_2}$ is given by
$$g=\left(\begin{matrix}g_{\mathrm{ImJ}}&0\\
0&g_{\mathrm{ReZ}}\end{matrix}\right)$$ with the matrics for the ImJ 
and ReZ directions being given by
$$g_{\mathrm{ImJ}}=g_{\mathrm{ReZ}}=
\mathrm{diag}\Big(\frac{1}{y_1^2},\,\cdots,\, 
\frac{1}{y_{r_1}^2},\
\Big(\frac{1}{v_1^2}\Big)^2,\,\cdots \,\Big(\frac{1}{v_{r_2}^2}\Big)^2\Big).$$
(Here we reminder the reader that the twist resulting from
$N_\tau=2$ for complex places $\tau$ is entered in the discussion: 
In the above matrix we used $\Big(\frac{1}{v_j^2}\Big)^2$ 
instead of a simple $\frac{1}{v_j^2}$.)

Since there is no changes in the ReZ directions, which 
play a totally independent role,  hence the ReZ part of the matrix
for the metric remains the same. Now let us turn to ImJ directions.
From above, the metric for the ImJ directions  is given by
$$(g_{ij}):=g_{\mathrm{ImJ}}=
\mathrm{diag}\Big(\frac{1}{y_1^2},\,\cdots,\, 
\frac{1}{y_{r_1}^2},\
\Big(\frac{1}{v_1^2}\Big)^2,\,\cdots \,\Big(\frac{1}{v_{r_2}^2}\Big)^2\Big).$$
In general term, to find the matrix $(\tilde g_{ij})$ obtained from 
$(g_{ij})$ by the change of variables, we first need to calculate the 
partial derivatives so as to get $\tilde g_{ij}$ from the formula
$$\tilde g_{ij}=\sum_{\alpha,\beta}\frac{\partial x^\alpha}
{\partial\tilde  x^i}\frac{\partial x^\beta}{\partial\tilde  x^j}
g_{\alpha\beta}.$$ (Here, in terms of $g_{ij}$ and $\tilde g_{ij}$, the 
variables are assumed to be renumbered as $x^1,x^2,\cdots, x^{r_1+r_2}$ 
and $\tilde x^1,\tilde x^2,\cdots, \tilde x^{r_1+r_2}$ respectively.)

By an obvious calculation,
 $$\begin{aligned}\frac{\partial y_{\sigma_i}}{\partial Y_0}
=&\frac{1}{(r_1+r_2)Y_0}y_{\sigma_i},\qquad i=1,\cdots, r_1,\\ 
\frac{\partial y_{\sigma_i}}{\partial Y_q}
=& 2\log\Big|\varepsilon_q^{(i)}\Big|\cdot y_{\sigma_i},
\qquad i=1,\cdots, r_1,\ q=1,\cdots,r_1+r_2-1\\
\frac{\partial t_{\tau_j}}{\partial Y_0}
=&\frac{1}{(r_1+r_2)Y_0}t_{\tau_j},\qquad j=1,\cdots, r_2,\\
\frac{\partial t_{\tau_j}}{\partial Y_q}=&2\log\Big|\varepsilon_q^{(j)}
\Big|^2\cdot t_{\tau_j},\qquad j=1,\cdots, r_2,\ q=1,\cdots,r_1+r_2-1.
\end{aligned}$$ 

Thus by the  formula $\tilde g_{ij}=\sum_{\alpha,\beta}\frac{\partial x^\alpha}
{\partial\tilde  x^i}\frac{\partial x^\beta}{\partial\tilde  x^j}
g_{\alpha\beta}$ for the change of variables
and the symmetry of the matric matrix, 
we are lead to calculate the following three types of products of matrices:
$$\begin{aligned}&X_0\cdot\mathrm{diag}\Big(\frac{1}{y_1^2},\,\cdots,\, 
\frac{1}{y_{r_1}^2},\
\Big(\frac{1}{v_1^2}\Big)^2,\,\cdots \,\Big(\frac{1}{v_{r_2}^2}\Big)^2\Big)
\cdot X_0^t,\\
&X_0\cdot\mathrm{diag}\Big(\frac{1}{y_1^2},\,\cdots,\, 
\frac{1}{y_{r_1}^2},\
\Big(\frac{1}{v_1^2}\Big)^2,\,\cdots \,\Big(\frac{1}{v_{r_2}^2}\Big)^2\Big)
\cdot X_q^t,\\
&X_p\cdot\mathrm{diag}\Big(\frac{1}{y_1^2},\,\cdots,\, 
\frac{1}{y_{r_1}^2},\
\Big(\frac{1}{v_1^2}\Big)^2,\,\cdots \,\Big(\frac{1}{v_{r_2}^2}\Big)^2\Big)
\cdot X_q^t\end{aligned}$$
where $$X_0:=
\Big(\frac{y_{\sigma_1}}{(r_1+r_2)Y_0},\cdots,\frac{y_{\sigma_{r_1}}}
{(r_1+r_2)Y_0},\frac{t_{\tau_1}}{(r_1+r_2)Y_0},\cdots, 
\frac{t_{\tau_{r_2}}}{(r_1+r_2)Y_0}\Big)$$
and $$X_p=\Bigg(2\log\Big|\varepsilon_p^{(1)}\Big|\cdot y_{\sigma_1},\cdots,
2\log\Big|\varepsilon_p^{(r_1)}\Big|\cdot y_{\sigma_{r_1}},
2\log\Big|\varepsilon_p^{(r_1+1)}\Big|^2\cdot v^2_{\tau_1},\cdots,
2\log\Big|\varepsilon_p^{(r_1+r_2)}\Big|^2\cdot v^2_{\tau_{r_2}}\Bigg),$$
for $p,\,q=1,2,\cdots,r_1+r_2-1$.

Hence, $$\tilde g_{11}=\frac{1}{(r_1+r_2)Y_0^2},\qquad 
\tilde g_{1j}=\tilde g_{j1}=0,\qquad j=2,\cdots, r_1+r_2,$$
since $$\sum_{i=1}^{r_1}\log\Big|\varepsilon_p^{(i)}\Big|
+\sum_{j=1}^{r_2}\log\Big|\varepsilon_p^{(r_1+j)}\Big|^2=0;$$
while
$$\tilde g_{(i+1)(j+1)}=4\sum_{p=1}^{r_1}\log\Big|\varepsilon_i^{(p)}\Big|
\log \Big|\varepsilon_j^{(p)}\Big|+4\sum_{p=r_1+1}^{r_1+r_2}
\log\Big|\varepsilon_i^{(p)}\Big|^2
\log \Big|\varepsilon_j^{(p)}\Big|^2$$ for $i,j=1,\cdots, r_1+r_2-1.$

Consequently, turning to the new volume element, we use 
$$\det(\tilde g_{ij})=\frac{4^{r_1+r_2-1}}{(r_1+r_2)Y_0^2}R^2,$$ where 
$$R:=(\log\Big\|\varepsilon_q^{(p)}\Big\|)_{p,q=1}^{r_1+r_2-1}$$ 
is the regulator of $K$. (See e.g. [Neu].) 
Thus, by taking also account of ReZ direction, we have
$$\begin{aligned}d\omega=&
\Big(\sqrt{\det(\tilde g_{ij})\cdot \frac{1}{Y_0^2}}\Big)\,dY_0\,dY_1\cdots 
dY_{r_1+r_2-1}\,
\prod_{\sigma:\mathbb R}dx_\sigma\prod_{\tau:\mathbb C} dz_\tau\\
=&\frac{2^{r_1+r_2-1}}{\sqrt{r_1+r_2}}
R\,\frac{dY_0}{Y_0^2}\,dY_1\cdots dY_{r_1+r_2-1}
\,\prod_{\sigma:\mathbb R}dx_\sigma\prod_{\tau:\mathbb C} dz_\tau.
\end{aligned}$$
Clearly, the boundary $\partial \mathcal D_T$ of 
$\mathcal D_T$ consists of 

\noindent
1) (the corresponding parts of) the
boundary of the fundamental domain of $\mathcal D$; and 

\noindent
2) the hyperplane of $\mathcal D$ defined by the condition 
$Y_0=T':=N(\frak a\frak b^{-2})\cdot T$.

\noindent
(Please note that the factor $N(\frak a\frak b^{-2})$ is added in 2). 
This is because in the definition of $\mathcal D_T$, what we used is the 
distance to cusp, not simply $Y_0$.)

\noindent
Consequently, if we set $d\mu$ to be the 
volume element of this hypersurface, then
$$d\mu=\frac{1}{\sqrt{\tilde g_{11}}}d\omega\Big|_{Y_0=T'}=
\frac{\sqrt{r_1+r_2}}{T'}
2^{r_1+r_2-1}R\,dY_1\cdots dY_{r_1+r_2-1}\prod_{\sigma:\mathbb R}
dx_\sigma\prod_{\tau:\mathbb C} dz_\tau.$$
Moreover, if we let $\nu$ be the unit normal to the hypesurface, since 
$$\Big\langle \frac{\partial}{\partial Y_0},\frac{\partial}{\partial Y_0}
\Big\rangle={\tilde g}^{11}\Big|_{Y_0=T'}=(r_1+r_2){T'}^2,$$ we have 
$$\nu=\Big(\frac{1}{\sqrt{r_1+r_2}T'},0,
\ldots,0\Big).$$ Thus the outer normal derivative of a function $f$ 
is given by 
$$\frac{\partial f}{\partial\nu}=\Big(\frac{1}{\sqrt{r_1+r_2}T'},0,\ldots,0
\Big)
\cdot\mathrm{grad}\, f=\sqrt {r_1+r_2}\cdot T'\,
\frac{\partial f}{\partial Y_0}.$$

Now by (the fact that the group $SL(\mathcal O_K\oplus\frak a)$ is 
finitely generated and) the concrete construction of our fundamental domain, 
we see that the boundary 
$\partial \mathcal D_T $ consists of finitely many of surfaces 
which are either parts of horospheres or parts $X_i(T)$
of planes cut out  by $Y_0=T_i'$, where $T_i'=N(\frak a\frak b_i^{-2})\cdot T$
(with $\frak b_i$ the fractional ideal associated to the cusp $\eta_i$).
Moreover, besides the  hyperplanes associated with 
$Y_0=T'$, the set of horospheres 
appeared on the boundary is divided into the sets of equivalent pairs 
for which 
the integral of the outer normal derivative along one surface in a pair is 
equal to the integral of the inner normal derivative along the other surface 
in the pair. (Say, in terms of $Y_{p\geq 1}$, they are given by
$Y_p=\pm\frac{1}{2}$ in pairs.)  
As such, we further conclude that
$$\begin{aligned}&\iiint_{\mathcal D_T}
\widehat E_{2,\frak a}(\tau,s)d\mu(\tau)\\
=&\frac{1}{r_1+4r_2}
\frac{1}{s(s-1)}\iiint_{\mathcal D_T}
\Delta_K\widehat E_{2,\frak a}(\tau,s)d\mu(\tau)\\
=&\iint_{\partial \mathcal D_T}
\frac{\partial \widehat E_{2,\frak a}(\tau,s)}{\partial \nu}d\mu\\
=&\frac{1}{r_1+4r_2}
\frac{1}{s(s-1)}\sum_{i=1}^h\iint_{X_i(T)}
\frac{\partial \widehat E_{2,\frak a}(\tau,s)}{\partial \nu}d s,\end{aligned}
$$ where $X_i(T)$ denotes the part of the boundary of $\mathcal D_T$ coming 
from the pull back of the intersection of the hypersurface $Y_0=T_i'$ with 
$F_{\eta_i}, i=1,2,\ldots,h$. 
(Here we used the fact that for $T\geq 1$, $X_i(T)$ are disjoint from each 
other. See e.g., the Lemma in 2.5.6.)

Now, we are ready to use Fourier expansion to do the calculation. Note 
that the average for $e^{2\pi i t}$ 
(together with its derivative) over an interval of length 1 is zero. 
Hence, in the above integration for $\widehat E$ over $\mathcal D_T$, 
we are in fact left with only the constant terms of the Fourier expansion for 
$\widehat E(s)$. Consequently, with 
$T_i'=N(\frak a\frak b_i^{-2})\cdot T$, then, 
up to constant factors depending only on $K$,
$$\begin{aligned}&\iiint_{\mathcal D_T}
\widehat E_{2,\frak a}(\tau,s)d\mu(\tau)\\
=&\frac{1}{r_1+4r_2}
\frac{1}{s(s-1)}\iiint_{\mathcal D_T}
\Delta_K\widehat E_{2,\frak a}(\tau,s)d\mu(\tau)\\
=&\frac{1}{r_1+4r_2}
\frac{1}{s(s-1)}\sum_{i=1}^h\iint_{X_i(T)}
\frac{\partial }{\partial \nu}\Big(A_{0i} Y_0^s+B_{0i}Y_0^{1-s}\Big)d \mu\\
=&\frac{1}{r_1+4r_2}
\frac{1}{s(s-1)}\sum_{i=1}^h\iint_{X_i(T)}\sqrt{r_1+r_2}\cdot T_i'
\frac{\partial }{\partial Y_0}\Big(A_{0i} Y_0^s+B_{0i}Y_0^{1-s}\Big)\\
&\qquad\cdot
\frac{\sqrt {r_1+r_2}}{T_i'}2^{r_1+r_2-1}R\cdot dY_1\ldots dY_{r_1+r_2-1}\cdot
\prod_{\sigma:\mathbb R}dx_\sigma\cdot\prod_{\tau:\mathbb C}dz_\tau\\
=&\frac{r_1+r_2}{r_1+4r_2}
\frac{1}{s(s-1)}\cdot 2^{r_1+r_2-1}R\cdot
\sum_{i=1}^h \iint_{X_i(T)} \Big(s\, A_{0i} {Y_0}^{s-1}-(s-1)\,
B_{0i}{Y_0}^{-s}\Big)\\
&\qquad dY_1\ldots dY_{r_1+r_2-1}\cdot
\prod_{\sigma:\mathbb R}dx_\sigma\cdot\prod_{\tau:\mathbb C}dz_\tau\\
=&\frac{r_1+r_2}{r_1+4r_2}
\frac{1}{s(s-1)}2^{r_1+r_2-1}R\\
&\cdot\sum_{i=1}^h\Big(s\, A_{0i} {T_i'}^{s-1}-(s-1)\,
B_{0i}{T_i'}^{-s}\Big)\cdot  \iint_{X_i(T)} 
dY_1\ldots dY_{r_1+r_2-1}\cdot
\prod_{\sigma:\mathbb R}dx_\sigma\cdot\prod_{\tau:\mathbb C}dz_\tau\\
=&\frac{r_1+r_2}{r_1+4r_2}
2^{r_1+r_2-1}R\cdot 
D_K^{\frac{1}{2}} \sum_{i=1}^h N(\frak a\frak b_i^{-2})\cdot
\Big(\frac{A_{0i}}{s-1}\cdot {T_i'}^{s-1}-
\frac{B_{0i}}{s}{T_i'}^{-s}\Big),\end{aligned}$$ due to the fact that the 
lattice corresponding 
to the cusp $\eta_i=\frac{\alpha_i}{\beta_i}$ is given by 
$\frak a\frak b_i^{-2}$ with $\frak b_i=\mathcal O_K\alpha_i+\frak a\beta_i$
and $Y_p\in[-\frac{1}{2},\frac{1}{2}]$.
Thus with the precisely formula we have for $A_{0i}(s)$ and the functional 
equation with the change $s\leftrightarrow 1-s$, we finally obtain 
the following
\vskip 0.30cm
\noindent
{\bf {\Large Theorem.}} {\it Up to a constant factor depending only on  $K$,
$$\iiint_{\mathcal D_T}\widehat E_{2,\frak a}(\tau,s)d\mu(\tau)=
\frac{\xi_K(2s)}{s-1}\Big(\Delta_K\cdot T\Big)^{s-1}-\frac{\xi_K(2-2s)}{s}
\Big(\Delta_K\cdot T\Big)^{-s}.$$}

\noindent
Proof. Indeed, by the functional equation, we only need to calculate the 
coefficient of $\frac{T^{s-1}}{s-1}$.
Note that, by Theorem in 3.3.2, the partial constant term $A_{0,i}$ 
for the completed Eisenstein series $\widehat E_{2,\frak a}(\tau,s)$ 
is given by $$\Bigg(\Big(\pi^{-s}\Gamma(s)\Big)^{r_1}
\Big(2\pi^{-2s}\Gamma(2s)\Big)^{r_2}\Big(N(\frak a)\Delta_K\Big)^s\Bigg)\cdot
\zeta([\frak a^{-1}\frak b_i], 2s)\cdot N(\frak a\frak b_i^{-1})^{-2s}.$$
Hence, up to a constant factor depending only on $K$, the coefficient
of $\frac{T^{s-1}}{s-1}$ in the integration
$\iiint_{\mathcal D_T}\widehat E_{2,\frak a}(\tau,s)d\mu(\tau)$
is simply the summation $\sum_{i=1}^h$ of
$N(\frak a\frak b_i^{-2})A_{0,i}$ timing with the factor
$N(\frak a\frak b_i^{-2})^{s-1}$
resulting from the discrepency between $T$ and $T_i'$. That is to say,
up to a constant factor depending only on  $K$,
the coefficient of $\frac{T^{s-1}}{s-1}$ is nothing but
$$\begin{aligned}&\sum_{i=1}^h\Big(\pi^{-s}\Gamma(s)\Big)^{r_1}
\Big(2\pi^{-2s}\Gamma(2s)\Big)^{r_2}\Big(N(\frak a)\Delta_K\Big)^s\cdot
\zeta([\frak a^{-1}\frak b_i], 2s)N(\frak a\frak b_i^{-1})^{-2s}
\cdot N(\frak a\frak b_i^{-2})\cdot N(\frak a\frak b_i^{-2})^{s-1}\\
=&\Delta_K^s\cdot\Big(\pi^{-s}\Gamma(s)\Big)^{r_1}
\Big(2\pi^{-2s}\Gamma(2s)\Big)^{r_2}\sum_{i=1}^h
\zeta([\frak a^{-1}\frak b_i],2s)\\
=&\Delta_K^s\cdot \Big(\pi^{-s}\Gamma(s)\Big)^{r_1}
\Big(2\pi^{-2s}\Gamma(2s)\Big)^{r_2}\zeta_K(2s)\\
=&\Delta_K^s\cdot\xi_K(2s),\end{aligned}$$
since $$\sum_{i=1}^h
\zeta([\frak a^{-1}\frak b_i],2s)=\zeta_K(2s),$$ 
resulting from the facts that 

\noindent
(i) the $h$ ideal classes $[\frak a^{-1}\frak b_i]$ for fixed 
$\frak a$ run over all elements of the class group of $K$; and that

\noindent
(ii) the total Dedekind zeta function decomposes into a summation of partial 
zeta functions associated to ideal classes.

\noindent 
This completes the proof.

Consequently, we have the following
\vskip 0.30cm
\noindent
{\bf{\Large Fact}} (VIII) {\it Up to a constant factor depending only on  $K$,
the rank two non-abelian zeta function $\xi_{K,2}(s)$ is given by
$$\xi_{K,2}(s)=
\frac{\xi_K(2s)}{s-1}\Delta_K^{s-1}-\frac{\xi_K(2s-1)}{s}
\Delta_K^{-s}\qquad \Re(s)\,>1.$$}

\noindent
Proof. This is because, by Fact ? in 2.5, we have the moduli space of 
rank two semi-stable lattices of volume $N(\frak a)\Delta_K$ with underlying
projective module $\mathcal O_K\oplus \frak a$ is given by 
$\mathcal D_1$. But from the Theorem above, up to a constant factor 
depending only on $K$, $$\iiint_{\mathcal D_1}
\widehat E_{2,\frak a}(\tau,s)d\mu(\tau)=
\frac{\xi_K(2s)}{s-1}\Delta_K^{s-1}-\frac{\xi_K(2-2s)}{s}\Delta_K^{-s}.$$
That is to say, up to a constant factor 
depending only on $K$,
$$\xi_{K,2;\frak a}(s)=
\frac{\xi_K(2s)}{s-1}\Delta_K^{s-1}-\frac{\xi_K(2-2s)}{s}\Delta_K^{-s}.$$ 
Therefore, by Fact IV in I.1.9, up to a constant factor 
depending only on $K$, $$\xi_{K,2}(s)=
\frac{\xi_K(2s)}{s-1}\Delta_K^{s-1}-\frac{\xi_K(2-2s)}{s}\Delta_K^{-s}.$$ This 
completes the proof.
\vskip 0.30cm
\noindent
{\bf Remark.} The reason for the degeneration of rank two non-abelian zeta 
functions are explained in my paper on
\lq Analytic truncation and Rankin-Selberg versus algebraic 
truncation and non-abelian zeta',
{\it Algebraic Number Theory and Related Topics}, RIMS Kokyuroku, 
No.1324 (2003). As pointed there, when rank is 3 or bigger,
essential non-abelian parts resulting from higher Fourier coefficients in the 
Fourier expansion of Eisenstein series do give their contributions 
(to non-abelian zeta functions.)

\chapter{Zeros of Rank Two Non-Abelian Zeta Functions for Number Fields}

\section{Zeros of Rank Two Non-Abelian Zeta Function of $\mathbb Q$}

Following what was happened in history, let me first start 
with  Suzuki's weak  result [Su] and then give Lagarias' unconditional result.
\vskip 0.20cm
\noindent
{\bf {\Large Theorem.}} {\it If the Riemann Hypothesis for the Riemann zeta 
function holds, then all zeros of 
$\xi_{\mathbb Q,2}(s)$ lie on the critical line $\Re(s)=\frac {1}{2}.$}
\vskip 0.20cm
This is a very clever observation, rooted back to 
Titchimashi's book on Riemann Zeta Functions.

\subsection{Product Formula for Entire Function of Order 1}

Let $f(z)$ be an entire function of order one on $\mathbb C$, that is, 

\noindent
(i) $f(z)$ is analytic over $\mathbb C$;

\noindent
(ii) $f(z)=O\Big(\exp(|z|^{1+\varepsilon})\Big),\ \ \forall \varepsilon>0$.

\noindent
Let $n(R)$ denote the number of zeros of $f(z)$ inside $C_R$, the circle 
of radius $R$ centered at the origin. Then 

\noindent
(1) $n(R)=O(R^\alpha),\ \ \forall\alpha>1$;

\noindent
(2) $\sum_{\rho_n,f(\rho_n)=0}|\rho_n|^{-\alpha}$ converges. 
In particular, $(1-\frac{z}{\rho_n})
\exp(\frac{z}{\rho_n})=1+O\Big((\frac{z}{\rho_n})^2\Big)$ as $n\to\infty$. 
As a direct consequnce,

\noindent
(3) $P(z):=\prod_{\rho_n}\Big(1-\frac{z}{\rho_n}\Big)\cdot 
\exp\Big(\frac{z}{\rho_n}\Big)$ 
converges and $\frac{f(z)}{P(z)}$ is an entire
function of order 1 without zeros, hence should be in the form 
$\exp(A+Bz)$ for certain constants $A,\,B$. 

That is to say, we have the following

\noindent
{\bf \Large  Hadamard Product Theorem}.
{\it Let $f(z)$ be an entire function of order one on $\mathbb C$, 
then  there exist constants $A,\,B$ such that 
$$f(z)=e^{A+Bz}\cdot 
\prod_\rho\Big(1-\frac{z}{\rho}\Big)\cdot\exp\Big(\frac{z}{\rho}\Big).$$}

\noindent 
{\bf Example.} (See e.g. [Ed]) We have  $$\frac{1}{2}s(s-1)\cdot
\xi(s)=e^{A+Bz}\cdot 
\prod_\rho\Big(1-\frac{z}{\rho}\Big)\cdot\exp\Big(\frac{z}{\rho}\Big)$$ 
with $A=-\log 2,\, 
B=-\frac{\gamma}{2}-1+\frac{1}{2}\log 4\pi$,
where $\gamma=\lim_{n\to\infty}\Big(1+\frac{1}{2}+\cdots+\frac{1}{n}
-\log n\Big)$ denotes the Euler constant.

\subsection{Proof}

Let $$F(z)=-Z(\frac{1}{2}+2 i z)\qquad\mathrm{with}\qquad Z(s)=s(1-s)\xi(s).$$

\noindent
{\bf{\Large Proposition}.} (Suzuki) (1) 
$F(z+\frac{i}{4})-F(z-\frac{i}{4})=iz(1+4z^2)\,
\xi_{\mathbb Q,2}(\frac{1}{2}+iz).$

\noindent
(2) {\it Assume the RH, then all zeros of $F(z+\frac{i}{4})-F(z-\frac{i}{4})$ 
are real.}

\noindent
{\it In particular, then the RH implies that 
$\xi_{\mathbb Q,2}(\frac{1}{2}+zi)$ admits only real zeros.}

\noindent
Proof. (1) Simple calculation. Indeed,
$$F\Big(z+\frac{i}{4}\Big)=-Z\Big(\frac{1}{2}+2i(\frac{i}{4}+z)\Big)
=-Z\Big(\frac{1}{2}-\frac{1}{2}+2iz\Big)=-Z(2iz).$$ So 
$$F\Big(z-\frac{i}{4}\Big)=-Z\Big(1+2iz\Big)$$ and 
$$\begin{aligned}F(z+\frac{i}{4})-F(z-\frac{i}{4})=&(1+2iz)(-2iz)\xi(1+2iz)-2iz(1-2iz)\xi(2iz)\\
=&2iz(1-2iz)(1+2iz)\cdot
\Big(\frac{\xi(1+2iz)}{2iz-1}-\frac{\xi(2iz)}{1+2iz}\Big)\\
=&iz(1+4z^2)\cdot \Big(\frac{\xi\big(2(\frac{1}{2}+iz)\big)}
{(\frac{1}{2}+iz)-1}-\frac{\xi(2(\frac{1}{2}+iz)-1)}{\frac{1}{2}+iz}\Big)\\
=&iz(1+4z^2)\,\xi_{\mathbb Q,2}(\frac{1}{2}+iz).\end{aligned}$$

\noindent
(2) Clearly, $F(z)$ is an entire function of order 1, so there are 
constants $A,\,B$ such that
$$F(z)=e^{A+Bz}\cdot 
\prod_{\rho:F(\rho)=0}\Big(1-\frac{z}{\rho}\Big)\cdot
\exp\Big(\frac{z}{\rho}\Big).$$
Note that essentially, $\rho$ are zeros of the completed Riemann 
zeta but transformed from $z$ to $\frac{1}{2}+2iz$.
Hence, {\it by the RH}, all $\rho$ are real.

 Moreover, since $F(z)=-Z(\frac{1}{2}+2iz)$ with $Z(s)=s(1-s)\,\xi(s)$, 
we have for $x\in\mathbb R$,
 $$\overline{F(x)}=\overline{-Z\Big({1}/{2}+2ix\Big)}
=-Z\Big(\overline{{1}/{2}+2ix}\Big)
 =-Z(\frac{1}{2}-2ix\Big)$$ which by the functional equation is simply
 $$-Z\Big(1-(\frac{1}{2}-2ix)\Big)=-Z\Big(\frac{1}{2}+2ix\Big)=F(x).$$ 
That is to say, for $x\in\mathbb R$, $F(x)$ takes only 
 real values. Hence, constants $A$ and $B$ are both real.

Now let $z_0=x_0+iy_0$ be a zero of $$F(z+\frac{i}{4})-F(z-\frac{i}{4})
=iz(1+4z^2)\,\xi_{\mathbb Q,2}
(\frac{1}{2}+iz).$$ Then $z_0=0$ and/or $z_0$ is a zero of 
$\xi_{\mathbb Q,2}(\frac{1}{2}+iz)$ since 
$\xi_{\mathbb Q,2}(\frac{1}{2}+iz)$ admits simple poles 
at $z=\pm\frac{1}{2}i$. 
 
In any case, $$F(z_0+\frac{i}{4})=F(z_0-\frac{i}{4}).$$
By taking absolute values on both sides,
$$\begin{aligned}&\Big|e^{A+B(z_0+\frac{i}{4})}\cdot 
\prod\Big(1-\frac{z_0+\frac{i}{4}}{\rho}\Big)\cdot
\exp\Big(\frac{z_0+
\frac{i}{4}}{\rho}\Big)\Big|\\
&\qquad=\Big|e^{A+B(z_0-\frac{i}{4})}
\cdot 
\prod\Big(1-\frac{z_0-\frac{i}{4}}{\rho}\Big)\cdot
\exp\Big(\frac{z_0-\frac{i}{4}}{\rho}\Big)\Big|.\end{aligned}$$
Since $B\in\mathbb R$ and $\rho_n\in\mathbb R$ (which is obtained by the RH 
as said above), we hence get
$$1=\prod_{n=1}^\infty\frac{(x_0-\rho_n)^2+(y_0-
\frac{1}{4})^2}{(x_0-\rho_n)^2+(y_0+\frac{1}{4})^2}.$$ Thus if 
$y_0>0$, then the right hand side is $<1$, while 
if $y_0<0$, then the right hand side is $>1$. 
Contradiction. This leads then $y_0=0$, hence completes the proof.
\vskip 0.20cm
With this in mind, note that in the proof above, the RH was 
used to ensure that $\rho$ are real, which have 
the effect that then in the calculation for the exponential 
factor $\exp\Big(\frac{z_0+\frac{i}{4}}{\rho}\Big)$, 
the ratio
$\frac{\Big|\exp\Big(\frac{z_0+\frac{i}{4}}{\rho}\Big)\Big|}
{\Big|\exp\Big(\frac{z_0-\frac{i}{4}}{\rho}\Big)\Big|}$ gives us
the exact value 1. That is to say, this factor of ratio of exp's
does not contribute.

However, one does not need such an argument from the very beginning
to eliminate the factors $\exp\Big(\frac{z_0\pm \frac{i}{4}}{\rho}\Big)$. 
In fact, this is the improvement of Lagarias, who gets his own unconditional
result totally independently, as a part of his understanding of 
de Branges's work  ([Lag]). 
The trick is very simple: Use the functional 
equation. So instead of working on individual $\rho_n$ in the product, 
we may equally use the functional equation 
to pair $\rho$ and $1-\rho$ for the zeros of
 the completed Riemann zeta function, or even to group $\rho$, $1-\rho$,
$\overline \rho$ and $1-\overline\rho$ together. Consequently, the exponential
factor appeared inside the infinite 
product  may be totally omitted. That is to say,
from the very beginning, we may simply assume that 
the Hadamard product involved takes the form $$F(z)=e^{A+Bz}\cdot 
\prod_{\rho:F(\rho)=0}'\Big(1-\frac{z}{\rho}\Big)$$ where $\prod'$ means 
that $\rho$'s are paired or grouped as above.
Form here, it is an easy exercise to deduce the following result of 
(Suzuki and) Lagarias.
\vskip 0.20cm
\noindent
{\Large Fact} (IX$)_{\mathbb Q}$ 
{\it All zeros of $\xi_{\mathbb Q,2}(s)$ lie on the line 
$\Re(s)=\frac {1}{2}$.}

\noindent
Proof. Alternatively, as above, we have $$\begin{aligned}&\Big|e^{A+B(z_0+\frac{i}{4})}\cdot 
\prod\Big(1-\frac{z_0+\frac{i}{4}}{\rho}\Big)\cdot
\exp\Big(\frac{z_0+
\frac{i}{4}}{\rho}\Big)\Big|\\
&\qquad=\Big|e^{A+B(z_0-\frac{i}{4})}
\cdot 
\prod\Big(1-\frac{z_0-\frac{i}{4}}{\rho}\Big)\cdot
\exp\Big(\frac{z_0-\frac{i}{4}}{\rho}\Big)\Big|.\end{aligned}$$
Since $B\in\mathbb R$, so if we can take care of
 the factors $\exp\Big(\frac{z_0+
\frac{i}{4}}{\rho}\Big)$ and  $\exp\Big(\frac{z_0-\frac{i}{4}}{\rho}\Big)$
in a nice way, we are done. For this, as said above, let us group $\rho,\,
\bar\rho,\,1-\rho,\,1-\bar\rho$ together, we see that
$\frac{1}{\rho}+\frac{1}{\bar\rho}=\frac{2\,\Re(\rho)}{|\rho|^2}$
and $\frac{1}{1-\rho}+\frac{1}{1-\bar\rho}=\frac{2-2\,\Re(\rho)}{|1-\rho|^2}$ 
are all reals, hence, the same prove as above works.

It is very beautiful. Is not it?!!!

\subsection{A Simple Generalization}

The above method works for the functions $\xi_{\mathbb Q,2}^T(s)$ as well, 
provided that $T\geq 1$. Indeed, first, recall that we have the precise
relation
$$\xi_{\mathbb Q,2}^T(s)=\frac{\xi(2s)}{s-1}\cdot T^{s-1}-
\frac{\xi(2s-1)}{s}\cdot T^{-s}.$$
Consequently,
$$F(z+\frac{i}{4})\cdot T^{-\frac{1}{2}-iz}-
F(z-\frac{i}{4})\cdot T^{-\frac{1}{2}+iz}
=iz(1+4z^2)\,\xi_{\mathbb Q,2}^T
(\frac{1}{2}+zi).$$
Therefore, using the same proof, we arrive at
the relation 
$$\begin{aligned}&\Big|e^{A+B(z_0+\frac{i}{4})}\cdot 
\prod\Big(1-\frac{z_0+\frac{i}{4}}{\rho}\Big)\cdot
\exp\Big(\frac{z_0+
\frac{i}{4}}{\rho}\Big)\Big|\cdot\Big|T^{-iz-\frac{1}{2}}\Big|\\
&\qquad=\Big|e^{A+B(z_0-\frac{i}{4})}
\cdot 
\prod\Big(1-\frac{z_0-\frac{i}{4}}{\rho}\Big)\cdot
\exp\Big(\frac{z_0-\frac{i}{4}}{\rho}\Big)\Big|\cdot 
\Big|T^{iz-\frac{1}{2}}\Big|.\end{aligned}$$
That is to say,
$$1=\prod_{n=1}^\infty\frac{(x_0-\rho_n)^2+(y_0-
\frac{1}{4})^2}{(x_0-\rho_n)^2+(y_0+\frac{1}{4})^2}
\cdot\frac{T^{-y_0}}{T^{y_0}}.$$
Or equivalently,
$$T^{2y_0}=\prod_{n=1}^\infty\frac{(x_0-\rho_n)^2+(y_0-
\frac{1}{4})^2}{(x_0-\rho_n)^2+(y_0+\frac{1}{4})^2}.$$

Thus with $T\geq 1$, we have

\noindent
(i) if $y_0>0$, the left hand side is $>1$, while
the right hand side is $<1$, contradiction; while

\noindent
(ii) if $y_0<0$, the left hand side is $<1$, while
the right hand side is $>1$, contradiction.
That is to say, we obtain the following
\vskip 0.20cm
\noindent
{\bf \Large Fact} (IX$')_{\mathbb Q}$ {\it For $T\geq 1$, all zeros of 
$\xi_{\mathbb Q,2}^T(s)$ lie on the critical line 
$\Re(s)=\frac {1}{2}$.}

Recall that 
$$\xi_{\mathbb Q,2}^T(s)=\frac{\xi(2s)}{s-1}\cdot T^{s-1}-
\frac{\xi(2s-1)}{s}\cdot T^{-s}.$$ Clearly 
$$\lim_{T\to 1^+}\xi_{\mathbb Q,2}^T(s)=\xi_{\mathbb Q,2}(s),$$
while  $\lim_{T\to +\infty}\xi_{\mathbb Q,2}^T(s)$ does not really 
make any sense. This says that  even when we have a family
of natural functions whose zeros all lie on the critical line,
in general, we cannot take limit for this family to preserve this property, 
no matter how careful we are. 

\section{Zeros of Rank Two Zetas for Number Fields: 
Generalized Riemann Hypothesis}

Finally, we are ready to state the following
\vskip 0.30cm
\noindent
{\bf \Large Fact} (IX) {\it All zeros of rank two non-abelian zeta functions
for number fields are on the critical line
$\Re(s)=\frac{1}{2}$.}

\noindent
Proof. This is a direct consequence of the following three facts:

\noindent
First, we know that, by the Rankin-Selberg \& Zagier method,
up to a constant factor depending only on $K$,
$$\xi_{K,2}(s)=\frac{\xi_K(2s)}{s-1}\Delta_K^{s-1}-
\frac{\xi_K(2s-1)}{s}\Delta_K^{-s}.$$ 
Secondly, $\Delta_K\geq 1$ for any number field; and Thridly, $s(s-1)\cdot
\xi_K(s)$ is also an entire function of order one [L1]. Consequently, 
the proof in the previous section, more precisely, that of 5.1.3, 
on the zeros works here as well
by a simple change from the Riemann $\xi$ for the field of rationals 
$\mathbb Q$ to the Dedekind $\xi_K$ for the number field $K$.
\vskip 0.30cm
\noindent
{\bf Remarks.} 1) For rank two non-abelian zeta of function fields, 
it is known that  Siegel type zeros do appear.
Moreover, it is expected that a modified Riemann Hypothesis holds as well. 
For example, by the precise formula given in [We4], one checks easily that
 rank two non-zbelian zeta of elliptic curves satisfy the modified 
Riemann Hypothesis. (The details for genus two curves will be given elsewhere.)

\noindent
2) One may wonder why the Riemann Hypothesis for rank two zeta can be proved.
One explanation is that there is an additional symmetry for them: while rank
two zetas
are supposed to be non-abelain by definition, they degenerate into 
combinations of abelian zetas (for the reason that there are not enough 
parabolic subgroups in $SL_2$).

\noindent 
3) Even though the proof above for the RH of rank two zeta is 
technically correct, 
it is not really  philosophically right. A geniune proof should 
work uniformly for all types of zetas: Both abelian and non-abelian zetas
are supposed to have Euler products resulting from 
abelian and non-abelain reciprocity law (see the discussion in [We1,2],) --
Abelian zetas admit commutative Euler product, while non-abelian zetas admit
non-commutative Euler product.
\newpage
\centerline{\Large REFERENCES}
\vskip 0.45cm
\noindent
\vskip 0.20cm
\noindent
[Bo1]  Borel, A. Some finiteness properties of adele groups over
number fields, Publ. Math., IHES, {\bf 16} (1963) 5-30
\vskip 0.20cm
\noindent
[Bo2] Borel, A. {\it Introduction aux groupes arithmetictiques}, Hermann, 1969
\vskip 0.30cm
\noindent 
[Bor] A. Borisov, Convolution structures and arithmetic
cohomology, to appear in Compositio Math.
\vskip 0.20cm
\noindent
[Bu2] Bump, D. The Rankin-Selberg method: a survey. {\it Number theory, 
trace formulas 
and discrete groups} (Oslo, 1987), 49--109, Academic Press, Boston, MA, 1989.
\vskip 0.20cm
\noindent 
[Ed] Edwards, H.M. {\it Riemann's Zeta Function}, Dover Publications, INC., 
1974
\vskip 0.20cm
\noindent 
[Ef] Efrat, I.Y. {\it The Selberg Trace Formula for $PSL_2(\mathbb R)^n$}, 
Memoirs of AMS, no. 359, 1987
\vskip 0.20cm
\noindent 
[EGM] Elstrodt, J. Grunewald, F. \& Mennicke, J.L. {\it Groups Acting on 
Hyperbolic Space: Harmonic Analysis and Number Theory},
 Springer Verlag(1997)
\vskip 0.20cm
\noindent 
[FT] Fr\"ohlich, A. \& Taylor, M.J. {\it Algebraic Number Theory}, Cambridge 
studies in advanced mathematics {\bf 27}, Cambridge Univ. Press (1991)
\vskip 0.20cm
\noindent 
[Ge] van der Geer, G. {\it Hilbert Modular Surfaces}, Ergebnisse der 
Mathematik und ihrer Grenzgebiete, 3. Folge Bd.16, Springer (1988) 
\vskip 0.20cm
\noindent 
[GS] van der Geer, G. \&  Schoof, R. Effectivity of Arakelov
Divisors and the Theta Divisor of a Number Field, Sel. Math., New ser.
{\bf 6} (2000), 377-398  
\vskip 0.20cm
\noindent 
[Gr] Grayson, D.R. Reduction theory using semistability. 
Comment. Math. Helv.  {\bf 59}  (1984),  no. 4, 600--634.
\vskip 0.20cm
\noindent 
[Gu] Gupta, Sh.D. On the Rankin-Selberg Method for functions not of rapid decay
on congruence subgroups, J Number Theory., no. 120, 1997
\vskip 0.20cm
\noindent 
[Ha] Hayashi, Tsukasa: The degeneration of new non-abelian zeta function 
of rank 2, 
priprint, 2003
\vskip 0.20cm
\noindent 
[Iw1] Iwasawa, K. Letter to Dieudonn\'e, April 8, 1952, in
 {\it Zeta Functions in Geometry}, 
Advanced Studies in Pure Math. {\bf 21} (1992), 445-450
\vskip 0.20cm
\noindent 
[Iw2] Iwasawa, K. {\it Lectures notes on Riemann(-Artin) Hypothesis,} 
noted by Kimura, Princeton Univ. 197?
\vskip 0.20cm
\noindent 
[Kub] Kubota, T. {\it Elementary theory of Eisenstein series.} Kodansha Ltd.,
Halsted Press, 1973.
\vskip 0.20cm
\noindent 
[Lag] Lagarias, J. private communication
\vskip 0.20cm
\noindent 
[L1]  Lang, S. {\it Algebraic Number Theory}, 
Springer-Verlag, 1986
\vskip 0.20cm
\noindent 
[L2] Lang, S. {\it Introduction to Arekelov Theory}, 
Springer-Verlag, 1988
\vskip 0.20cm
\noindent 
[Mi] Miyake, M. {\it Modular Forms}, Springer (1989)
\vskip 0.20cm
\noindent 
[Mo] Moreno, C. {\it  Algebraic curves over finite fields.}
Cambridge Tracts in Mathematics, {\bf 97}, Cambridge University Press, 1991
\vskip 0.20cm
\noindent 
[Neu] Neukirch, J. {\it Algebraic Number Theory}, Grundlehren der
Math. Wissenschaften, Vol. {\bf 322}, Springer-Verlag, 1999
\vskip 0.20cm
\noindent 
[Se] Serre, J.-P. {\it Algebraic Groups and Class Fields}, GTM
{\bf 117}, Springer-Verlag (1988)
\vskip 0.20cm
\noindent 
[S] Siegel, C.L. {\it Lectures on advanced analytic number theory.} Notes by 
S. Raghavan. Tata Institute of Fundamental Research Lectures on 
Mathematics, No. 23  1965
\vskip 0.30cm
\noindent 
[St] Stuhler, U. Eine Bemerkung zur Reduktionstheorie
quadratischer Formen.  Arch. Math.  {\bf 27} (1976), no. 6,
604--610.
\vskip 0.20cm
\noindent 
[Su] Suzuki, Masatoshi: priviate communication
\vskip 0.20cm
\noindent 
[Ta] Tate, J. Fourier analysis in number fields and Hecke's
zeta functions, Thesis, Princeton University, 1950 
\vskip 0.20cm
\noindent
[Te] Terras, A. {\it Harmonic analysis on symmetric spaces and applications II}, Springer-Verlag, 1988
\vskip 0.20cm
\noindent 
[Ti] Titchmarsh, E.C. {\it The Theory of Riemann Zeta-Function}, Oxford Univ. 
Press, 1951
\vskip 0.20cm
\noindent
[W] Weil, A. {\it Basic Number Theory}, Springer-Verlag, 1973
\vskip 0.20cm
\noindent 
[We1] Weng, L.  A Program for Geometric Arithmetic, at math.AG/0111241
\vskip 0.20cm
\noindent 
[We2] Weng, L. Non-Abelian $L$-Functions for Number Fields, submitted, 2003
\vskip 0.20cm
\noindent 
[We3] Weng, L. Non-Abelian Zeta Functions for Function Fields, Amer J Math., 
to appear
\vskip 0.20cm
\noindent 
[We4] Weng, L. Refined Brill-Noether locus and non-abelian Zeta 
functions for elliptic curves, {\it Algebraic geometry in East Asia} 
(Kyoto, 2001), 245--262, World Sci. 2002
\vskip 0.20cm
\noindent
[Z] Zagier, D.  The Rankin-Selberg method for automorphic functions which are not of rapid decay. 
J. Fac. Sci. Univ. Tokyo Sect. IA Math. 28(3), 415--437 (1982). 
\vskip 2.0cm
This is a corrected version of the one distributed on October 23, 2004.
\enddocument